\documentclass{amsart}

\usepackage{color}
\usepackage[dvips]{graphicx}
\usepackage{amscd}
\usepackage{amsmath}
\usepackage{ascmac}
\usepackage{cases}
\usepackage{amssymb}
\usepackage{amsfonts}
\usepackage{amsthm}
\usepackage[all]{xy}
\usepackage{enumerate}
\usepackage{wrapfig}

\theoremstyle{definition}
\newtheorem{Thm}{{\bf Theorem}}[section]

\newtheorem{Lem}[Thm]{{\bf Lemma}}
\newtheorem{Prop}[Thm]{{\bf Proposition}}
\newtheorem{Cor}[Thm]{{\bf Corollary}}

\newtheorem{Def}[Thm]{{\bf Definition}}
\newtheorem{Rem}[Thm]{{\bf Remark}}

\newcounter{Exami}

\numberwithin{equation}{section}

\title[Description of generalized isomonodromic deformations]{
Description of generalized isomonodromic deformations 
of rank two linear differential equations
using apparent singularities}
\author{Arata Komyo}
\address{Center for Mathematical and Data Sciences, Kobe Univ.,
1-1 Rokkodai-cho, Nada-ku, Kobe, 657-8501, Japan}
\email{akomyo@math.kobe-u.ac.jp}

\begin{document}

\maketitle
\thispagestyle{empty}

\begin{abstract}
In this paper, we consider the generalized isomonodromic deformations 
of rank two irregular connections on the Riemann sphere. 
We introduce Darboux coordinates on the parameter space of a family of rank two irregular connections
by apparent singularities.
By the Darboux coordinates,
we describe the generalized isomonodromic deformations
as Hamiltonian systems.
\end{abstract}

\section{Introduction}

For connections on the trivial bundle on $\mathbb{P}^1$,
the
regular singular isomonodromic deformation is the Schlesinger equation 
and the unramified irregular singular generalized isomonodromic deformation is 
the Jimbo--Miwa--Ueno
equation which is completely given in \cite{JMU-1}, \cite{JM-1}, and \cite{JM-2}.
Bertola--Mo and
Bremer--Sage have generalized 
the Jimbo--Miwa--Ueno equation (see \cite{BM}, \cite{BS-1}, and \cite{BS-2}).
Boalch \cite{Boalch} has given the symplectic geometry of 
Jimbo--Miwa--Ueno equation.
That is, the Jimbo--Miwa--Ueno equations 
are equivalent to a flat symplectic Ehresmann connection 
on a certain symplectic fiber bundle.
The fibers of the symplectic fiber bundle 
are certain
moduli spaces of meromorphic connections over $\mathbb{P}^1$.
In this paper, we consider 
the generalized isomonodromic deformation
from this point of view.

As in \cite{JMU-1} 
the monodromy data 
for certain families of irregular singular differential equations
involve the asymptotic behavior of solutions along Stokes sectors at each singular point.
Here we impose that the singularities of irregular singular differential equations 
satisfy some generic condition. 
More precisely, these singularities are regular or unramified irregular. 
If we have such a generic family of irregular singular differential equations,
then we can locally define the monodromy map (in other words, Riemann--Hilbert map) from 
the space of parameters of this family to the moduli space of irregular monodromy representations
(for details, for example, see \cite{Boalch} and \cite{IS}).
The fibers of the monodromy map are a foliation of the space of parameters of the family
(see \cite{Boalch}, \cite{IS} and \cite{Malgrange-1}). 
The foliation is called the ({\it generalized}) {\it isomonodromic foliation}.
(The corresponding vector field is called ({\it generalized}) {\it isomonodromic deformations}).
On the other hand, there exists another approach to 
generalized isomonodromic deformations.
As in \cite[Appendix]{Boalch}, 
a submanifold $\mathcal{L}$ in the space of parameters of the family
is contained in a leaf of this foliation if and only if the
family of connections corresponding to $\mathcal{L}$ is integrable.

When irregular singular differential equations have special singularities
(so-called (generic) ramified irregular singularities), 
the formulation of the Riemann--Hilbert map is still not clear. 
But by using the integrable condition (which is the second point of view),
we can define the generalized isomonodromic deformation for 
such a special family of irregular singular differential equations.
Moreover, this generalized isomonodromic deformation is integrable.
So we have the generalized isomonodromic foliation
(see \cite{BM}, \cite{BS-2}, and \cite{I2}).

In this paper, 
we consider generalized isomonodromic deformations only
from the viewpoint of the integrability condition.
More specifically, 
we construct 
a {\it horizontal lift} of the family of connections as in \cite[Theorem 6.2]{IS} 
and \cite[Section 9]{I2}.
Here the horizontal lift 
is a first order infinitesimal extension of the relative connection with an integrability condition.

Let $D$ be the effective divisor on $\mathbb{P}^1$ defined as 
\begin{equation*}
D= \sum_{i=1}^{\nu}n_i \cdot t_i+n_{\infty} \cdot \infty
\quad \text{and} \quad
n:=\deg(D) =\sum_{i=1}^{\nu}n_i +n_{\infty}.
\end{equation*}
Let $E$ be a rank $2$ vector bundle on $\mathbb{P}^1$.
Let $\nabla \colon E \rightarrow E \otimes \Omega_{\mathbb{P}^1}^1(D)$ 
be a connection on $E$ with the polar divisor $D$.
We call such pairs $(E,\nabla)$ {\it connections}.
Remark that 
we can shift the degree of the vector bundle $E$ by arbitrary integers
by applying some natural operations
(twisting by rank one meromorphic connections and 
birational bundle modifications, 
called {\it canonical transformations}, {\it elementary transformations},
or {\it Hecke modifications}).
In this paper, we assume that the degree of the vector bundle is $1$.
If $H^1(\mathbb{P}^1, E^{\vee})=0$, then we have 
$E\cong \mathcal{O}_{\mathbb{P}^1}\oplus \mathcal{O}_{\mathbb{P}^1}(1)$.
Here $E^{\vee}$ is the dual of the vector bundle $E$.
If there exists a family of connections of degree $1$, 
then there exists a Zariski open subset of the parameter space
such that this open subset parametrizes connections with 
the bundle type $\mathcal{O}_{\mathbb{P}^1}\oplus \mathcal{O}_{\mathbb{P}^1}(1)$.
For connections with the bundle type 
$\mathcal{O}_{\mathbb{P}^1}\oplus \mathcal{O}_{\mathbb{P}^1}(1)$,
we can define {\it apparent singularities} of the connections 
and also we can define dual parameters for the apparent singularities.
By the apparent singularities and dual parameters,
we may give a map from a moduli space of connections
with 
the bundle type $\mathcal{O}_{\mathbb{P}^1}\oplus \mathcal{O}_{\mathbb{P}^1}(1)$ 
to $\mathrm{Sym}^{(n-3)}(\mathbb{C}^2)$.
On the other hand,
Diarra--Loray \cite[Section 6]{DF} gave global normal forms of 
the connections with 
the bundle type $\mathcal{O}_{\mathbb{P}^1}\oplus \mathcal{O}_{\mathbb{P}^1}(n-2)$,
whose connection matrices are companion matrices.
By this normal form, we may construct a family of connections with 
bundle type $\mathcal{O}_{\mathbb{P}^1}\oplus \mathcal{O}_{\mathbb{P}^1}(1)$
parametrized by a Zariski open subset of 
$\mathrm{Sym}^{(n-3)}(\mathbb{C}^2)$.
By this family, we have a map from
the Zariski open subset of 
$\mathrm{Sym}^{(n-3)}(\mathbb{C}^2)$ to the moduli space of connections.
Finally, we have a birational correspondence 
between the moduli space of connections 
and $\mathrm{Sym}^{(n-3)}(\mathbb{C}^2)$.
Our point of view is that this birational correspondence 
gives coordinates on (a Zariski open set of) the moduli space of connections.
In this paper, we consider the generalized isomonodromic deformations 
(integrable deformations) of connections. 
We may regard the generalized isomonodromic deformations
as vector fields on the moduli space of connections.
The main purpose of this paper is giving explicit description of 
generalized isomonodromic deformations 
by using these coordinates.
Here the eigenvalues of the leading coefficients 
of the Laurent expansions of the connections
at each irregular singular point 
 are not necessarily distinct.
(If any leading coefficients have distinct eigenvalues respectively,
then the generalized isomonodromic deformations of 
this family of connections corresponds to the Jimbo--Miwa--Ueno equations).
That is,
we will consider not only \textit{unramified irregular singular points}
 (Definition \ref{2021.11.12.15.16_2} below)
but also 
\textit{ramified irregular singular points} (Definition \ref{2021.11.12.15.16} below).

There exist many studies on Hamiltonians 
of the Jimbo--Miwa--Ueno equation 
(\cite{Fedorov}, \cite{Hurtubise}, \cite{Woodhouse-1}, \cite{Woodhouse-2}, 
\cite{Wong-1},
and \cite{Yamakawa-1}).
The main subject of this paper is to give explicit descriptions of 
the symplectic structure and 
the Hamiltonians of the generalized isomonodromic deformations 
by using apparent singularities.
For the regular singular isomonodromic deformations,
Dubrovin--Mazzocco \cite{DM}
have introduced 
isomonodromic Darboux coordinates on the moduli space of Fuchsian systems,
which are connections on the trivial bundle over $\mathbb{P}^1$. 
They have described the isomonodromic deformations of Fuchsian systems 
as Hamiltonian systems by using the isomonodromic Darboux coordinates.
Roughly speaking, 
we extend their argument for the regular singular case
to the irregular singular (rank $2$) case.
In our calculation, 
Krichever's formula of the symplectic form \cite[Section 5]{Krichever}
is used as in \cite{DM}.
On the other hand, Kimura \cite{Kimura} has studied 
the degeneration of the two-dimensional Garnier systems.
By the confluence procedure, 
Hamiltonian systems of (generalized) isomonodromic deformations 
of certain rank $2$ liner differential equations are described explicitly.
We try to compare our Hamiltonian systems and 
Kimura's Hamiltonian systems by an example.
Our Hamiltonian systems are not necessarily commuting Hamiltonian systems
in \cite[Definition 3.5]{DM}.

\subsection{Space of deformation parameters}

Now we describe {\it the space of deformation parameters}
for our generalized isomonodromic deformations.
Put $I:= \{ 1,2,\ldots,\nu,\infty\}$, $t_1:=0$, $t_2:=1$, and  
$t_{\infty}:= \infty \in \mathbb{P}^1$.
We take a decomposition $I= I_{\text{reg}} \cup I_{\text{un}} \cup I_{\text{ra}}$ 
such that $I_{\text{reg}}$, $I_{\text{un}}$, and $I_{\text{ra}}$ are disjoint each other.
We assume that $n_i=1 $ for $i\in I_{\text{reg}}$ and
$n_i >1$ for $ i \in I_{\text{un}} \cup I_{\text{ra}}$.
We set
\begin{equation*}
T_{\boldsymbol{t}}:= \left\{ (t_3,\ldots, t_{\nu})  \in \mathbb{C}^{\nu-2}
\ \middle| \  
\begin{array}{l}
t_i \neq t_j \ ( i\neq j),\text{ and}\\
t_i \notin \{ 0,1 \} \ (i=3,\ldots,\nu)  \\
\end{array}
\right\}.
\end{equation*}
Moreover, put 
\begin{equation*}
\begin{aligned}
&T_{\boldsymbol{\theta}}^{\text{res}}:= 
\left\{ \boldsymbol{\theta}_0 \in \mathbb{C}^{2(\nu+1)} 
\ \middle|\  
\begin{array}{l}
\sum_{i\in I_{\text{reg}}\cup I_{\text{un}}}(\theta^+_{n_i-1,t_i}+\theta^-_{n_i-1,t_i} )
+\sum_{i\in I_{\text{ra}}} ( \theta_{2n_i-2,t_i} - \frac{1}{2}) =-1 \\
\theta_{0,t_i}^{+}-\theta_{0,t_i}^{-} \notin \mathbb{Z} \text{ for } i  \in I_{\text{reg}}
\end{array}
\right\}
\text{ and } \\
&T_{\boldsymbol{\theta}}:= 
\left\{ \boldsymbol{\theta} \in 
\prod_{i \in I_{\text{un}}} \mathbb{C}^{2(n_i-1)}
\times 
\prod_{i \in I_{\text{ra}}} \mathbb{C}^{2n_i-2} \ \middle| \ 
\begin{array}{l}
\theta^+_{0,t_i}- \theta^-_{0,t_i} \neq 0 \text{ for $i \in I_{\mathrm{un}}$} \\
 \theta_{1,t_i} \neq 0 \text{ for $i \in I_{\mathrm{ra}}$}
\end{array}
 \right\}.
\end{aligned}
\end{equation*}
Here we set 
\begin{equation*}
\boldsymbol{\theta}_0:= ((\theta^+_{n_i-1,t_i},\theta^-_{n_i-1,t_i} )_{i\in I_{\text{reg}}\cup I_{\text{un}}},
(\theta_{2n_i-2,t_i} )_{i\in I_{\text{ra}}})
\end{equation*}
and we denote by
$\boldsymbol{\theta} = (\boldsymbol{\theta}_{\text{un}},\boldsymbol{\theta}_{\text{ra}})$
an element of $T_{\boldsymbol{\theta}}$, where
\begin{equation*}
\begin{aligned}
\boldsymbol{\theta}_{\text{un}}&=
 ((\theta^+_{0,t_i},\theta^-_{0,t_i} ),\ldots,
 (\theta^+_{n_i-2,t_i},\theta^-_{n_i-2,t_i} ))_{i \in I_{\text{un}}} \text{ and} \\
 \boldsymbol{\theta}_{\text{ra}}&= 
 (\theta_{0,t_i},\ldots, \theta_{2n_i-3,t_i} )_{i \in I_{\text{ra}}}.
\end{aligned}
 \end{equation*}
The relation in the definition of $T_{\boldsymbol{\theta}}^{\text{res}}$ 
is called the {\it Fuchs relation}.
Fix a tuple of complex numbers
$\boldsymbol{t}_{\text{ra}} = (t_i)_{i \in \{3,4,\ldots,\nu\} \cap I_{\text{ra}}}$,
where 
$t_i \neq t_j \ ( i\neq j)$ and
$t_i \notin \{ 0,1 \}$.
We denote by $(T_{\boldsymbol{t}})_{\boldsymbol{t}_{\text{ra}}}$
the fiber of $\boldsymbol{t}_{\text{ra}}$ under the projection
\begin{equation*}
T_{\boldsymbol{t}} \longrightarrow
\prod_{i \in \{3,4,\ldots,\nu\} \cap I_{\text{ra}}} \mathbb{C}.
\end{equation*}

\begin{Def}
We define
{\it the space of deformation parameters} as
$(T_{\boldsymbol{t}})_{\boldsymbol{t}_{\text{ra}}} \times T_{\boldsymbol{\theta}}$.
\end{Def} 
 
Hence we consider the positions of the points $t_i$ for 
$i \in I_{\mathrm{reg}} \cup I_{\mathrm{un}}$ as 
deformation parameters.
On the other hand,
we do not consider 
the positions of the points $t_i$ for $i \in I_{\mathrm{ra}}$ as 
deformation parameters,
since the integrable deformations whose deformation parameters 
are the positions of the ramified irregular points are more complicated.

\subsection{Symplectic fiber bundle}

Next we define an algebraic variety over the space of deformation parameters 
$(T_{\boldsymbol{t}})_{\boldsymbol{t}_{\text{ra}}} \times T_{\boldsymbol{\theta}}$
such that this algebraic variety parametrizes connections and
there exists a symplectic form on each fiber.
This algebraic variety is considered 
as {\it the phase space} of our generalized isomonodromic deformations.
We set
\begin{equation*}
\widehat{\mathcal{M}}_{\boldsymbol{t}_{\text{ra}}} :=
\left\{ 
\begin{array}{l}
(\{ (q_1,p_1),\ldots ,(q_{n-3},p_{n-3}) \} ,(t_3,\ldots,t_{\nu}) )\\
\in \mathrm{Sym}^{(n-3)}(\mathbb{C}^2) \times 
(T_{\boldsymbol{t}})_{\boldsymbol{t}_{\text{ra}}}
\end{array}
\ \middle|\ 
\begin{array}{l}
q_i \neq q_j \ ( i\neq j) \text{ and }\\
q_j \notin \{ 0,1,t_3,\ldots,t_{\nu},\infty\} \ (j=1,\ldots,n-3)
\end{array}
\right\}.
\end{equation*}
If we take a point $\boldsymbol{t}_0=(t_3,\ldots, t_{\nu}) 
$ of $(T_{\boldsymbol{t}})_{\boldsymbol{t}_{\text{ra}}}$,
we put
\begin{equation*}
\mathcal{M}_{\boldsymbol{t}_{0},\boldsymbol{t}_{\text{ra}}} :=
\left\{ 
\begin{array}{l}
\{ (q_1,p_1),\ldots ,(q_{n-3},p_{n-3}) \} \\
\in \mathrm{Sym}^{(n-3)}(\mathbb{C}^2)
\end{array}
\ \middle|\ 
\begin{array}{l}
q_i \neq q_j \ ( i\neq j) \text{ and }\\
q_j \notin \{0,1,\infty\} \cup \boldsymbol{t}_{0}  \ (j=1,\ldots,n-3)
\end{array}
\right\}.
\end{equation*}

\begin{Def}
We define a symplectic fiber bundle
 $\pi_{\boldsymbol{t}_{\mathrm{ra}}, \boldsymbol{\theta}_0}$ as
 the natural projection
\begin{equation}\label{2020.1.20.21.15}
\pi_{\boldsymbol{t}_{\mathrm{ra}}, \boldsymbol{\theta}_0} \colon 
\widehat{\mathcal{M}}_{\boldsymbol{t}_{\text{ra}}}\times T_{\boldsymbol{\theta}} \longrightarrow
(T_{\boldsymbol{t}})_{\boldsymbol{t}_{\text{ra}}} \times T_{\boldsymbol{\theta}}.
\end{equation}
Here the symplectic structure on the fiber 
$\mathcal{M}_{\boldsymbol{t}_{0},\boldsymbol{t}_{\text{ra}}} \times \{ \boldsymbol{\theta}\}$
of 
$(\boldsymbol{t}_0, \boldsymbol{\theta})\in 
(T_{\boldsymbol{t}})_{\boldsymbol{t}_{\text{ra}}} \times T_{\boldsymbol{\theta}}$ 
is defined by
\begin{equation}\label{2020.2.7.11.24}
\sum_{j=1}^{n-3}d \left(\frac{p_j}{\prod_{i=1}^{\nu} (q_j-t_i)^{n_i}}  \right)\wedge dq_j.
\end{equation}
\end{Def}

Set $D(\boldsymbol{t}_0):= n_0 \cdot 0 + n_1\cdot 1 +\sum_{i=3}^{\nu} n_i \cdot t_i +
n_{\infty} \cdot \infty$
for $\boldsymbol{t}_0=(t_3,\ldots, t_{\nu}) 
\in (T_{\boldsymbol{t}})_{\boldsymbol{t}_{\text{ra}}}$.
We may check that
the fiber 
$\mathcal{M}_{\boldsymbol{t}_{0},\boldsymbol{t}_{\text{ra}}} \times \{ \boldsymbol{\theta}\}$
is isomorphic to the moduli space $\mathfrak{Conn}_{(\boldsymbol{t}_0,\boldsymbol{\theta},\boldsymbol{\theta}_0)}$.
Here $\mathfrak{Conn}_{(\boldsymbol{t}_0,\boldsymbol{\theta},\boldsymbol{\theta}_0)}$
is the moduli space of $(\boldsymbol{\theta},\boldsymbol{\theta}_0)$-connections on 
$\mathcal{O}_{\mathbb{P}^1}\oplus \mathcal{O}_{\mathbb{P}^1}(1)$
such that the polar divisors of connections are $D(\boldsymbol{t}_0)$ and
connections satisfy some generic conditions
(see \eqref{2021.4.29.11.42} below).
The correspondence between 
$\mathcal{M}_{\boldsymbol{t}_{0},\boldsymbol{t}_{\text{ra}}} \times \{ \boldsymbol{\theta}\}$ 
and 
$\mathfrak{Conn}_{(\boldsymbol{t}_0,\boldsymbol{\theta},\boldsymbol{\theta}_0)}$
is given by the theory of apparent singularities (Section \ref{2021.4.29.11.57} below)
and construction of a family of connections parametrized by 
$\mathcal{M}_{\boldsymbol{t}_{0},\boldsymbol{t}_{\text{ra}}} \times \{ \boldsymbol{\theta}\}$
(Section \ref{2021.4.29.11.58} below).
For the construction of a family, we will use
Diarra--Loray's global normal form.
We will call $\widehat{\mathcal{M}}_{\boldsymbol{t}_{\text{ra}}}\times T_{\boldsymbol{\theta}}$
an {\it extended moduli space of connections}.

\subsection{Main results}

For the vector fields 
$\partial / \partial \theta_{l,t_i}^{\pm}$ ($i \in I_{\mathrm{un}}$, $l=0,1,\ldots,n_i-2$),  
$\partial / \partial \theta_{l',t_i}$ ($i \in I_{\mathrm{ra}}$, $l'=0,1,\ldots,2n_i-3$), and
$\partial / \partial t_i$ ($i\in \{3,4\ldots,\nu\} \cap (I_{\mathrm{reg}} \cup I_{\mathrm{un}})$),
we define the vector fields 
$\delta^{\mathrm{IMD}}_{\theta_{l,t_i}^{\pm}}$,
$\delta^{\mathrm{IMD}}_{\theta_{l',t_i}}$, and
$\delta^{\mathrm{IMD}}_{t_i}$
on $\widehat{\mathcal{M}}_{\boldsymbol{t}_{\text{ra}}}\times T_{\boldsymbol{\theta}}$
by the integrable deformations of the family of connections 
parametrized by $\widehat{\mathcal{M}}_{\boldsymbol{t}_{\text{ra}}}\times T_{\boldsymbol{\theta}}$.
(in Section \ref{2020.1.24.22.25}, Section \ref{2020.2.7.11.34}, and Section \ref{2020.1.25.8.11}). 
We define a $2$-form $\hat{\omega}$
on $\widehat{\mathcal{M}}_{\boldsymbol{t}_{\text{ra}}}\times T_{\boldsymbol{\theta}}$
such that 
the restriction of $\hat{\omega}$ to each fiber of 
$\pi_{\boldsymbol{t}_{\mathrm{ra}}, \boldsymbol{\theta}_0}$
coincides with the symplectic form \eqref{2020.2.7.11.24} and 
the interior products with the vector fields 
determined by the integrable deformations vanish:
\begin{equation*}
\iota(\delta^{\mathrm{IMD}}_{\theta_{l,t_i}^{\pm}})\hat{\omega}=
\iota(\delta^{\mathrm{IMD}}_{\theta_{l',t_i}})\hat{\omega}=
\iota(\delta^{\mathrm{IMD}}_{t_i})\hat{\omega}=0.
\end{equation*}
We call the $2$-form $\hat{\omega}$ the \textit{isomonodromy $2$-form} 
as in \cite{Yamakawa-1}.
The main result of this paper is explicit description of 
the isomonodromy $2$-form using apparent singularities.
Our description of $\hat{\omega}$ is as follows:
\begin{equation*}
\begin{aligned}
\hat{\omega} =&\ \sum_{j=1}^{n-3} 
d \left(\frac{p_j}{\prod_{i=1}^{\nu} (q_j-t_i)^{n_i}} 
 - \sum_{i=1}^{\nu} \frac{D_i(q_j;\boldsymbol{t}, \boldsymbol{\theta})}{(q_j -t_i)^{n_i}} 
 - D_{\infty}(q_j;\boldsymbol{t}, \boldsymbol{\theta})\right) \wedge dq_j \\
&+ \sum_{i\in I_{\mathrm{un}}} 
\sum_{l=0}^{n_{i}-2} \left( d  H_{\theta_{l,t_i}^{+}} \wedge d \theta_{l,t_i}^{+}
+ d H_{\theta_{l,t_i}^{-}} \wedge d \theta_{l,t_i}^{-} \right) \\
&+ \sum_{i \in I_{\mathrm{ra}}}\sum_{l'=0}^{2n_i-3} d H_{\theta_{l' , t_i}} \wedge d\theta_{l' , t_i}
+\sum_{i\in \{3,4,\ldots,\nu \}\cap (I_{\mathrm{reg}} \cup I_{\mathrm{un}} )} d H_{t_i} \wedge dt_i\\
&+ \text{[ a section of $\pi_{\boldsymbol{t}_{\mathrm{ra}}, \boldsymbol{\theta}_0}^* 
(\Omega^2_{(T_{\boldsymbol{t}})_{\boldsymbol{t}_{\text{ra}}} \times T_{\boldsymbol{\theta}}} )$ ]}
\end{aligned}
\end{equation*}
(Theorem \ref{2020.1.21.13.48} and Theorem \ref{2020.1.29.11.46}).
Here $D_i(q_j;\boldsymbol{t}, \boldsymbol{\theta})$ where $i\in I$ are
defined in Lemma \ref{2019.12.30.22.00}
as $D_i$.
The Hamiltonians $H_{\theta_{l,t_i}^{\pm}}$, $H_{\theta_{l' , t_i}}$, and $H_{t_i}$
are defined in Definition \ref{2020.2.7.11.59},
Definition \ref{2020.2.7.12.01}, and
Definition \ref{2020.2.7.12.00}, respectively.
Roughly speaking, the Hamiltonians $H_{\theta_{l,t_i}^{\pm}}$ and $H_{\theta_{l' , t_i}}$ 
appear in 
the holomorphic parts of the diagonalizations of connections at each singular point $t_i$.
By this description of the isomonodromy $2$-form,
we obtain hamiltonian descriptions of the vector fields determined by 
the integrable deformations: 
\begin{equation*}
\begin{aligned}
\delta_{\theta_{l,t_i}^{\pm}}^{\mathrm{IMD}}&= 
\frac{\partial}{\partial \theta_{l,t_i}^{\pm}}
-\sum_{j=1}^{n-3} \left( 
 \frac{\partial H_{\theta_{l,t_i}^{\pm}}}{\partial \eta_j} 
\frac{\partial}{\partial q_j} 
- \frac{\partial H_{\theta_{l,t_i}^{\pm}}}{\partial q_j} 
\frac{\partial}{\partial \eta_j}
\right), \\
\delta_{\theta_{l',t_i}}^{\mathrm{IMD}}&= 
\frac{\partial}{\partial \theta_{l',t_i}}
-\sum_{j=1}^{n-3} \left( 
 \frac{\partial H_{\theta_{l',t_i}}}{\partial \eta_j} 
\frac{\partial}{\partial q_j} 
- \frac{\partial H_{\theta_{l',t_i}}}{\partial q_j} 
\frac{\partial}{\partial \eta_j}
\right), \text{ and}\\
\delta_{t_i}^{\mathrm{IMD}}&= 
\frac{\partial}{\partial t_i}
-\sum_{j=1}^{n-3} \left( 
 \frac{\partial  H_{t_i}}{\partial \eta_j} 
\frac{\partial}{\partial q_j} 
- \frac{\partial H_{t_i}}{\partial q_j} 
\frac{\partial}{\partial \eta_j}
\right)
\end{aligned}
\end{equation*}
(Corollary \ref{2020.1.21.15.11} and Corollay \ref{2020.2.2.19.01}).
Here we put $\eta_j:= \frac{p_j}{\prod_{i=1}^{\nu} (q_j-t_i)^{n_i}}
 - \sum_{i=1}^{\nu} \frac{D_i(q_j;\boldsymbol{t}, \boldsymbol{\theta})}{(q_j -t_i)^{n_i}} 
 - D_{\infty}(q_j;\boldsymbol{t}, \boldsymbol{\theta})$.

The organization of this paper is as follows.
In Section 2, we recall the definition of the apparent singularities 
and Diarra--Loray's global normal form. 
In Section 3, we consider the integrable deformations of connections
which have only regular singularities and unramified irregular singularities.
We define a $2$-form on the fiber
$\mathcal{M}_{\boldsymbol{t}_{0},\boldsymbol{t}_{\text{ra}}}$
by Krichever's formula \cite[Section 5]{Krichever}.
We show that this $2$-form coincides with the symplectic form \eqref{2020.2.7.11.24}.
Also by Krichever's formula,
we define a $2$-form on $\widehat{\mathcal{M}}_{\boldsymbol{t}_{\text{ra}}}\times T_{\boldsymbol{\theta}}$.
We show that this $2$-form is the isomonodromy $2$-form.
By calculation of this $2$-form on 
$\widehat{\mathcal{M}}_{\boldsymbol{t}_{\text{ra}}}\times T_{\boldsymbol{\theta}}$
by using Diarra--Loray's global normal form,
we have an explicit formula of this $2$-form.
In Section 4, we extend the argument of Section 3 to the integrable deformations of connections
which have ramified irregular singularities.
In Section 5, we consider two examples.
The first example is the case where $D=2\cdot [0] + 2 \cdot [1]+2 \cdot [\infty]$.
We assume that $0,1,\infty \in \mathbb{P}^1$ are unramified irregular singular points.
The dimension of the space of deformation parameters is $6$
and the dimension of the fiber $\widehat{\mathcal{M}}_{\boldsymbol{t}_{0},\boldsymbol{t}_{\text{ra}}}$ 
is $6$.
The second example is the case where $D=5 \cdot [\infty]$.
We assume that $\infty \in \mathbb{P}^1$ is ramified irregular singular point.
This example corresponds to Kimura's $H(9/2)$ in \cite{Kimura}.
We consider the family of connections corresponding to Kimura's family $L(9/2;2)$.
We reproduce the Hamiltonian system $H(9/2)$.

\section{Normal forms for rank two linear irregular differential equations}\label{2021.4.29.12.46}

In the first part of this section, we will give a correspondence 
between the moduli space of connections 
and $\mathcal{M}_{\boldsymbol{t}_{0},\boldsymbol{t}_{\text{ra}}} 
\subset \mathrm{Sym}^{(n-3)}(\mathbb{C}^2)$.
First, we recall the theory of apparent singularities in Section \ref{2021.4.29.11.57}.
This theory gives a map from
the moduli space of connections 
to $\mathcal{M}_{\boldsymbol{t}_{0},\boldsymbol{t}_{\text{ra}}}$.
Second, we recall Diarra--Loray's global normal form in Section \ref{2021.4.29.11.58}.
This normal form
gives a map from $\mathcal{M}_{\boldsymbol{t}_{0},\boldsymbol{t}_{\text{ra}}}$
to the moduli space of connections.
In the second part of this section (Section \ref{2021.4.29.12.39} and Section \ref{2021.4.29.12.40}),
we will consider infinitesimal deformations of connections 
and define horizontal lifts of connections.
If we construct a horizontal lift of a connection,
then we have an integrable deformation of a connection.
After Section \ref{2021.4.29.12.46}, we will discuss construction of horizontal lifts.
In the third part of this section
(Section \ref{2021.4.29.12.41}), we will discuss local solutions of 
the differential equations with respect to the connections
at the apparent singularities.
We will use these solutions for the definition of the $2$-form $\omega$ 
on $\mathcal{M}_{\boldsymbol{t}_{0},\boldsymbol{t}_{\text{ra}}}$ (in Section \ref{2020.2.8.10.24} below). 

We take a natural affine open covering $\{ U_0 , U_{\infty}\}$ of $\mathbb{P}^1$.
Denote by $x$ a coordinate on $U_0$ and 
by $w$ a coordinate on $U_{\infty}$.
That is, $w=x^{-1}$ on $U_0 \cap U_{\infty}$.
Let $\infty$ be the point $w=0$ on $U_{\infty}$.
Set
$E_k := \mathcal{O}_{\mathbb{P}^1} \oplus \mathcal{O}_{\mathbb{P}^1}(k)$.
Here we define the vector bundle
 $E_k$
by two trivializations $\{ \varphi^{(k)}_{U_0} \colon E_k|_{U_0} 
\xrightarrow{\cong}\mathcal{O}_{U_0}^{\oplus 2} ,
\varphi^{(k)}_{U_\infty} \colon 
E_k|_{U_\infty} 
\xrightarrow{\cong}\mathcal{O}_{U_\infty}^{\oplus 2} \}$
such that
\begin{equation}\label{2022.4.8.9.42}
\xymatrix{
E_k|_{U_0\cap U_\infty }  \ar[d]^-= \ar[r]^-{\varphi^{(k)}_{U_\infty}}& 
\mathcal{O}_{U_0\cap U_\infty}^{\oplus 2} \ar[d]^-{G_k} \\
E_k|_{U_0\cap U_\infty} \ar[r]^-{\varphi^{(k)}_{U_0}  }&
 \mathcal{O}_{U_0\cap U_\infty}^{\oplus 2}, 
 }
\end{equation}
where $G_k =
 \begin{pmatrix}
 1 & 0 \\
 0 & x^{k}
 \end{pmatrix}$.
Fix a tuple of complex numbers
$\boldsymbol{t}_{\text{ra}} = (t_i)_{i \in \{3,4,\ldots,\nu\} \cap I_{\text{ra}}}$,
where 
$t_i \neq t_j \ ( i\neq j)$ and
$t_i \notin \{ 0,1 \}$.

\subsection{Apparent singularities}\label{2021.4.29.11.57}

Take $\boldsymbol{t}_0= (t_i)_{i \in \{3,4,\ldots,\nu\} }
 \in (T_{\boldsymbol{t}})_{\boldsymbol{t}_{\text{ra}}}$
and
set $D=n_1 \cdot 0+n_2 \cdot 1+ n_3 \cdot t_3 + \cdots +
n_\nu \cdot t_{\nu}+n_\infty \cdot \infty$.
For a connection $(E_1,\nabla\colon E_1 \rightarrow
E_1 \otimes \Omega^1_{\mathbb{P}^1}(D))$, 
we define the \textit{apparent singularities of $(E_1,\nabla)$} as follows. 
Consider the sequence of maps
\begin{equation*}
\xymatrix{
\mathcal{O}_{\mathbb{P}^1}(1) \ar[r]^-{\subset} & E_1 \ar[r]^-{\nabla} 
& E_1 \otimes \Omega^1_{\mathbb{P}^1}(D)
\ar[r]^-{\text{quotient}} & (E_1/\mathcal{O}_{\mathbb{P}^1}(1) )\otimes \Omega^1_{\mathbb{P}^1}(D) 
\cong \mathcal{O}_{\mathbb{P}^1}(n-2).
}
\end{equation*}
This composition is an $\mathcal{O}_{\mathbb{P}^1}$-morphism.
We denote by 
$\varphi_{\nabla}\colon \mathcal{O}_{\mathbb{P}^1}(1)
 \rightarrow \mathcal{O}_{\mathbb{P}^1}(n-2)$
this composition.  
We assume that
the subbundle $\mathcal{O}_{\mathbb{P}^1}(1) \subset E_1$ is not $\nabla$-invariant.
Then $\varphi_{\nabla}$ is not the zero morphism.
The $\mathcal{O}_{\mathbb{P}^1}$-morphism $\varphi_{\nabla}$ 
has $n - 3$ zeroes counted with multiplicity.

\begin{Def}
We define {\it apparent singularities of $(E_1,\nabla)$}
as $\mathrm{div}(\varphi_{\nabla}) \in |\mathcal{O}_{\mathbb{P}^1}(n-3)| 
\cong \mathrm{Sym}^{(n-3)}(\mathbb{P}^1)$.
\end{Def}

By the trivialization $\varphi^{(1)}_{U_0} \colon E_1|_{U_0} 
\xrightarrow{\cong}\mathcal{O}_{U_0}^{\oplus 2}$,
we have the following description:
$$
\nabla|_{U_0}=
d+ 
\begin{pmatrix}
A(x) & B(x) \\
C(x) & D(x)
\end{pmatrix}\frac{dx}{P(x)},
$$
where $P(x):=\prod_{i=1}^{\nu} (x-t_i)^{n_i}$ and 
$A,B,C,D$ are polynomials such that 
$\deg(A) \leq n-2, \deg(B) \leq n-3, \deg(C) \leq n-1,  \deg (D) \leq n-2$. 
Then the apparent singularities of $(E_1,\nabla)$ 
are zeros of the polynomial $B(x)$.

Assume that the apparent singularities of $(E_1,\nabla)$ consist of distinct points and 
all of them are distinct from the poles $t_1, \ldots , t_{\nu},\infty$ of the connection $\nabla$.
We can define a birational bundle transformation
\begin{equation*}
\phi_{\nabla} := \mathrm{id} \oplus \varphi_{\nabla} \colon
\mathcal{O}_{\mathbb{P}^1} \oplus \mathcal{O}_{\mathbb{P}^1}(1) 
\dashrightarrow \mathcal{O}_{\mathbb{P}^1} \oplus \mathcal{O}_{\mathbb{P}^1}(n-2)
\end{equation*}
and consider the pushed-forward connection $(\phi_{\nabla})_*\nabla$ on 
$\mathcal{O}_{\mathbb{P}^1}\oplus \mathcal{O}_{\mathbb{P}^1}(n-2)$.
Then we have a transformation of a connection with bundle type 
$\mathcal{O}_{\mathbb{P}^1}\oplus \mathcal{O}_{\mathbb{P}^1}(1)$:
\begin{equation}\label{2020.1.11.22.49}
(E_1, \nabla) \longmapsto 
(E_{n-2}, (\phi_{\nabla})_*\nabla).
\end{equation}
The connection $(\phi_{\nabla})_*\nabla$ has simple poles $q_1,\ldots,q_{n-3}$
with residual eigenvalues $0$ and $-1$ at each pole.
Let $D_{\text{App}}$ be the effective divisor $q_1+\cdots+q_{n-3}$.
We may decompose $(\phi_{\nabla})_*\nabla$ as follows:
\begin{equation}\label{2021.3.27.18.48}
\begin{pmatrix}
\nabla_{11} & \Phi_{12} \\
\Phi_{21} & \nabla_{22}
\end{pmatrix},
\end{equation}
where $\nabla_{11}\colon \mathcal{O}_{\mathbb{P}^1} \rightarrow 
\mathcal{O}_{\mathbb{P}^1} \otimes \Omega^1_{\mathbb{P}^1}(D+D_{\text{App}})$
and 
$\nabla_{22}\colon \mathcal{O}_{\mathbb{P}^1}(n-2) \rightarrow 
\mathcal{O}_{\mathbb{P}^1} (n-2)\otimes \Omega^1_{\mathbb{P}^1}(D+D_{\text{App}})$
are connections.
Moreover $\Phi_{12}\colon \mathcal{O}_{\mathbb{P}^1} (n-2)\rightarrow 
\mathcal{O}_{\mathbb{P}^1} \otimes \Omega^1_{\mathbb{P}^1}(D+D_{\text{App}})$
and 
$\Phi_{21}\colon \mathcal{O}_{\mathbb{P}^1} \rightarrow 
\mathcal{O}_{\mathbb{P}^1} (n-2)\otimes \Omega^1_{\mathbb{P}^1}(D+D_{\text{App}})$
are $\mathcal{O}_{\mathbb{P}^1}$-morphisms.
Since the birational bundle transformation $\phi_{\nabla}$
is given by 
$ \begin{pmatrix}
 1 & 0 \\
 0 & B(x)
 \end{pmatrix}$,
 the connection \eqref{2021.3.27.18.48} has the following description:
$$
\left.
\begin{pmatrix}
\nabla_{11} & \Phi_{12} \\
\Phi_{21} & \nabla_{22}
\end{pmatrix} \right|_{U_0} 
= 
\begin{pmatrix}
d + \frac{A(x)dx}{P(x)} & \frac{dx}{P(x)} \\
\frac{C(x)B(x)dx}{P(x)} & d+ \frac{D(x)dx}{P(x)} - \sum_{j=1}^{n-3} \frac{dx}{x-q_j}
\end{pmatrix}.
$$
By an automorphism of the bundle $E_{n-2}$,
we may normalize the connection $(\phi_{\nabla})_*\nabla$ 
so that the normalized connection has the following conditions
(for details, see \cite[Proposition 3]{DF}):
\begin{itemize}
\item the connection $\nabla_{11}$ is the trivial connection; and
\item the $\mathcal{O}_{\mathbb{P}^1}$-morphism $\Phi_{12}$
corresponds to the section
\begin{equation}\label{2021.3.31.15.26}
\frac{dx}{\prod_{i=1}^{\nu}(x-t_i)^{n_i}}  \quad \text{(on $U_0$) \quad and \quad }
\frac{-w^{-n+2}dw}{w^2\prod_{i=1}^{\nu}(1/w-t_i)^{n_i}}  \quad \text{(on $U_{\infty}$)}
\end{equation}
under the isomorphism 
$\mathrm{Hom}_{\mathcal{O}_{\mathbb{P}^1}}( \mathcal{O}_{\mathbb{P}^1} 
 (n-2),
\mathcal{O}_{\mathbb{P}^1}\otimes \Omega^1_{\mathbb{P}^1}(D+D_{\text{App}}))
\cong H^0( \mathcal{O}_{\mathbb{P}^1} (-n+2)
\otimes \Omega^1_{\mathbb{P}^1}(D+D_{\text{App}}))$.
\end{itemize}
Automorphisms of $E_{n-2}$
preserving these conditions of $(\phi_{\nabla})_*\nabla$
are just scalars (see \cite[Section 3]{DF}).
For each $j=1,2,\ldots,n-3$,
the $0$-eigendirection of
the residue matrix of
the normalized connection at $q_j$ 
corresponds to a point 
$\boldsymbol{p}_j \in \mathbb{P}(E_{n-2})|_{q_j} \cong \mathbb{P}^1$.
Here this identification is given by the trivialization $\varphi^{(n-2)}_{U_0}$.

\begin{Def} 
Since the $(-1)$-eigendirection is contained in the second factor of 
$E_{n-2} = \mathcal{O}_{\mathbb{P}^1}\oplus \mathcal{O}_{\mathbb{P}^1}(n-2)$,
we have $p_j \in \mathbb{C}$ such that $\boldsymbol{p}_j =[1:p_j] $.
We call $p_j$ a {\it dual parameter with respect to an apparent singularity $q_j$}.
\end{Def}

\subsection{Global normal form for rank two linear irregular differential equations}\label{2021.4.29.11.58}

In the previous section, we assign a point on 
$\mathrm{Sym}^{(n-3)}(\mathbb{C}^2)$
to a connection $E_1 \rightarrow E_1 \otimes \Omega^1_{\mathbb{P}^1} (D)$.
Conversely, take a point 
$\{ (q_1,p_1),\ldots , (q_{n-3},p_{n-3}) \} $
on $\mathcal{M}_{\boldsymbol{t}_0 ,\boldsymbol{t}_{\text{ra}}}
\subset \mathrm{Sym}^{(n-3)}(\mathbb{C}^2)$.
Then we may construct a connection 
$E_1 \rightarrow E_1 \otimes \Omega^1_{\mathbb{P}^1} (D)$
such that $q_1+\cdots+q_{n-3}$ is the apparent singularities 
and $p_j$ is the dual parameter of $q_j$, $j=1,2,\ldots,n-3$
(Proposition \ref{2021.4.11.16.40} below).
Now we discuss this construction.
We define an effective divisor $D_{\text{App}}$ on $\mathbb{P}^1$ as
$D_{\text{App}}=q_1+\cdots+q_{n-3}$.

\begin{Def}\label{2021.4.11.20.42}
For the point on $\mathcal{M}_{\boldsymbol{t}_0 ,\boldsymbol{t}_{\text{ra}}}$,
let 
$$
\nabla^{(n-2)}_{\mathrm{DL}}\colon E_{n-2} \longrightarrow 
E_{n-2} \otimes \Omega_{\mathbb{P}^1}^1 (D+D_{\mathrm{App}})
$$
be a connection with the following connection matrix on $U_0$:
\begin{equation}\label{2021.4.18.13.06}
\Omega^{(n-2)}=
\begin{pmatrix}
0& \frac{1}{P(x)} \\
c_0(x) & d_0(x)
\end{pmatrix}dx.
\end{equation}
Here we put
$P(x):=\prod_{i=1}^{\nu} (x-t_i)^{n_i}$,
\begin{equation}\label{2020.12.26.21.03}
\begin{aligned}
c_0(x)&:=\sum_{i=1}^{\nu} \frac{C_i(x)}{(x-t_i)^{n_i}}
+\sum_{j=1}^{n-3}\frac{p_j}{x-q_j} +\tilde{C}(x) +x^{n-3}C_{\infty} (x)   
, \text{ and} \\
d_0(x)&:=\sum_{i=1}^{\nu} \frac{D_i(x)}{(x-t_i)^{n_i}}
+\sum_{j=1}^{n-3}\frac{-1}{x-q_j} +D_{\infty} (x),
\end{aligned}
\end{equation}
where $C_i$, $D_i$, $(i=1,\ldots,\nu)$, $C_{\infty}$, $D_{\infty}$, and $\tilde{C}$ are 
polynomials in $x$ such that 
\begin{itemize}
\item $\deg (C_i),\deg (D_i) \le n_i - 1$ for $i=1,...,\nu,$
\item $\deg (C_\infty) \le n_{\infty}-1$, $\deg (D_\infty) \le n_\infty - 2$,
\item $\deg (\tilde{C})\le n-4$.
\end{itemize}
We assume that $q_1,\ldots,q_{n-3}$ are apparent singularities,
that is, the elementary transformation 
$$
\left((\tilde{\Phi}_{q_j})^{-1}d\tilde{\Phi}_{q_j}
+(\tilde{\Phi}_{q_j})^{-1}\Omega^{(n-2)} \tilde{\Phi}_{q_j}\right),
\quad \text{where }\tilde{\Phi}_{q_j}:= 
\begin{pmatrix}
1 & 0 \\
p_j & x-q_j
\end{pmatrix},
$$
of $\Omega^{(n-2)}$ by $\tilde{\Phi}_{q_j}$
has no pole at $q_j$.
We call such a connection $\nabla^{(n-2)}_{\mathrm{DL}}$
{\it Diarra--Loray's global normal form}.
\end{Def}

The corresponding
connection matrix $\Omega^{(n-2)}_{\infty}$ on $U_{\infty}$ 
of $\nabla^{(n-2)}_{\mathrm{DL}}$
is
$$
\begin{aligned}
\Omega^{(n-2)}_{\infty}&=
G_{n-2}^{-1} dG_{n-2}+G_{n-2}^{-1} \Omega^{(n-2)} G_{n-2} \\
&= \begin{pmatrix}
0& 0 \\
0& -n+2
\end{pmatrix}\frac{dw}{w}+
\begin{pmatrix}
0& \frac{1}{w^{n-2}P(1/w)} \\
w^{n-2}c_0(1/w) & d_0(1/w)
\end{pmatrix}\frac{-dw}{w^2}.
\end{aligned}
$$
We may check that
$\Omega^{(n-2)}_{\infty}$ has a pole of order $n_{\infty}$ at $\infty$.
We decompose the connection $\nabla^{(n-2)}_{\mathrm{DL}}$
as in \eqref{2021.3.27.18.48}.
Since the $(1,1)$-entry of $\Omega^{(n-2)}$ is zero 
and the $(1,2)$-entry of $\Omega^{(n-2)}$ is $\frac{dx}{P(x)}$,
the connection $\nabla_{11}$ is the trivial connection 
and the $\mathcal{O}_{\mathbb{P}^1}$-morphism $\Phi_{12}$
corresponds to the section \eqref{2021.3.31.15.26}.
The vector $(1,p_j)$ is
a $0$-eigenvector of the residue matrix of $\Omega^{(n-2)}$ at $q_j$.

 Now we consider a transformation of
 the connection $\nabla_{\mathrm{DL}}^{(n-2)}$ on $E_{n-2}$
 into a connection on $E_1$.
Set $Q_1 (x)=\prod^{n-3}_{j=1} (x -q_j)$. 
Let $Q_2(x)$ 
be the unique polynomial of degree $n-4$ such that 
$Q_2(q_j) = p_j$ for $j=1,2,\ldots,n-3$.
Set 
\begin{equation}\label{2021.4.2.20.06}
\begin{aligned}
&\tilde{G} := \begin{pmatrix} 1 & 0 \\ Q_2(x) & Q_1(x) \end{pmatrix}
\colon  \mathcal{O}_{U_0}^{\oplus 2} \dashrightarrow 
 \mathcal{O}_{U_0}^{\oplus 2} \quad 
\text{ and }\\
&\tilde{G}_{\infty} := G_{n-2}^{-1} \begin{pmatrix} 1 & 0 \\ Q_2(1/w) & Q_1(1/w) \end{pmatrix} G_{1}
\colon  \mathcal{O}_{U_\infty}^{\oplus 2} \dashrightarrow 
 \mathcal{O}_{U_\infty}^{\oplus 2}.
\end{aligned}
\end{equation}
Let $\Omega^{(1)}$ be the transformation of $\Omega^{(n-2)}$ by $\tilde{G}$:
\begin{equation}\label{2021.4.17.10.00}
\begin{aligned}
\Omega^{(1)}&=
\tilde{G}^{-1} d\tilde{G}+\tilde{G}^{-1} \Omega^{(n-2)} \tilde{G} \\
&=
\begin{pmatrix}\frac{Q_2(x)}{P(x)}&\frac{Q_1(x)}{P(x)}\\ 
\frac{c_0(x)+Q_2 (x)d_0(x)+(Q_2(x))'}{Q_1(x)} 
- \frac{(Q_2(x))^2}{P(x) Q_1(x)}&d_0(x)+\frac{(Q_1(x))'}{Q_1(x)} 
- \frac{Q_2(x)}{P(x)}\end{pmatrix} dx.
\end{aligned}
\end{equation}

\begin{Prop}\label{2021.4.11.16.40}
{\it We set
$$
\nabla_{\mathrm{DL}}^{(1)} := 
\begin{cases}
 d+ \Omega^{(1)} & \text{on $U_{0}$} \\
 d+ G_1^{-1} dG_1 + G_1^{-1}
\Omega^{(1)} \, G_1 & \text{on $U_{\infty}$}.
\end{cases}
$$
Then $\nabla_{\mathrm{DL}}^{(1)}$ is a 
connection
$$
\mathcal{O}_{\mathbb{P}^1} \oplus \mathcal{O}_{\mathbb{P}^1} (1)
\longrightarrow (\mathcal{O}_{\mathbb{P}^1} \oplus \mathcal{O}_{\mathbb{P}^1} (1)) \otimes 
\Omega_{\mathbb{P}^1}^1 (D).$$
That is, the pole divisor of $\nabla_{\mathrm{DL}}^{(1)} $ 
is $D$.
Moreover, the apparent singularities of $\nabla_{\mathrm{DL}}^{(1)}$
is $q_1+\cdots+q_{n-3}$ and 
the dual parameter with respect to $q_j$ is $p_j$.
}
\end{Prop}

\begin{proof}
First we will show that $\nabla_{\mathrm{DL}}^{(1)}$ has no pole at 
$q_1,\ldots,q_{n-3}$ by induction. 
Let $s \in \{1,\ldots,n-3\}$.
We define $Q^{(s)}_1 (x)=\prod^{s}_{j=1} (x -q_j)$ and $Q^{(s)}_2(x)$ 
is the unique polynomial of degree $s-1$ such that 
$Q^{(s)}_2(q_j) = p_j$ for $j=1,2,\ldots,s$.
Set 
$$
\tilde{G}^{(s)} := \begin{pmatrix} 1 & 0 \\ Q^{(s)}_2(x) & Q^{(s)}_1(x) \end{pmatrix}.
$$
Assume that $d+ (\tilde{G}^{(s)})^{-1} d \tilde{G}^{(s)} + 
(\tilde{G}^{(s)})^{-1}\Omega^{(n-2)} \, \tilde{G}^{(s)}$
has no pole at $q_1,q_2,\ldots,q_{s}$.
We will show that 
$d+ (\tilde{G}^{(s+1)})^{-1} d \tilde{G}^{(s+1)} + 
(\tilde{G}^{(s+1)})^{-1}\Omega^{(n-2)} \, \tilde{G}^{(s+1)}$
has no pole at $q_1,q_2,\ldots,q_{s+1}$.
We may check the following equalities:
$$
\begin{aligned}
& (\tilde{G}^{(s)})^{-1} d \tilde{G}^{(s)} + 
(\tilde{G}^{(s)})^{-1}\Omega^{(n-2)} \, \tilde{G}^{(s)}\\
&=\begin{pmatrix}\frac{Q^{(s)}_2(x)}{P(x)}&\frac{Q^{(s)}_1(x)}{P(x)}\\ 
\frac{c_0(x)+Q^{(s)}_2 (x)d_0(x)+(Q^{(s)}_2(x))'}{Q^{(s)}_1(x)} 
- \frac{(Q^{(s)}_2(x))^2}{P(x) Q^{(s)}_1(x)}&d_0(x)+\frac{(Q_1^{(s)}(x))'}{Q^{(s)}_1(x)} 
- \frac{Q^{(s)}_2(x)}{P(x)}\end{pmatrix} dx\\
&=\begin{pmatrix}
0 & 0 \\
\frac{p_{s+1} -Q_2^{(s)} (q_{s+1})}{Q_1^{(s)}(q_{s+1})}  &-1 
\end{pmatrix}\frac{dx}{x-q_{s+1}} + [\text{ holomorphic parts }].
\end{aligned}
$$
Here the last equality is the expansion at $x=q_{s+1}$.
Since $q_{s+1}$ is an apparent singularity,
we can transform the connection $d+ (\tilde{G}^{(s)})^{-1} d \tilde{G}^{(s)} + 
(\tilde{G}^{(s)})^{-1}\Omega^{(n-2)} \, \tilde{G}^{(s)}$ 
into a connection which is holomorphic at $q_{s+1}$ by the matrix
$$
\begin{pmatrix}
1 & 0 \\
\frac{p_{s+1} -Q_2^{(s)} (q_{s+1})}{Q_1^{(s)}(q_{s+1})}  &x-q_{s+1} 
\end{pmatrix}.
$$
We may check that 
$$
\tilde{G}^{(s+1)}_{0}=
\tilde{G}^{(s)}_{0}
\begin{pmatrix}
1 & 0 \\
\frac{p_{s+1} -Q_2^{(s)} (q_{s+1})}{Q_1^{(s)}(q_{s+1})}  &x-q_{s+1} 
\end{pmatrix}.
$$
Then we have that $d+ (\tilde{G}_0^{(s+1)})^{-1} d \tilde{G}_0^{(s+1)} + 
(\tilde{G}_0^{(s+1)})^{-1}\Omega^{(n-2)} \, \tilde{G}_0^{(s+1)}$
has no pole at $q_1,q_2,\ldots,q_{s+1}$.
So $\nabla_{\mathrm{DL}}^{(1)}$ has no pole at 
$q_1,\ldots,q_{n-3}$ by induction. 
Since $\tilde{G}_{\infty}$ is holomorphic at $w=0$ 
and
the determinant of $\tilde{G}_{\infty}$ does not vanish at $w=0$,
$$
 d+ G_1^{-1} dG_1 + G_1^{-1}
\Omega^{(1)} \, G_1
= d+ \tilde{G}_{\infty}^{-1} d\tilde{G}_{\infty} + \tilde{G}_{\infty}^{-1}
\Omega^{(n-2)}_{\infty} \, \tilde{G}_{\infty},
$$
has a pole of order $n_{\infty}$ at $\infty$.
Then the polar divisor of $\nabla_{\mathrm{DL}}^{(1)}$ is $D$.

By \eqref{2021.4.17.10.00}, 
the $(1,2)$-term of $\Omega^{(1)}$ is
$\frac{Q_1(x)}{P(x)}$.
The apparent singularities of $\nabla_{\mathrm{DL}}^{(1)}$
are the zeros of $\frac{Q_1(x)}{P(x)}$.
Then 
the apparent singularities of $\nabla_{\mathrm{DL}}^{(1)}$
are $q_1+\cdots+q_{n-3}$,
since the birational bundle transformation $\phi_{\nabla} $ 
is given by 
$\begin{pmatrix}
1 & 0 \\
0 &Q_1 (x) 
\end{pmatrix}$. 
Moreover,
$$
\begin{pmatrix}
1 & 0 \\
Q_2(x) &Q_1(x) 
\end{pmatrix}=
\begin{pmatrix}
1 & 0 \\
Q_2(x) &1 
\end{pmatrix}
\begin{pmatrix}
1 & 0 \\
0 &Q_1 (x) 
\end{pmatrix}.
$$
Here $\begin{pmatrix}
1 & 0 \\
Q_2(x) &1 
\end{pmatrix}$ is an automorphism of $E_{n-2}$.
Since 
the vector $(1,p_j)$ is
a $0$-eigenvector of the residue matrix of $\Omega^{(n-2)}$ at $q_j$,
the dual parameter with respect to $q_j$ is $p_j$.
\end{proof}

\subsection{Local formal data}\label{2020.1.21.14.24}

We put $x_{t_i} = (x-t_i)$ for $i=1,\ldots,\nu$ and $x_{\infty}= w$.
Put $I:= \{ 1,2,\ldots,\nu,\infty\}$, $t_1:=0$, $t_2:=1$, and  
$t_{\infty}:= \infty \in \mathbb{P}^1$.
We take a decomposition $I= I_{\text{reg}} \cup I_{\text{un}} \cup I_{\text{ra}}$ 
such that $I_{\text{reg}}$, $I_{\text{un}}$, and $I_{\text{ra}}$ are disjoint each other.
We assume that $n_i =1$ for $i \in I_{\text{reg}}$ and 
$n_i >1$ for $i \in I_{\text{un}} \cup I_{\text{ra}}$.

Let $\nabla$ be a connection on $E_1$:
$$
\nabla \colon E_1 \longrightarrow E_1 \otimes \Omega_{\mathbb{P}^1}^1(D),
$$
where $D=\sum_{i \in I} n_i \cdot t_i$.
For each $i \in I$, we take an affine open subset $U_{i} \subset \mathbb{P}^1$ 
such that $t_i \in U_i$.
We take a trivialization 
$E_1|_{U_i}\cong \mathcal{O}_{U_i}^{\oplus 2}$
and choose the coordinate $x_{t_i}$ on $U_i$ 
such that the point $t_i$ is defined by $x_{t_i}=0$.
Let
$\Omega$ be the connection matrix
of $\nabla$ associated to this trivialization.
We may describe $\Omega$ as follows:
\begin{equation*}
\Omega =
\Omega_{t_i}(0) \frac{dx_{t_i}}{x_{t_i}^{n_i}} + [\text{ higher order terms }], \quad 
\Omega_{t_i}(0)\in \mathfrak{gl}(2,\mathbb{C})
\end{equation*}
for each $i\in I$.

\begin{Def}\label{2021.11.12.15.16_2}
We say $t_i$ is an {\it unramified irregular singular point of $\nabla$}
if $n_{i} >1$ and $\Omega_{t_i}(0)$ has distinct eigenvalues.
\end{Def}

Let $( \mathrm{Det}(E_1) ,\mathrm{Tr} (\nabla) )$ be 
the determinant bundle of $E_1$ with the induced connection $\mathrm{Tr} (\nabla)$,
that is, $\mathrm{Tr} (\nabla) =\nabla \wedge \mathrm{id} + \mathrm{id} \wedge \nabla$.
We consider the trivialization $\mathrm{Det}(E_1) |_{U_i} \cong \mathcal{O}_{U_i}$
induced by the trivialization of $E_1|_{U_i}$.
Let $\alpha' \in \mathcal{O}_{U_i} dx_{t_i}/x_{t_i}^{n_i}$ be the connection matrix of 
$\mathrm{Tr} (\nabla)|_{U_i}$ associated to this trivialization (if necessary, $U_i$ shrinks).
We consider the tensor product $(\mathcal{O}_{U_i}^{\oplus 2}, d+\Omega)
\otimes (\mathcal{O}_{U_i}, d -\frac{1}{2} \alpha')$.
Let $N(x_{t_i}) dx_{t_i}/x_{t_i}^{n_i}$ be the connection matrix of this tensor product.
Remark that $N(x_{t_i}) \in \mathrm{End} (\mathcal{O}_{U_i}^{\oplus 2})$.

\begin{Def}\label{2021.11.12.15.16}
We say $t_i$ is a {\it ramified irregular singular point of $\nabla$}
if $n_{i} >1$, $N(0)$ is a nonzero nilpotent matrix, and $N(x_{t_i})^2 \not\equiv 0$ (mod $x_{t_i}^2$).
\end{Def}

We assume that 
\begin{itemize}
\item the differences of the eigenvalues of $\Omega_{t_i}(0)$ are not integers 
for any $i \in I_{\text{reg}}$, 
\item $t_i$ are 
unramified irregular singular points for any $i\in I_{\text{un}}$, and
\item $t_i$ are 
ramified irregular singular points for any $i\in I_{\text{ra}}$.
\end{itemize}

\begin{Lem}[For example {\cite[Proposition 9 and Proposition 10]{DF}}]
\begin{itemize}
\item[(1)] {\it If $i \in I_{\mathrm{un}}$,
then there exists a matrix $M \in \mathrm{GL}_2(\mathbb{C}[[x_{t_i}]])$ such 
that 
$$
M^{-1} dM +M^{-1} \Omega M
=
\frac{
\begin{pmatrix}
\theta_{0,t_i}^+ & 0 \\
0 & \theta_{0,t_i}^-
\end{pmatrix}
}{x_{t_i}^{n_i}}dx_{t_i} + \cdots +
\frac{
\begin{pmatrix}
\theta_{n_i-1,t_i}^+ & 0 \\
0 & \theta_{n_i-1,t_i}^-
\end{pmatrix}
}{x_{t_i}} dx_{t_i}.
$$
We call the tuple $((\theta^+_{0,t_i},\theta^-_{0,t_i} ),\ldots,
 (\theta^+_{n_i-1,t_i},\theta^-_{n_i-1,t_i} ))$ the 
 local formal data of $\nabla$ at $t_i$.}
\item[(2)] {\it If $i \in I_{\mathrm{ra}}$,
then there exists a matrix $M \in \mathrm{GL}_2(\mathbb{C}[[x_{t_i}]])$ such 
that 
$$
M^{-1} dM +M^{-1} \Omega M =
\begin{pmatrix}
\alpha_i & \beta_i \\
x_{t_i}\beta_i & \alpha_i -\frac{dx_{t_i}}{2x_{t_i}}
\end{pmatrix},
$$
where
\begin{equation*}
\begin{cases}
\alpha_i:= \frac{\theta_{0,t_i}}{2}\frac{dx_{t_i}}{x_{t_i}^{n_i}}+ \cdots 
+\frac{\theta_{2l,t_i}}{2} \frac{dx_{t_i}}{x_{t_i}^{n_i-l}}+\cdots 
+\frac{\theta_{2n_i-2,t_i}}{2} \frac{dx_{t_i}}{x_{t_i}} \\
\beta_i:= \frac{\theta_{1,t_i}}{2}\frac{dx_{t_i}}{x_{t_i}^{n_i}}+ \cdots 
+ \frac{\theta_{2l+1,t_i}}{2}\frac{dx_{t_i}}{x_{t_i}^{n_i-l}}+\cdots 
+\frac{\theta_{2n_i-3,t_i}}{2}\frac{dx_{t_i}}{x^2_{t_i}}.
\end{cases}
\end{equation*}
We call the tuple $(\theta_{0,t_i},\ldots, \theta_{2n_i-2,t_i})$ the 
 local formal data of $\nabla$ at $t_i$.}
\end{itemize}
\end{Lem}

If we define $\zeta_i$ as $x_{t_i}=\zeta_i^2$ and put 
 \begin{equation}\label{2020.1.10.13.58}
 M_{\zeta_i}
: =\begin{pmatrix}
 1 & 1 \\
 \zeta_i & -\zeta_i
 \end{pmatrix},
 \end{equation}
then we have the following diagonalization:
\begin{equation*}
\begin{aligned}
& M_{\zeta_i}^{-1} d  M_{\zeta_i} +
 M_{\zeta_i}^{-1}
\begin{pmatrix}
\alpha_i & \beta_i \\
x_{t_i}\beta_i & \alpha_i -\frac{dx_{t_i}}{2x_{t_i}}
\end{pmatrix} M_{\zeta_i} \\
&=\sum_{l=0,1,\ldots, n_i-1 }
\begin{pmatrix}
\frac{\theta_{2l,t_i}d\zeta_i}{\zeta_i^{2(n_i-l) -1}} & 0 \\
0 & \frac{\theta_{2l,t_i}d\zeta_i}{\zeta_i^{2(n_i-l) -1}} 
\end{pmatrix}+
\sum_{l=0,1,\ldots, n_i-2}
\begin{pmatrix}
\frac{\theta_{2l+1,t_i}d\zeta_i}{\zeta_i^{2(n_i-l) -2}} & 0 \\
0 &- \frac{\theta_{2l+1,t_i}d\zeta_i}{\zeta_i^{2(n_i-l) -2}} 
\end{pmatrix}.
\end{aligned}
\end{equation*}

\begin{Def}
\begin{itemize}
\item Let $(\boldsymbol{\theta},\boldsymbol{\theta}_0) \in
T_{\boldsymbol{\theta}} \times T_{\boldsymbol{\theta}}^{\text{res}}$
and let
$\nabla \colon E_1 \rightarrow E_1 \otimes \Omega^1_{\mathbb{P}^1}(D)$
be a connection.
If the tuple of the local formal data of $\nabla$ is $(\boldsymbol{\theta},\boldsymbol{\theta}_0)$,
we call this connection a 
{\it $(\boldsymbol{\theta},\boldsymbol{\theta}_0)$-connection on $E_1$}.

\item We say that Diarra--Loray's normal form 
 $\nabla_{\mathrm{DL}}^{(n-2)} 
\colon E_{n-2} \rightarrow E_{n-2} \otimes \Omega^1_{\mathbb{P}^1}(D+D_{\mathrm{App}})$
is {\it $(\boldsymbol{\theta},\boldsymbol{\theta}_0)$-connection}
if the corresponding $\nabla_{\mathrm{DL}}^{(1)} 
\colon E_{1} \rightarrow E_{1} \otimes \Omega^1_{\mathbb{P}^1}(D)$
is  a 
$(\boldsymbol{\theta},\boldsymbol{\theta}_0)$-connection on $E_1$.
\end{itemize}
\end{Def}

\begin{Lem}[{\cite[Lemma 16 and Lemma 19]{{DF}}}]\label{2019.12.30.22.00}
{\it 
Let
$(\{ (q_1,p_1),\ldots ,(q_{n-3},p_{n-3}) \} ,(t_3,\ldots,t_{\nu}) )
\in \widehat{\mathcal{M}}_{\boldsymbol{t}_{\text{ra}}}$ and 
 $(\boldsymbol{\theta},\boldsymbol{\theta}_0) \in
T_{\boldsymbol{\theta}} \times T_{\boldsymbol{\theta}}^{\text{res}}$.
Set 
$D=n_1 \cdot 0+n_2 \cdot 1+ \sum_{i=3}^{\nu} n_i \cdot t_i+n_\infty \cdot \infty$.
There exists
a unique tuple of the polynomials $((C_i, D_i)_{i \in I}, \tilde{C})$ in \eqref{2020.12.26.21.03}
such that 
\begin{itemize}
\item the polar divisor of $\nabla_{\mathrm{DL}}^{(n-2)}$
is $D+q_1+\cdots +q_{n-3}$;
\item $q_1,\ldots,q_{n-3}$ are apparent singularities;
\item the dual parameter with respect to $q_j$ is $p_j$ ($j=1,2,\ldots,n-3$);
\item $\nabla_{\mathrm{DL}}^{(n-2)}$ 
is a $(\boldsymbol{\theta},\boldsymbol{\theta}_0)$-connection.
\end{itemize}
}
\end{Lem}

Let $((C_i, D_i)_{i \in I}, \tilde{C})$ be 
the tuple of the polynomials in Lemma \ref{2019.12.30.22.00}.
The polynomials $C_i$ and $D_i$ ($i \in I$) have
simple descriptions.
Now we give explicit descriptions of $C_i$ and $D_i$ ($i \in I$).
For $i \in I_{\mathrm{reg}} \cup I_{\mathrm{un}}$,
we define a polynomial $\Theta^{\pm}_{i}$ in $x$ as
$\sum^{n_i-1}_{l=0}  \frac{\theta_{l,t_i}^{\pm}}{(x-t_i)^{n_i-l}} 
= \frac{\Theta^{\pm}_{i}}{(x-t_i)^{n_i}}$.
For $i \in I_{\mathrm{ra}}$, 
we define polynomials $A_{i}$ and $B_{i}$ in $x$ as
$\sum^{n_i-1}_{l=0} \frac{\theta_{2l,t_i}}{2(x-t_i)^{n_i-l}} = \frac{A_{i}}{(x-t_i)^{n_i}}$
and 
$\sum^{n_i-2}_{l=0} \frac{\theta_{2l+1,t_i}}{2(x-t_i)^{n_i-l}} 
= \frac{B_{i}}{(x-t_i)^{n_i}}$,
respectively.
For $i \in I_{\mathrm{reg}} \cup I_{\mathrm{un}}$,
the polynomials $C_i$ and $D_i$ have following description:
\begin{equation}\label{2021.4.30.13.59}
\begin{cases}
C_i = - \left( \Theta^+_{i}
\cdot \Theta^-_{i} \cdot
\prod_{j\neq i } ( x- t_j)^{n_j} \right)  \quad \mathrm{mod}  \  (x-t_{i})^{n_i} \\
D_i =\Theta^+_{i}+\Theta^-_{i}.
\end{cases}
\end{equation}
For $i \in I_{\mathrm{ra}}$,
the polynomials $C_i$ and $D_i$ have following description:
\begin{equation}\label{2021.4.30.13.59.2}
\begin{cases}
C_i = - \left( 
\left( A_i^2- \frac{(x-t_i)^{n_i-1}}{2}A_i-(x-t_i)B_i^2 \right)\cdot
\prod_{j\neq i} (x-t_j)^{n_j}
\right)  \quad \mathrm{mod}  \  (x-t_{i})^{n_i} \\
D_i =2A_i -\frac{(x-t_i)^{n_i-1}}{2}.
\end{cases}
\end{equation}

\subsection{Family of connections}

Let $(\boldsymbol{\theta},\boldsymbol{\theta}_0) \in
T_{\boldsymbol{\theta}} \times T_{\boldsymbol{\theta}}^{\text{res}}$
and $\boldsymbol{t}_0 =(t_3 ,\ldots,t_{\nu})
\in (T_{\boldsymbol{t}})_{\boldsymbol{t}_0}$.
Set $D=n_1 \cdot 0+n_2 \cdot 1+ n_3 \cdot t_3 + \cdots +
n_\nu \cdot t_{\nu}+n_\infty \cdot \infty$.
Let $\mathfrak{Conn}_{(\boldsymbol{t}_0,\boldsymbol{\theta},\boldsymbol{\theta}_0)}$
be the moduli space of $(\boldsymbol{\theta},\boldsymbol{\theta}_0)$-connections
satisfying some generic conditions:
\begin{equation}\label{2021.4.29.11.42}
\mathfrak{Conn}_{(\boldsymbol{t}_0,\boldsymbol{\theta},\boldsymbol{\theta}_0)}:=
\left\{ (E,\nabla)
\ \middle| \ 
\begin{array}{l}
\text{$E \cong \mathcal{O}_{\mathbb{P}^1} \oplus \mathcal{O}_{\mathbb{P}^1}(1)$ and }
\text{$\nabla$ is a $(\boldsymbol{\theta},\boldsymbol{\theta}_0)$-connection such that}\\
\text{$\mathcal{O}_{\mathbb{P}^1}(1) \subset E_1$ is not $\nabla$-invariant,
$D_{\mathrm{App}}$ is reduced, and} \\
\text{$D_{\mathrm{App}}$ has disjoint support with $D$} 
\end{array}
\right\} / \sim.
\end{equation}
Here $(E,\nabla) \sim (E',\nabla')$ means that 
there exists an isomorphism $\varphi \colon E \rightarrow E'$ such that the following diagram 
is commutative:
$$
\xymatrix{
E  \ar[r]^-{\nabla} \ar[d]^-{\varphi} & E \times \Omega_{\mathbb{P}^1}^1(D) \ar[d]^-{\varphi \otimes \mathrm{id}} \\
E' \ar[r]^-{\nabla'}  & E' \times \Omega_{\mathbb{P}^1}^1(D).
}
$$

By taking apparent singularities and the dual parameters from a connection 
$(E_1,\nabla) \in \mathfrak{Conn}_{(\boldsymbol{t}_0,\boldsymbol{\theta},\boldsymbol{\theta}_0)}$,
we may define a map
\begin{equation*}
\begin{aligned}
\mathrm{App} \colon 
\mathfrak{Conn}_{(\boldsymbol{t}_0,\boldsymbol{\theta},\boldsymbol{\theta}_0)} 
&\longrightarrow
\mathcal{M}_{\boldsymbol{t}_{0},\boldsymbol{t}_{\text{ra}}} 
\subset \mathrm{Sym}^{(n-3)}(\mathbb{C}^2) \\
(E,\nabla)
&\longmapsto 
\{ (q_1,p_1),\ldots,(q_{n-3},p_{n-3}) \} .
\end{aligned}
\end{equation*}
Now we construct an inverse map of $\mathrm{App}$ as follows.
Let $d$ be the relative exterior derivative of 
$\mathbb{P}^1 \times  \mathcal{M}_{\boldsymbol{t}_{0},\boldsymbol{t}_{\text{ra}}}
\rightarrow  \mathcal{M}_{\boldsymbol{t}_{0},\boldsymbol{t}_{\text{ra}}}$.
By Definition \ref{2021.4.11.20.42} and 
Lemma \ref{2019.12.30.22.00},
we may construct an algebraic family 
$$
\tilde{\nabla}_{\mathrm{DL}}^{(n-2)}=
\begin{cases}
d+ \Omega_{(\boldsymbol{t}_0,\boldsymbol{\theta},\boldsymbol{\theta}_0)}^{(n-2)} 
& \text{ on $U_0 \times \mathcal{M}_{\boldsymbol{t}_{0},\boldsymbol{t}_{\text{ra}}} $} \\
d + G_{n-2}^{-1} dG_{n-2}+G_{n-2}^{-1} 
\Omega_{(\boldsymbol{t}_0,\boldsymbol{\theta},\boldsymbol{\theta}_0)}^{(n-2)} \, G_{n-2}
 & \text{ on $U_\infty\times \mathcal{M}_{\boldsymbol{t}_{0},\boldsymbol{t}_{\text{ra}}} $}
\end{cases}
$$
of $(\boldsymbol{\theta},\boldsymbol{\theta}_0)$-connections 
on $E_{n-2}$
parametrized by $\mathcal{M}_{\boldsymbol{t}_{0},\boldsymbol{t}_{\text{ra}}} $.
We set 
$\Omega_{(\boldsymbol{t}_0,\boldsymbol{\theta},\boldsymbol{\theta}_0)}^{(1)}=
\tilde{G}^{-1} d\tilde{G}+\tilde{G}^{-1} 
\Omega_{(\boldsymbol{t}_0,\boldsymbol{\theta},\boldsymbol{\theta}_0)}^{(n-2)}
 \tilde{G}$.
Then we have an algebraic family
$$
\tilde{\nabla}_{\mathrm{DL}}^{(1)}=
\begin{cases}
d+ \Omega_{(\boldsymbol{t}_0,\boldsymbol{\theta},\boldsymbol{\theta}_0)}^{(1)}
 & \text{ on $U_0 \times \mathcal{M}_{\boldsymbol{t}_{0},\boldsymbol{t}_{\text{ra}}} $} \\
d + G_{1}^{-1} dG_{1}+G_{1}^{-1} 
\Omega_{(\boldsymbol{t}_0,\boldsymbol{\theta},\boldsymbol{\theta}_0)}^{(1)} \, G_{1}
 & \text{ on $U_\infty\times \mathcal{M}_{\boldsymbol{t}_{0},\boldsymbol{t}_{\text{ra}}} $}
\end{cases}
$$
of $(\boldsymbol{\theta},\boldsymbol{\theta}_0)$-connections
$E_1 \rightarrow E_1 \otimes \Omega_{\mathbb{P}^1}^1(D)$ parametrized by 
$\mathcal{M}_{\boldsymbol{t}_{0},\boldsymbol{t}_{\text{ra}}} $
by Proposition \ref{2021.4.11.16.40}.
The algebraic family $\tilde{\nabla}_{\mathrm{DL}}^{(1)}$
parametrized by $\mathcal{M}_{\boldsymbol{t}_{0},\boldsymbol{t}_{\text{ra}}} $
gives the inverse map of $\mathrm{App}$:
\begin{equation*}
\begin{aligned}
\mathrm{App}^{-1} \colon \mathcal{M}_{\boldsymbol{t}_{0},\boldsymbol{t}_{\text{ra}}} 
&\longrightarrow
\mathfrak{Conn}_{(\boldsymbol{t}_0,\boldsymbol{\theta},\boldsymbol{\theta}_0)}  \\
\boldsymbol{p} =\{ (q_1,p_1),\ldots,(q_{n-3},p_{n-3}) \}
&\longmapsto 
(E_1,\tilde\nabla^{(1)}_{\mathrm{DL}}|_{\mathbb{P}^1 \times \{\boldsymbol{p} \}} ).
\end{aligned}
\end{equation*}

Next we consider the extended moduli space 
$\widehat{\mathfrak{Conn}}_{(\boldsymbol{t}_{\mathrm{ra}},\boldsymbol{\theta}_0)}$
of 
$\mathfrak{Conn}_{(\boldsymbol{t}_0,\boldsymbol{\theta},\boldsymbol{\theta}_0)} $.
We set
$D(\boldsymbol{t}_0) :=n_1 \cdot 0+n_2 \cdot 1+ \sum_{i=3}^{\nu} n_i \cdot t_i 
+n_\infty \cdot \infty$ for $\boldsymbol{t}_0 =(t_3 ,\ldots,t_{\nu})
\in (T_{\boldsymbol{t}})_{\boldsymbol{t}_0}$.
This extended moduli space 
$\widehat{\mathfrak{Conn}}_{(\boldsymbol{t}_{\mathrm{ra}},\boldsymbol{\theta}_0)}$
is defined by
$$
\widehat{\mathfrak{Conn}}_{(\boldsymbol{t}_{\mathrm{ra}},\boldsymbol{\theta}_0)}:=
\left\{ (E,\nabla,\boldsymbol{t}_0,\boldsymbol{\theta})
\ \middle| \ 
\begin{array}{l}
\text{$\boldsymbol{t}_0 \in (T_{\boldsymbol{t}})_{\boldsymbol{t}_{\mathrm{ra}}}$, 
$\boldsymbol{\theta} \in T_{\boldsymbol{\theta}}$},  
\text{$E \cong \mathcal{O}_{\mathbb{P}^1} \oplus \mathcal{O}_{\mathbb{P}^1}(1)$, and}\\
\text{$\nabla$ is a $(\boldsymbol{\theta},\boldsymbol{\theta}_0)$-connection such that
the polar divisor of $\nabla$ is $D(\boldsymbol{t}_0)$,} \\
\text{$\mathcal{O}_{\mathbb{P}^1}(1) \subset E_1$ is not $\nabla$-invariant,
$D_{\mathrm{App}}$ is reduced, and}\\
\text{$D_{\mathrm{App}}$ has disjoint support with $D(\boldsymbol{t}_0)$} 
\end{array}
\right\} / \sim.
$$
Here $(E,\nabla,\boldsymbol{t}_0,\boldsymbol{\theta}) 
\sim (E',\nabla',\boldsymbol{t}_0',\boldsymbol{\theta}')$ means that 
$\boldsymbol{t}_0=\boldsymbol{t}_0'$, $\boldsymbol{\theta}=\boldsymbol{\theta}'$, and
there exists an isomorphism $\varphi \colon E \rightarrow E'$ such that the following diagram 
is commutative:
$$
\xymatrix{
E  \ar[r]^-{\nabla} \ar[d]^-{\varphi} & E \times \Omega_{\mathbb{P}^1}^1(D) \ar[d]^-{\varphi \otimes \mathrm{id}} \\
E' \ar[r]^-{\nabla'}  & E' \times \Omega_{\mathbb{P}^1}^1(D).
}
$$
By taking apparent singularities, the dual parameters,
the position of singular points, and the local formal data from a connection 
$(E_1,\nabla) \in \widehat{\mathfrak{Conn}}_{(\boldsymbol{t}_{\mathrm{ra}},\boldsymbol{\theta}_0)}$,
we may define a map
\begin{equation}\label{2020.1.6.13.00}
\begin{aligned}
\widehat{\mathrm{App}} \colon 
\widehat{\mathfrak{Conn}}_{(\boldsymbol{t}_{\mathrm{ra}},\boldsymbol{\theta}_0)}
&\longrightarrow
\widehat{\mathcal{M}}_{\boldsymbol{t}_{\text{ra}}} \times T_{\boldsymbol{\theta}}\subset
(\mathrm{Sym}^{(n-3)}(\mathbb{C}^2)\times T_{\boldsymbol{t}} )
\times T_{\boldsymbol{\theta}} \\
(E,\nabla,\boldsymbol{t}_0,\boldsymbol{\theta})
&\longmapsto 
((\{ (q_1,p_1),\ldots,(q_{n-3},p_{n-3}) \}, \boldsymbol{t}_0),\boldsymbol{\theta}) .
\end{aligned}
\end{equation}
Now we may also construct an inverse map of $\widehat{\mathrm{App}}$ as follows.
Here let $d$ be the relative exterior derivative of 
$\mathbb{P}^1 \times  (\widehat{\mathcal{M}}_{\boldsymbol{t}_{\text{ra}}} \times T_{\boldsymbol{\theta}})
\rightarrow  \widehat{\mathcal{M}}_{\boldsymbol{t}_{\text{ra}}} \times T_{\boldsymbol{\theta}}$.
By Definition \ref{2021.4.11.20.42} and 
Lemma \ref{2019.12.30.22.00},
we may construct an algebraic family 
$$
\tilde{\nabla}_{\mathrm{DL, ext}}^{(n-2)}=
\begin{cases}
d+ \widehat\Omega_{(\boldsymbol{t}_{\mathrm{ra}},\boldsymbol{\theta}_0)}^{(n-2)} 
& \text{ on $U_0 \times( \widehat{\mathcal{M}}_{\boldsymbol{t}_{\text{ra}}} 
\times T_{\boldsymbol{\theta}})$} \\
d + G_{n-2}^{-1} dG_{n-2}+G_{n-2}^{-1} 
\widehat\Omega_{(\boldsymbol{t}_{\mathrm{ra}},\boldsymbol{\theta}_0)}^{(n-2)} \, G_{n-2}
 & \text{ on $U_\infty
 \times (\widehat{\mathcal{M}}_{\boldsymbol{t}_{\text{ra}}} \times T_{\boldsymbol{\theta}})$}
\end{cases}
$$
of connections 
on $E_{n-2}$
parametrized by $\widehat{\mathcal{M}}_{\boldsymbol{t}_{\text{ra}}} \times T_{\boldsymbol{\theta}}$.
We set 
$\widehat\Omega_{(\boldsymbol{t}_{\mathrm{ra}},\boldsymbol{\theta}_0)}^{(1)}=
\tilde{G}^{-1} d\tilde{G}+\tilde{G}^{-1} 
\widehat\Omega_{(\boldsymbol{t}_{\mathrm{ra}},\boldsymbol{\theta}_0)}^{(n-2)}
 \tilde{G}$.
Then we have an algebraic family
\begin{equation}\label{2021.4.14.10.36}
\tilde{\nabla}_{\mathrm{DL, ext}}^{(1)}=
\begin{cases}
d+ \widehat\Omega_{(\boldsymbol{t}_{\mathrm{ra}},\boldsymbol{\theta}_0)}^{(1)} 
& \text{ on $U_0\times(\widehat{\mathcal{M}}_{\boldsymbol{t}_{\text{ra}}} 
\times T_{\boldsymbol{\theta}})$} \\
d + G_{1}^{-1} dG_{1}+G_{1}^{-1} 
\widehat\Omega_{(\boldsymbol{t}_{\mathrm{ra}},\boldsymbol{\theta}_0)}^{(1)} \, G_{1}
 & \text{ on $U_\infty \times (\widehat{\mathcal{M}}_{\boldsymbol{t}_{\text{ra}}} 
 \times T_{\boldsymbol{\theta}})$}
\end{cases}
\end{equation}
of connections on 
$E_1$ parametrized by 
$\widehat{\mathcal{M}}_{\boldsymbol{t}_{\text{ra}}} \times T_{\boldsymbol{\theta}}$
by Proposition \ref{2021.4.11.16.40}.
Let $\tilde{\boldsymbol{t}}_0=(\tilde{t}_1,\ldots,\tilde{t}_{\nu},\tilde{t}_{\infty})$ 
be a family of $(\nu+1)$-points on $\mathbb{P}^1$
parametrized by $(T_{\boldsymbol{t}})_{\boldsymbol{t}_{\text{ra}}}$
and $\tilde{\boldsymbol{\theta}}$ be 
a family of tuples of complex numbers 
parametrized by $T_{\boldsymbol{\theta}}$.
We denote by the same characters $\tilde{\boldsymbol{t}}_0$ and
$\tilde{\boldsymbol{\theta}}$
the pull-backs of $\tilde{\boldsymbol{t}}_0$ and $\tilde{\boldsymbol{\theta}}$
under the compositions
\begin{equation*}
\widehat{\mathcal{M}}_{\boldsymbol{t}_{\text{ra}}}\times T_{\boldsymbol{\theta}} 
\xrightarrow{\ \mathrm{id} \times \pi_{\boldsymbol{t}_{\mathrm{ra}}, \boldsymbol{\theta}_0} \ }
(T_{\boldsymbol{t}})_{\boldsymbol{t}_{\text{ra}}} \times T_{\boldsymbol{\theta}}
\xrightarrow{\text{projection}}
(T_{\boldsymbol{t}})_{\boldsymbol{t}_{\text{ra}}}
\end{equation*}
and
\begin{equation*}
\widehat{\mathcal{M}}_{\boldsymbol{t}_{\text{ra}}}\times T_{\boldsymbol{\theta}} 
\xrightarrow{\ \mathrm{id} \times \pi_{\boldsymbol{t}_{\mathrm{ra}}, \boldsymbol{\theta}_0} \ }
(T_{\boldsymbol{t}})_{\boldsymbol{t}_{\text{ra}}} \times T_{\boldsymbol{\theta}}
\xrightarrow{\text{projection}}
T_{\boldsymbol{\theta}},
\end{equation*}
respectively.
The algebraic family $\tilde{\nabla}_{\mathrm{DL, ext}}^{(1)}$
parametrized by $\mathcal{M}_{\boldsymbol{t}_{\text{ra}}} \times T_{\boldsymbol{\theta}}$
gives the inverse map of $\widehat{\mathrm{App}}$:
\begin{equation*}
\begin{aligned}
\widehat{\mathrm{App}}^{-1} \colon
\widehat{\mathcal{M}}_{\boldsymbol{t}_{\text{ra}}} \times T_{\boldsymbol{\theta}}
&\longrightarrow
\widehat{\mathfrak{Conn}}_{(\boldsymbol{t}_{\mathrm{ra}},\boldsymbol{\theta}_0)}  \\
\hat{\boldsymbol{p}} =
((\{ (q_1,p_1),\ldots,(q_{n-3},p_{n-3}) \}, \boldsymbol{t}_0),\boldsymbol{\theta})
&\longmapsto 
(E_1,\tilde\nabla^{(1)}_{\mathrm{DL, ext}}|_{\mathbb{P}^1 \times \{\hat{\boldsymbol{p}} \}},
\tilde{\boldsymbol{t}}_0|_{\hat{\boldsymbol{p}} }, 
\tilde{\boldsymbol{\theta}}|_{\hat{\boldsymbol{p}} }
) .
\end{aligned}
\end{equation*}

\subsection{Infinitesimal deformations of connections}\label{2021.4.29.12.39}

Let $U$ be an open subset of 
$\mathcal{M}_{\boldsymbol{t}_{0},\boldsymbol{t}_{\text{ra}}}
 \subset \mathrm{Sym}^{(n-3)}(\mathbb{C}^2)$.
Let $\delta$ be a vector field on $U$.
By the vector field $\delta$, we have a map
\begin{equation}\label{2021.4.12.13.52}
\mathbb{P}^1 \times \mathrm{Spec}\,
 \mathcal{O}_{U}[\epsilon]
 \longrightarrow \mathbb{P}^1 \times U,
\end{equation}
where $\epsilon^2=0$. 
We take 
the pull-back of 
the family $\tilde{\nabla}_{\mathrm{DL}}^{(n-2)}|_{\mathbb{P}^1 \times U}$
under the map \eqref{2021.4.12.13.52}.
We denote by  
$$
\begin{cases}
d+ \Omega_{(\boldsymbol{t}_0,\boldsymbol{\theta},\boldsymbol{\theta}_0)}^{(n-2)} 
+\epsilon \delta(\Omega_{(\boldsymbol{t}_0,\boldsymbol{\theta},\boldsymbol{\theta}_0)}^{(n-2)})
& \text{ on $U_0 \times \mathrm{Spec}\,
 \mathcal{O}_{U}[\epsilon]$} \\
d + G_{n-2}^{-1} dG_{n-2}+G_{n-2}^{-1} 
\Omega_{(\boldsymbol{t}_0,\boldsymbol{\theta},\boldsymbol{\theta}_0)}^{(n-2)} \, G_{n-2}
+\epsilon G_{n-2}^{-1} 
\delta(\Omega_{(\boldsymbol{t}_0,\boldsymbol{\theta},\boldsymbol{\theta}_0)}^{(n-2)})
G_{n-2}
 & \text{ on $U_\infty \times \mathrm{Spec}\,
 \mathcal{O}_{U}[\epsilon]$}
\end{cases}
$$
the expansion of this pull-back of $\tilde{\nabla}_{\mathrm{DL}}^{(n-2)}|_{\mathbb{P}^1 \times U}$
with respect to $\epsilon$.
We also denote by  
$$
\begin{cases}
d+ \Omega_{(\boldsymbol{t}_0,\boldsymbol{\theta},\boldsymbol{\theta}_0)}^{(1)} 
+\epsilon \delta(\Omega_{(\boldsymbol{t}_0,\boldsymbol{\theta},\boldsymbol{\theta}_0)}^{(1)})
& \text{ on $U_0 \times \mathrm{Spec}\,
 \mathcal{O}_{U}[\epsilon]$} \\
d + G_{1}^{-1} dG_{1}+G_{1}^{-1} 
\Omega_{(\boldsymbol{t}_0,\boldsymbol{\theta},\boldsymbol{\theta}_0)}^{(1)} \, G_{1}
+\epsilon G_{1}^{-1} 
\delta(\Omega_{(\boldsymbol{t}_0,\boldsymbol{\theta},\boldsymbol{\theta}_0)}^{(1)})
G_{1}
 & \text{ on $U_\infty \times \mathrm{Spec}\,
 \mathcal{O}_{U}[\epsilon]$}
\end{cases}
$$
the expansion of
the pull-back of 
the family $\tilde{\nabla}_{\mathrm{DL}}^{(1)}|_{\mathbb{P}^1 \times U}$
under the map \eqref{2021.4.12.13.52}.

Let $\hat{U}$ be an open subset of 
$\widehat{\mathcal{M}}_{\boldsymbol{t}_{\text{ra}}}\times T_{\boldsymbol{\theta}}$.
Let $\hat\delta$ be a vector field on $\hat{U}$.
By the vector field $\hat\delta$, we have a map
\begin{equation}\label{2021.4.12.13.55}
\mathbb{P}^1 \times \mathrm{Spec}\,
 \mathcal{O}_{\hat{U}}[\epsilon]
 \longrightarrow \mathbb{P}^1 \times \hat{U},
\end{equation}
where $\epsilon^2=0$. 
We take
the pull-back of 
the family $\tilde{\nabla}_{\mathrm{DL, ext}}^{(n-2)}|_{\mathbb{P}^1 \times \hat{U}}$
under the map \eqref{2021.4.12.13.55}.
We denote by  
$$
\begin{cases}
d+ \widehat\Omega_{(\boldsymbol{t}_{\mathrm{ra}},\boldsymbol{\theta}_0)}^{(n-2)} 
+\epsilon \delta(\widehat\Omega_{(\boldsymbol{t}_{\mathrm{ra}},\boldsymbol{\theta}_0)}^{(n-2)})
& \text{ on $U_0 \times \mathrm{Spec}\,
 \mathcal{O}_{\hat{U}}[\epsilon]$} \\
d + G_{n-2}^{-1} dG_{n-2}+G_{n-2}^{-1} 
\widehat\Omega_{(\boldsymbol{t}_{\mathrm{ra}},\boldsymbol{\theta}_0)}^{(n-2)} \, G_{n-2}
+\epsilon G_{n-2}^{-1} 
\delta(\widehat\Omega_{(\boldsymbol{t}_{\mathrm{ra}},\boldsymbol{\theta}_0)}^{(n-2)})
G_{n-2}
 & \text{ on $U_\infty \times \mathrm{Spec}\,
 \mathcal{O}_{\hat{U}}[\epsilon]$}
\end{cases}
$$
the expansion of this pull-back of $\tilde{\nabla}_{\mathrm{DL, ext}}^{(n-2)}|_{\mathbb{P}^1 \times \hat{U}}$
with respect to $\epsilon$.
Here this $\widehat\Omega_{(\boldsymbol{t}_{\mathrm{ra}},\boldsymbol{\theta}_0)}^{(n-2)}$
means the pull-back of $\widehat\Omega_{(\boldsymbol{t}_{\mathrm{ra}},\boldsymbol{\theta}_0)}^{(n-2)}$
on $U_0 \times \hat{U}$ by the trivial projection 
$U_0 \times \mathrm{Spec}\,
 \mathcal{O}_{\hat{U}}[\epsilon]
 \rightarrow U_0 \times \hat{U}$.
 Remark that there is a difference between
 this $\widehat\Omega_{(\boldsymbol{t}_{\mathrm{ra}},\boldsymbol{\theta}_0)}^{(n-2)}$
 and the pull-back of $\widehat\Omega_{(\boldsymbol{t}_{\mathrm{ra}},\boldsymbol{\theta}_0)}^{(n-2)}$
on $U_0 \times \hat{U}$ by \eqref{2021.4.12.13.55}.
The $\epsilon$-part $ \delta(\widehat\Omega_{(\boldsymbol{t}_{\mathrm{ra}},\boldsymbol{\theta}_0)}^{(n-2)})$
adjusts this difference.
We also denote by  
$$
\begin{cases}
d+ \widehat\Omega_{(\boldsymbol{t}_{\mathrm{ra}},\boldsymbol{\theta}_0)}^{(1)} 
+\epsilon \delta(\widehat\Omega_{(\boldsymbol{t}_{\mathrm{ra}},\boldsymbol{\theta}_0)}^{(1)})
& \text{ on $U_0 \times \mathrm{Spec}\,
 \mathcal{O}_{\tilde{U}}[\epsilon]$} \\
d + G_{1}^{-1} dG_{1}+G_{1}^{-1} 
\widehat\Omega_{(\boldsymbol{t}_{\mathrm{ra}},\boldsymbol{\theta}_0)}^{(1)} \, G_{1}
+\epsilon G_{1}^{-1} 
\delta(\widehat\Omega_{(\boldsymbol{t}_{\mathrm{ra}},\boldsymbol{\theta}_0)}^{(1)})
G_{1}
 & \text{ on $U_\infty \times \mathrm{Spec}\,
 \mathcal{O}_{\tilde{U}}[\epsilon]$}
\end{cases}
$$
the expansion of the pull-back of 
the family $\tilde{\nabla}_{\mathrm{DL, ext}}^{(1)}|_{\mathbb{P}^1 \times \hat{U}}$
under the map \eqref{2021.4.12.13.55}.

\subsection{Horizontal lifts of a family of connections}\label{2021.4.29.12.40}

Let $\widehat{E}_{1}$ and $\widehat{E}_{n-2}$ be 
the pull-backs of $E_1$ and $E_{n-2}$ under the projection
$\mathbb{P}^1 \times(\widehat{\mathcal{M}}_{\boldsymbol{t}_{\text{ra}}} \times T_{\boldsymbol{\theta}})
\rightarrow \mathbb{P}^1$, respectively.
Set
$D(\tilde{\boldsymbol{t}}_0) :=
\sum_{i=1}^{\nu} n_i \cdot \tilde{t}_i+ n_{\infty} \cdot \tilde{t}_{\infty},$
which is a Cartier divisor 
on $\mathbb{P}^1 \times \widehat{\mathcal{M}}_{\boldsymbol{t}_{\text{ra}}}
\times T_{\boldsymbol{\theta}} $,
which is flat over 
$\widehat{\mathcal{M}}_{\boldsymbol{t}_{\text{ra}}}\times T_{\boldsymbol{\theta}}$.
Let 
$$
\tilde{\nabla}_{\mathrm{DL, ext}}^{(1)}=
\begin{cases}
d+ \widehat\Omega_{(\boldsymbol{t}_{\mathrm{ra}},\boldsymbol{\theta}_0)}^{(1)} 
& \text{ on $U_0\times(\widehat{\mathcal{M}}_{\boldsymbol{t}_{\text{ra}}} 
\times T_{\boldsymbol{\theta}})$} \\
d + G_{1}^{-1} dG_{1}+G_{1}^{-1} 
\widehat\Omega_{(\boldsymbol{t}_{\mathrm{ra}},\boldsymbol{\theta}_0)}^{(1)} \, G_{1}
 & \text{ on $U_\infty \times (\widehat{\mathcal{M}}_{\boldsymbol{t}_{\text{ra}}} 
 \times T_{\boldsymbol{\theta}})$}
\end{cases}
$$
be the family \eqref{2021.4.14.10.36}.
This family is a relative connection on $\widehat{E}_{1}$:
$$
\tilde{\nabla}_{\mathrm{DL, ext}}^{(1)} 
\colon \widehat{E}_{1}
\longrightarrow 
\widehat{E}_{1} \otimes 
\Omega^1_{\mathbb{P}^1 \times (\widehat{\mathcal{M}}_{\boldsymbol{t}_{\text{ra}}}
\times T_{\boldsymbol{\theta}} )
/ \widehat{\mathcal{M}}_{\boldsymbol{t}_{\text{ra}}}\times T_{\boldsymbol{\theta}} }
(D(\tilde{\boldsymbol{t}}_0)).
$$
We will consider an infinitesimal deformation of $\tilde{\nabla}_{\mathrm{DL, ext}}^{(1)}$,
which means an ``integrable deformation''.

Let $\widehat{E}_1^{\epsilon}$ be the pull-back of $\widehat{E}_1$
under the trivial projection  
\begin{equation}\label{2021.4.14.10.57}
\mathbb{P}^1 \times (\widehat{\mathcal{M}}_{\boldsymbol{t}_{\text{ra}}}\times T_{\boldsymbol{\theta}} )
\times \mathrm{Spec}\, \mathbb{C} [\epsilon] \longrightarrow
\mathbb{P}^1 \times (\widehat{\mathcal{M}}_{\boldsymbol{t}_{\text{ra}}}\times T_{\boldsymbol{\theta}} ).
\end{equation}
Let $\delta_{\text{time}}$ be 
a vector field on $(T_{\boldsymbol{t}})_{\boldsymbol{t}_{\text{ra}}} \times T_{\boldsymbol{\theta}}$.
The vector field $\delta_{\text{time}}$ gives a map 
$$
\pi_{\delta_{\text{time}}} \colon 
((T_{\boldsymbol{t}})_{\boldsymbol{t}_{\text{ra}}} \times T_{\boldsymbol{\theta}})
\times \mathrm{Spec}\, \mathbb{C} [\epsilon]
 \longrightarrow 
(T_{\boldsymbol{t}})_{\boldsymbol{t}_{\text{ra}}} \times T_{\boldsymbol{\theta}}.
$$
Set
$$
D_{\text{red}}(\tilde{\boldsymbol{t}}_0) :=
\sum_{i \in I}  \tilde{t}_i .
$$
We consider $D_{\text{red}}(\tilde{\boldsymbol{t}}_0)$
as a Cartier divisor on 
$\mathbb{P}^1 \times ((T_{\boldsymbol{t}})_{\boldsymbol{t}_{\text{ra}}}
 \times T_{\boldsymbol{\theta}})$.
We denote by 
$D(\tilde{\boldsymbol{t}}_0)_{\epsilon}$ and
$D_{\text{red}}(\tilde{\boldsymbol{t}}_0)_{\epsilon}$
the pull-backs of 
$D(\tilde{\boldsymbol{t}}_0)$ and
$D_{\text{red}}(\tilde{\boldsymbol{t}}_0)$
under the composition 
\begin{equation}\label{2021.4.20.16.11}
\mathbb{P}^1 \times (\widehat{\mathcal{M}}_{\boldsymbol{t}_{\text{ra}}}\times T_{\boldsymbol{\theta}} )
\times \mathrm{Spec}\, \mathbb{C} [\epsilon]
\xrightarrow{\mathrm{id} \times  \pi_{\boldsymbol{t}_{\mathrm{ra}}, \boldsymbol{\theta}_0} 
\times \mathrm{id}}
\mathbb{P}^1 \times((T_{\boldsymbol{t}})_{\boldsymbol{t}_{\text{ra}}} \times T_{\boldsymbol{\theta}})
\times \mathrm{Spec}\, \mathbb{C} [\epsilon]
 \xrightarrow{\mathrm{id}\times \pi_{\delta_{\text{time}}}}
\mathbb{P}^1 \times (T_{\boldsymbol{t}})_{\boldsymbol{t}_{\text{ra}}} \times T_{\boldsymbol{\theta}},
\end{equation}
respectively.
Take local defining equations $\tilde{x}_{t_i}$ of 
the Cartier divisor $D_{\text{red}}(\tilde{\boldsymbol{t}}_0)_{\epsilon}$.
Let $\tilde{\Omega}^1_{\delta_{\text{time}}}$ be a coherent subsheaf 
of $\Omega^1_{\mathbb{P}^1 \times
 (\widehat{\mathcal{M}}_{\boldsymbol{t}_{\text{ra}}}\times T_{\boldsymbol{\theta}} )
 \times \mathrm{Spec}\, \mathbb{C} [\epsilon]
/ \widehat{\mathcal{M}}_{\boldsymbol{t}_{\text{ra}}}\times T_{\boldsymbol{\theta}} }
(D(\tilde{\boldsymbol{t}}_0)_{\epsilon})$ which is locally defined by
\begin{equation}\label{2021.4.14.23.20}
\begin{aligned}
\tilde{\Omega}^1_{\delta_{\text{time}}} &=
\mathcal{O}_{\mathbb{P}^1 \times
 (\widehat{\mathcal{M}}_{\boldsymbol{t}_{\text{ra}}}\times T_{\boldsymbol{\theta}} )
 \times \mathrm{Spec}\, \mathbb{C} [\epsilon]} \frac{d\tilde{x}_{t_i}}{\tilde{x}_{t_i}^{n_i}} 
 + \mathcal{O}_{\mathbb{P}^1 \times
 (\widehat{\mathcal{M}}_{\boldsymbol{t}_{\text{ra}}}\times T_{\boldsymbol{\theta}} )}
  \frac{d\epsilon}{\tilde{x}_{t_i}^{n_i-1}}.
\end{aligned}
\end{equation}
Let 
$\nabla_{\delta_{\text{time}}} \colon \widehat{E}_1^{\epsilon}
\rightarrow \widehat{E}_1^{\epsilon} \otimes 
\tilde{\Omega}^1_{\delta_{\text{time}}}$
be a morphism with the Leibniz rule.
That is, 
$\nabla_{\delta_{\text{time}}}(fa) =a \otimes \hat{d}f + f \nabla_{\delta_{\text{time}}}(a)$
for $f \in \mathcal{O}_{\mathbb{P}^1 
\times (\widehat{\mathcal{M}}_{\boldsymbol{t}_{\text{ra}}}\times T_{\boldsymbol{\theta}} )
\times \mathrm{Spec}\, \mathbb{C} [\epsilon]}$ and 
$a \in \widehat{E}_1^{\epsilon}$.
Here $\hat{d}$ is the relative exterior derivative of 
$\mathbb{P}^1 \times 
(\widehat{\mathcal{M}}_{\boldsymbol{t}_{\text{ra}}}\times T_{\boldsymbol{\theta}} )
\times \mathrm{Spec}\, \mathbb{C} [\epsilon]
\rightarrow  \widehat{\mathcal{M}}_{\boldsymbol{t}_{\text{ra}}}\times T_{\boldsymbol{\theta}} $.
We denote by
$$
\nabla_{\delta_{\text{time}}}=
\begin{cases}
\hat{d}+ \widehat\Omega_{\delta_{\text{time}}}^{(1)} 
+\epsilon \delta(\widehat\Omega_{\delta_{\text{time}}}^{(1)}) + \Upsilon_{\delta_{\text{time}}} d\epsilon\\
\qquad \qquad 
\text{ on $U_0 \times(\widehat{\mathcal{M}}_{\boldsymbol{t}_{\text{ra}}}\times T_{\boldsymbol{\theta}} )
\times \mathrm{Spec}\, \mathbb{C} [\epsilon]$ } \\
\hat{d} + G_{1}^{-1} dG_{1}+G_{1}^{-1} 
\widehat\Omega_{\delta_{\text{time}}}^{(1)} \, G_{1}
+\epsilon G_{1}^{-1} 
\delta(\widehat\Omega_{\delta_{\text{time}}}^{(1)})
G_{1}
+ G_1^{-1} \Upsilon_{\delta_{\text{time}}} G_1 d\epsilon \\
\qquad\qquad 
\text{ on $U_\infty \times(\widehat{\mathcal{M}}_{\boldsymbol{t}_{\text{ra}}}\times T_{\boldsymbol{\theta}} )
\times \mathrm{Spec}\, \mathbb{C} [\epsilon]$ }.
\end{cases}
$$
the expansion of the morphism $\nabla_{\delta_{\text{time}}} $
with respect to $\epsilon$.
Here the connection matrices are decomposed into $dx_{t_i}$-terms and $d \epsilon$-terms.
Remark that $dx_{t_i}$ is the pull-back of $dx_{t_i}$ under the morphism
$\mathbb{P}^1 \times (\widehat{\mathcal{M}}_{\boldsymbol{t}_{\text{ra}}}\times T_{\boldsymbol{\theta}} )
\times \mathrm{Spec}\, \mathbb{C} [\epsilon]
\rightarrow
\mathbb{P}^1 \times (T_{\boldsymbol{t}})_{\boldsymbol{t}_{\text{ra}}} \times T_{\boldsymbol{\theta}}$
defined by $  \pi_{\boldsymbol{t}_{\mathrm{ra}}, \boldsymbol{\theta}_0} $ and 
the trivial projection.
On the other hand,
$d\tilde{x}_{t_i}$ is the pull-back of $dx_{t_i}$
under the morphism 
$\mathbb{P}^1 \times (\widehat{\mathcal{M}}_{\boldsymbol{t}_{\text{ra}}}\times T_{\boldsymbol{\theta}} )
\times \mathrm{Spec}\, \mathbb{C} [\epsilon]
\rightarrow
\mathbb{P}^1 \times (T_{\boldsymbol{t}})_{\boldsymbol{t}_{\text{ra}}} \times T_{\boldsymbol{\theta}}$
defined in \eqref{2021.4.20.16.11}.
Moreover, remark that $\widehat\Omega_{\delta_{\text{time}}}^{(1)}$
has a pole on the divisor $D(\tilde{\boldsymbol{t}}_0) \times \mathrm{Spec}\, \mathbb{C} [\epsilon]$,
which is different from the divisor $D(\tilde{\boldsymbol{t}}_0)_{\epsilon}$.
So $\widehat\Omega_{\delta_{\text{time}}}^{(1)}$ does not belong to 
$\mathcal{E}nd (\widehat{E}_1^{\epsilon}) \otimes 
\tilde{\Omega}^1_{\delta_{\text{time}}}$.
The $\epsilon$-term $\epsilon \delta(\widehat\Omega_{\delta_{\text{time}}}^{(1)}) + \Upsilon_{\delta_{\text{time}}} d\epsilon$
adjust the condition that the image of $\nabla_{\delta_{\text{time}}}$
is contained in $\widehat{E}_1^{\epsilon} \otimes 
\tilde{\Omega}^1_{\delta_{\text{time}}}$.

Let $\overline{\nabla}_{\delta_{\text{time}}}$
be the relative connection induced by $\nabla_{\delta_{\text{time}}}$:
$$
\overline{\nabla}_{\delta_{\text{time}}} \colon 
\widehat{E}_1^{\epsilon}
\longrightarrow \widehat{E}_1^{\epsilon} \otimes 
\Omega^1_{\mathbb{P}^1 \times (\widehat{\mathcal{M}}_{\boldsymbol{t}_{\text{ra}}}
\times T_{\boldsymbol{\theta}} ) \times \mathrm{Spec}\, \mathbb{C} [\epsilon]
/ (\widehat{\mathcal{M}}_{\boldsymbol{t}_{\text{ra}}}\times T_{\boldsymbol{\theta}})
\times \mathrm{Spec}\, \mathbb{C} [\epsilon]} 
(D(\tilde{\boldsymbol{t}}_0)_{\epsilon}).
$$
That is, 
$$
\overline{\nabla}_{\delta_{\text{time}}}=
\begin{cases}
d+ \widehat\Omega_{\delta_{\text{time}}}^{(1)} 
+\epsilon \delta(\widehat\Omega_{\delta_{\text{time}}}^{(1)}) 
& \text{ on $U_0 \times(\widehat{\mathcal{M}}_{\boldsymbol{t}_{\text{ra}}}\times T_{\boldsymbol{\theta}} )
\times \mathrm{Spec}\, \mathbb{C} [\epsilon]$ } \\
d + G_{1}^{-1} dG_{1}+G_{1}^{-1} 
\widehat\Omega_{\delta_{\text{time}}}^{(1)} \, G_{1}
+\epsilon G_{1}^{-1} 
\delta(\widehat\Omega_{\delta_{\text{time}}}^{(1)})
G_{1}
 & \text{ on $U_\infty \times(\widehat{\mathcal{M}}_{\boldsymbol{t}_{\text{ra}}}\times 
 T_{\boldsymbol{\theta}} )
\times \mathrm{Spec}\, \mathbb{C} [\epsilon]$ }.
\end{cases}
$$
We consider $\widehat\Omega_{(\boldsymbol{t}_{\mathrm{ra}},\boldsymbol{\theta}_0)}^{(1)}$ as
a matrix
with values in $\Omega^1_{\mathbb{P}^1\times (\widehat{\mathcal{M}}_{\boldsymbol{t}_{\text{ra}}}\times 
 T_{\boldsymbol{\theta}} )}( D(\tilde{\boldsymbol{t}}_0))$.
We take a pull-back of the matrix 
$\widehat\Omega_{(\boldsymbol{t}_{\mathrm{ra}},\boldsymbol{\theta}_0)}^{(1)}$
by the trivial projection 
$U_0 \times (\widehat{\mathcal{M}}_{\boldsymbol{t}_{\text{ra}}}\times T_{\boldsymbol{\theta}} )
\times \mathrm{Spec}\, \mathbb{C} [\epsilon]
\rightarrow
U_0 \times (\widehat{\mathcal{M}}_{\boldsymbol{t}_{\text{ra}}}\times T_{\boldsymbol{\theta}} )$.
This pull-back induces a matrix with 
values in $\Omega^1_{\mathbb{P}^1\times (\widehat{\mathcal{M}}_{\boldsymbol{t}_{\text{ra}}}\times 
 T_{\boldsymbol{\theta}} )\times \mathrm{Spec}\, \mathbb{C} [\epsilon]/
(\widehat{\mathcal{M}}_{\boldsymbol{t}_{\text{ra}}}\times 
 T_{\boldsymbol{\theta}} )  }( D(\tilde{\boldsymbol{t}}_0) \times\mathrm{Spec}\, \mathbb{C} [\epsilon] )$.
 We also denote by $\widehat\Omega_{(\boldsymbol{t}_{\mathrm{ra}},\boldsymbol{\theta}_0)}^{(1)}$
 this induced matrix.

\begin{Def}
We say $\nabla_{\delta_{\text{time}}}$ is a \textit{horizontal lift of $\tilde{\nabla}_{\mathrm{DL, ext}}^{(1)}$} 
if $\nabla_{\delta_{\text{time}}}$ satisfies $ \widehat\Omega_{\delta_{\text{time}}}^{(1)}  
=  \widehat\Omega_{(\boldsymbol{t}_{\mathrm{ra}},\boldsymbol{\theta}_0)}^{(1)} $
and the integrable condition
$$
\delta(\widehat\Omega_{\delta_{\text{time}}}^{(1)}) \wedge d \epsilon = 
d  \Upsilon_{\delta_{\text{time}}}   \wedge d \epsilon
+[\widehat\Omega_{\delta_{\text{time}}}^{(1)}  ,\Upsilon_{\delta_{\text{time}}}]
 \wedge d \epsilon.
$$ 
If $\nabla_{\delta_{\text{time}}}$ is a horizontal lift of $\tilde{\nabla}_{\mathrm{DL, ext}}^{(1)}$,
the relative connection 
$\overline{\nabla}_{\delta_{\text{time}}}$ means an integrable deformation 
of $\tilde{\nabla}_{\mathrm{DL, ext}}^{(1)}$.
\end{Def}
Construction of horizontal lifts of $\tilde{\nabla}_{\mathrm{DL, ext}}^{(1)}$
is discussed in 
Section \ref{2020.1.24.22.25}, 
Section \ref{2020.1.25.8.11}, 
and Section \ref{2020.2.7.11.34}.

\subsection{Solutions of 
$d  +\Omega_{(\boldsymbol{t}_{0},\boldsymbol{\theta},\boldsymbol{\theta}_0)}^{(n-2)} =0$ 
at the apparent singularities}\label{2021.4.29.12.41}

Since $q_j$ ($j=1,2,\ldots,n-3$) are apparent singularities, 
then we have the following lemma:
\begin{Lem}\label{2019.12.30.21.53}
\textit{
For each $j \in \{ 1,\ldots,n-3 \}$,
the equation $(d  +
 \Omega_{(\boldsymbol{t}_{0},\boldsymbol{\theta},\boldsymbol{\theta}_0)}^{(n-2)}) \Psi =0$ 
has a solution
$\psi_{q_j}= \Phi_{q_j} \Xi_{q_j}(x)\Lambda_{q_j}(x)$
at $q_j$.
Here
\begin{equation}\label{2021.4.27.16,01}
\begin{aligned}
&\Phi_{q_j}:=
\begin{pmatrix}
1 &0 \\
p_j &1 
\end{pmatrix}, \qquad 
\Lambda_{q_j}(x):=
\begin{pmatrix}
1 & 0 \\
0 &x-q_j
\end{pmatrix},  \text{ and}, \\
&\Xi_{q_j}(x) := 
\begin{pmatrix}
1 & 0\\
0 & 1
\end{pmatrix}+
\sum_{s=1}^{\infty} 
\begin{pmatrix}
(\xi^{q_j}_s)_{11} & (\xi^{q_j}_s)_{12} \\
(\xi^{q_j}_s)_{21} & (\xi^{q_j}_s)_{22}
\end{pmatrix}
(x-q_j)^s,
\end{aligned}
\end{equation}
where $(\xi^{q_j}_1)_{11} = -\frac{p_j}{P(q_j)}$,
$(\xi^{q_j}_1)_{12}= -\frac{1}{2P(q_j)}$,
$(\xi^{q_j}_1)_{21} = 0$, and
\begin{equation*}
(\xi^{q_j}_1)_{22}= \frac{p_j}{P(q_j)} - \sum_{i=1}^{\nu} \frac{D_i(q_j)}{(q_j-t_i)^{n_i}}
+ \sum_k \frac{1}{q_j - q_k}- D_{\infty} (q_j) .
\end{equation*}
The solution 
$\psi_{q_j}$
has converging entries.
}
\end{Lem}

\begin{proof}

The connection matrix 
$\Omega^{(n-2)}
_{(\boldsymbol{t}_{0},\boldsymbol{\theta},\boldsymbol{\theta}_0)}$
has the following description at $q_j$ by \eqref{2021.4.18.13.06}:
\begin{equation*}
\Omega^{(n-2)}
_{(\boldsymbol{t}_{0},\boldsymbol{\theta},\boldsymbol{\theta}_0)}
= \begin{pmatrix}
 0 & 0 \\
 \frac{p_j}{x-q_j} & \frac{-1}{x-q_j} 
\end{pmatrix} dx
+\begin{pmatrix}
 0 & \frac{1}{P(q_j)} \\
c^{(0)}_{q_j}
& d^{(0)}_{q_j}
\end{pmatrix}dx +O(x-q_j).
\end{equation*}
Here we set
$$
\begin{aligned}
c^{(0)}_{q_j}:=&\ 
\sum_{i=1}^{\nu} \frac{C_i(q_j)}{(q_j-t_i)^{n_i}}
+\sum_{k\neq j}\frac{p_j}{q_j-q_k} +\tilde{C}(q_j) +q_j^{n-3}C_{\infty} (q_j)  \text{ and} \\
d^{(0)}_{q_j}:=&\ 
\sum_{i=1}^{\nu} \frac{D_i(q_j)}{(q_j-t_i)^{n_i}}
+\sum_{k\neq j}\frac{-1}{q_j-q_k} +D_{\infty} (q_j).
\end{aligned}
$$
Let $\Phi_{q_j}$ be the matrix in \eqref{2021.4.27.16,01}.
We may check the following equality:
\begin{equation*}
\Phi_{q_j}^{-1} \Omega^{(n-2)}
_{(\boldsymbol{t}_{0},\boldsymbol{\theta},\boldsymbol{\theta}_0)}
\Phi_{q_j}
= \begin{pmatrix}
 0 & 0 \\
0 & \frac{-1}{x-q_j} 
\end{pmatrix} dx
+\begin{pmatrix}
 \frac{p_j}{P(q_j)} & \frac{1}{P(q_j)} \\
0 & d^{(0)}_{q_j}- \frac{p_j}{P(q_j)}
\end{pmatrix}dx +O(x-q_j).
\end{equation*}
Here the $(2,1)$-entry of this constant term is zero,
since $q_j$ is an apparent singular point.
Since $\Phi_{q_j}^{-1} \Omega^{(n-2)}
_{(\boldsymbol{t}_{0},\boldsymbol{\theta},\boldsymbol{\theta}_0)}
\Phi_{q_j}$ has simple pole at $q_j$ and $q_j$ is an apparent singular point,
there exists a convergent power series $\Xi_{q_j}(x)=\mathrm{id}
+ \Xi_{q_j}^{(1)} (x-q_j) + \cdots$ such that
\begin{equation}\label{2021.4.27.16.15}
(\Phi_{q_j} \Xi_{q_j}(x))^{-1} d(\Phi_{q_j} \Xi_{q_j}(x))+
(\Phi_{q_j} \Xi_{q_j}(x))^{-1} 
\Omega_{(\boldsymbol{t}_{0},\boldsymbol{\theta},\boldsymbol{\theta}_0)}^{(n-2)}
(\Phi_{q_j} \Xi_{q_j}(x)) =
\begin{pmatrix}
 0 & 0 \\
 0 & -1 
\end{pmatrix}\frac{dx}{x-q_j}.
\end{equation}
We calculate the left hand side of \eqref{2021.4.27.16.15}.
The constant term of this left hand side is
$$
 \Xi_{q_j}^{(1)} + \left[ \begin{pmatrix}
 0 & 0 \\
0 & -1 
\end{pmatrix},
  \Xi_{q_j}^{(1)}  \right]+
 \begin{pmatrix}
 \frac{p_j}{P(q_j)} & \frac{1}{P(q_j)} \\
0 & d^{(0)}_{q_j}- \frac{p_j}{P(q_j)}
\end{pmatrix}.
$$
This matrix is a zero matrix.
So we may check that $\Xi_{q_j}^{(1)}$ is determined as follows:
$$
\Xi_{q_j}^{(1)} = 
\begin{pmatrix}
 -\frac{p_j}{P(q_j)} & -\frac{1}{2P(q_j)} \\
(\xi_1^{q_j})_{21} & -d^{(0)}_{q_j}+ \frac{p_j}{P(q_j)}
\end{pmatrix}.
$$
We may determine $(\xi_1^{q_j})_{21}$ freely.
Here we set $(\xi_1^{q_j})_{21}=0$.
\end{proof}

\section{Unramified irregular singularities}

In this section,
we assume that $ I_{\text{ra}} =\emptyset$.
In Section \ref{2020.2.8.10.24}, we define a $2$-form on the fiber
$\mathcal{M}_{\boldsymbol{t}_{0},\boldsymbol{t}_{\text{ra}}}$
by Krichever's formula \cite[Section 5]{Krichever}.
Remark that $\mathcal{M}_{\boldsymbol{t}_{0},\boldsymbol{t}_{\text{ra}}}$
is isomorphic to the moduli space
$\mathfrak{Conn}_{(\boldsymbol{t}_0,\boldsymbol{\theta},\boldsymbol{\theta}_0)}$.
We show that this $2$-form coincides with the symplectic form \eqref{2020.2.7.11.24}.
In Section \ref{2020.1.24.22.25} and Section \ref{2020.1.25.8.11},
we will construct horizontal lifts of $\tilde{\nabla}_{\mathrm{DL, ext}}^{(1)}$.
Let
$\partial / \partial \theta_{l,t_i}^{\pm}$ ($i \in I_{\mathrm{un}}$ and $l=0,1,\ldots,n_i-2$) and
$\partial / \partial t_i$ ($i=3,4\ldots,\nu$)
be the vector fields on $(T_{\boldsymbol{t}})_{\boldsymbol{t}_{\text{ra}}} \times T_{\boldsymbol{\theta}}$.
By the construction of the horizontal lifts,
we have the vector fields 
$\delta^{\mathrm{IMD}}_{\theta_{l,t_i}^{\pm}}$ and 
$\delta^{\mathrm{IMD}}_{t_i}$
on
$\widehat{\mathcal{M}}_{\boldsymbol{t}_{\text{ra}}}\times T_{\boldsymbol{\theta}}$
determined by the integrable deformations
with respect to $\partial / \partial \theta_{l,t_i}^{\pm}$
and $\partial / \partial t_i$,
respectively.
Remark that 
$\widehat{\mathcal{M}}_{\boldsymbol{t}_{\text{ra}}}\times T_{\boldsymbol{\theta}}$ 
is isomorphic to the extended moduli space
$\widehat{\mathfrak{Conn}}_{(\boldsymbol{t}_{\mathrm{ra}},\boldsymbol{\theta}_0)}$.
In Section \ref{2020.2.8.10.28},
we define a $2$-form on 
$\widehat{\mathcal{M}}_{\boldsymbol{t}_{\text{ra}}}\times T_{\boldsymbol{\theta}}$
by Krichever's formula.
We show that this $2$-form is the isomonodromy $2$-form.
In Section \ref{2020.2.8.10.30},
we calculate this $2$-form on 
$\widehat{\mathcal{M}}_{\boldsymbol{t}_{\text{ra}}}\times T_{\boldsymbol{\theta}}$
by using Diarra--Loray's global normal form.
Then we obtian an explicit formula of this $2$-form.

We consider the leading coefficient of 
$\Omega^{(n-2)}
_{(\boldsymbol{t}_{0},\boldsymbol{\theta},\boldsymbol{\theta}_0)} $
at $t_i$:
\begin{equation*}
\Omega^{(n-2)}
_{(\boldsymbol{t}_{0},\boldsymbol{\theta},\boldsymbol{\theta}_0)}
= \begin{pmatrix}
0 & \frac{1}{\prod_{j \neq i}(t_i-t_j)^{n_j}} \\
\theta_{0,t_i}^+ \theta_{0,t_i}^-\prod_{j \neq i}(t_i-t_j)^{n_j} & \theta_{0,t_i}^++ \theta_{0,t_i}^-
\end{pmatrix}\frac{dx_{t_i}}{x_{t_i}^{n_i}}
 + [\text{ higher order terms }].
\end{equation*}
Remark that this leading coefficient at $t_i$
is independent of $\{ (q_j ,p_j) \}_{j=1,2,\ldots,n-3}$.
We fix $\Phi_i \in \mathrm{GL}(2,\mathbb{C})$ so that 
\begin{equation*}
\Phi_i^{-1} 
\Omega^{(n-2)}
_{(\boldsymbol{t}_{0},\boldsymbol{\theta},\boldsymbol{\theta}_0)} 
\Phi_i
= \begin{pmatrix}
\theta_{0,t_i}^+ & 0 \\
0 & \theta_{0,t_i}^-
\end{pmatrix}\frac{dx_{t_i}}{x_{t_i}^{n_i}}
 + [\text{ higher order terms }].
\end{equation*}
We call the matrix $\Phi_i$ a {\it compatible framing at }$t_i$.
If we have another $\Phi_i'$ 
such that the leading coefficient matrix of
 $(\Phi_i')^{-1} \Omega^{(n-2)}
_{(\boldsymbol{t}_{0},\boldsymbol{\theta},\boldsymbol{\theta}_0)} \Phi'_i$
is diagonal as above,
then there exists a diagonal matrix $C_{t_i}$ such that
$\Phi_i' =\Phi_i C_{t_i}$, since $\theta_{0,t_i}^+-\theta_{0,t_i}^-\neq0$.

\begin{Lem}[{For example \cite[Lemma 3.1]{Krichever}}]\label{2019.12.30.22.09}
\textit{
Assume that $\theta_{0,t_i}^{+}-\theta_{0,t_i}^{-} \neq 0$ if $n_i >1$ 
and $\theta_{0,t_i}^{+}-\theta_{0,t_i}^{-} \notin \mathbb{Z}$ if $n_i =1$.
For a compatible framing $\Phi_i$,
there exist unique
\begin{itemize}
\item $\theta_{l,t_i}^{\pm} \in \Gamma (\mathcal{M}_{\boldsymbol{t}_0,\boldsymbol{t}_{\mathrm{ra}}},
\mathcal{O}_{\mathcal{M}_{\boldsymbol{t}_0,\boldsymbol{t}_{\mathrm{ra}}}})$ 
$(l \geq n_i$ and $i\in I)$,
\item $\sum_{s=1}^{\infty} (\xi^{(i)}_s)_{12} x_{t_i}^s
\in \Gamma (\mathcal{M}_{\boldsymbol{t}_0,\boldsymbol{t}_{\mathrm{ra}}},
\mathcal{O}_{\mathcal{M}_{\boldsymbol{t}_0,\boldsymbol{t}_{\mathrm{ra}}}}) [[x_{t_i}]]$, and
\item $\sum_{s=1}^{\infty} (\xi^{(i)}_s)_{21} x_{t_i}^s 
\in \Gamma (\mathcal{M}_{\boldsymbol{t}_0,\boldsymbol{t}_{\mathrm{ra}}},
\mathcal{O}_{\mathcal{M}_{\boldsymbol{t}_0,\boldsymbol{t}_{\mathrm{ra}}}}) [[x_{t_i}]]$
\end{itemize}
such that 
$\psi_i:= \Phi_i \Xi_i(x_{t_i})
\mathrm{exp} (-\Lambda_i(x_{t_i})) $ 
satisfies the equation $(d  + \Omega^{(n-2)}
_{(\boldsymbol{t}_{0},\boldsymbol{\theta},\boldsymbol{\theta}_0)}) \psi_i =0$ 
 formally at $t_i$.
Here we put
\begin{equation*}
\begin{aligned}
\Lambda_i(x_{t_i}):=&\  
\begin{pmatrix}
 \hat{\lambda}_{i}^+(x_{t_i}) & 0 \\
0 &\hat{\lambda}_{i}^-(x_{t_i})
\end{pmatrix} \\
\Xi_i(x_{t_i}) := &\ 
\begin{pmatrix}
1 & 0\\
0 & 1
\end{pmatrix}+
\sum_{s=1}^{\infty} 
\begin{pmatrix}
0 & (\xi^{(i)}_s)_{12} \\
(\xi^{(i)}_s)_{21} & 0
\end{pmatrix}
x_{t_i}^s 
\end{aligned}
\end{equation*}
where $\hat{\lambda}_i^{\pm}(x_{t_i}):=
\sum_{l=0}^{\infty}\theta_{l,t_i}^{\pm} \int x_{t_i}^{-n_i+l} dx_{t_i}$.
That is, $\psi_i$ is a
formal fundamental matrix solution at $t_i$.
}
\end{Lem}

\begin{Rem}\label{2021.4.30.17.15}
By the equations in \eqref{2021.4.30.13.59},
we have that
the polynomials $C_i$ and $D_i$ ($i=1,2,\ldots, \nu, \infty$) 
in \eqref{2021.4.18.13.06} are independent of 
the parameters $\{ (q_j,p_j) \}_{j=1,2,\ldots,n-3}$ of 
$\mathcal{M}_{\boldsymbol{t}_{0},\boldsymbol{t}_{\text{ra}}}$.
If we take a compatible framing $\Phi_i$
so that $\Phi_i$ is independent of 
$\{ (q_j,p_j) \}_{j=1,2,\ldots,n-3}$, 
then the coefficients of the formal power series $\Phi_i \Xi_i(x_{t_i})$
until the $(x-t_i)^{n_i-1}$-term are independency of $\{ (q_j,p_j) \}_{j=1,2,\ldots,n-3}$.
This independency is the assumption of Lemma \ref{2020.1.20.21.24} (below).
We will use this fact for the calculation of Hamiltonians and the isomonodromy 2-form.
\end{Rem}

\subsection{Symplectic structure}\label{2020.2.8.10.24}

\begin{Def}[{\cite[Section 5]{Krichever} and \cite[Formula (3.16), p.306]{DM}}]\label{2022.4.6.16.3}
Let $\delta_1$ and $\delta_2$ be vector fields on 
$\mathcal{M}_{\boldsymbol{t}_{0},\boldsymbol{t}_{\text{ra}}}\subset
\mathrm{Sym}^{(n-3)}(\mathbb{C}^2)$,
which is isomorphic to the moduli space
$\mathfrak{Conn}_{(\boldsymbol{t}_0,\boldsymbol{\theta},\boldsymbol{\theta}_0)}$.
We fix a formal fundamental matrix solution $\psi_i$ of 
$(d+\Omega^{(n-2)}
_{(\boldsymbol{t}_{0},\boldsymbol{\theta},\boldsymbol{\theta}_0)})\psi_i=0$ at $x=t_i$
as in Lemma \ref{2019.12.30.22.09}.
Moreover,
we fix a fundamental matrix solution $\psi_{q_j}$ of 
$(d+\Omega^{(n-2)}
_{(\boldsymbol{t}_{0},\boldsymbol{\theta},\boldsymbol{\theta}_0)})\psi_{q_j}=0$ at $x=q_j$
as in Lemma \ref{2019.12.30.21.53}.
We set 
$$
  \delta (\Omega^{(n-2)}
_{(\boldsymbol{t}_{0},\boldsymbol{\theta},\boldsymbol{\theta}_0)}) 
\wedge \delta(\psi_i)\psi_i^{-1} :=
 \delta_1 (\Omega^{(n-2)}
_{(\boldsymbol{t}_{0},\boldsymbol{\theta},\boldsymbol{\theta}_0)}) 
\delta_2(\psi_i)\psi_i^{-1} - \delta_1(\psi_i) \psi_i^{-1}
\delta_2 (\Omega^{(n-2)}
_{(\boldsymbol{t}_{0},\boldsymbol{\theta},\boldsymbol{\theta}_0)})
$$
and 
$$
  \delta (\Omega^{(n-2)}
_{(\boldsymbol{t}_{0},\boldsymbol{\theta},\boldsymbol{\theta}_0)}) 
\wedge \delta(\psi_{q_j})\psi_{q_j}^{-1} :=
 \delta_1 (\Omega^{(n-2)}
_{(\boldsymbol{t}_{0},\boldsymbol{\theta},\boldsymbol{\theta}_0)}) 
\delta_2(\psi_{q_j})\psi_{q_j}^{-1} 
- \delta_1(\psi_{q_j}) \psi_{q_j}^{-1}
\delta_2 (\Omega^{(n-2)}
_{(\boldsymbol{t}_{0},\boldsymbol{\theta},\boldsymbol{\theta}_0)}).
$$
We define a $2$-form $\omega$ on
$\mathcal{M}_{\boldsymbol{t}_{0},\boldsymbol{t}_{\text{ra}}}$ as
\begin{equation}\label{2020.1.7.16.55}
\begin{aligned}
\omega (\delta_1, \delta_2) :=\ &  
 \frac{1}{2} \sum_{i \in I}   \mathrm{res}_{x=t_i} 
\mathrm{Tr} \left(   \delta (\Omega^{(n-2)}
_{(\boldsymbol{t}_{0},\boldsymbol{\theta},\boldsymbol{\theta}_0)}) 
\wedge \delta(\psi_i)\psi_i^{-1}  \right)
+\frac{1}{2} \sum_{j=1}^{n-3} \mathrm{res}_{x=q_j} 
\mathrm{Tr} \left(   \delta (\Omega^{(n-2)}
_{(\boldsymbol{t}_{0},\boldsymbol{\theta},\boldsymbol{\theta}_0)}) 
\wedge \delta(\psi_{q_j})\psi_{q_j}^{-1}  \right),
\end{aligned}
\end{equation}
where $I:= \{ 1,2,\ldots, \nu, \infty \}$.
\end{Def}

In \cite[Section 5]{Krichever}, it is discussed that this definition is well-defined.
We recall this argument in \cite[Section 5]{Krichever}.
First, we show that
the right hand side of (\ref{2020.1.7.16.55}) is independent of the choice of 
$\psi_{q_j}$.
If we have another solution $\psi'_{q_j}$,
then we have matrix $C_{q_j}$
such that this matrix is independent of parameters on $\mathbb{P}^1$
and $\psi'_{q_j}= \psi_{q_j} C_{q_j}$.
By the Leibniz rule, we have
$$
\begin{aligned}
\mathrm{Tr} (\delta_1 (\Omega^{(n-2)}
_{(\boldsymbol{t}_{0},\boldsymbol{\theta},\boldsymbol{\theta}_0)}) 
\delta_2(\psi'_{q_j})(\psi'_{q_j})^{-1} )
&=\mathrm{Tr}  (\delta_1 (\Omega^{(n-2)}
_{(\boldsymbol{t}_{0},\boldsymbol{\theta},\boldsymbol{\theta}_0)}) 
\delta_2(\psi_{q_j}C_{q_j})(\psi_{q_j}C_{q_j})^{-1} ) \\
&=\mathrm{Tr}( \delta_1 (\Omega^{(n-2)}
_{(\boldsymbol{t}_{0},\boldsymbol{\theta},\boldsymbol{\theta}_0)}) 
\delta_2(\psi_{q_j})\psi_{q_j}^{-1})
+\mathrm{Tr} (\psi_{q_j}^{-1} \delta_1 (\Omega^{(n-2)}
_{(\boldsymbol{t}_{0},\boldsymbol{\theta},\boldsymbol{\theta}_0)}) 
\psi_{q_j}\delta_2(C_{q_j})C_{q_j}^{-1}).
\end{aligned}
$$
We take variations of the both hand side of the equation $d \psi_{q_j} =- \Omega^{(n-2)}
_{(\boldsymbol{t}_{0},\boldsymbol{\theta},\boldsymbol{\theta}_0)} \psi_{q_j} $
with respect of $\delta_1$.
Then we have equalities
\begin{equation}\label{2021.4.18.17.39}
\begin{aligned}
\psi_{q_j}^{-1} \delta_1 (\Omega^{(n-2)}
_{(\boldsymbol{t}_{0},\boldsymbol{\theta},\boldsymbol{\theta}_0)}) 
\psi_{q_j}
&= - \psi_{q_j}^{-1}d (\delta_1(\psi_{q_j}) ) 
-\psi_{q_j}^{-1}\Omega^{(n-2)}
_{(\boldsymbol{t}_{0},\boldsymbol{\theta},\boldsymbol{\theta}_0)}
\delta_1(\psi_{q_j}) \\
&=- \psi_{q_j}^{-1}d (\delta_1(\psi_{q_j}) ) 
+\psi_{q_j}^{-1} d (\psi_{q_j}) \psi_{q_j}^{-1}
\delta_1(\psi_{q_j})\\
&=  - d(\psi_{q_j}^{-1} \delta_1(\psi_{q_j}) ) .
\end{aligned}
\end{equation}
Here the second equality is given by $d \psi_{q_j} =-
\Omega^{(n-2)}
_{(\boldsymbol{t}_{0},\boldsymbol{\theta},\boldsymbol{\theta}_0)} \psi_{q_j} $.
So we have 
\begin{equation*}
\mathrm{Tr} (\delta_1 (\Omega^{(n-2)}
_{(\boldsymbol{t}_{0},\boldsymbol{\theta},\boldsymbol{\theta}_0)}) 
\delta_2(\psi'_{q_j})(\psi'_{q_j})^{-1} )
= \mathrm{Tr}( \delta_1 (\Omega^{(n-2)}
_{(\boldsymbol{t}_{0},\boldsymbol{\theta},\boldsymbol{\theta}_0)}) 
\delta_2(\psi_{q_j})\psi_{q_j}^{-1})
- \mathrm{Tr} (d(\psi_{q_j}^{-1} \delta_1(\psi_{q_j}) )\delta_2(C_{q_j})C_{q_j}^{-1}).
\end{equation*}
Since the solution $\psi_{q_j}$ is holomorphic at $q_j$ and 
$C_{q_j}$ is independent of parameters on $\mathbb{P}^1$,
the residue of the second term of the right hand side is zero.
This fact means that 
the right hand side of (\ref{2020.1.7.16.55}) is independent of the choice of 
$\psi_{q_j}$.

We may check that
the residue of 
$\mathrm{Tr} (\delta (\Omega^{(n-2)}_{(\boldsymbol{t}_{0},\boldsymbol{\theta},\boldsymbol{\theta}_0)}) 
 \wedge \delta(\psi_i)\psi_i^{-1})$ at $\tilde{t}_i$
is well-defined as follows.
We have the following equality:
\begin{equation}\label{2021.4.30.8.17}
\begin{aligned}
\delta(\psi_i)\psi_i^{-1} &= \delta \big(\Phi_i \Xi_i(x_{t_i})
\mathrm{exp} (-\Lambda_i(x_{t_i})) \big) \big(\Phi_i \Xi_i(x_{t_i})
\mathrm{exp} (-\Lambda_i(x_{t_i})) \big)^{-1} \\
&=
\delta \big(\Phi_i \Xi_i(x_{t_i}) \big) \big(\Phi_i \Xi_i(x_{t_i}) \big)^{-1}
-(\Phi_i \Xi_i(x_{t_i}))
\delta \big( \Lambda_i(x_{t_i}) \big)
(\Phi_i \Xi_i(x_{t_i}))^{-1}.
\end{aligned}
\end{equation}
Since $\theta^{\pm}_{n_i-1,t_i}$ is constant on $\mathcal{M}_{\boldsymbol{t}_{0},\boldsymbol{t}_{\text{ra}}}$,
$\delta(\theta^{\pm}_{n_i-1,t_i})=0$. 
Then $\delta(\theta^{\pm}_{n_i-1,t_i} \int x_{t_i}^{-1} dx_{t_i} )=\delta(c)$.
Here $c$ is an integration constant.
If we fix integration constants on $\Lambda_i(x_{t_i})$,
then we can take 
the residue of 
$\mathrm{Tr} (   \delta (\Omega^{(n-2)}
_{(\boldsymbol{t}_{0},\boldsymbol{\theta},\boldsymbol{\theta}_0)}) 
\wedge \delta(\psi_i)\psi_i^{-1}  )$ at $\tilde{t}_i$.
We may check that $\mathrm{res}_{x=t_i} 
\mathrm{Tr} (   \delta (\Omega^{(n-2)}
_{(\boldsymbol{t}_{0},\boldsymbol{\theta},\boldsymbol{\theta}_0)}) 
\wedge \delta(\psi_i)\psi_i^{-1}  )$
is independent of the choice of the integration constant as follows.
We take other integration constants
and a formal solution $\psi_i'$ is given for the integration constants.
There exists a diagonal matrix $C_{t_i}$ such that 
$\psi_i'=\psi_i C_{t_i}$
and
$C_{t_i}$ is independent of parameters on $\mathbb{P}^1$.
By the same argument as above,
we have the equality 
\begin{equation}\label{2021.4.15.13.34}
\mathrm{Tr} (\delta_1 (\Omega^{(n-2)}
_{(\boldsymbol{t}_{0},\boldsymbol{\theta},\boldsymbol{\theta}_0)}) 
\delta_2(\psi'_{i})(\psi'_{i})^{-1} )
= \mathrm{Tr}( \delta_1 (\Omega^{(n-2)}
_{(\boldsymbol{t}_{0},\boldsymbol{\theta},\boldsymbol{\theta}_0)}) 
\delta_2(\psi_{i})\psi_{i}^{-1})
- \mathrm{Tr} (d(\psi_{i}^{-1} \delta_1(\psi_{i}) )\delta_2(C_{t_i})C_{t_i}^{-1}).
\end{equation} 
Since $\Lambda_i(x_{t_i})$ and $C_{t_i}$ are diagonal,
we have 
$\mathrm{exp} (-\Lambda_i(x_{t_i}))\delta_2(C_{t_i})C_{t_i}^{-1}
 \mathrm{exp} (-\Lambda_i(x_{t_i}))^{-1} =\delta_2(C_{t_i})C_{t_i}^{-1}$.
We calculate the second term of the left hand side of \eqref{2021.4.15.13.34} as follows.
\begin{equation}\label{2021.4.17.19.11}
\begin{aligned}
&\mathrm{Tr} \Big(d(\psi_{i}^{-1} \delta_1(\psi_{i}) )\delta_2(C_{t_i})C_{t_i}^{-1} \Big)\\
&=\mathrm{Tr} \Big(
d\Big(\mathrm{exp} (-\Lambda_i(x_{t_i}))^{-1}
(\Phi_i \Xi_i(x_{t_i}))^{-1} \delta_1(\Phi_i \Xi_i(x_{t_i})) \mathrm{exp} (-\Lambda_i(x_{t_i})) \Big)
\delta_2(C_{t_i})C_{t_i}^{-1}\Big) \\
& \qquad + \mathrm{Tr} \Big(
d\Big(\mathrm{exp} (-\Lambda_i(x_{t_i}))^{-1}
 \delta_1(  \mathrm{exp} (-\Lambda_i(x_{t_i}))) \Big)
\delta_2(C_{t_i})C_{t_i}^{-1}\Big)\\
&=\mathrm{Tr} \Big(
d\Lambda_i(x_{t_i})
(\Phi_i \Xi_i(x_{t_i}))^{-1} \delta_1(\Phi_i \Xi_i(x_{t_i}))  
\delta_2(C_{t_i})C_{t_i}^{-1}\Big)
- \mathrm{Tr} \Big(
(\Phi_i \Xi_i(x_{t_i}))^{-1} \delta_1(\Phi_i \Xi_i(x_{t_i})) 
d\Lambda_i(x_{t_i})  
\delta_2(C_{t_i})C_{t_i}^{-1}\Big) \\
& \qquad + 
\mathrm{Tr} \Big(
d\Big(
(\Phi_i \Xi_i(x_{t_i}))^{-1} \delta_1(\Phi_i \Xi_i(x_{t_i}))  \Big)
\delta_2(C_{t_i})C_{t_i}^{-1}\Big)
+\mathrm{Tr} \Big(
d\Big( \delta_1(  -\Lambda_i(x_{t_i}))\Big)
\delta_2(C_{t_i})C_{t_i}^{-1}\Big) \\
&=\mathrm{Tr} \Big(
d\Big(
(\Phi_i \Xi_i(x_{t_i}))^{-1} \delta_1(\Phi_i \Xi_i(x_{t_i}))  \Big)
\delta_2(C_{t_i})C_{t_i}^{-1}\Big)
+\mathrm{Tr} \Big(
d \Big(  \delta_1(  -\Lambda_i(x_{t_i}) ) \Big)
\delta_2(C_{t_i})C_{t_i}^{-1}\Big).
\end{aligned}
\end{equation}
The residue parts of
 $d( (\Phi_i \Xi_i(x_{t_i}))^{-1} \delta_1(\Phi_i \Xi_i(x_{t_i})))$ 
and
$d(\delta_1(  -\Lambda_i(x_{t_i})))$ vanish.
Since $\delta_2(C_{t_i})C_{t_i}^{-1}$ is independent of parameters on $\mathbb{P}^1$,
the residues of the formal meromorphic differentials 
of the last line of \eqref{2021.4.17.19.11} at $t_i$ are zero.
Then we have that $\mathrm{res}_{x=t_i} 
\mathrm{Tr} (   \delta (\Omega^{(n-2)}
_{(\boldsymbol{t}_{0},\boldsymbol{\theta},\boldsymbol{\theta}_0)}) 
\wedge \delta(\psi_i)\psi_i^{-1}  )$
is independent of the choice of the integration constant.
Finally, the residue of 
$\mathrm{Tr} \left(   \delta (\Omega^{(n-2)}
_{(\boldsymbol{t}_{0},\boldsymbol{\theta},\boldsymbol{\theta}_0)}) 
\wedge \delta(\psi_i)\psi_i^{-1}  \right)$ at $\tilde{t}_i$
is well-defined.

Next we show that the right hand side of (\ref{2020.1.7.16.55}) is independent of the choice of 
a formal solution $\psi_i$.
Let $C_{t_i} (x_{t_i})$ be the following diagonal matrix:
$$ 
\begin{pmatrix}
c_{t_i,11}(x_{t_i}) & 0 \\
0 & c_{t_i,22}(x_{t_i})
\end{pmatrix}
=
\begin{pmatrix}
c_{t_i,11}^{(0)} & 0 \\
0 & c_{t_i,22}^{(0)}
\end{pmatrix}
 + \begin{pmatrix}
c_{t_i,11}^{(1)} & 0 \\
0 & c_{t_i,22}^{(1)}
\end{pmatrix}
 x_{t_i} + 
 \begin{pmatrix}
c_{t_i,11}^{(2)} & 0 \\
0 & c_{t_i,22}^{(2)}
\end{pmatrix} x_{t_i}^2 + \cdots.
$$
We define $\Xi'(x_{t_i})$ and $\Lambda_i' (x_{t_i})$ by
$$
\begin{aligned}
 \Xi'(x_{t_i}) =&\  \Xi(x_{t_i})C_{t_i} (x_{t_i}),\\
\Lambda_i' (x_{t_i}) =&\  \Lambda_i (x_{t_i}) + 
\begin{pmatrix}
\int c_{t_i,11}(x_{t_i})^{-1}d(c_{t_i,11}(x_{t_i})) & 0 \\
0 & \int c_{t_i,22} (x_{t_i})^{-1}d(c_{t_i,22}(x_{t_i}))
\end{pmatrix}.
\end{aligned}
$$
Then we have another formal fundamental matrix solution
$\psi_i' = \Phi_i \Xi'(x_{t_i}) \exp (-  \Lambda_i' (x_{t_i}))$.
There exists a diagonal matrix $C_{t_i}$ such that 
$\psi_i'=\psi_i C_{t_i}$
and
$C_{t_i}$ is independent of parameters on $\mathbb{P}^1$.
By the same argument as above, we have
\begin{equation*}
\mathrm{Tr} (\delta_1 (\Omega^{(n-2)}
_{(\boldsymbol{t}_{0},\boldsymbol{\theta},\boldsymbol{\theta}_0)}) 
\delta_2(\psi'_{i})(\psi'_{i})^{-1} )
= \mathrm{Tr}( \delta_1 (\Omega^{(n-2)}
_{(\boldsymbol{t}_{0},\boldsymbol{\theta},\boldsymbol{\theta}_0)}) 
\delta_2(\psi_{i})\psi_{i}^{-1}).
\end{equation*} 
Then we obtain that
the right hand side of (\ref{2020.1.7.16.55}) is independent of the choice of 
$\psi_i$.

\begin{Thm}\label{2020.1.21.16.17}
\textit{
Let $\omega$ be the $2$-form
on $\mathcal{M}_{\boldsymbol{t}_{0},\boldsymbol{t}_{\text{ra}}}$ 
defined by \eqref{2020.1.7.16.55} in Definition \ref{2022.4.6.16.3}.
The $2$-form $\omega$ coincides with
\begin{equation*}
\sum_{j=1}^{n-3} d \left(\frac{p_j}{P(q_j)} \right) \wedge dq_j.
\end{equation*}
}
\end{Thm}

\begin{proof}
Recall that $\omega (\delta_1, \delta_2)$ is
\begin{equation*}
 \frac{1}{2} \sum_{i \in I}   \mathrm{res}_{x=t_i} 
\mathrm{Tr} \left(   \delta (\Omega^{(n-2)}
_{(\boldsymbol{t}_{0},\boldsymbol{\theta},\boldsymbol{\theta}_0)}) 
\wedge \delta(\psi_i)\psi_i^{-1}  \right)
+\frac{1}{2} \sum_{j=1}^{n-3} \mathrm{res}_{x=q_j} 
\mathrm{Tr} \left(   \delta (\Omega^{(n-2)}
_{(\boldsymbol{t}_{0},\boldsymbol{\theta},\boldsymbol{\theta}_0)}) 
\wedge \delta(\psi_{q_j})\psi_{q_j}^{-1}  \right).
\end{equation*}
We calculate the residue of
$\mathrm{Tr} \left( \delta (\Omega^{(n-2)}
_{(\boldsymbol{t}_{0},\boldsymbol{\theta},\boldsymbol{\theta}_0)}) 
\wedge \delta(\psi_{q_j})\psi_{q_j}^{-1}    \right)$ at $x=q_j$.
For this purpose, 
first, we calculate $\delta (\Omega^{(n-2)}
_{(\boldsymbol{t}_{0},\boldsymbol{\theta},\boldsymbol{\theta}_0)})$ 
around $x=q_j$ as follows.
The connection matrix 
$\Omega^{(n-2)}
_{(\boldsymbol{t}_{0},\boldsymbol{\theta},\boldsymbol{\theta}_0)}$
has the following description at $q_j$ by \eqref{2021.4.18.13.06}:
\begin{equation}\label{2021.4.23.9.54}
\Omega^{(n-2)}
_{(\boldsymbol{t}_{0},\boldsymbol{\theta},\boldsymbol{\theta}_0)}
= \begin{pmatrix}
 0 & 0 \\
 \frac{p_j}{x-q_j} & \frac{-1}{x-q_j} 
\end{pmatrix} dx
+\begin{pmatrix}
 0 & b_0' \\
c_0' & d_0' 
\end{pmatrix}dx.
\end{equation}
Here $b_0'$, $c_0'$, and $d_0'$ are holomorphic at $x=q_j$.
Since $t_i$ and $\theta_{l,t_i}^{\pm}$ ($i \in I, 0\leq l \leq n_i -1$)
are constants on $\mathcal{M}_{\boldsymbol{t}_{0},\boldsymbol{t}_{\text{ra}}}$,
we have $\delta(t_i)=0$ and $\delta(\theta_{l,t_i}^{\pm})=0$.
By $\delta(t_i)=0$, we have $\delta(b_0' )=0$.
By \eqref{2021.4.30.13.59},
we have $\delta(D_i)=0$ for $i=1,2,\ldots,\nu, \infty$.
We take the variation
$\delta (\Omega^{(n-2)}
_{(\boldsymbol{t}_{0},\boldsymbol{\theta},\boldsymbol{\theta}_0)})$
 of $\Omega^{(n-2)}
_{(\boldsymbol{t}_{0},\boldsymbol{\theta},\boldsymbol{\theta}_0)}$
associated to $\delta$:
\begin{equation*}
\delta (\Omega^{(n-2)}
_{(\boldsymbol{t}_{0},\boldsymbol{\theta},\boldsymbol{\theta}_0)})=
\begin{pmatrix}
0& 0 \\
\delta(c_0) & \delta(d_0)
\end{pmatrix} dx,
\end{equation*}
where
\begin{equation*}
\begin{aligned}
\delta(c_0)&=
\frac{p_j \delta (q_j)}{(x-q_j)^2}
+\frac{\delta(p_j)}{(x-q_j)}
 + O(x-q_j)^0 \text{ and }\\
\delta(d_0)&=
-\frac{\delta(q_j)}{(x-q_j)^2}
-\sum_{k\neq j}\frac{\delta(q_k)}{(q_j-q_k)^2} + O(x-q_j). 
\end{aligned}
\end{equation*}
Second, we consider $\delta(\psi_{q_j})\psi_{q_j}^{-1}$. 
We have
\begin{equation}\label{2021.4.23.9.37}
\delta(\psi_{q_j})\psi_{q_j}^{-1} 
= \delta(\Phi_{q_j} \Xi_{q_j}(x) )(\Phi_{q_j} \Xi_{q_j}(x))^{-1}+ 
(\Phi_{q_j} \Xi_{q_j}(x) ) 
\begin{pmatrix}
0 & 0 \\
0 &\frac{-\delta(q_j)}{x-q_j}
\end{pmatrix}
 (\Phi_{q_j} \Xi_{q_j}(x))^{-1}.
\end{equation}
By using Lemma \ref{2019.12.30.21.53},
we have the following equality:
$$
\begin{aligned}
\Phi_{q_j}\Xi_{q_j}(x)
&= \begin{pmatrix}
1 & 0 \\
p_j & 1 
\end{pmatrix}
+
\begin{pmatrix}
-\frac{p_j}{P(q_j)} & -\frac{1}{2P(q_j)} \\
-\frac{p_j^2}{P(q_j)} & \frac{p_j}{2P(q_j)} - \sum_{i=1}^{\nu} \frac{D_i(q_j)}{(q_j-t_i)^{n_i}}
+ \sum_{k\neq j} \frac{1}{q_j - q_k}- D_{\infty} (q_j)
\end{pmatrix}(x-q_j) \\
&\qquad + O(x-q_j)^2.
\end{aligned}
$$
By this description of $\Phi_{q_j}\Xi_{q_j}(x)$, we may check that 
the constant term of the expansion of $\delta(\Phi_{q_j} \Xi_{q_j}(x))$ at $q_j$ 
has following description:
$$
\begin{pmatrix}
0 & 0 \\
\delta(p_j) & 0
\end{pmatrix} - \delta(q_j)
\begin{pmatrix}
-\frac{p_j}{P(q_j)} & -\frac{1}{2P(q_j)} \\
-\frac{p_j^2}{P(q_j)} &*
\end{pmatrix} 
$$
and the coefficient of
the $(x-q_j)$-term of the expansion of $\delta(\Phi_{q_j} \Xi_{q_j}(x))$ 
has following description:
$$
\begin{pmatrix}
* & 0 \\
* & \frac{\delta(p_j)}{2P(q_j)} + \sum_{k \neq j} \frac{\delta(q_k)}{(q_j - q_k)^2}
\end{pmatrix} - \delta(q_j)
\begin{pmatrix}
* & * \\
* &*
\end{pmatrix} .
$$
Here we put the entries having $\delta(q_j)$
together in the second matrices.
Moreover,
we may check that 
the constant term of the expansion of $(\Phi_{q_j} \Xi_{q_j}(x))^{-1}$ at $q_j$ 
is $\begin{pmatrix}
1 & 0 \\
-p_j & 1
\end{pmatrix}$
and the coefficient of
the $(x-q_j)$-term of the expansion of $(\Phi_{q_j} \Xi_{q_j}(x))^{-1}$ 
has following description:
$$
-
\begin{pmatrix}
-\frac{p_j}{2P(q_j)} & -\frac{1}{2P(q_j)} \\
* & \frac{p_j}{P(q_j)} - \sum_{i=1}^{\nu} \frac{D_i(q_j)}{(q_j-t_i)^{n_i}}
+ \sum_{k\neq j} \frac{1}{q_j - q_k}- D_{\infty} (q_j)
\end{pmatrix} .
$$
By the calculation of $\delta(\Phi_{q_j} \Xi_{q_j}(x))$ and
$(\Phi_{q_j} \Xi_{q_j}(x))^{-1}$,
we may show that 
$\delta(\Phi_{q_j} \Xi_{q_j}(x)) (\Phi_{q_j} \Xi_{q_j}(x))^{-1}$ is
\begin{equation}\label{2021.4.23.9.52}
\begin{pmatrix}
* & \frac{\delta(q_j)}{2P(q_j)} \\
*& *
\end{pmatrix}
+ 
\begin{pmatrix}
* &  f^{(1)}_{12} \delta(q_j) \\
* &  \frac{\delta(p_j)}{P(q_j)} 
+ \sum_k \frac{ \delta(q_k)}{(q_j - q_k)^2}
- f^{(1)}_{22} \delta(q_j)
\end{pmatrix}(x-q_j)+ O(x-q_j)^2,
\end{equation}
where $f_{12}^{(1)}$ and $f_{22}^{(1)}$ are rational functions 
on $\mathcal{M}_{\boldsymbol{t}_{0},\boldsymbol{t}_{\text{ra}}}$.
We consider the second term of \eqref{2021.4.23.9.37}.
We may show that 
\begin{equation}\label{2021.4.23.9.51}
\begin{aligned}
&(\Phi_{q_j} \Xi_{q_j}(x) ) 
\begin{pmatrix}
0 & 0 \\
0 &\frac{-\delta(q_j)}{x-q_j}
\end{pmatrix}
 (\Phi_{q_j} \Xi_{q_j}(x))^{-1}\\
 &=
\frac{\begin{pmatrix}
0 & 0 \\
*& - \delta(q_j)
\end{pmatrix}}{x-q_j}
+ 
\begin{pmatrix}
* &  \frac{ \delta(q_j)}{2P(q_j)} \\
* & *
\end{pmatrix}+ 
\begin{pmatrix}
* &  g^{(1)}_{12} \delta(q_j) \\
* &  g^{(1)}_{22} \delta(q_j)
\end{pmatrix}(x-q_j)
+O(x-q_j)^2,
\end{aligned}
\end{equation}
where $g^{(1)}_{12}$ and $g^{(1)}_{22}$ are rational functions 
on $\mathcal{M}_{\boldsymbol{t}_{0},\boldsymbol{t}_{\text{ra}}}$.
By \eqref{2021.4.23.9.52} and \eqref{2021.4.23.9.51}, we have 
\begin{equation*}
\begin{aligned}
\delta(\psi_{q_j})\psi_{q_j}^{-1} 
&= 
\frac{\begin{pmatrix}
0 & 0 \\
*& -\delta(q_j)
\end{pmatrix}}{x-q_j}
+ 
\begin{pmatrix}
* &  \frac{ \delta(q_j)}{P(q_j)} \\
* &  *
\end{pmatrix}+ 
\begin{pmatrix}
* &  (g^{(1)}_{12} +f^{(1)}_{12}) \delta(q_j) \\
* & \frac{\delta(p_j)}{P(q_j)} 
+ \sum_k \frac{ \delta(q_k)}{(q_j - q_k)^2}
+(g^{(1)}_{22} - f^{(1)}_{22}) \delta(q_j)
\end{pmatrix}(x-q_j)\\
&\qquad +O(x-q_j)^2.
\end{aligned}
\end{equation*}
By \eqref{2021.4.23.9.54} and this equality,
we have 
$$
\begin{aligned}
&\mathrm{res}_{x=q_j}\mathrm{Tr} \left( \delta_1 (\Omega^{(n-2)}
_{(\boldsymbol{t}_{0},\boldsymbol{\theta},\boldsymbol{\theta}_0)})
\delta_2(\psi_{q_j})\psi_{q_j}^{-1}  \right)\\
&= \frac{ \delta_1(p_j) \delta_2(q_j)}{P(q_j)}
- \frac{ \delta_1(q_j) \delta_2(p_j)}{P(q_j)}
+\sum_{k\neq j}\frac{\delta_1(q_k)\delta_2(q_j)}{(q_j-q_k)^2}
-\sum_{k\neq j}\frac{\delta_1(q_j)\delta_2(q_k)}{(q_j-q_k)^2} \\
&\qquad+ \left( p_j(g^{(1)}_{12} +f^{(1)}_{12})
-(g^{(1)}_{22} - f^{(1)}_{22}) \right) \delta_1(q_j)\delta_2(q_j).
\end{aligned}
$$
Since $\sum_{j=1}^{n-3}
\sum_{k\neq j}\frac{\delta_1(q_k)\delta_2(q_j)-\delta_1(q_j)\delta_2(q_k)}{(q_j-q_k)^2}=0$,
we have
\begin{equation*}
\mathrm{res}_{x=q_j} 
\mathrm{Tr} \left( \delta (\Omega^{(n-2)}
_{(\boldsymbol{t}_{0},\boldsymbol{\theta},\boldsymbol{\theta}_0)})
\wedge \delta(\psi_{q_j})\psi_{q_j}^{-1}   
\right) 
= 
\frac{2\delta_1(p_j)\delta_2(q_j)}{P(q_j)}
-\frac{2\delta_2(p_j)\delta_1(q_j)}{P(q_j)} .
\end{equation*}

Next we calculate the residue of
$\mathrm{Tr} \left( \delta (\Omega^{(n-2)}
_{(\boldsymbol{t}_{0},\boldsymbol{\theta},\boldsymbol{\theta}_0)}) 
\wedge\delta(\psi_{i})\psi_{i}^{-1} 
 \right)$ at $x=t_i$.
First, we consider the expansion of
$\delta (\Omega^{(n-2)}
_{(\boldsymbol{t}_{0},\boldsymbol{\theta},\boldsymbol{\theta}_0)})$ 
at $x=t_i$.
Since 
$\delta(C_i)=\delta(D_i)=0$ 
for $i=1,2,\ldots,\nu, \infty$, we have
$\delta(c_0)= O(x_{t_i}^0)$ and 
$\delta(d_0)= O(x_{t_i}^0)$.
Second, we consider $\delta(\psi_{i})\psi_{i}^{-1}$. 
By Lemma \ref{2019.12.30.22.09}, we have
\begin{equation*}
\delta(\psi_{i})\psi_{i}^{-1} 
= \delta(\Phi_{i} \Xi_{i}(x_{t_i}) )(\Phi_{i} \Xi_{i}(x_{t_i}))^{-1}+ 
(\Phi_{i} \Xi_{i}(x_{t_i}) ) 
\begin{pmatrix}
-\delta (\hat{\lambda}_i^+(x_{t_i})) & 0 \\
0 &-\delta (\hat{\lambda}_i^-(x_{t_i}))
\end{pmatrix}
 (\Phi_{i} \Xi_{i}(x))^{-1}.
\end{equation*}
Since $\delta (\hat{\lambda}_i^{\pm}(x_{t_i}))=O(x_{t_i})$,
we have that the residue of 
$ \mathrm{Tr} \left( \delta (\Omega^{(n-2)}
_{(\boldsymbol{t}_{0},\boldsymbol{\theta},\boldsymbol{\theta}_0)}) 
\wedge\delta(\psi_{i})\psi_{i}^{-1}   \right)$ at $t_i$ is
zero.
Then 
we obtain
\begin{equation*}
\omega (\delta_1, \delta_2) = 
\sum_{j=1}^{n-3}
\left( \frac{\delta_1(p_j)\delta_2(q_j)}{P(q_j)}
-\frac{\delta_2(p_j)\delta_1(q_j)}{P(q_j)} \right),
\end{equation*}
which means that 
$\omega$ 
coincides with $\sum_{j=1}^{n-3} d\left( \frac{p_j}{P(q_j)} \right) \wedge dq_j$.
\end{proof}

\subsection{Note on the relation to the symplectic structure of the coadjoint orbits}
We apply the argument in \cite[the proof of Theorem 3.3]{DM} 
for our $\omega$.
Let $d+\Omega^0$ be a connection on 
$E_1=\mathcal{O}_{\mathbb{P}^1} \oplus \mathcal{O}_{\mathbb{P}^1}(1)$,
whose polar divisor is $D$.
Remark that
the connection $d+\Omega^0$
 is related to a connection on 
 $E_{n-2}=\mathcal{O}_{\mathbb{P}^1}\oplus \mathcal{O}_{\mathbb{P}^1}(n-2)$
via the transformation (\ref{2020.1.11.22.49}). 
Let $t$ be a component of the divisor $D$.
Choosing a formal coordinate $x_{t}$ near $t$ 
and a trivialization of $E$ on the formal neighborhood of $t$,
we describe $\nabla$ near $t$ by
\begin{equation*}
d + \Omega_0^0 \frac{dx_{t}}{x_{t}^{n_t}} + [\text{ higher order terms }], \quad 
\Omega_0^0 \in \mathfrak{gl}(2,\mathbb{C}).
\end{equation*}
Let $\psi$ be a formal solution at $t$, that is $(d+\Omega^0) \psi=0$.
For $j=1,2$, let $\delta_j(\Omega^0)$ and $\delta_j(\psi)$ be 
the variations of $\Omega^0$ and $\psi$, respectively.
Here $\delta_j$ ($j=1,2$) mean vector fields on 
$\mathcal{M}_{\boldsymbol{t}_{0},\boldsymbol{t}_{\text{ra}}}$.

We define $G_{n_t}$ as $G_{n_i}= \mathrm{GL}(2,\mathbb{C}[x_{t}]/(x_{t}^{n_t}))$.
Let $\mathfrak{g}_{n_t}$ be the Lie algebra of $G_{n_t}$ and
$\mathfrak{g}_{n_t}^*$ be the dual of $\mathfrak{g}_{n_t}$.
We define $\Omega^0_{\le n_t-1}$ and $U^{(j)}_{\le n_t-1}$ ($j=1,2$) as
$\Omega^0 = (\Omega^0_{\le n_t-1} )\cdot x_{t}^{-n_t}  + O(x_{t}^0)$ and
$\delta_{j}(\psi) \psi^{-1} = U^{(j)}_{\le n_t-1} + O(x_{t}^{n_t})$,
respectively.
We identify $\Omega^0_{\le n_t-1}$ and $U^{(j)}_{\le n_t-1}$
as elements of $\mathfrak{g}_{n_t}^*$
by the pairing $\langle X,Y \rangle= \sum_{k=0}^{n_t-1} (X_k Y_{n_t-1-k})$ where 
$X=X_0 + X_1 x_{t}+\cdots+X_{n_t-1}x_{t}^{n_t-1}\in \mathfrak{g}_{n_t}$ and 
$Y=Y_0 + Y_1 x_{t}+\cdots+Y_{n_t-1}x_{t}^{n_t-1}\in \mathfrak{g}_{n_t}$.
Since $\delta_{j} (\Omega^0) = -[\Omega^0, \delta_{j}(\psi) \psi^{-1}] 
- \frac{d}{dx_{t}} (\delta_j(\psi) \psi^{-1}) $
for $j=1,2$, 
we have
\begin{equation}\label{2020.1.11.22.54}
\begin{aligned}
\delta_j(\Omega^0_{\le n_t-1}) = - [ \Omega^0_{\le n_t-1} , U^{(j)}_{\le n_t-1} ].
\end{aligned}
\end{equation}
By this equality, we have the following equality:
\begin{equation}\label{2020.1.11.22.55}
\begin{aligned}
&\frac{1}{2}\mathrm{res}_{x_{t}=0}\mathrm{Tr} \left( \delta_1 (\Omega^0)
\delta_2(\psi) (\psi)^{-1}  
-\delta_1(\psi) (\psi)^{-1} \delta_2 (\Omega^0) \right) \\
&=-\mathrm{Tr} \left\langle \Omega^0_{\le n_t-1} , [ U^{(1)}_{\le n_t-1} ,
U^{(2)}_{\le n_t-1}] \right\rangle.
\end{aligned}
\end{equation}
If we consider the elementary transformation (in other words, the Hecke modification),
we have a connection on the rank $2$ trivial bundle from
the connection $d+\Omega^0$ on $E_1$.
By the equalities (\ref{2020.1.11.22.54}) and (\ref{2020.1.11.22.55}),
we have a relation between $\omega$ and the symplectic form on 
the product of the coadjoint orbits of $G_{n_t}$ for each component $t$ of $D$
(see \cite[Proposition 2.1]{Boalch}).

\subsection{Integrable deformations associated to $T_{\boldsymbol{\theta}}$}\label{2020.1.24.22.25}

First we fix $i \in I$ and $l \in \{ 0,1,\ldots, n_i -2\}$.
Let $\widehat{E}_{1}$ be 
the pull-back of $E_1$ under the projection
$\mathbb{P}^1 \times(\widehat{\mathcal{M}}_{\boldsymbol{t}_{\text{ra}}} \times T_{\boldsymbol{\theta}})
\rightarrow \mathbb{P}^1$.
Let 
$$
\tilde{\nabla}_{\mathrm{DL, ext}}^{(1)}=
\begin{cases}
d+ \widehat\Omega_{(\boldsymbol{t}_{\mathrm{ra}},\boldsymbol{\theta}_0)}^{(1)} 
& \text{ on $U_0\times(\widehat{\mathcal{M}}_{\boldsymbol{t}_{\text{ra}}} \times
 T_{\boldsymbol{\theta}})$} \\
d + G_{1}^{-1} dG_{1}+G_{1}^{-1} 
\widehat\Omega_{(\boldsymbol{t}_{\mathrm{ra}},\boldsymbol{\theta}_0)}^{(1)} \, G_{1}
 & \text{ on $U_\infty \times (\widehat{\mathcal{M}}_{\boldsymbol{t}_{\text{ra}}} 
 \times T_{\boldsymbol{\theta}})$}
\end{cases}
$$
be the family \eqref{2021.4.14.10.36} 
of connections on $\widehat{E}_{1}$.
Let $\theta_{l,t_i}^{\pm}$ be the natural coordinate of 
$(T_{\boldsymbol{t}})_{\boldsymbol{t}_{\text{ra}}} \times T_{\boldsymbol{\theta}}$
and $\partial/\partial \theta_{l,t_i}^{\pm}$ be the vector field on
$(T_{\boldsymbol{t}})_{\boldsymbol{t}_{\text{ra}}} \times T_{\boldsymbol{\theta}}$
associated to $\theta_{l,t_i}^{\pm}$.
We will construct a horizontal lift of $\tilde{\nabla}_{\mathrm{DL, ext}}^{(1)}$
with respect to $\partial/\partial \theta_{l,t_i}^{\pm}$.

We consider diagonalizations of $\tilde{\nabla}_{\mathrm{DL, ext}}^{(1)}$ until 
some degree term at each $\tilde{t}_{i'}$ ($i' \in I$).
 By using the explicit form of $d+ \widehat\Omega_{(\boldsymbol{t}_{\mathrm{ra}},\boldsymbol{\theta}_0)}^{(n-2)} $,
we take a family of compatible framings of
$d+ \widehat\Omega_{(\boldsymbol{t}_{\mathrm{ra}},\boldsymbol{\theta}_0)}^{(n-2)} $
 at $\tilde{t}_{i'}$ for each $i' \in I$.
We denote by $\Phi_{i'}$ this family of compatible framings 
 at $\tilde{t}_{i'}$ for each $i' \in I$.
Let $\Xi_{i'}(x_{t_{i'}})$ be the formal transformation 
of
$d+ \widehat\Omega_{(\boldsymbol{t}_{\mathrm{ra}},\boldsymbol{\theta}_0)}^{(n-2)} $ 
at $\tilde{t}_{i'}$
with respect to $\Phi_{i'}$
appeared in Lemma \ref{2019.12.30.22.09}.
Let $\tilde{G}$ be the matrix defined in \eqref{2021.4.2.20.06}.
We denote 
formal expansion of $\tilde{G}^{-1}\Phi_{i'}\Xi_{i'}(x_{t_{i'}})$ at $x_{t_{i'}}=0$ by\
\begin{equation}\label{2021.4.18.21.46}
\tilde{G}^{-1}\Phi_{i'}\Xi_{i'}(x_{t_{i'}})
= P_{i',0} + P_{i',1} x_{t_{i'}} + P_{i',2} x^2_{t_{i'}} +\cdots.
\end{equation}
Set 
\begin{equation*}
\begin{aligned}
P_{i'} :=&\  P_{i',0} + P_{i',1} x_{t_{i'}} + P_{i',2} x^2_{t_{i'}} +\cdots +
 P_{i',2n_{i'}-1} x^{2n_{i'}-1}_{t_{i'}}
\ \  \text{(for $i'\in I$)}, \ \   \text{and} \\ 
P_{\nu+1} := &\ \mathrm{id}.
\end{aligned}
\end{equation*} 
We take an affine open covering $\{ \hat{U}_{i'}\}_{i' \in I \cup \{ \nu+1 \}}$ of 
$\mathbb{P}^1\times (\widehat{\mathcal{M}}_{\boldsymbol{t}_{\text{ra}}}
 \times T_{\boldsymbol{\theta}})$ 
such that 
\begin{itemize}
\item for $i'\in I$, $\tilde{t}_{i'} \subset \hat{U}_{i'}$,
$\tilde{t}_j \cap \hat{U}_{i'} =\emptyset$ (for any $j\neq i', j \in I$),
and $\sum_{s=1}^{2n_{i'}-1} P_{i',s} x^s_{t_{i'}} $ 
is invertible on each point of $\hat{U}_{i'}$ 
\item for $i'=\nu+1$, $\tilde{t}_j \cap \hat{U}_{i'} =\emptyset$ (for any $j \in I$).
\end{itemize}
Set $\hat{U}_{i'_1,i'_2}:=\hat{U}_{i'_1}\cap \hat{U}_{i'_2}$.

Now we define new trivializations $\{ ( \hat{U}_{i'} ,\hat{\varphi}_{i'}) \}_{i' \in I \cup \{ \nu+1\}}$ of $\widehat{E}_1$.
We denote also by $(U_0\times (\widehat{\mathcal{M}}_{\boldsymbol{t}_{\text{ra}}} \times T_{\boldsymbol{\theta}})
, \varphi^{(1)}_{U_0})$ and 
$(U_{\infty} \times (\widehat{\mathcal{M}}_{\boldsymbol{t}_{\text{ra}}} \times T_{\boldsymbol{\theta}}) ,
\varphi^{(1)}_{U_\infty})$
the trivializations of $\widehat{E}_1$ induced by
\eqref{2022.4.8.9.42}.
We define $\hat{\varphi}_{i'}$ (for $i' \in (I \cup \{ \nu+1\})\setminus \{ \infty\}$) 
by the composition
$$
\hat{\varphi}_{i'} \colon \widehat{E}_1 |_{\hat{U}_{i'}} \xrightarrow{\ \varphi^{(1)}_{U_0}|_{\hat{U}_{i'}}  \ }
\mathcal{O}^{\oplus 2}_{\hat{U}_{i'}} \xrightarrow{\ P^{-1}_{i'}  \ }\mathcal{O}^{\oplus 2}_{\hat{U}_{i'}}.
$$
We define $\hat{\varphi}_{\infty}$ by the composition
$$
\hat{\varphi}_{\infty} \colon \widehat{E}_1 |_{\hat{U}_{\infty}} \xrightarrow{\ \varphi^{(1)}_{U_\infty}|_{\hat{U}_{\infty}}  \ }
\mathcal{O}^{\oplus 2}_{\hat{U}_{\infty}} \xrightarrow{\ P_{\infty}^{-1}  \ }\mathcal{O}^{\oplus 2}_{\hat{U}_{\infty}}.
$$
Then we have new trivializations $\{ ( \hat{U}_{i'} ,\hat{\varphi}_{i'}) \}_{i' \in I \cup \{ \nu+1\}}$ of $\widehat{E}_1$. 
Let $\hat{\Omega}_{i'}$ be
the connection matrix 
of $ \tilde{\nabla}_{\mathrm{DL, ext}}^{(1)}$
under the new trivialization $\hat{\varphi}_{i'}$:
\begin{equation}\label{2021.11.11.23.26}
\begin{aligned}
\hat{\Omega}_{i'}&=
P_{i'}^{-1} dP_{i'}+P_{i'}^{-1} 
\widehat\Omega_{(\boldsymbol{t}_{\mathrm{ra}},\boldsymbol{\theta}_0)}^{(1)} |_{\hat{U}_{i'}}
P_{i'} \quad \text{ for $i' \in (I \cup \{ \nu+1\})\setminus \{ \infty\}$, and}\\
\hat{\Omega}_{\infty} &=
(G_1P_{\infty})^{-1} d (G_1P_{\infty})+(G_1P_{\infty})^{-1} 
\widehat\Omega_{(\boldsymbol{t}_{\mathrm{ra}},\boldsymbol{\theta}_0)}^{(1)} |_{\hat{U}_{\infty}}
(G_1P_{\infty}).
\end{aligned}
\end{equation}
Remark that
$\hat{\Omega}_{i'}$ is diagonal 
until the $x_{t_{i'}}^{n_{i'}-1}$-term for each $i' \in I$.

Now we construct an integrable deformations of 
$  \tilde{\nabla}_{\mathrm{DL, ext}}^{(1)}$. 
For the fixed $i \in I$ and $l$ ($0\le l \le n_i-2$),
we define matrices $B_{\theta_{l,t_i}^{\pm}}(x_{t_i})$ 
by
\begin{equation*}
B_{\theta_{l,t_i}^+}(x_{t_i}):=  \frac{\begin{pmatrix}
\frac{\delta(\theta_{l,t_i}^+)}{-n_i+l+1} & 0 \\
0 & 0
\end{pmatrix}}{x_{t_i}^{n_i -l-1}}  \text{ and }
B_{\theta_{l,t_i}^-}(x_{t_i}):= \frac{ \begin{pmatrix}
0& 0 \\
0 & \frac{\delta(\theta_{l,t_i}^-)}{-n_i+l+1}
\end{pmatrix}}{x_{t_i}^{n_i -l-1}}.
\end{equation*}
For each $i' \in I\cup\{ \nu+1\}$,
we set $(\widehat{E}_1)_{i',\epsilon} =
 \widehat{E}_1|_{\hat{U}_{i'}} \otimes_{\mathbb{C}} \mathbb{C}[\epsilon]$,
 $\hat{U}^{\epsilon}_{i'} = \hat{U}_{i'} \times \mathrm{Spec}\, \mathbb{C}[\epsilon]$, 
 and $\hat{U}^{\epsilon}_{i'_1,i'_2}=
 \hat{U}^{\epsilon}_{i'_1} \cap \hat{U}^{\epsilon}_{i'_2}$.
We define matrices $P^{\epsilon}_i$ and $P^{\epsilon}_{i'}$ by
\begin{equation}\label{2021.4.16.15.17}
P^{\epsilon}_i = P_i ( \mathrm{id}+ \epsilon B_{\theta_{l,t_i}^{\pm}}(x_{t_i}))
\quad \text{and} \quad 
P^{\epsilon}_{i'} = P_{i'} \otimes \mathrm{id}
\quad
 (\text{where $i' \in (I\setminus \{i\} ) \cup \{ \nu+1\}$}),
\end{equation}
respectively.
In the argument below, 
we will replace $P^{\epsilon}_{\infty}$ with $G_1P^{\epsilon}_{\infty}$.
The matrices give isomorphisms
$$ 
 \mathcal{O}^{\oplus 2}_{\hat{U}_{i'_1,i'_2} }
 \otimes_{\mathbb{C}} \mathbb{C}[\epsilon]
\xrightarrow{\  P_{i'_1}^{\epsilon} \ }
\mathcal{O}^{\oplus 2}_{\hat{U}_{i'_1,i'_2}}
 \otimes_{\mathbb{C}} \mathbb{C}[\epsilon]
$$
for each $i_1', i_2' \in I \cup \{ \nu+1\}$.
First, we define
a vector bundle $(\widehat{E}_1)_{\theta_{l,t_i}^{\pm}}^{\epsilon}$ on 
$\mathbb{P}^1\times(\widehat{\mathcal{M}}_{\boldsymbol{t}_{\text{ra}}} \times
 T_{\boldsymbol{\theta} } ) \times \mathrm{Spec}\, \mathbb{C}[\epsilon]$ by
gluing $\{ (\widehat{E}_1)_{i',\epsilon}\}_{{i'} \in I\cup \{ \nu+1\}}$ as follows:
We glue $(\widehat{E}_1)_{i'_1,\epsilon}$ and $(\widehat{E}_1)_{i'_2,\epsilon}$
($i'_1,i'_2 \in I \cup \{\nu+1 \}$) by the composition
$$
(\widehat{E}_1)_{i'_1,\epsilon}|_{\hat{U}^{\epsilon}_{i'_1,i'_2}}
 \xrightarrow{ \hat\varphi_{i'_1}|_{\hat{U}_{i'_1,i'_2}} \otimes 1 } 
 \mathcal{O}^{\oplus 2}_{\hat{U}_{i'_1,i'_2}}
 \otimes_{\mathbb{C}} \mathbb{C}[\epsilon]
\xrightarrow{ \ (P^{\epsilon}_{i'_2})^{-1}P^{\epsilon}_{i'_1}  \  }
\mathcal{O}^{\oplus 2}_{\hat{U}_{i'_1,i'_2}}
 \otimes_{\mathbb{C}} \mathbb{C}[\epsilon]
\xrightarrow{\hat \varphi^{-1}_{i'_2}|_{\hat{U}_{i'_1,i'_2}} \otimes 1 }  
(\widehat{E}_1)_{i'_2,\epsilon}|_{\hat{U}^{\epsilon}_{i'_1,i'_2}} .
$$
By construction, we have 
$(\widehat{E}_1)_{\theta_{l,t_i}^{\pm}}^{\epsilon}
\otimes \mathbb{C}[\epsilon]/(\epsilon) = \widehat{E}_1$.
Second, we define a morphism 
$$
\nabla^{\epsilon}_{\partial/\partial \theta_{l,t_i}^{\pm}}  \colon 
(\widehat{E}_1)_{\theta_{l,t_i}^{\pm}}^{\epsilon}
 \longrightarrow 
 (\widehat{E}_1)_{\theta_{l,t_i}^{\pm}}^{\epsilon}
   \otimes
\tilde\Omega^1_{\partial/\partial \theta_{l,t_i}^{\pm}}
$$
with the Leibniz rule.
Here $\tilde\Omega^1_{\partial/\partial \theta_{l,t_i}^{\pm}}$
is the coherent subsheaf \eqref{2021.4.14.23.20}.
We define $\nabla^{\epsilon}_{i'} $ ($i' \in I \cup \{ \nu+1 \}$) as follows:
\begin{equation}\label{2021.4.16.15.34}
\left\{
\begin{array}{ll}
\nabla^{\epsilon}_{i'}  = \hat{d}+ \hat{\Omega}_{i'}  \quad \text{
for $i' \in (I\setminus \{ i \})\cup \{ \nu+1 \}$,} \\
\nabla^{\epsilon}_{i}=  \hat{d}  + \hat{\Omega}_i + \epsilon\left( \frac{\partial}{\partial x_{t_i}}
 (B_{\theta_{l,t_i}^{\pm}} )dx_{t_i}
+ [ \hat{\Omega}_i, B_{\theta_{l,t_i}^{\pm}}]  \right) 
+B_{\theta_{l,t_i}^{\pm}} d \epsilon.
\end{array}
\right.
\end{equation}
We can consider $ \nabla^{\epsilon}_{i'} $ ($i' \in I \cup \{ \nu+1 \}$) as a morphism
$$
\nabla^{\epsilon}_{i'} \colon 
(\widehat{E}_1)_{\theta_{l,t_i}^{\pm}}^{\epsilon}|_{\hat{U}^{\epsilon}_{i'}} \longrightarrow 
(\widehat{E}_1)_{\theta_{l,t_i}^{\pm}}^{\epsilon}|_{\hat{U}^{\epsilon}_{i'}} \otimes
\tilde\Omega^1_{\partial/\partial \theta_{l,t_i}^{\pm}}|_{\hat{U}^{\epsilon}_{i'}}
$$
by using the trivialization $\hat\varphi_{i'}\otimes 1
 \colon (\widehat{E}_1)_{\theta_{l,t_i}^{\pm}}^{\epsilon}|_{\hat{U}^{\epsilon}_{i'}}
=\widehat{E}_1|_{\hat{U}_{i'}} \otimes_{\mathbb{C}}\mathbb{C}[\epsilon] 
\rightarrow \mathcal{O}^{\oplus 2}_{\hat{U}_{i'}}
 \otimes_{\mathbb{C}} \mathbb{C}[\epsilon] $.
We may glue $ \{ \nabla^{\epsilon}_{i'} \}_{i' \in I \cup \{ \nu+1 \}}$.
Finally, we obtain $\nabla^{\epsilon}_{\partial/\partial \theta_{l,t_i}^{\pm}}$ by
this gluing.
Since $\hat{\Omega}_i$ and $B_{\theta_{l,t_i}^{\pm}}$ are diagonal 
until the $x_{t_i}^{n_i-1}$-terms, 
the negative parts of the relative connections
$\overline{\nabla}^{\epsilon}_{\partial/\partial \theta_{l,t_i}^{+}}$
and
$\overline{\nabla}^{\epsilon}_{\partial/\partial \theta_{l,t_i}^{-}}$
 along the divisor $\tilde{t}_i$ are
\begin{equation*}
\begin{aligned}
&\begin{pmatrix}
\theta_{0,t_i}^{+}& 0 \\
0 & \theta_{0,t_i}^{-}
\end{pmatrix}\frac{dx_{t_i}}{x_{t_i}^{n_i}}+\cdots 
+\begin{pmatrix}
\theta_{l,t_i}^{+} + \epsilon \delta(\theta_{l,t_i}^{+}) & 0 \\
0 & \theta_{l,t_i}^{-}
\end{pmatrix}\frac{dx_{t_i}}{x_{t_i}^{n_i- l}}+\cdots 
+\begin{pmatrix}
 \theta_{n_i-1,t_i}^+ & 0 \\
0 & \theta_{n_i-1,t_i}^-
\end{pmatrix}\frac{dx_{t_i}}{x_{t_i}}, \text{ and }\\
&\begin{pmatrix}
\theta_{0,t_i}^{+}& 0 \\
0 & \theta_{0,t_i}^{-}
\end{pmatrix}\frac{dx_{t_i}}{x_{t_i}^{n_i}}+\cdots 
+\begin{pmatrix}
\theta_{l,t_i}^{+}  & 0 \\
0 & \theta_{l,t_i}^{-}+ \epsilon \delta(\theta_{l,t_i}^{-})
\end{pmatrix}\frac{dx_{t_i}}{x_{t_i}^{n_i- l}}+\cdots 
+\begin{pmatrix}
 \theta_{n_i-1,t_i}^+ & 0 \\
0 & \theta_{n_i-1,t_i}^-
\end{pmatrix}\frac{dx_{t_i}}{x_{t_i}},
\end{aligned}
\end{equation*}
respectively.

Let 
$\widehat{E}^{\epsilon}_1$
be the pull-back of $\widehat{E}_1$ under the projection \eqref{2021.4.14.10.57}.
We consider a short exact sequence 
$$
0 \longrightarrow 
\epsilon \mathcal{H} om (  (\widehat{E}_1)_{\theta_{l,t_i}^{\pm}}^{\epsilon},
\widehat{E}^{\epsilon}_1)
\longrightarrow  \mathcal{H} om (  (\widehat{E}_1)_{\theta_{l,t_i}^{\pm}}^{\epsilon},
\widehat{E}^{\epsilon}_1)
\longrightarrow  \mathcal{E} nd (\hat{E}_1) 
\longrightarrow 0.
$$
Note that 
$\epsilon \mathcal{H} om (  (\widehat{E}_1)_{\theta_{l,t_i}^{\pm}}^{\epsilon},
\widehat{E}^{\epsilon}_1) \cong
 (\epsilon ) \otimes \mathcal{E} nd (\hat{E}_1)$.
Since the bundle type is $\mathcal{O}_{\mathbb{P}^1}\oplus \mathcal{O}_{\mathbb{P}^1}(1)$, 
we have $R^{1} \pi_* ((\epsilon ) \otimes \mathcal{E} nd (\hat{E}_1)) =0$,
which means the rigidity of $\mathcal{O}_{\mathbb{P}^1}\oplus \mathcal{O}_{\mathbb{P}^1}(1)$.
Here $\pi$ is the projection 
$\mathbb{P}^1 \times(\widehat{\mathcal{M}}_{\boldsymbol{t}_{\text{ra}}} \times T_{\boldsymbol{\theta}})
\rightarrow \widehat{\mathcal{M}}_{\boldsymbol{t}_{\text{ra}}} \times T_{\boldsymbol{\theta}}$.
So we have a short exact sequence 
$$
0 \longrightarrow 
\pi_*\left( \epsilon \mathcal{H} om (  (\widehat{E}_1)_{\theta_{l,t_i}^{\pm}}^{\epsilon} ,
\widehat{E}^{\epsilon}_1)\right)
\longrightarrow \pi_*\left( \mathcal{H} om (  (\widehat{E}_1)_{\theta_{l,t_i}^{\pm}}^{\epsilon},
\widehat{E}^{\epsilon}_1) \right)
\longrightarrow \pi_*\left( \mathcal{E} nd (\hat{E}_1) \right)
\longrightarrow 0.
$$
Since $\widehat{\mathcal{M}}_{\boldsymbol{t}_{\text{ra}}} \times T_{\boldsymbol{\theta}}$
is affine, we have that 
$$
\mathrm{Hom} (  (\widehat{E}_1)_{\theta_{l,t_i}^{\pm}}^{\epsilon},
\widehat{E}^{\epsilon}_1)
=\Gamma\left( 
\widehat{\mathcal{M}}_{\boldsymbol{t}_{\text{ra}}} \times T_{\boldsymbol{\theta}},
 \pi_*\left( \mathcal{H} om (  (\widehat{E}_1)_{\theta_{l,t_i}^{\pm}}^{\epsilon},
\widehat{E}^{\epsilon}_1) \right) \right)
\longrightarrow 
\Gamma\left( 
\widehat{\mathcal{M}}_{\boldsymbol{t}_{\text{ra}}} \times T_{\boldsymbol{\theta}},
\pi_*\left( \mathcal{E} nd (\hat{E}_1) \right)
\right) =\mathrm{End} (\hat{E}_1)
$$
is surjective.
Then we have a lift $\varphi_{\Upsilon}
\in \mathrm{Hom} (  (\widehat{E}_1)_{\theta_{l,t_i}^{\pm}}^{\epsilon},
\widehat{E}^{\epsilon}_1)
$ of $\mathrm{id} \in \mathrm{End} (\hat{E}_1)$.
This lift $\varphi_{\Upsilon}$ is an isomorphism
$$
\varphi_{\Upsilon} \colon
  (\widehat{E}_1)_{\theta_{l,t_i}^{\pm}}^{\epsilon}
   \xrightarrow{ \ \cong \ } \widehat{E}^{\epsilon}_1.
$$

We consider a pair 
$(\widehat{E}^{\epsilon}_1 ,(\varphi_{\Upsilon}^{-1})^*\nabla^{\epsilon}_{\partial/\partial \theta_{l,t_i}^{\pm}}  )$
induced by
$((\widehat{E}_1)_{\theta_{l,t_i}^{\pm}}^{\epsilon},
\nabla^{\epsilon}_{\partial/\partial \theta_{l,t_i}^{\pm}})$.
By construction, the pair
$(\widehat{E}^{\epsilon}_1 ,(\varphi_{\Upsilon}^{-1})^*\nabla^{\epsilon}_{\partial/\partial \theta_{l,t_i}^{\pm}}  )$
is a horizontal lift of $\tilde{\nabla}_{\mathrm{DL, ext}}^{(1)}$
with respect to $\partial / \partial \theta_{l,t_i}^{\pm}$.
Let $\overline{(\varphi_{\Upsilon}^{-1})^*\nabla^{\epsilon}_{\partial/\partial \theta_{l,t_i}^{\pm}}}$
be the relative connection induced by
$(\varphi_{\Upsilon}^{-1})^*\nabla^{\epsilon}_{\partial/\partial \theta_{l,t_i}^{\pm}}$.
Since $\varphi_{\Upsilon}$ is holomorphic and invertible
along the pole divisor of $\nabla^{\epsilon}_{\partial/\partial \theta_{l,t_i}^{\pm}}$,
the local formal data of 
$\overline{(\varphi_{\Upsilon}^{-1})^*\nabla^{\epsilon}_{\partial/\partial \theta_{l,t_i}^{\pm}}}$
is the same as in $\nabla^{\epsilon}_{\partial/\partial \theta_{l,t_i}^{\pm}}$.
The family of connections 
$(\widehat{E}^{\epsilon}_1,
\overline{(\varphi_{\Upsilon}^{-1})^*\nabla^{\epsilon}_{\partial/\partial \theta_{l,t_i}^{\pm}}})$
parametrized by $( \widehat{\mathcal{M}}_{\boldsymbol{t}_{\text{ra}}}\times T_{\boldsymbol{\theta}})
\times \mathrm{Spec}\, \mathbb{C}[\epsilon]$ gives a map 
from the base space
 $( \widehat{\mathcal{M}}_{\boldsymbol{t}_{\text{ra}}}\times T_{\boldsymbol{\theta}})
\times \mathrm{Spec}\, \mathbb{C}[\epsilon]$ 
to the moduli space $\widehat{\mathfrak{Conn}}_{(\boldsymbol{t}_{\mathrm{ra}},\boldsymbol{\theta}_0)}$.
By taking composition with $\widehat{\mathrm{App}}$ defined in \eqref{2020.1.6.13.00},
we have a map 
\begin{equation}\label{2021.4.15.11.04}
( \widehat{\mathcal{M}}_{\boldsymbol{t}_{\text{ra}}}\times T_{\boldsymbol{\theta}})
\times \mathrm{Spec}\, \mathbb{C}[\epsilon]
\longrightarrow
\widehat{\mathcal{M}}_{\boldsymbol{t}_{\text{ra}}} \times T_{\boldsymbol{\theta}}.
\end{equation}

\begin{Def}
Then 
we may define the vector field 
on $\widehat{\mathcal{M}}_{\boldsymbol{t}_{\text{ra}}}\times T_{\boldsymbol{\theta}}$
associated 
to the map \eqref{2021.4.15.11.04}.
We denote by $\delta_{\theta_{l,t_i}^{\pm}}^{\mathrm{IMD}}$ this vector field 
on $\widehat{\mathcal{M}}_{\boldsymbol{t}_{\text{ra}}}\times T_{\boldsymbol{\theta}}$.
\end{Def}

Let $f_{\theta_{l,t_i}^{\pm}}^{\mathrm{IMD}} \colon 
( \widehat{\mathcal{M}}_{\boldsymbol{t}_{\text{ra}}}\times T_{\boldsymbol{\theta}})
\times \mathrm{Spec}\, \mathbb{C}[\epsilon]
\longrightarrow
\widehat{\mathcal{M}}_{\boldsymbol{t}_{\text{ra}}} \times T_{\boldsymbol{\theta}}$
be the map induced by the vector field
$\delta_{\theta_{l,t_i}^{\pm}}^{\mathrm{IMD}}$.
We have $\widehat{E}^{\epsilon}_1=
(\mathrm{id} \times f_{\theta_{l,t_i}^{\pm}}^{\mathrm{IMD}})^* \widehat{E}_1$.
We denote by
\begin{equation}\label{2021.4.16.12.11}
\begin{cases}
d+ \widehat\Omega_{(\boldsymbol{t}_{\mathrm{ra}},\boldsymbol{\theta}_0)}^{(1)} 
+\epsilon \delta_{\theta_{l,t_i}^{\pm}}^{\mathrm{IMD}}
(\widehat\Omega_{(\boldsymbol{t}_{\mathrm{ra}},\boldsymbol{\theta}_0)}^{(1)} ) \\
\qquad \qquad 
\text{ on $U_0 \times(\widehat{\mathcal{M}}_{\boldsymbol{t}_{\text{ra}}}\times T_{\boldsymbol{\theta}} )
\times \mathrm{Spec}\, \mathbb{C} [\epsilon]$ } \\
d + G_{1}^{-1} dG_{1}+G_{1}^{-1} 
\widehat\Omega_{(\boldsymbol{t}_{\mathrm{ra}},\boldsymbol{\theta}_0)}^{(1)} 
 \, G_{1}
+\epsilon G_{1}^{-1} 
\delta_{\theta_{l,t_i}^{\pm}}^{\mathrm{IMD}}
(\widehat\Omega_{(\boldsymbol{t}_{\mathrm{ra}},\boldsymbol{\theta}_0)}^{(1)} )
G_{1}  \\
\qquad\qquad 
\text{ on $U_\infty \times(\widehat{\mathcal{M}}_{\boldsymbol{t}_{\text{ra}}}\times T_{\boldsymbol{\theta}} )
\times \mathrm{Spec}\, \mathbb{C} [\epsilon]$ }.
\end{cases}
\end{equation}
the pull-back of $\tilde{\nabla}_{\mathrm{DL, ext}}^{(1)}$ 
under the map $\mathrm{id} \times f_{\theta_{l,t_i}^{\pm}}^{\mathrm{IMD}}$.
Since $(\widehat{E}^{\epsilon}_1 , 
(\mathrm{id} \times f_{\theta_{l,t_i}^{\pm}}^{\mathrm{IMD}})^*\tilde{\nabla}_{\mathrm{DL, ext}}^{(1)} )$
is isomorphic to 
$(\widehat{E}^{\epsilon}_1,
\overline{(\varphi_{\Upsilon}^{-1})^*\nabla^{\epsilon}_{\partial/\partial \theta_{l,t_i}^{\pm}}})$,
we have a lift of 
$(\mathrm{id} \times f_{\theta_{l,t_i}^{\pm}}^{\mathrm{IMD}})^*\tilde{\nabla}_{\mathrm{DL, ext}}^{(1)}$:
\begin{equation}\label{2021.4.16.12.12}
\begin{cases}
\hat{d}+ \widehat\Omega_{(\boldsymbol{t}_{\mathrm{ra}},\boldsymbol{\theta}_0)}^{(1)} 
+\epsilon \delta_{\theta_{l,t_i}^{\pm}}^{\mathrm{IMD}}
(\widehat\Omega_{(\boldsymbol{t}_{\mathrm{ra}},\boldsymbol{\theta}_0)}^{(1)} )
+ \Upsilon_{\theta_{l,t_i}^{\pm}}^{\mathrm{IMD}}
 d\epsilon\\
\qquad \qquad 
\text{ on $U_0 \times(\widehat{\mathcal{M}}_{\boldsymbol{t}_{\text{ra}}}\times T_{\boldsymbol{\theta}} )
\times \mathrm{Spec}\, \mathbb{C} [\epsilon]$ } \\
\hat{d} + G_{1}^{-1} dG_{1}+G_{1}^{-1} 
\widehat\Omega_{(\boldsymbol{t}_{\mathrm{ra}},\boldsymbol{\theta}_0)}^{(1)} 
 \, G_{1}
+\epsilon G_{1}^{-1} 
\delta_{\theta_{l,t_i}^{\pm}}^{\mathrm{IMD}}
(\widehat\Omega_{(\boldsymbol{t}_{\mathrm{ra}},\boldsymbol{\theta}_0)}^{(1)} )
G_{1} 
+ G_1^{-1} \Upsilon_{\theta_{l,t_i}^{\pm}}^{\mathrm{IMD}}
 G_1 d\epsilon \\
\qquad\qquad 
\text{ on $U_\infty \times(\widehat{\mathcal{M}}_{\boldsymbol{t}_{\text{ra}}}\times T_{\boldsymbol{\theta}} )
\times \mathrm{Spec}\, \mathbb{C} [\epsilon]$ },
\end{cases}
\end{equation}
which is 
a morphism $\widehat{E}^{\epsilon}_1 \rightarrow 
\widehat{E}^{\epsilon}_1 \otimes
\tilde\Omega^1_{\partial/\partial \theta_{l,t_i}^{\pm}}$
with the Leibniz rule.
Since $(\widehat{E}^{\epsilon}_1,
\overline{(\varphi_{\Upsilon}^{-1})^*\nabla^{\epsilon}_{\partial/\partial \theta_{l,t_i}^{\pm}}})$
is relativization of the horizontal lift,
we have the following equality
\begin{equation}\label{2021.4.16.14.37}
 \delta_{\theta_{l,t_i}^{\pm}}^{\mathrm{IMD}}
(\widehat\Omega_{(\boldsymbol{t}_{\mathrm{ra}},\boldsymbol{\theta}_0)}^{(1)} )
 = 
d  \Upsilon_{\theta_{l,t_i}^{\pm}}^{\mathrm{IMD}} 
+[\widehat\Omega_{(\boldsymbol{t}_{\mathrm{ra}},\boldsymbol{\theta}_0)}^{(1)}  ,
\Upsilon_{\theta_{l,t_i}^{\pm}}^{\mathrm{IMD}}],
\end{equation}
which means the integrable condition.

\subsection{Integrable deformations associated to 
$(T_{\boldsymbol{t}})_{\boldsymbol{t}_{\mathrm{ra}}}$}\label{2020.1.25.8.11}

First we fix $i \in \{3,4,\ldots,\nu\} $.
Let 
$$
\tilde{\nabla}_{\mathrm{DL, ext}}^{(1)}=
\begin{cases}
d+ \widehat\Omega_{(\boldsymbol{t}_{\mathrm{ra}},\boldsymbol{\theta}_0)}^{(1)} 
& \text{ on $U_0\times(\widehat{\mathcal{M}}_{\boldsymbol{t}_{\text{ra}}} \times
 T_{\boldsymbol{\theta}})$} \\
d + G_{1}^{-1} dG_{1}+G_{1}^{-1} 
\widehat\Omega_{(\boldsymbol{t}_{\mathrm{ra}},\boldsymbol{\theta}_0)}^{(1)} \, G_{1}
 & \text{ on $U_\infty \times (\widehat{\mathcal{M}}_{\boldsymbol{t}_{\text{ra}}} 
 \times T_{\boldsymbol{\theta}})$}
\end{cases}
$$
be the family \eqref{2021.4.14.10.36}.
Let $t_i$ be the natural coordinate of 
$(T_{\boldsymbol{t}})_{\boldsymbol{t}_{\text{ra}}} \times T_{\boldsymbol{\theta}}$
and $\partial/\partial t_i$ be the vector field on
$(T_{\boldsymbol{t}})_{\boldsymbol{t}_{\text{ra}}} \times T_{\boldsymbol{\theta}}$
associated to $t_i$.
We will construct a horizontal lift of $\tilde{\nabla}_{\mathrm{DL, ext}}^{(1)}$
with respect to $\partial/\partial t_i$.

For the fixed index $i$,
we define a matrix $B_{t_i}(x)$ 
by
\begin{equation*}
B_{t_i}(x):=-\sum^{n_i-1}_{l=0}   \begin{pmatrix}
\theta_{l,t_i}^+ & 0 \\
0 & \theta_{l,t_i}^-
\end{pmatrix}\frac{\hat{\delta}(t_i)}{(x-t_i)^{n_i -l}} .
\end{equation*}
We define
a vector bundle $(\widehat{E}_1)_{t_i}^{\epsilon}$ on 
$\mathbb{P}^1\times(\widehat{\mathcal{M}}_{\boldsymbol{t}_{\text{ra}}} \times
 T_{\boldsymbol{\theta} } ) \times \mathrm{Spec}\, \mathbb{C}[\epsilon]$
by the same argument as in the construction of $(\widehat{E}_1)_{\theta_{l,t_i}^{\pm}}^{\epsilon}$.
That is, we replace $B_{\theta_{l,t_i}^{\pm}}(x_{t_i})$ in \eqref{2021.4.16.15.17}
with $B_{t_i}(x)$.
We define a morphism
$$
\nabla^{\epsilon}_{\partial/\partial t_i}  \colon 
(\widehat{E}_1)_{t_i}^{\epsilon} 
\longrightarrow 
(\widehat{E}_1)_{t_i}^{\epsilon} \otimes
\tilde\Omega^1_{\partial/\partial t_i}
$$
by gluing the following connections:
\begin{equation*}
\left\{
\begin{array}{ll}
\nabla^{\epsilon}_{i'}  = \hat{d}+ \hat{\Omega}_{i'}  \quad \text{
for $i' \in (I\setminus \{ i \})\cup \{ \nu+1 \}$,} \\
\nabla^{\epsilon}_{i}=\hat{d}+ \hat{\Omega}_i + \epsilon\left( \frac{\partial}{\partial x_{t_i}}
 (B_{t_i} )dx_{t_i}
+ [ \hat{\Omega}_i, B_{t_i}]  \right) 
+B_{t_i} d \epsilon
\end{array}
\right.
\end{equation*}
as in the construction of $\nabla^{\epsilon}_{\partial/\partial \theta^{\pm}_{l,t_i}}$ in the previous section.
Here $\hat{\Omega}_{i'}$ ($i' \in I \cup \{ \nu+1 \}$) is defined 
in \eqref{2021.11.11.23.26}.
Now we check that the connection matrix of $\nabla^{\epsilon}_{i}$
is a section of $\tilde{\Omega}^1_{\partial/\partial t_i}$ defined in
\eqref{2021.4.14.23.20}.
We set $\tilde{x}_{t_i} := x-(t_i +\epsilon \hat{\delta}(t_i))=x_{t_i} -\epsilon \hat{\delta}(t_i)$.
Since $\epsilon^2=0$ and $\epsilon d \epsilon =0$, we may check the following equalities:
\begin{equation*}
\begin{aligned}
&\begin{pmatrix}
\theta_{0,t_i}^{+}& 0 \\
0 & \theta_{0,t_i}^{-}
\end{pmatrix}\frac{\hat{d}
\tilde{x}_{t_i}}{\tilde{x}_{t_i}^{n_i}  }+\cdots 
+\begin{pmatrix}
 \theta_{n_i-1,t_i}^+ & 0 \\
0 & \theta_{n_i-1,t_i}^-
\end{pmatrix}\frac{\hat{d}
\tilde{x}_{t_i} }{\tilde{x}_{t_i}} \\
&=\begin{pmatrix}
\theta_{0,t_i}^{+}& 0 \\
0 & \theta_{0,t_i}^{-}
\end{pmatrix}\frac{d
x_{t_i}}{(x_{t_i} -\epsilon \hat{\delta}(t_i))^{n_i}  }+\cdots 
+\begin{pmatrix}
 \theta_{n_i-1,t_i}^+ & 0 \\
0 & \theta_{n_i-1,t_i}^-
\end{pmatrix}\frac{d
x_{t_i}}{x_{t_i} -\epsilon \hat{\delta}(t_i)} \\
&\quad \qquad -\left( \begin{pmatrix}
\theta_{0,t_i}^{+}& 0 \\
0 & \theta_{0,t_i}^{-}
\end{pmatrix}\frac{
 \hat{\delta}(t_i) d\epsilon }{x_{t_i}^{n_i}  }+\cdots 
+\begin{pmatrix}
 \theta_{n_i-1,t_i}^+ & 0 \\
0 & \theta_{n_i-1,t_i}^-
\end{pmatrix}\frac{
\hat{\delta}(t_i)d \epsilon }{x_{t_i} } \right)\\
&=\begin{pmatrix}
\theta_{0,t_i}^{+}& 0 \\
0 & \theta_{0,t_i}^{-}
\end{pmatrix}\frac{d
x_{t_i}}{x_{t_i} ^{n_i}  }+\cdots 
+\begin{pmatrix}
 \theta_{n_i-1,t_i}^+ & 0 \\
0 & \theta_{n_i-1,t_i}^-
\end{pmatrix}\frac{d
x_{t_i}}{x_{t_i} } \\
&\qquad \quad -\epsilon \hat{\delta}(t_i) \cdot \frac{\partial }{\partial x_{t_i}} \left(
\begin{pmatrix}
\theta_{0,t_i}^{+}& 0 \\
0 & \theta_{0,t_i}^{-}
\end{pmatrix}\frac{1}{x_{t_i} ^{n_i}  }+\cdots 
+\begin{pmatrix}
 \theta_{n_i-1,t_i}^+ & 0 \\
0 & \theta_{n_i-1,t_i}^-
\end{pmatrix}\frac{1}{x_{t_i} } \right)d
x_{t_i}+B_{t_i} d \epsilon\\
&=\begin{pmatrix}
\theta_{0,t_i}^{+}& 0 \\
0 & \theta_{0,t_i}^{-}
\end{pmatrix}\frac{d
x_{t_i}}{x_{t_i} ^{n_i}  }+\cdots 
+\begin{pmatrix}
 \theta_{n_i-1,t_i}^+ & 0 \\
0 & \theta_{n_i-1,t_i}^-
\end{pmatrix}\frac{d
x_{t_i}}{x_{t_i} } 
 +\epsilon  \frac{\partial }{\partial x_{t_i}} \left(
B_{t_i} \right)d
x_{t_i}+B_{t_i} d \epsilon.
\end{aligned}
\end{equation*}
Moreover $\hat{\Omega}_i$ and $B_{t_i}$ are diagonal 
until the $x_{t_i}^{n_i-1}$-terms.
So we have that 
$$
\hat{\Omega}_i + \epsilon\left( \frac{\partial}{\partial x_{t_i}}
 (B_{t_i} )dx_{t_i}
+ [ \hat{\Omega}_i, B_{t_i}]  \right) 
+B_{t_i} d \epsilon
$$
is a section of $\tilde{\Omega}^1_{\partial/\partial t_i}$.

As in the previous section, 
$\widehat{E}_1^{\epsilon} \cong (\widehat{E}_1)_{t_i}^{\epsilon}$.
If we consider the pull-back of $\nabla^{\epsilon}_{\partial/\partial t_i}$
under this isomorphism, 
then we have a horizontal lift of $\tilde{\nabla}_{\mathrm{DL, ext}}^{(1)}$
with respect to $\partial/\partial t_i$.
If we take a relativization of this horizontal lift,
we have a family of connections parametrized by
$( \widehat{\mathcal{M}}_{\boldsymbol{t}_{\text{ra}}}\times T_{\boldsymbol{\theta}})
\times \mathrm{Spec}\, \mathbb{C}[\epsilon]$.
This family gives a map from the base space
$( \widehat{\mathcal{M}}_{\boldsymbol{t}_{\text{ra}}}\times T_{\boldsymbol{\theta}})
\times \mathrm{Spec}\, \mathbb{C}[\epsilon]$
to the moduli space $\widehat{\mathfrak{Conn}}_{(\boldsymbol{t}_{\mathrm{ra}},\boldsymbol{\theta}_0)}$.
By taking composition with $\widehat{\mathrm{App}}$ defined in \eqref{2020.1.6.13.00},
we have a map 
\begin{equation}\label{2021.4.16.15.40}
( \widehat{\mathcal{M}}_{\boldsymbol{t}_{\text{ra}}}\times T_{\boldsymbol{\theta}})
\times \mathrm{Spec}\, \mathbb{C}[\epsilon]
\longrightarrow
\widehat{\mathcal{M}}_{\boldsymbol{t}_{\text{ra}}} \times T_{\boldsymbol{\theta}}.
\end{equation}

\begin{Def}
Then 
we may define the vector field 
on $\widehat{\mathcal{M}}_{\boldsymbol{t}_{\text{ra}}}\times T_{\boldsymbol{\theta}}$
associated 
to the map \eqref{2021.4.16.15.40}.
We denote by $\delta_{t_i}^{\mathrm{IMD}}$ this vector field 
on $\widehat{\mathcal{M}}_{\boldsymbol{t}_{\text{ra}}}\times T_{\boldsymbol{\theta}}$.
\end{Def}

Let $f_{t_i}^{\mathrm{IMD}} \colon 
( \widehat{\mathcal{M}}_{\boldsymbol{t}_{\text{ra}}}\times T_{\boldsymbol{\theta}})
\times \mathrm{Spec}\, \mathbb{C}[\epsilon]
\longrightarrow
\widehat{\mathcal{M}}_{\boldsymbol{t}_{\text{ra}}} \times T_{\boldsymbol{\theta}}$
be the map induced by the vector field
$\delta_{t_i}^{\mathrm{IMD}}$.
We have $\widehat{E}^{\epsilon}_1=
(\mathrm{id} \times f_{t_i}^{\mathrm{IMD}})^* \widehat{E}_1$.
We denote by
\begin{equation}\label{2021.4.16.12.12_y}
\begin{cases}
d+ \widehat\Omega_{(\boldsymbol{t}_{\mathrm{ra}},\boldsymbol{\theta}_0)}^{(1)} 
+\epsilon \delta_{t_i}^{\mathrm{IMD}}
(\widehat\Omega_{(\boldsymbol{t}_{\mathrm{ra}},\boldsymbol{\theta}_0)}^{(1)} ) \\
\qquad \qquad 
\text{ on $U_0 \times(\widehat{\mathcal{M}}_{\boldsymbol{t}_{\text{ra}}}\times T_{\boldsymbol{\theta}} )
\times \mathrm{Spec}\, \mathbb{C} [\epsilon]$ } \\
d + G_{1}^{-1} dG_{1}+G_{1}^{-1} 
\widehat\Omega_{(\boldsymbol{t}_{\mathrm{ra}},\boldsymbol{\theta}_0)}^{(1)} 
 \, G_{1}
+\epsilon G_{1}^{-1} 
\delta_{t_i}^{\mathrm{IMD}}
(\widehat\Omega_{(\boldsymbol{t}_{\mathrm{ra}},\boldsymbol{\theta}_0)}^{(1)} )
G_{1}  \\
\qquad\qquad 
\text{ on $U_\infty \times(\widehat{\mathcal{M}}_{\boldsymbol{t}_{\text{ra}}}\times T_{\boldsymbol{\theta}} )
\times \mathrm{Spec}\, \mathbb{C} [\epsilon]$ }
\end{cases}
\end{equation}
the pull-back of $\tilde{\nabla}_{\mathrm{DL, ext}}^{(1)}$ 
under the map $\mathrm{id} \times f_{t_i}^{\mathrm{IMD}}$.
As in the previous section,
we have a lift of 
$(\mathrm{id} \times f_{t_i}^{\mathrm{IMD}})^*\tilde{\nabla}_{\mathrm{DL, ext}}^{(1)}$:
\begin{equation}\label{2021.4.16.12.12_x}
\begin{cases}
\hat{d}+ \widehat\Omega_{(\boldsymbol{t}_{\mathrm{ra}},\boldsymbol{\theta}_0)}^{(1)} 
+\epsilon \delta_{t_i}^{\mathrm{IMD}}
(\widehat\Omega_{(\boldsymbol{t}_{\mathrm{ra}},\boldsymbol{\theta}_0)}^{(1)} )
+ \Upsilon_{t_i}^{\mathrm{IMD}}
 d\epsilon\\
\qquad \qquad 
\text{ on $U_0 \times(\widehat{\mathcal{M}}_{\boldsymbol{t}_{\text{ra}}}\times T_{\boldsymbol{\theta}} )
\times \mathrm{Spec}\, \mathbb{C} [\epsilon]$ } \\
\hat{d} + G_{1}^{-1} dG_{1}+G_{1}^{-1} 
\widehat\Omega_{(\boldsymbol{t}_{\mathrm{ra}},\boldsymbol{\theta}_0)}^{(1)} 
 \, G_{1}
+\epsilon G_{1}^{-1} 
\delta_{t_i}^{\mathrm{IMD}}
(\widehat\Omega_{(\boldsymbol{t}_{\mathrm{ra}},\boldsymbol{\theta}_0)}^{(1)} )
G_{1} 
+ G_1^{-1} \Upsilon_{t_i}^{\mathrm{IMD}}
 G_1 d\epsilon \\
\qquad\qquad 
\text{ on $U_\infty \times(\widehat{\mathcal{M}}_{\boldsymbol{t}_{\text{ra}}}\times T_{\boldsymbol{\theta}} )
\times \mathrm{Spec}\, \mathbb{C} [\epsilon]$ },
\end{cases}
\end{equation}
which is 
a morphism $\widehat{E}^{\epsilon}_1 \rightarrow 
\widehat{E}^{\epsilon}_1 \otimes
\tilde\Omega^1_{\partial/\partial t_i}$
with the Leibniz rule and the following equality
\begin{equation}\label{2021.11.16.17.17}
 \delta_{t_i}^{\mathrm{IMD}}
(\widehat\Omega_{(\boldsymbol{t}_{\mathrm{ra}},\boldsymbol{\theta}_0)}^{(1)} )
 = 
d  \Upsilon_{t_i}^{\mathrm{IMD}} 
+[\widehat\Omega_{(\boldsymbol{t}_{\mathrm{ra}},\boldsymbol{\theta}_0)}^{(1)}  ,
\Upsilon_{t_i}^{\mathrm{IMD}}],
\end{equation}
which means the integrable condition.

\subsection{Isomonodromy $2$-form}\label{2020.2.8.10.28}

\begin{Def}\label{2021.4.23.11.26}
Let $\hat{\delta}_1$ and $\hat{\delta}_2$ be vector fields on 
$\widehat{\mathcal{M}}_{\boldsymbol{t}_{\text{ra}}}\times T_{\boldsymbol{\theta}}$,
which is isomorphic to the extended moduli space
$\widehat{\mathfrak{Conn}}_{(\boldsymbol{t}_{\mathrm{ra}},\boldsymbol{\theta}_0)}$.
We fix a formal fundamental matrix solution $\psi_i$ of 
$(d+\widehat\Omega_{(\boldsymbol{t}_{\mathrm{ra}},\boldsymbol{\theta}_0)}^{(n-2)})
\psi_i=0$ at $x=t_i$
as in Lemma \ref{2019.12.30.22.09}.
We take a fundamental matrix solution $\psi_{q_j}$ of 
$(d+\widehat\Omega_{(\boldsymbol{t}_{\mathrm{ra}},\boldsymbol{\theta}_0)}^{(n-2)})
\psi_{q_j}=0$ at $x=q_j$
as in Lemma \ref{2019.12.30.21.53}.
We define a $2$-form $\hat{\omega}$ 
on $\widehat{\mathcal{M}}_{\boldsymbol{t}_{\text{ra}}}\times T_{\boldsymbol{\theta}}$ as
\begin{equation}\label{2020.1.7.17.07}
\begin{aligned}
\hat{\omega} (\hat{\delta}_1,\hat{\delta}_2) :=\ &  \frac{1}{2} \sum_{i\in I} 
\mathrm{res}_{x=t_i} 
\mathrm{Tr} \left( \hat{\delta} (
\widehat\Omega_{(\boldsymbol{t}_{\mathrm{ra}},\boldsymbol{\theta}_0)}^{(n-2)} ) 
\wedge \hat{\delta}(\psi_i)\psi_i^{-1} )   \right)
+\frac{1}{2} \sum_{j=1}^{n-3} \mathrm{res}_{x=q_j} 
\mathrm{Tr} \left( \hat{\delta} 
(\widehat\Omega_{(\boldsymbol{t}_{\mathrm{ra}},\boldsymbol{\theta}_0)}^{(n-2)} )
\wedge \hat{\delta}(\psi_{q_j})\psi_{q_j}^{-1} 
  \right),
\end{aligned}
\end{equation}
where $I:= \{ 1,2,\ldots, \nu, \infty \}$.
Here we set $\hat\delta(A) \wedge \hat\delta(\psi) \psi^{-1}:= 
\hat\delta_1 (A)
\hat\delta_2(\psi) (\psi)^{-1}  
-\hat\delta_1(\psi) (\psi)^{-1} \hat\delta_2 (A)$.
\end{Def}

Since $\theta^{\pm}_{n_i-1,t_i}$ is constant on 
$\widehat{\mathcal{M}}_{\boldsymbol{t}_{\text{ra}}}\times T_{\boldsymbol{\theta}}$,
$\hat{\delta}(\theta^{\pm}_{n_i-1,t_i})=0$. 
Then $\hat{\delta}(\theta^{\pm}_{n_i-1,t_i} \int x_{t_i}^{-1} dx_{t_i} )=
\frac{-\theta^{\pm}_{n_i-1,t_i}\hat\delta(t_i)}{x-{t_i}}+\hat\delta(c)$.
Here $c$ is an integration constant.
By the same argument as in Section \ref{2020.2.8.10.24},
we have that the residue of 
$\hat{\delta} (\widehat\Omega_{(\boldsymbol{t}_{\mathrm{ra}},\boldsymbol{\theta}_0)}^{(n-2)}) 
 \wedge \hat{\delta}(\psi_i)\psi_i^{-1}$ at $\tilde{t}_i$
is well-defined.
By the same argument as in Section \ref{2020.2.8.10.24},
We may check that 
the right hand side of \eqref{2020.1.7.17.07} is independent of the choice of 
$\psi_{q_j}$ and $\psi_{i}$.

We will show a transformation formula (Lemma \ref{2020.12.24.17.31} below).
We will use this transformation formula 
for calculation of $\hat{\omega} (\hat{\delta}_1,\hat{\delta}_2)$.
We show this transformation formula
for general situations.
Let $C$ be a smooth projective curve over $\mathbb{C}$
and let $M$ be an algebraic variety over $\mathbb{C}$.
Let $U$ be an analytic open subset of $C$.
Let $x$ be a parameter on $U$.
Let $ d+ A dx$ be a family of connections on $\mathcal{O}^{\oplus 2}_U$ 
parameterized by $M$.
Assume we can take a (formal) fundamental matrix solution $\psi$ 
of $(d + A dx) \psi = 0 $,
that is, there exists  
$\psi \in \mathcal{E}nd (\mathcal{O}_M^{\oplus 2} ) \otimes \widehat{\mathcal{O}}_{U,0}$
such that $d \psi + A \psi dx  = 0 $.
Here $d$ means the relative exterior derivative on
the projection $U\times M \rightarrow M$.

\begin{Lem}\label{2020.12.24.17.31}
\textit{
Let $\delta_1$ and $\delta_2$ be vector fields on $M$.
Let $g$
be a family of matrices parameterized by $M$ 
such that the entries of the matrix $g_m$ for $m \in M$
 are meromorphic functions on $U$.
 We assume that we can define $g^{-1}$, which is 
 the family of matrices parameterized by $M$ 
such that for $m\in M$, $g_m (g^{-1})_m =\mathrm{id}$ 
and the entries of $(g^{-1})_m$ are meromorphic functions on $U$.
Set $A':= g^{-1} dg +g^{-1} A g$ and $\psi':=g^{-1} \psi$.
Moreover set $u^{(l)}:= \delta_l(g)g^{-1}$ and
 $\tilde{u}^{(l)}:= g^{-1}\delta_l(g)$ for $l \in \{ 1,2\}$.
Then we have the following equality:
\begin{equation}\label{2020.1.8.12.55}
\begin{split}
&\mathrm{Tr} \left( \delta(A') \wedge \delta(\psi') (\psi')^{-1}\right)
- \mathrm{Tr} \left( \delta(A) \wedge \delta(\psi) \psi^{-1} \right)\\
& =-\mathrm{Tr} \left( \delta_1(A') \tilde{u}^{(2)} -
\tilde{u}^{(1)}\delta_2(A')  \right)
-\mathrm{Tr} \left( \delta_1(A)   u^{(2)} -
u^{(1)}\delta_2(A)   \right)  \\
&\quad + \mathrm{Tr} \left(  d (\psi^{-1} u^{(1)} \delta_2(\psi) 
-\psi^{-1} u^{(2)} \delta_1(\psi) ) \right).
\end{split}
\end{equation}}
\end{Lem}

\begin{proof}
Since $\psi'=g^{-1} \psi$, we have the following equalities:
\begin{equation}\label{eq:2020.11.19.15.37}
\begin{aligned}
\mathrm{Tr} \left( \delta(A') \wedge \delta(\psi') (\psi')^{-1}\right)
&=\mathrm{Tr} \left( \delta(A') \wedge \delta(g^{-1}\psi) \psi^{-1}g\right)\\
&=\mathrm{Tr} \left( \delta(A') \wedge  ( - g^{-1}\delta(g) g^{-1}\psi \psi^{-1} g 
+  g^{-1}\delta( \psi) \psi^{-1}g )\right) \\
&=\mathrm{Tr} \left( \delta(A') \wedge  ( - g^{-1}\delta(g) 
+  g^{-1}\delta( \psi) \psi^{-1}g )\right).
\end{aligned}
\end{equation}
We calculate $\mathrm{Tr} \left( \delta(A') \wedge 
  ( g^{-1}\delta( \psi) \psi^{-1}g ) \right)$ as follows:
\begin{equation}\label{eq:2020.11.19.15.38}
\begin{aligned}
\mathrm{Tr} \left( \delta(A') \wedge   ( g^{-1}\delta( \psi) \psi^{-1}g ) \right)
&=\mathrm{Tr} \left( \delta(g^{-1} dg + g^{-1}A g ) \wedge  
(  g^{-1}\delta( \psi) \psi^{-1}g ) \right)\\
&=\mathrm{Tr} \left( ( -g^{-1} \delta(g) g^{-1} dg+ g^{-1} \delta( dg )) \wedge 
 (  g^{-1}\delta( \psi) \psi^{-1}g ) \right) \\
&\qquad +\mathrm{Tr} \left( (-g^{-1} \delta( g ) g^{-1} A g 
+g^{-1} \delta(A) g + g^{-1} A \delta(g) )
\wedge  (  g^{-1}\delta( \psi) \psi^{-1}g ) \right)\\
&=\mathrm{Tr} \left( d (  \delta(g)   g^{-1} ) \wedge 
 ( \delta( \psi) \psi^{-1}) \right) 
 -\mathrm{Tr} \left( ( \delta( g ) g^{-1} A   )
\wedge  (  \delta( \psi) \psi^{-1} ) \right)\\
& \qquad+\mathrm{Tr} \left(  \delta(A)   
\wedge  ( \delta( \psi) \psi^{-1} ) \right)
+\mathrm{Tr} \left( (A \delta(g)  g^{-1} )
\wedge  (  \delta( \psi) \psi^{-1} ) \right).
\end{aligned}
\end{equation}
By the equalities \eqref{eq:2020.11.19.15.37} 
and \eqref{eq:2020.11.19.15.38}, we have the following equality:
\begin{equation}\label{eq:2020.11.19.15.43}
\begin{aligned}
&\mathrm{Tr} \left( \delta(A') \wedge \delta(\psi') (\psi')^{-1}\right)
+\mathrm{Tr} \left( \delta(A') \wedge  ( g^{-1}\delta(g) )\right)
-\mathrm{Tr} \left(  \delta(A)   
\wedge  ( \delta( \psi) \psi^{-1} ) \right) \\
&=\mathrm{Tr} \left( d(  \delta(g) g^{-1} ) \wedge  (  \delta( \psi) \psi^{-1}) \right) \\
&\qquad -\mathrm{Tr} \left( (\delta( g ) g^{-1} A   )
\wedge  (  \delta( \psi) \psi^{-1} ) \right)
+\mathrm{Tr} \left( (A \delta(g)  g^{-1})
\wedge  (  \delta( \psi) \psi^{-1} ) \right).
\end{aligned}
\end{equation}
We calculate $\mathrm{Tr} \left(  d (\psi^{-1} u^{(1)} \delta_2(\psi) 
-\psi^{-1} u^{(2)} \delta_1(\psi) ) \right)$ as follows:
\begin{equation}\label{eq:2020.11.19.15.44}
\begin{aligned}
& \mathrm{Tr} \left(  d (\psi^{-1} u^{(1)} \delta_2(\psi) 
-\psi^{-1} u^{(2)} \delta_1(\psi) ) \right) \\
&=\mathrm{Tr} \left( - \psi^{-1} d (\psi) \psi^{-1} u^{(1)} \delta_2(\psi)
+\psi^{-1}   d (u^{(1)})  \delta_2(\psi) 
+\psi^{-1}   u^{(1)}  \delta_2( d \psi)  \right) \\ 
&\qquad  - \mathrm{Tr} \left(  -\psi^{-1} d (\psi) \psi^{-1} u^{(2)} \delta_1(\psi)
+\psi^{-1}   d (u^{(2)})  \delta_1(\psi) 
+\psi^{-1}   u^{(2)}  \delta_1( d \psi)  \right)  \\
&=\mathrm{Tr} \left(  A u^{(1)} \delta_2(\psi)  \psi^{-1} 
+   d (u^{(1)})  \delta_2(\psi) \psi^{-1}
-   u^{(1)}  \delta_2( A \psi) \psi^{-1}  \right)  \\ 
&\qquad - \mathrm{Tr} \left(   A u^{(2)} \delta_1(\psi)\psi^{-1} 
+   d (u^{(2)})  \delta_1(\psi)\psi^{-1} 
-   u^{(2)}  \delta_1( A \psi) \psi^{-1}   \right)  \\
&=\mathrm{Tr} \left( (A \delta(g)  g^{-1})
\wedge  (  \delta( \psi) \psi^{-1} ) \right)
+\mathrm{Tr} \left( d(  \delta(g) g^{-1} ) \wedge  (  \delta( \psi) \psi^{-1}) \right) \\ 
&\qquad  - \mathrm{Tr} \left(     u^{(1)}  \delta_2( A) 
-   u^{(2)}  \delta_1( A)  
+  u^{(1)}   A \delta_2(\psi) \psi^{-1}
-   u^{(2)}   A \delta_1( \psi) \psi^{-1} \right)  \\
&=\mathrm{Tr} \left( (A \delta(g)  g^{-1})
\wedge  (  \delta( \psi) \psi^{-1} ) \right)
+\mathrm{Tr} \left( d(  \delta(g) g^{-1} ) \wedge  (  \delta( \psi) \psi^{-1}) \right) \\
&\qquad+ \mathrm{Tr} \left( \delta ( A   )
\wedge  (  \delta( g) g^{-1} ) \right)dx
 -\mathrm{Tr} \left( (\delta( g ) g^{-1} A  )
\wedge  (  \delta( \psi) \psi^{-1} ) \right).
\end{aligned}
\end{equation}
Here the second equality follows from $d \psi = - A \psi $.
The equality \eqref{2020.1.8.12.55} follows
from the equalities \eqref{eq:2020.11.19.15.43} and \eqref{eq:2020.11.19.15.44}.
\end{proof}

\begin{Prop}\label{2021.4.26.22.19}
{\it 
Let
$\tilde{G}$ and
$\tilde{G}_{\infty}$
be the matrices defined in
\eqref{2021.4.2.20.06}.
Set $\tilde\psi_i := \tilde{G}^{-1} \psi_i$ for any $i \in I_{\mathrm{un}} \setminus \{\infty \}$
and
$\tilde\psi_{\infty} := \tilde{G}^{-1}_{\infty} \psi_\infty$.
We have the following equality: 
\begin{equation}\label{eq:2020.11.19.16.08}
\begin{split}
\hat\omega (\hat\delta_1, \hat\delta_2) &=
 \frac{1}{2} \sum_{i \in I}   \mathrm{res}_{x=t_i} 
\mathrm{Tr} \left( \hat\delta
(\widehat\Omega_{(\boldsymbol{t}_{\mathrm{ra}},\boldsymbol{\theta}_0)}^{(1)}) 
\wedge\hat\delta(\tilde\psi_i)(\tilde\psi_i)^{-1}  \right).
\end{split}
\end{equation}}
\end{Prop}

\begin{proof}
By Proposition \ref{2021.4.11.16.40}, we have 
$\widehat\Omega_{(\boldsymbol{t}_{\mathrm{ra}},\boldsymbol{\theta}_0)}^{(1)}
= \tilde{G}^{-1} d \tilde{G}+
\tilde{G}^{-1}
\widehat\Omega_{(\boldsymbol{t}_{\mathrm{ra}},\boldsymbol{\theta}_0)}^{(n-2)} \tilde{G}$.
Set $u_{\tilde{G}}^{(l)}:= \delta_l(\tilde{G})\tilde{G}^{-1}$ 
and $\tilde{u}_{\tilde{G}}^{(l)}:= \tilde{G}^{-1}\delta_l(\tilde{G})$ for $l \in \{ 1,2\}$.
Set $\tilde\psi_{q_j} := \tilde{G}^{-1} \psi_{q_j}$ for any $j \in  \{1,2,\ldots,n-3 \}$.
We calculate the difference between the right hand side
and the left hand side of \eqref{eq:2020.11.19.16.08}.
By Lemma \ref{2019.12.30.21.53}, we have 
$\tilde{\psi}_{q_j}= \tilde{G}^{-1} \Phi_{q_j} \Xi_{q_j}(x) \Lambda_{q_j}(x) $.
We calculate $ \tilde{G}^{-1} \Phi_{q_j}$ as follows.
$$
\begin{aligned}
 \tilde{G}^{-1} \Phi_{q_j}
&=\begin{pmatrix} 1 & 0 \\ -\frac{Q_2(x)}{Q_1(x)} & \frac{1}{Q_1(x)} \end{pmatrix}
\begin{pmatrix} 1 & 0 \\ p_j & 1 \end{pmatrix}
= \begin{pmatrix} 1 & 0 \\ \frac{-Q_2(x) +p_j}{Q_1(x)} & 1 \end{pmatrix}
\begin{pmatrix} 1 & 0 \\ 0 & \frac{1}{Q_1(x)} \end{pmatrix}.
\end{aligned}
$$
Since $Q_2(q_j)=p_j$, we may remove a pole of $\frac{-Q_2(x) +p_j}{Q_1(x)}$ at $q_j$.
By Lemma \ref{2019.12.30.21.53}, we may check that the pole of
$\begin{pmatrix} 1 & 0 \\ 0 & \frac{1}{Q_1(x)} \end{pmatrix}
\Xi_{q_j}(x) \Lambda_{q_j}(x)$
at $q_j$ is removable.
So we have that 
$\mathrm{Tr} \left( \delta(
\widehat\Omega_{(\boldsymbol{t}_{\mathrm{ra}},\boldsymbol{\theta}_0)}^{(1)})
\wedge \delta(\tilde\psi_{q_j}) 
(\tilde\psi_{q_j})^{-1}\right)$
has no pole at $q_j$ ($j=1,2,\ldots,n-3$).
Then the difference between the right hand side
and the left hand side of \eqref{eq:2020.11.19.16.08}
is equal to
\begin{equation}\label{2021.4.20.13.28}
\begin{split}
& \frac{1}{2} \sum_{i \in I}   \mathrm{res}_{x=t_i} \left( 
\mathrm{Tr} \left( \delta
(\widehat\Omega_{(\boldsymbol{t}_{\mathrm{ra}},\boldsymbol{\theta}_0)}^{(1)}) 
\wedge \delta(\tilde\psi_i) (\tilde\psi_i)^{-1}\right)- 
\mathrm{Tr} \left( \delta
(\widehat\Omega_{(\boldsymbol{t}_{\mathrm{ra}},\boldsymbol{\theta}_0)}^{(n-2)})
 \wedge \delta(\psi_i) \psi_i^{-1} \right) \right)\\
& \quad  +\frac{1}{2} \sum_{j =1}^{n-3}   \mathrm{res}_{x=q_j} \left( 
\mathrm{Tr} \left( \delta
(\widehat\Omega_{(\boldsymbol{t}_{\mathrm{ra}},\boldsymbol{\theta}_0)}^{(1)}) 
\wedge \delta(\tilde\psi_{q_j}) (\tilde\psi_{q_j})^{-1}\right)- 
\mathrm{Tr} \left( \delta
(\widehat\Omega_{(\boldsymbol{t}_{\mathrm{ra}},\boldsymbol{\theta}_0)}^{(n-2)})
 \wedge \delta(\psi_{q_j}) \psi_{q_j}^{-1} \right) \right).
\end{split}
\end{equation}
By the equation \eqref{2020.1.8.12.55}, the difference \eqref{2021.4.20.13.28} is equal to 
\begin{equation}\label{eq:2020.11.19.16.18}
\begin{split}
 &-\frac{1}{2} \sum_{i \in I}   \mathrm{res}_{x=t_i} \left( 
\mathrm{Tr} \left( \delta_1(
\widehat\Omega_{(\boldsymbol{t}_{\mathrm{ra}},\boldsymbol{\theta}_0)}^{(1)})
 \tilde{u}_{\tilde{G}}^{(2)} -
\tilde{u}^{(1)}_{\tilde{G}}\delta_2(
\widehat\Omega_{(\boldsymbol{t}_{\mathrm{ra}},\boldsymbol{\theta}_0)}^{(1)})
  \right)
-\mathrm{Tr} \left( \delta_1(
\widehat\Omega_{(\boldsymbol{t}_{\mathrm{ra}},\boldsymbol{\theta}_0)}^{(n-2)})
   u^{(2)}_{\tilde{G}} -
u^{(1)}_{\tilde{G}}\delta_2(
\widehat\Omega_{(\boldsymbol{t}_{\mathrm{ra}},\boldsymbol{\theta}_0)}^{(n-2)})   \right)   \right) \\
& \quad - \frac{1}{2} \sum_{j =1}^{n-3}   \mathrm{res}_{x=q_j}  \left( 
\mathrm{Tr} \left( \delta_1(
\widehat\Omega_{(\boldsymbol{t}_{\mathrm{ra}},\boldsymbol{\theta}_0)}^{(1)})
 \tilde{u}_{\tilde{G}}^{(2)} -
\tilde{u}^{(1)}_{\tilde{G}}\delta_2(
\widehat\Omega_{(\boldsymbol{t}_{\mathrm{ra}},\boldsymbol{\theta}_0)}^{(1)})
  \right)
-\mathrm{Tr} \left( \delta_1
(\widehat\Omega_{(\boldsymbol{t}_{\mathrm{ra}},\boldsymbol{\theta}_0)}^{(n-2)})
   u_{\tilde{G}}^{(2)} -
u^{(1)}_{\tilde{G}}\delta_2(
\widehat\Omega_{(\boldsymbol{t}_{\mathrm{ra}},\boldsymbol{\theta}_0)}^{(n-2)}) 
  \right)   \right).
\end{split}
\end{equation}
Here note that the third term of the right hand side of \eqref{2020.1.8.12.55}
is an exact form. 
Then the residue of this third term vanishes.
Claim that \eqref{eq:2020.11.19.16.18} vanishes.
We show this claim as follows.
Set
$u_{\tilde{G},\infty}^{(l)}:= \delta_l(\tilde{G}_{\infty})\tilde{G}_{\infty}^{-1}$ 
and $\tilde{u}_{\tilde{G},\infty}^{(l)}:= \tilde{G}_{\infty}^{-1}\delta_l(\tilde{G}_{\infty})$ 
for $l \in \{ 1,2\}$.
Since $\tilde{G}_{\infty} = G_{n-2}^{-1} \tilde{G} G_1$,
we have $\tilde{u}^{(l)}_{\tilde{G},\infty} = G_1^{-1} \tilde{u}^{(l)}_{\tilde{G}} G_1$
and $u^{(l)}_{\tilde{G},\infty} = G_{n-2}^{-1} u^{(l)}_{\tilde{G}} G_{n-2}$ for $l=1,2$.
The meromorphic differential form
\begin{equation}\label{eq:globalmerodiff_unrami}
\begin{cases}
\frac{1}{2} \mathrm{Tr} \left( \delta_1(
\widehat\Omega_{(\boldsymbol{t}_{\mathrm{ra}},\boldsymbol{\theta}_0)}^{(1)})
 \tilde{u}^{(2)}_{\tilde{G}} -
\tilde{u}_{\tilde{G}}^{(1)}\delta_2(
\widehat\Omega_{(\boldsymbol{t}_{\mathrm{ra}},\boldsymbol{\theta}_0)}^{(1)})
  \right)
-\frac{1}{2}\mathrm{Tr} \left( \delta_1(
\widehat\Omega_{(\boldsymbol{t}_{\mathrm{ra}},\boldsymbol{\theta}_0)}^{(n-2)})
   u^{(2)}_{\tilde{G}} -
u^{(1)}_{\tilde{G}}\delta_2(
\widehat\Omega_{(\boldsymbol{t}_{\mathrm{ra}},\boldsymbol{\theta}_0)}^{(n-2)}) 
  \right)  \\
\qquad\qquad 
\text{ on $U_0\times ( \widehat{\mathcal{M}}_{\boldsymbol{t}_{\text{ra}}}
\times T_{\boldsymbol{\theta}})$} \\
\frac{1}{2} \mathrm{Tr} \left( G_1^{-1}\delta_1(
\widehat\Omega_{(\boldsymbol{t}_{\mathrm{ra}},\boldsymbol{\theta}_0)}^{(1)})G_1
 \tilde{u}^{(2)}_{\tilde{G},\infty} -
\tilde{u}_{\tilde{G},\infty}^{(1)} G_1^{-1} \delta_2(
(\widehat\Omega_{(\boldsymbol{t}_{\mathrm{ra}},\boldsymbol{\theta}_0)}^{(1)}) G_1
 \right) \\
\quad -\frac{1}{2}\mathrm{Tr} \left( G_{n-2}^{-1} \delta_1(
\widehat\Omega_{(\boldsymbol{t}_{\mathrm{ra}},\boldsymbol{\theta}_0)}^{(n-2)})G_{n-2}
   u^{(2)}_{\tilde{G},\infty} -
u^{(1)}_{\tilde{G},\infty}  G_{n-2}^{-1} \delta_2((
\widehat\Omega_{(\boldsymbol{t}_{\mathrm{ra}},\boldsymbol{\theta}_0)}^{(n-2)})G_{n-2}
   \right)  \\
\qquad\qquad \text{ on $U_{\infty} \times 
( \widehat{\mathcal{M}}_{\boldsymbol{t}_{\text{ra}}}\times T_{\boldsymbol{\theta}})$}
\end{cases}
\end{equation}
is a family of global meromorphic differential forms on $\mathbb{P}^1$
parametrized by 
$ \widehat{\mathcal{M}}_{\boldsymbol{t}_{\text{ra}}}\times T_{\boldsymbol{\theta}}$.
The differential forms
have poles at only $t_i$ and $q_j$ ($i\in I$ and $j=1,2,\ldots , n-3$).
The sums of residues \eqref{eq:2020.11.19.16.18}
are just the sums of all residues of
the global meromorphic differential forms
 \eqref{eq:globalmerodiff_unrami} on $\mathbb{P}^1$
 parametrized by $\widehat{\mathcal{M}}_{\boldsymbol{t}_{\text{ra}}}\times T_{\boldsymbol{\theta}}$. 
By the residue theorem, we have 
that \eqref{eq:2020.11.19.16.18}
is zero.
Finally we obtain the equality \eqref{eq:2020.11.19.16.08}.
\end{proof}

\begin{Thm}\label{2020.1.21.13.51}
\textit{
For the vector field $\delta_{\theta_{l,t_i}^{\pm}}^{\mathrm{IMD}}$, we have
$\hat{\omega}(\delta_{\theta_{l,t_i}^{\pm}}^{\mathrm{IMD}},\hat{\delta}) =0$ for any
vector field $\hat{\delta} 
\in \Theta_{\widehat{\mathcal{M}}_{\boldsymbol{t}_{\text{ra}}}\times T_{\boldsymbol{\theta}}}$.
Moreover, 
for the vector field $\delta_{t_i}^{\mathrm{IMD}}$, we have
$\hat{\omega}(\delta_{t_i}^{\mathrm{IMD}}, \hat{\delta}) =0$ for any
vector field $\hat{\delta} \in 
\Theta_{\widehat{\mathcal{M}}_{\boldsymbol{t}_{\text{ra}}}\times T_{\boldsymbol{\theta}}}$.
}
\end{Thm}

\begin{proof}

By the equality \eqref{eq:2020.11.19.16.08}, we have 
\begin{equation}\label{2021.4.16.12.45}
\begin{split}
\hat\omega (\delta_{\theta_{l,t_i}^{\pm}}^{\mathrm{IMD}}, \hat\delta) &=
 \frac{1}{2} \sum_{i' \in I}   \mathrm{res}_{x=t_{i'}} 
\mathrm{Tr} \left( \delta_{\theta_{l,t_i}^{\pm}}^{\mathrm{IMD}}
(\widehat\Omega_{(\boldsymbol{t}_{\mathrm{ra}},\boldsymbol{\theta}_0)}^{(1)}) 
\hat\delta(\tilde\psi_{i'})(\tilde\psi_{i'})^{-1} - 
\delta_{\theta_{l,t_i}^{\pm}}^{\mathrm{IMD}}(\tilde\psi_{i'}) (\tilde\psi_{i'})^{-1}
\hat\delta (\widehat\Omega_{(\boldsymbol{t}_{\mathrm{ra}},\boldsymbol{\theta}_0)}^{(1)}) 
  \right).
\end{split}
\end{equation}
Here $\delta_{\theta_{l,t_i}^{\pm}}^{\mathrm{IMD}}
(\widehat\Omega_{(\boldsymbol{t}_{\mathrm{ra}},\boldsymbol{\theta}_0)}^{(1)}) $
appears in the $\epsilon$-term of the morphism
\eqref{2021.4.16.12.11}.

Now we consider 
 replacement of 
 $\delta_{\theta_{l,t_i}^{\pm}}^{\mathrm{IMD}}(\tilde\psi_{i'}) (\tilde\psi_{i'})^{-1}$
 in \eqref{2021.4.16.12.45} for each $i' \in I$.
 We will show that 
 we may replace
 $\delta_{\theta_{l,t_i}^{\pm}}^{\mathrm{IMD}}(\tilde\psi_{i'}) (\tilde\psi_{i'})^{-1}$
with $\Upsilon_{\theta_{l,t_i}^{\pm}}^{\mathrm{IMD}}$ 
as follows.
Here $\Upsilon_{\theta_{l,t_i}^{\pm}}^{\mathrm{IMD}}$ 
appeared in \eqref{2021.4.16.12.12}.
We take an analytic open subset $V$ of 
$\widehat{\mathcal{M}}_{\boldsymbol{t}_{\text{ra}}}\times T_{\boldsymbol{\theta}}$.
We take an inverse image of $V$ under the projection
 $$
p_{\mathbb{P}^1} \colon
\mathbb{P}^1 \times
( \widehat{\mathcal{M}}_{\boldsymbol{t}_{\text{ra}}}\times T_{\boldsymbol{\theta}})
\longrightarrow
\widehat{\mathcal{M}}_{\boldsymbol{t}_{\text{ra}}}\times T_{\boldsymbol{\theta}}.
 $$
Let $\widehat{\Delta}_{i'}^{\text{an}}$ ($i' \in I$) be an analytic open
subset of the inverse image $p_{\mathbb{P}^1}^{-1} (V)$
such that $\tilde{t}_{i'} \cap p_{\mathbb{P}^1}^{-1} (V) \subset \widehat{\Delta}_{i'}^{\text{an}}$ 
and the fibers of 
$p_{\mathbb{P}^1} |_{\widehat{\Delta}_{i'}^{\text{an}}} \colon 
\widehat{\Delta}_{i'}^{\text{an}} \rightarrow V$ 
for each point of $V$
are unit disks such that $x_{t_i}$ gives a coordinate 
of the unit disks.
Let $U_t^{\text{an}}$ be an analytic open subset of $\mathbb{C}^1 = 
\mathrm{Spec}\,  \mathbb{C}[t]$ such that $ 0 \in U_t^{\text{an}} $ and 
$U_t^{\text{an}}$ is small enough.
We consider the restriction 
of \eqref{2021.4.16.12.11} to 
$\widehat{\Delta}_{i'}^{\text{an}}
\times \mathrm{Spec}\, \mathbb{C}[\epsilon] $.
This is a morphism $\widehat{E}_1^{\epsilon}|_{\widehat{\Delta}_{i'}^{\text{an}} 
\times \mathrm{Spec}\, \mathbb{C}[\epsilon]} 
\rightarrow (\widehat{E}_1^{\epsilon}
 \otimes \tilde\Omega^1_{\partial/\partial \theta_{l,t_i}^{\pm}})
 |_{\widehat{\Delta}_{i'}^{\text{an}}
\times \mathrm{Spec}\, \mathbb{C}[\epsilon]}$.
Let $\widehat{E}_1^{t}$
be the pull-back of $E_1$ under the first projection
$$
\widehat{\Delta}_{i'}^{\text{an}}  \times U_t^{\text{an}}
\hookrightarrow 
\mathbb{P}^1 \times
( \widehat{\mathcal{M}}_{\boldsymbol{t}_{\text{ra}}}\times T_{\boldsymbol{\theta}}) \times U_t^{\text{an}}
\rightarrow \mathbb{P}^1.
$$
Let
$D(\tilde{\boldsymbol{t}}_0)_{\epsilon}$
be the pull-back of $D(\tilde{\boldsymbol{t}}_0)$
under the composition \eqref{2021.4.20.16.11}.
We take divisor
$D(\tilde{\boldsymbol{t}}_0)_{t}$ 
on $\widehat{\Delta}_{i'}^{\text{an}}  \times U_t^{\text{an}}$
such that  the pull-back of $D(\tilde{\boldsymbol{t}}_0)_t$
under the map
$\widehat{\Delta}_{i'}^{\text{an}} \times \mathrm{Spec}\, \mathbb{C}[\epsilon ]
\rightarrow  \widehat{\Delta}_{i'}^{\text{an}}  \times U_t^{\text{an}}$
is
$D(\tilde{\boldsymbol{t}}_0)_{\epsilon}
|_{\widehat{\Delta}_{i'}^{\text{an}} \times \mathrm{Spec}\, \mathbb{C}[\epsilon ]}$.
Here this map 
$\widehat{\Delta}_{i'}^{\text{an}} \times \mathrm{Spec}\, \mathbb{C}[\epsilon]
\rightarrow  
\widehat{\Delta}_{i'}^{\text{an}}
\times U_t^{\text{an}}$
is given by the substitution $t= \epsilon$.
We take a relative connection on $\widehat{E}_1^t$:
$$
\widehat{E}_1^t
\longrightarrow
\widehat{E}_1^t \otimes
\Omega^1_{\widehat{\Delta}_{i'}^{\text{an}}  \times U_t^{\text{an}}/
\widehat{\mathcal{M}}_{\boldsymbol{t}_{\text{ra}}}\times T_{\boldsymbol{\theta}} \times U_t^{\text{an}}}
(D(\tilde{\boldsymbol{t}}_0)_{t})
$$
such that the pull-back of this relative connection on $\widehat{E}_1^t$ under the map 
$\widehat{\Delta}_{i'}^{\text{an}} \times \mathrm{Spec}\, \mathbb{C}[\epsilon]
\rightarrow  
\widehat{\Delta}_{i'}^{\text{an}}
\times U_t^{\text{an}}$
is just the restriction 
of \eqref{2021.4.16.12.11}.
We denote by
$\widehat\Omega_{(\boldsymbol{t}_{\mathrm{ra}},\boldsymbol{\theta}_0)}^{(1)}
(x_{t_{i'}},t)$
the connection matrix of the relative connection,
where 
$\widehat\Omega_{(\boldsymbol{t}_{\mathrm{ra}},\boldsymbol{\theta}_0)}^{(1)}
(x_{t_{i'}},0)=
\widehat\Omega_{(\boldsymbol{t}_{\mathrm{ra}},\boldsymbol{\theta}_0)}^{(1)}
|_{\widehat{\Delta}_{i'}^{\text{an}} }$
and
$\frac{\partial}{\partial t}
\widehat\Omega_{(\boldsymbol{t}_{\mathrm{ra}},\boldsymbol{\theta}_0)}^{(1)}
(x_{t_{i'}},t)|_{t=0}
= \delta_{\theta_{l,t_i}^{\pm}}^{\mathrm{IMD}}
(\widehat\Omega_{(\boldsymbol{t}_{\mathrm{ra}},\boldsymbol{\theta}_0)}^{(1)}) 
|_{\widehat{\Delta}_{i'}^{\text{an}}}$.
Let $\widehat{\Sigma} \subset \widehat{\Delta}^{\text{an}}_{i'}$ 
be a family of sufficiently small sectors in 
$\widehat{\Delta}^{\text{an}}_{i'}
\rightarrow 
\widehat{\mathcal{M}}_{\boldsymbol{t}_{\text{ra}}}\times T_{\boldsymbol{\theta}}$. 
We take a fundamental matrix
solution $\Psi_{\widehat{\Sigma}}(x_{t_{i'}},t)$ 
on $ \widehat{\Sigma}  \times U_t^{\text{an}}$
of the connection 
$$
d\Psi_{\widehat\Sigma}(x_{t_{i'}},t)
+\widehat\Omega_{(\boldsymbol{t}_{\mathrm{ra}},\boldsymbol{\theta}_0)}^{(1)}
(x_{t_{i'}},t)
\Psi_{\widehat\Sigma}(x_{t_{i'}},t) =0
$$
with uniform asymptotic relation 
$$
\Psi_{\widehat\Sigma}(x_{t_{i'}},t) \exp\left( \Lambda^-_{i'}(x_{t_{i'}} ,t )\right) \sim
\widehat{P}_{i'}(x_{t_{i'}},t)
\quad (x_{t_{i'}} \rightarrow 0, x_{t_{i'}} \in \widehat\Sigma).
$$
Here we set
$$
\Lambda^-_{i'}(x_{t_{i'}}) :=\sum_{l=0}^{n_i-1}
\begin{pmatrix}
\theta_{l,t_{i'}}^{+}  \int x_{t_{i'}}^{-n_{i'}+l} dx_{t_{i'}}  & 0 \\
0 & \theta_{l,t_{i'}}^{-}  \int x_{t_{i'}}^{-n_{i'}+l} dx_{t_{i'}}
\end{pmatrix}
$$
and we take
$$
\begin{aligned}
&\widehat{P}_{i'}(x_{t_{i'}},t) =\widehat{P}_{i',0}(t) + \widehat{P}_{i',1}(t) x_{t_{i'}}+\cdots
\quad \text{and} \quad \\
&\Lambda^-_{i'}(x_{t_{i'}} ,t )= \sum_{l=0}^{n_i-1}
\begin{pmatrix}
\theta_{l,t_{i'}}^{+} (t) \int x_{t_{i'}}^{-n_{i'}+l} dx_{t_{i'}}  & 0 \\
0 & \theta_{l,t_{i'}}^{-} (t) \int x_{t_{i'}}^{-n_{i'}+l} dx_{t_{i'}}
\end{pmatrix}
\end{aligned}
$$
so that the expansions of $\Lambda^-_{i'}(x_{t_{i'}} ,\epsilon )$
and $\widehat{P}_{i'}(x_{t_{i'}},\epsilon)$ 
with respect to $\epsilon$ are the following:
$$
\begin{aligned}
\Lambda^-_{i'}(x_{t_{i'}} ,\epsilon ) &= \Lambda^-_{i'}(x_{t_{i'}})
+\epsilon \cdot \delta_{\theta_{l,t_i}^{\pm}}^{\mathrm{IMD}} ( \Lambda^-_{i'}(x_{t_{i'}}) ) \\
\widehat{P}_{i'}(x_{t_{i'}},\epsilon) 
&= \tilde\psi_{i'} \exp\left(\Lambda^-_{i'}(x_{t_{i'}})\right) 
+\epsilon \cdot \delta_{\theta_{l,t_i}^{\pm}}^{\mathrm{IMD}}
\left(
\tilde\psi_{i'} \exp\left( \Lambda^-_{i'}(x_{t_{i'}})\right)
\right).
\end{aligned}
$$
The uniform asymptotic relation means that
$$
\lim_{x_{t_{i'}} \rightarrow 0, x_{t_{i'}} \in 
\Gamma_{\widehat\Sigma}} 
\frac{\| \Psi_{\widehat\Sigma}(x_{t_{i'}},t) \exp\left(\Lambda^-_{i'}(x_{t_{i'}},t)\right)
- \sum_{j=0}^{N}\widehat{P}_{i',j}x_{t_{i'}}^j
 \|}{|x_{t_{i'}}|^N}=0, \quad (\text{uniformally})
$$
for any positive integer $N$.
We may check that 
\begin{equation}\label{2021.4.17.19.02}
\left.\frac{\partial \Psi_{\widehat\Sigma}(x_{t_{i'}},t)}{\partial t} 
\Psi_{\widehat\Sigma}(x_{t_{i'}},t)^{-1}  \right|_{t=0} 
\sim 
\delta_{\theta_{l,t_{i}}^{\pm}}^{\mathrm{IMD}}(\tilde\psi_{i'}) (\tilde\psi_{i'})^{-1}.
\end{equation}
By the integrable condition \eqref{2021.4.16.14.37},
we may take a fundamental matrix
solution $\Psi^{\text{flat}}_{\widehat\Sigma}(x_{t_{i'}},t)$ 
on $\widehat\Sigma \times U_t^{\text{an}}$
of the connection 
$d+\widehat\Omega_{(\boldsymbol{t}_{\mathrm{ra}},\boldsymbol{\theta}_0)}^{(1)}
(x_{t_{i'}},t)=0$ 
such that $\Psi^{\text{flat}}_{\widehat\Sigma}(x_{t_{i'}},t)$ satisfies
$\Psi^{\text{flat}}_{\widehat\Sigma}(x_{t_{i'}},0)
=\Psi_{\widehat\Sigma}(x_{t_{i'}},0)$ and
$$
\left.\frac{\partial \Psi^{\text{flat}}_{\widehat\Sigma}(x_{t_{i'}},t)}{\partial t} 
\Psi^{\text{flat}}_{\widehat\Sigma}(x_{t_{i'}},t)^{-1}  \right|_{t=0} 
=\Upsilon_{\theta_{l,t_{i}}^{\pm}}^{\mathrm{IMD}}|_{\widehat{\Sigma}}.
$$
There exists a matrix $C_{t_{i'}} (t)$ 
such that $C_{t_{i'}} (t)$ is independent of $x_{t_{i'}}$ and
$\Psi^{\text{flat}}_{\widehat\Sigma}(x_{t_{i'}},t) 
= \Psi_{\widehat\Sigma}(x_{t_{i'}},t)C_{t_{i'}} (t)$.
We calculate 
$\Upsilon_{\theta_{l,t_{i}}^{\pm}}^{\mathrm{IMD}}|_{\widehat{\Sigma}}$ as follows:
\begin{equation}\label{2021.4.17.19.03}
\begin{aligned}
\Upsilon_{\theta_{l,t_{i}}^{\pm}}^{\mathrm{IMD}}|_{\widehat{\Sigma}} &=\left.
\frac{\partial \Psi^{\text{flat}}_{\widehat{\Sigma}}(x_{t_{i'}},t)}{\partial t} 
\Psi^{\text{flat}}_{\widehat{\Sigma}}(x_{t_{i'}},t)^{-1}  \right|_{t=0} \\
&= 
\left.\frac{\partial \Psi_{\widehat{\Sigma}}(x_{t_{i'}},t)}{\partial t} 
\Psi_{\widehat{\Sigma}}(x_{t_{i'}},t)^{-1}  \right|_{t=0} 
+\Psi_{\widehat{\Sigma}}(x_{t_{i'}},0) 
\left.
\left( \frac{\partial C_{t_{i'}}(t)}{\partial t} 
C_{t_{i'}}(t)^{-1} \right) \right|_{t=0} 
\Psi_{\widehat{\Sigma}}(x_{t_{i'}},0)^{-1}  .
\end{aligned}
\end{equation}
Set
$$
\tilde{C}_{t_{i'}}(x_{t_{i'}}) := 
\left.
\exp\left( -  \Lambda^-_{i'}(x_{t_{i'}},t)\right)
\left( \frac{\partial C_{t_{i'}}(t)}{\partial t} 
C_{t_{i'}}(t)^{-1} \right) \exp\left( \Lambda^-_{i'}(x_{t_{i'}},t)\right) \right|_{t=0}.
$$
By the equalities \eqref{2021.4.17.19.02} and \eqref{2021.4.17.19.03},
we have 
\begin{equation*}
\tilde{C}_{t_{i'}}(x_{t_{i'}}) \sim 
\widehat{P}_{i'}(x_{t_{i'}},0)^{-1} 
\Upsilon_{\theta_{l,t_{i}}^{\pm}}^{\mathrm{IMD}}
\widehat{P}_{i'}(x_{t_{i'}},0)
- \widehat{P}_{i'}(x_{t_{i'}},0)^{-1}
\delta_{\theta_{l,t_{i}}^{\pm}}^{\mathrm{IMD}}(\psi'_{i'}) (\psi'_{i'})^{-1}
\widehat{P}_{i'}(x_{t_{i'}},0).
\end{equation*}
By this asymptotic relation, we have that
$x_{t_{i'}}^{n_{i'}}\tilde{C}_{t_{i'}}(x_{t_{i'}})$ is 
bounded on $\widehat\Sigma \times U_{t}^{\text{an}}$.
Then we may check that
$\left(\frac{\partial C_{t_{i'}} (t)}{\partial t}C_{t_{i'}} (t)^{-1} \right)|_{t=0}$ is 
a triangular matrix and $\tilde{C}_{t_{i'}}(x_{t_{i'}})
\sim \tilde{C}^{\text{diag}}_{t_{i'}}$,
where $\tilde{C}^{\text{diag}}_{t_{i'}}$ 
is a diagonal matrix
and $\tilde{C}^{\text{diag}}_{t_{i'}}$ is independent of $x_{t_{i'}}$.
By combining this asymptotic relation,
the asymptotic relation \eqref{2021.4.17.19.02}
and the equality \eqref{2021.4.17.19.03},
we have the following asymptotic relation:
\begin{equation}\label{2021.4.18.17.41_1}
\begin{aligned}
\delta_{\theta_{l,t_i}^{\pm}}^{\mathrm{IMD}}(\tilde\psi_{i'}) (\tilde\psi_{i'})^{-1}
&\sim 
\Upsilon_{\theta_{l,t_i}^{\pm}}^{\mathrm{IMD}}
-
\widehat{P}_{i'}(x_{t_{i'}},0)
 \tilde{C}^{\text{diag}}_{t_{i'}}
\widehat{P}_{i'}(x_{t_{i'}},0)^{-1} .
\end{aligned}
\end{equation}
So we have 
\begin{equation}\label{2021.4.18.12.36}
\begin{split}
 &  \mathrm{res}_{x=t_{i'}} 
\mathrm{Tr} \left( \delta_{\theta_{l,t_i}^{\pm}}^{\mathrm{IMD}}
(\widehat\Omega_{(\boldsymbol{t}_{\mathrm{ra}},\boldsymbol{\theta}_0)}^{(1)}) 
\hat\delta(\tilde\psi_{i'})(\tilde\psi_{i'})^{-1} - 
\delta_{\theta_{l,t_i}^{\pm}}^{\mathrm{IMD}}(\tilde\psi_{i'}) (\tilde\psi_{i'})^{-1}
\hat\delta (\widehat\Omega_{(\boldsymbol{t}_{\mathrm{ra}},\boldsymbol{\theta}_0)}^{(1)}) 
  \right)\\
&= \mathrm{res}_{x=t_{i'}} 
\mathrm{Tr} \left( \delta_{\theta_{l,t_i}^{\pm}}^{\mathrm{IMD}}
(\widehat\Omega_{(\boldsymbol{t}_{\mathrm{ra}},\boldsymbol{\theta}_0)}^{(1)}) 
\hat\delta(\tilde\psi_{i'})(\tilde\psi_{i'})^{-1} - 
\Upsilon_{\theta_{l,t_i}^{\pm}}^{\mathrm{IMD}}
 \hat\delta (\widehat\Omega_{(\boldsymbol{t}_{\mathrm{ra}},\boldsymbol{\theta}_0)}^{(1)})  
  \right) \\
  &\qquad  +\mathrm{res}_{x=t_{i'}} 
\mathrm{Tr} \left( 
\widehat{P}_{i'}(x_{t_{i'}},0)
 \tilde{C}^{\text{diag}}_{t_{i'}}
\widehat{P}_{i'}(x_{t_{i'}},0)^{-1} 
 \hat\delta (\widehat\Omega_{(\boldsymbol{t}_{\mathrm{ra}},\boldsymbol{\theta}_0)}^{(1)})  
\right) \\
&= \mathrm{res}_{x=t_{i'}} 
\mathrm{Tr} \left( \delta_{\theta_{l,t_i}^{\pm}}^{\mathrm{IMD}}
(\widehat\Omega_{(\boldsymbol{t}_{\mathrm{ra}},\boldsymbol{\theta}_0)}^{(1)}) 
\hat\delta(\tilde\psi_{i'})(\tilde\psi_{i'})^{-1} 
\right)- 
\mathrm{res}_{x=t_{i'}} 
\mathrm{Tr} \left( 
\Upsilon_{\theta_{l,t_i}^{\pm}}^{\mathrm{IMD}}
 \hat\delta (\widehat\Omega_{(\boldsymbol{t}_{\mathrm{ra}},\boldsymbol{\theta}_0)}^{(1)})  
  \right) .
\end{split}
\end{equation}
Here we may check 
the last equality as follows.
By the same calculations as in
\eqref{2021.4.18.17.39},
we have
\begin{equation*}
\begin{aligned}
\mathrm{res}_{x=t_{i'}} 
\mathrm{Tr} \left( 
\widehat{P}_{i'}(x_{t_{i'}},0)
 \tilde{C}^{\text{diag}}_{t_{i'}}
\widehat{P}_{i'}(x_{t_{i'}},0)^{-1} 
 \hat\delta (\widehat\Omega_{(\boldsymbol{t}_{\mathrm{ra}},\boldsymbol{\theta}_0)}^{(1)})  
\right)
&=\mathrm{res}_{x=t_{i'}} \mathrm{Tr} \left( 
\tilde\psi_{i'} 
 \tilde{C}^{\text{diag}}_{t_{i'}}
 \tilde\psi_{i'}^{-1} 
\cdot
 \hat\delta (\widehat\Omega_{(\boldsymbol{t}_{\mathrm{ra}},\boldsymbol{\theta}_0)}^{(1)})  
\right)\\
& =\mathrm{res}_{x=t_{i'}} \mathrm{Tr} \left( 
 \tilde{C}^{\text{diag}}_{t_{i'}}
\cdot d \Big((\tilde\psi_{i'} )^{-1}
 \hat\delta (\tilde\psi_{i'}) \Big)  
\right) .
\end{aligned}
\end{equation*}
Here remark that
$\tilde{C}^{\text{diag}}_{t_{i'}}$ and $\mathrm{exp} (-\Lambda^-_{i'}(x_{t_{i'}}))$ are diagonal,
and $\tilde\psi_{i'}=\widehat{P}_{i'}(x_{t_{i'}},0) \exp(-\Lambda^-_{i'}(x_{t_{i'}}))$.
This
is equal to 
\begin{equation}\label{2021.4.18.18.02}
\begin{aligned}
&\mathrm{res}_{x=t_{i'}} \mathrm{Tr} \Big(
 \tilde{C}^{\text{diag}}_{t_{i'}}
d\Big(\mathrm{exp} (-\Lambda^-_{i'}(x_{t_{i'}}))^{-1}
\widehat{P}_{i'}(x_{t_{i'}},0)^{-1} \hat{\delta}(\widehat{P}_{i'}(x_{t_{i'}},0)) 
\mathrm{exp} (-\Lambda^-_{i'}(x_{t_{i'}})) \Big)
\Big) \\
& \qquad +\mathrm{res}_{x=t_{i'}}\mathrm{Tr} \Big(
 \tilde{C}^{\text{diag}}_{t_{i'}}
d\Big(\mathrm{exp} (-\Lambda^-_{i'}(x_{t_{i'}}))^{-1}
 \hat{\delta}(  \mathrm{exp} (-\Lambda^-_{i'}(x_{t_{i'}}))) \Big)
\Big)\\
&=\mathrm{res}_{x=t_{i'}}\mathrm{Tr} \Big(
 \tilde{C}^{\text{diag}}_{t_{i'}}
d\Lambda^-_{i'}(x_{t_{i'}})
\widehat{P}_{i'}(x_{t_{i'}},0)^{-1} \hat{\delta}(\widehat{P}_{i'}(x_{t_{i'}},0))  
\Big) \\
&\quad- \mathrm{res}_{x=t_{i'}}\mathrm{Tr} \Big(
 \tilde{C}^{\text{diag}}_{t_{i'}}
\widehat{P}_{i'}(x_{t_{i'}},0)^{-1} \hat{\delta}(\widehat{P}_{i'}(x_{t_{i'}},0))
d\Lambda^-_{i'}(x_{t_{i'}})  
\Big) \\
& \qquad + \mathrm{res}_{x=t_{i'}}
\mathrm{Tr} \Big(
 \tilde{C}^{\text{diag}}_{t_{i'}}
d\Big(
\widehat{P}_{i'}(x_{t_{i'}},0)^{-1} \hat{\delta}(\widehat{P}_{i'}(x_{t_{i'}},0))  \Big)
 \Big) 
+\mathrm{res}_{x=t_{i'}}\mathrm{Tr} \Big(
\tilde{C}^{\text{diag}}_{t_{i'}}
d\Big( \hat{\delta}(  -\Lambda^-_{i'}(x_{t_{i'}}))\Big)
\Big) \\
&=\mathrm{res}_{x=t_{i'}}\mathrm{Tr} \Big(
\tilde{C}^{\text{diag}}_{t_{i'}} d\Big(
\widehat{P}_{i'}(x_{t_{i'}},0)^{-1} \hat{\delta}(\widehat{P}_{i'}(x_{t_{i'}},0))  \Big)
\Big)  
+\mathrm{res}_{x=t_{i'}}\mathrm{Tr} \Big(\tilde{C}^{\text{diag}}_{t_{i'}}
d\Big( \hat{\delta}(  -\Lambda^-_{i'}(x_{t_{i'}}))\Big)
\Big).
\end{aligned}
\end{equation}
Since $\tilde{C}^{\text{diag}}_{t_{i'}}$
is independent of $x_{t_{i'}}$,
then the last line of 
\eqref{2021.4.18.18.02} is zero.
Then we have
the last equality of \eqref{2021.4.18.12.36}.
This means that
we may replace
 $\delta_{\theta_{l,t_i}^{\pm}}^{\mathrm{IMD}}(\tilde\psi_{i'}) (\tilde\psi_{i'})^{-1}$
with $\Upsilon_{\theta_{l,t_i}^{\pm}}^{\mathrm{IMD}}$.

Next we will calculate $\mathrm{Tr} (\delta_{\theta_{l,t_i}^{\pm}}^{\mathrm{IMD}}
(\widehat\Omega_{(\boldsymbol{t}_{\mathrm{ra}},\boldsymbol{\theta}_0)}^{(1)}) 
\hat\delta(\tilde\psi_{i'})(\tilde\psi_{i'})^{-1})$.
By taking variations of the both hand sides of
$d\tilde\psi_{i'} =- \widehat\Omega_{(\boldsymbol{t}_{\mathrm{ra}},\boldsymbol{\theta}_0)}^{(1)} 
\tilde\psi_{i'}$,
we have
\begin{equation}\label{2021.4.18.12.40}
\widehat\Omega_{(\boldsymbol{t}_{\mathrm{ra}},\boldsymbol{\theta}_0)}^{(1)}
\hat\delta(\tilde\psi_{i'})(\tilde\psi_{i'})^{-1}
=  -d (\hat\delta(\tilde\psi_{i'}) ) (\tilde\psi_{i'})^{-1} -
\hat\delta(\widehat\Omega_{(\boldsymbol{t}_{\mathrm{ra}},\boldsymbol{\theta}_0)}^{(1)}).
\end{equation}
By the integrable condition \eqref{2021.4.16.14.37}, we have
\begin{equation}\label{2021.4.18.12.38}
\begin{aligned}
&\mathrm{Tr} \left( \delta_{\theta_{l,t_i}^{\pm}}^{\mathrm{IMD}}
(\widehat\Omega_{(\boldsymbol{t}_{\mathrm{ra}},\boldsymbol{\theta}_0)}^{(1)}) 
\hat\delta(\tilde\psi_{i'})(\tilde\psi_{i'})^{-1} \right)\\
&=\mathrm{Tr} \left( \left( 
d  \Upsilon_{\theta_{l,t_i}^{\pm}}^{\mathrm{IMD}}
+[\widehat\Omega_{(\boldsymbol{t}_{\mathrm{ra}},\boldsymbol{\theta}_0)}^{(1)}  ,
\Upsilon_{\theta_{l,t_i}^{\pm}}^{\mathrm{IMD}}]
\right)\hat\delta(\tilde\psi_{i'})(\tilde\psi_{i'})^{-1} \right) \\
&=  \mathrm{Tr} \left(
d  (\Upsilon_{\theta_{l,t_i}^{\pm}}^{\mathrm{IMD}})
\hat\delta(\tilde\psi_{i'})(\tilde\psi_{i'})^{-1}
+\widehat\Omega_{(\boldsymbol{t}_{\mathrm{ra}},\boldsymbol{\theta}_0)}^{(1)}
 \Upsilon_{\theta_{l,t_i}^{\pm}}^{\mathrm{IMD}}
\hat\delta(\tilde\psi_{i'})(\tilde\psi_{i'})^{-1}
-\Upsilon_{\theta_{l,t_i}^{\pm}}^{\mathrm{IMD}} 
\widehat\Omega_{(\boldsymbol{t}_{\mathrm{ra}},\boldsymbol{\theta}_0)}^{(1)}
\hat\delta(\tilde\psi_{i'})(\tilde\psi_{i'})^{-1} \right) \\
&=  \mathrm{Tr} \left(
d  (\Upsilon_{\theta_{l,t_i}^{\pm}}^{\mathrm{IMD}})
\hat\delta(\tilde\psi_{i'})(\tilde\psi_{i'})^{-1}
- d\tilde\psi_{i'} (\tilde\psi_{i'})^{-1} 
 \Upsilon_{\theta_{l,t_i}^{\pm}}^{\mathrm{IMD}}
\hat\delta(\tilde\psi_{i'})(\tilde\psi_{i'})^{-1}
+\Upsilon_{\theta_{l,t_i}^{\pm}}^{\mathrm{IMD}} 
 d (\hat\delta(\tilde\psi_{i'}) ) (\tilde\psi_{i'})^{-1} 
+\Upsilon_{\theta_{l,t_i}^{\pm}}^{\mathrm{IMD}}
\hat\delta(\widehat\Omega_{(\boldsymbol{t}_{\mathrm{ra}},\boldsymbol{\theta}_0)}^{(1)}) 
\right)\\
&=   \mathrm{Tr} \left(
d  (\Upsilon_{\theta_{l,t_i}^{\pm}}^{\mathrm{IMD}})
\hat\delta(\tilde\psi_{i'})(\tilde\psi_{i'})^{-1}
+\Upsilon_{\theta_{l,t_i}^{\pm}}^{\mathrm{IMD}} 
 d (\hat\delta(\tilde\psi_{i'}) ) (\tilde\psi_{i'})^{-1}  
 + \Upsilon_{\theta_{l,t_i}^{\pm}}^{\mathrm{IMD}}
\hat\delta(\tilde\psi_{i'})d\left( (\tilde\psi_{i'})^{-1} \right)
 \right)
+ \mathrm{Tr} \left(
\Upsilon_{\theta_{l,t_i}^{\pm}}^{\mathrm{IMD}}
\hat\delta(\widehat\Omega_{(\boldsymbol{t}_{\mathrm{ra}},\boldsymbol{\theta}_0)}^{(1)}) 
\right)\\
&= \mathrm{Tr} \left(
d  \left(\Upsilon_{\theta_{l,t_i}^{\pm}}^{\mathrm{IMD}}
\hat\delta(\tilde\psi_{i'})(\tilde\psi_{i'})^{-1}  \right)\right)
+ \mathrm{Tr} \left(
\Upsilon_{\theta_{l,t_i}^{\pm}}^{\mathrm{IMD}}
\hat\delta(\widehat\Omega_{(\boldsymbol{t}_{\mathrm{ra}},\boldsymbol{\theta}_0)}^{(1)}) 
\right).
\end{aligned}
\end{equation}
Here the third equality of \eqref{2021.4.18.12.38} is given by
\eqref{2021.4.18.12.40}.
By combining \eqref{2021.4.18.12.36} and \eqref{2021.4.18.12.38},
we have the following equalities:
\begin{equation}\label{2021.4.18.12.45}
\begin{split}
 &  \mathrm{res}_{x=t_{i'}} 
\mathrm{Tr} \left( \delta_{\theta_{l,t_i}^{\pm}}^{\mathrm{IMD}}
(\widehat\Omega_{(\boldsymbol{t}_{\mathrm{ra}},\boldsymbol{\theta}_0)}^{(1)}) 
\hat\delta(\tilde\psi_{i'})(\tilde\psi_{i'})^{-1} - 
\delta_{\theta_{l,t_i}^{\pm}}^{\mathrm{IMD}}(\tilde\psi_{i'}) (\tilde\psi_{i'})^{-1}
\hat\delta (\widehat\Omega_{(\boldsymbol{t}_{\mathrm{ra}},\boldsymbol{\theta}_0)}^{(1)}) 
  \right)\\
&= \mathrm{res}_{x=t_{i'}} 
\mathrm{Tr} \left( \delta_{\theta_{l,t_i}^{\pm}}^{\mathrm{IMD}}
(\widehat\Omega_{(\boldsymbol{t}_{\mathrm{ra}},\boldsymbol{\theta}_0)}^{(1)}) 
\hat\delta(\tilde\psi_{i'})(\tilde\psi_{i'})^{-1} 
\right)- 
\mathrm{res}_{x=t_{i'}} 
\mathrm{Tr} \left( 
\Upsilon_{\theta_{l,t_i}^{\pm}}^{\mathrm{IMD}}
 \hat\delta (\widehat\Omega_{(\boldsymbol{t}_{\mathrm{ra}},\boldsymbol{\theta}_0)}^{(1)})  
  \right) \\
  &=\mathrm{res}_{x=t_{i'}}  \mathrm{Tr} \left(
d  \left(\Upsilon_{\theta_{l,t_i}^{\pm}}^{\mathrm{IMD}}
\hat\delta(\tilde\psi_{i'})(\tilde\psi_{i'})^{-1}  \right)\right)=0.
\end{split}
\end{equation}
By combining \eqref{2021.4.16.12.45} and \eqref{2021.4.18.12.45},
we have 
\begin{equation*}
\hat\omega (\delta_{\theta_{l,t_i}^{\pm}}^{\mathrm{IMD}}, \hat\delta) =
 \frac{1}{2} \sum_{i' \in I}   \mathrm{res}_{x=t_{i'}} 
\mathrm{Tr} \left( \delta_{\theta_{l,t_i}^{\pm}}^{\mathrm{IMD}}
(\widehat\Omega_{(\boldsymbol{t}_{\mathrm{ra}},\boldsymbol{\theta}_0)}^{(1)}) 
\hat\delta(\tilde\psi_{i'})(\tilde\psi_{i'})^{-1} - 
\delta_{\theta_{l,t_i}^{\pm}}^{\mathrm{IMD}}(\tilde\psi_{i'}) (\tilde\psi_{i'})^{-1}
\hat\delta (\widehat\Omega_{(\boldsymbol{t}_{\mathrm{ra}},\boldsymbol{\theta}_0)}^{(1)}) 
  \right) =0.
\end{equation*}

Next we will show that 
$\hat\omega (\delta_{t_i}^{\mathrm{IMD}}, \hat\delta) =0$
for any $\hat{\delta}$.
For this purpose,
we show that 
\begin{equation}\label{2021.11.16.17.14}
\hat\omega (\delta_{t_i}^{\mathrm{IMD}}, \hat\delta) =
 \frac{1}{2} \sum_{i' \in I}   \mathrm{res}_{x=t_{i'}} 
\mathrm{Tr} \left( \delta_{t_i}^{\mathrm{IMD}}
(\widehat\Omega_{(\boldsymbol{t}_{\mathrm{ra}},\boldsymbol{\theta}_0)}^{(1)}) 
\hat\delta(\tilde\psi_{i'})(\tilde\psi_{i'})^{-1} - 
\Upsilon_{t_i}^{\mathrm{IMD}}
\hat\delta (\widehat\Omega_{(\boldsymbol{t}_{\mathrm{ra}},\boldsymbol{\theta}_0)}^{(1)}) 
  \right) 
\end{equation}
for any $\hat{\delta}$.
Here $\Upsilon_{t_i}^{\mathrm{IMD}}$ 
appeared in \eqref{2021.4.16.12.12_x}.
We define $\widehat{\Delta}_{i'}^{\text{an}}$ and $U_t^{\text{an}}$ as above.
Let
$D(\tilde{\boldsymbol{t}}_0)_{\epsilon}$
be the pull-back of $D(\tilde{\boldsymbol{t}}_0)$
under the composition \eqref{2021.4.20.16.11} with respect to
the vector field $\partial / \partial t_i$.
We define a divisor $D(\tilde{\boldsymbol{t}}_0)_{t}$ on 
$\widehat{\Delta}_{i'}^{\text{an}}\times U_t^{\text{an}}$ 
for this divisor $D(\tilde{\boldsymbol{t}}_0)_{\epsilon}$.
We take a vector bundle $\widehat{E}_1^t$ on $\widehat{\Delta}_{i'}^{\text{an}}\times U_t^{\text{an}}$
and a relative connection
$$
\widehat{E}_1^t
\longrightarrow
\widehat{E}_1^t \otimes
\Omega^1_{\widehat{\Delta}_{i'}^{\text{an}}  \times U_t^{\text{an}}/
\widehat{\mathcal{M}}_{\boldsymbol{t}_{\text{ra}}}\times T_{\boldsymbol{\theta}} \times U_t^{\text{an}}}
(D(\tilde{\boldsymbol{t}}_0)_{t})
$$
corresponding to \eqref{2021.4.16.12.12_y}.
We denote by
$\widehat\Omega_{(\boldsymbol{t}_{\mathrm{ra}},\boldsymbol{\theta}_0)}^{(1)}
(x_{t_{i'}},t)$
the connection matrix of the relative connection,
where the expansion of $\widehat\Omega_{(\boldsymbol{t}_{\mathrm{ra}},\boldsymbol{\theta}_0)}^{(1)}
(x_{t_{i'}}-\epsilon \delta_{t_i}^{\mathrm{IMD}} (t_{i'}) ,\epsilon)$ with respect to $\epsilon$
is equal to 
$$
\widehat\Omega_{(\boldsymbol{t}_{\mathrm{ra}},\boldsymbol{\theta}_0)}^{(1)}
|_{\widehat{\Delta}_{i'}^{\text{an}} }
+ \epsilon \cdot \delta_{t_i}^{\mathrm{IMD}}
(\widehat\Omega_{(\boldsymbol{t}_{\mathrm{ra}},\boldsymbol{\theta}_0)}^{(1)}) 
|_{\widehat{\Delta}_{i'}^{\text{an}}}.
$$
Let $\widehat{\Sigma} \subset \widehat{\Delta}^{\text{an}}_{i'}$ 
be a family of sufficiently small sectors in 
$\widehat{\Delta}^{\text{an}}_{i'}
\rightarrow 
\widehat{\mathcal{M}}_{\boldsymbol{t}_{\text{ra}}}\times T_{\boldsymbol{\theta}}$. 
We take a fundamental matrix
solution $\Psi_{\widehat{\Sigma}}(x_{t_{i'}},t)$ 
on $ \widehat{\Sigma}  \times U_t^{\text{an}}$
of the connection 
$$
d\Psi_{\widehat\Sigma}(x_{t_{i'}},t)
+\widehat\Omega_{(\boldsymbol{t}_{\mathrm{ra}},\boldsymbol{\theta}_0)}^{(1)}
(x_{t_{i'}},t)
\Psi_{\widehat\Sigma}(x_{t_{i'}},t) =0
$$
with uniform asymptotic relation 
\begin{equation}\label{2021.11.16.17.13}
\Psi_{\widehat\Sigma}(x_{t_{i'}},t) \exp\left( \Lambda^-_{i'}(x_{t_{i'}} ,t )\right) \sim
\widehat{P}_{i'}(x_{t_{i'}},t)
\quad (x_{t_{i'}} \rightarrow 0, x_{t_{i'}} \in \widehat\Sigma).
\end{equation}
Here we define $ \Lambda^-_{i'}(x_{t_{i'}} ,t )$ and
$\widehat{P}_{i'}(x_{t_{i'}},t)$ 
so that
the expansions of $\Lambda^-_{i'}(x_{t_{i'}}-\epsilon \delta_{t_i}^{\mathrm{IMD}} (t_{i'}) ,\epsilon)$
and $\widehat{P}_{i'}(x_{t_{i'}}-\epsilon \delta_{t_i}^{\mathrm{IMD}} (t_{i'}) ,\epsilon)$ 
with respect to $\epsilon$ are the following:
$$
\begin{aligned}
\Lambda^-_{i'}(x_{t_{i'}}-\epsilon \delta_{t_i}^{\mathrm{IMD}} (t_{i'}) ,\epsilon)
 &= \Lambda^-_{i'}(x_{t_{i'}})
+\epsilon \cdot \delta_{t_i}^{\mathrm{IMD}} ( \Lambda^-_{i'}(x_{t_{i'}}) ) \\
\widehat{P}_{i'}
(x_{t_{i'}}-\epsilon \delta_{t_i}^{\mathrm{IMD}} (t_{i'}) ,\epsilon)
&= \tilde\psi_{i'} \exp\left(\Lambda^-_{i'}(x_{t_{i'}})\right) 
+\epsilon \cdot \delta_{t_i}^{\mathrm{IMD}}
\left(
\tilde\psi_{i'} \exp\left( \Lambda^-_{i'}(x_{t_{i'}})\right)
\right).
\end{aligned}
$$
By the asymptotic relation \eqref{2021.11.16.17.13} and the same argument as above, we have 
\begin{equation*}
\delta_{t_i}^{\mathrm{IMD}}(\tilde\psi_{i'}) (\tilde\psi_{i'})^{-1}
\sim 
\Upsilon_{t_i}^{\mathrm{IMD}}
-
\widehat{P}_{i'}(x_{t_{i'}},0)
 \tilde{C}^{\text{diag}}_{t_{i'}}
\widehat{P}_{i'}(x_{t_{i'}},0)^{-1} 
\end{equation*}
as \eqref{2021.4.18.17.41_1}.
Here $\tilde{C}^{\text{diag}}_{t_{i'}}$ 
is a diagonal matrix
and $\tilde{C}^{\text{diag}}_{t_{i'}}$ is independent of $x_{t_{i'}}$.
By this asymptotic relation, 
we may check that \eqref{2021.11.16.17.14} as above.
By the integrable condition \eqref{2021.11.16.17.17},
we may check that 
\begin{equation*}
\hat\omega (\delta_{t_i}^{\mathrm{IMD}}, \hat\delta) =
 \frac{1}{2} \sum_{i' \in I}   \mathrm{res}_{x=t_{i'}} 
\mathrm{Tr} \left( \delta_{t_i}^{\mathrm{IMD}}
(\widehat\Omega_{(\boldsymbol{t}_{\mathrm{ra}},\boldsymbol{\theta}_0)}^{(1)}) 
\hat\delta(\tilde\psi_{i'})(\tilde\psi_{i'})^{-1} - 
\Upsilon_{t_i}^{\mathrm{IMD}}
\hat\delta (\widehat\Omega_{(\boldsymbol{t}_{\mathrm{ra}},\boldsymbol{\theta}_0)}^{(1)}) 
  \right) =0
\end{equation*}
for any $\hat{\delta}$.
\end{proof}

By Theorem \ref{2020.1.21.16.17} and Theorem \ref{2020.1.21.13.51},
we have that the $2$-form $\hat{\omega}$ is the isomonodromy $2$-form.

\subsection{Hamiltonian systems}\label{2020.2.8.10.30}

First, we define {\it Hamiltonians} 
on $\widehat{\mathcal{M}}_{\boldsymbol{t}_{\text{ra}}}\times T_{\boldsymbol{\theta}}$
as follows.
By Lemma \ref{2019.12.30.22.09},
we have the following diagonalization:
$$
\begin{aligned}
\Omega_{t_i}^{\text{diag}}&=
\frac{
\begin{pmatrix}
\theta_{0,t_i}^+ & 0 \\
0 & \theta_{0,t_i}^-
\end{pmatrix}
}{x_{t_i}^{n_i}}dx_{t_i} + \cdots +
\frac{
\begin{pmatrix}
\theta_{n_i-1,t_i}^+ & 0 \\
0 & \theta_{n_i-1,t_i}^-
\end{pmatrix}
}{x_{t_i}} dx_{t_i}\\
&\qquad + \begin{pmatrix}
\theta_{n_i,t_i}^+ & 0 \\
0 & \theta_{n_i,t_i}^-
\end{pmatrix}
 dx_{t_i}+ \cdots +
 \begin{pmatrix}
\theta_{2n_i-1,t_i}^+ & 0 \\
0 & \theta_{2n_i-1,t_i}^-
\end{pmatrix}
x_{t_i}^{n_i-1} dx_{t_i}+\cdots.
\end{aligned}
$$
Here we set $\Omega_{t_i}^{\text{diag}}:=
(\Phi_i \Xi_i)^{-1} d(\Phi_i \Xi_i) +(\Phi_i \Xi_i)^{-1} 
\widehat\Omega_{(\boldsymbol{t}_{\mathrm{ra}},\boldsymbol{\theta}_0)}^{(n-2)} 
(\Phi_i \Xi_i)$.
 Remark that
 we have an equation $(d +\Omega_{t_i}^{\text{diag}}) \exp(-\Lambda_i(x_{t_i}))=0$.
We set 
\begin{equation*}
\Lambda^+_i(x_{t_i})= 
\begin{pmatrix}
\sum_{k=n_i}^{\infty}\theta_{k,t_i}^{+} \int x_{t_i}^{-n_i+k} dx_{t_i} & 0 \\
0 &\sum_{k=n_i}^{\infty}\theta_{k,t_i}^{-} \int x_{t_i}^{-n_i+k} dx_{t_i}
\end{pmatrix} .
\end{equation*}

\begin{Def}\label{2020.2.7.11.59}
For each $i \in I_{\text{un}}$ and 
each $l$ ($0\le l \le n_i-2$),
we define rational functions $H_{\theta_{l,t_i}^{\pm}}$
on 
$\widehat{\mathcal{M}}_{\boldsymbol{t}_{\text{ra}}}\times T_{\boldsymbol{\theta}}$
as
\begin{equation*}
\begin{aligned}
H_{\theta_{l,t_i}^{+}}&=- [\text{ the coefficient 
of the $x_{t_i}^{n_{i}-l-1}$-term of 
the $(1,1)$-entry of $\Lambda^+_i(x_{t_i})$ }] \\
&=-\frac{ \theta_{2n_i-l-2,t_i}^{+}}{n_i-l-1}, \text{ and}\\
H_{\theta_{l,t_i}^{-}}
&=- [\text{ the coefficient 
of the $x_{t_i}^{n_{i}-l-1}$-term of 
the $(2,2)$-entry of $\Lambda^+_i(x_{t_i})$ }] \\
&=-\frac{ \theta_{2n_i-l-2,t_i}^{-}}{n_i-l-1}.
\end{aligned}
\end{equation*}
We call $H_{\theta_{l,t_i}^{\pm}}$ the {\it Hamiltonian associated to $\theta_{l,t_i}^{\pm}$}.
\end{Def}

\begin{Def}\label{2020.2.7.12.00}
For each $i\in \{ 3,4, \ldots,\nu\}$, 
put $\hat{\lambda}_{i,\pm}^{\le 2n_i-1}(x):=
\sum_{l=0}^{2n_i-1}\theta_{l,t_i}^{\pm} \int (x-t_i)^{-n_i+l} dx$.
We define a rational function $H_{t_i}$
on 
$\widehat{\mathcal{M}}_{\boldsymbol{t}_{\text{ra}}}\times T_{\boldsymbol{\theta}}$
as
\begin{equation*}
\begin{aligned}
H_{t_i}:=&\ 
- \frac{1}{2} \cdot \mathrm{res}_{x=t_i} \left( 
\mathrm{Tr}( \Omega_{t_i}^{\text{diag}} )^2
 \right) \\
=&\ 
-\sum_{l=0}^{n_i-1} \theta_{l,t_i}^+ \cdot \theta_{2n_i-l-1,t_i}^+
-\sum_{l=0}^{n_i-1} \theta_{l,t_i}^- \cdot \theta_{2n_i-l-1,t_i}^-.
\end{aligned}
\end{equation*}
We call $H_{t_i}$ the {\it Hamiltonian associated to $t_i$}.
\end{Def}

We will give description of integrable deformations by using these Hamiltonians.
This description is derived by calculation of the
isomonodromy $2$-form $\hat{\omega}$.
Now, we prepare a lemma for calculation of the
isomonodromy $2$-form $\hat{\omega}$.
Let $t$ be a component of the divisor $D$ and
$( U_t , x_t) $ be a couple of an affine open subset such that $t \in U_t$ 
and $x_t=x-t$.
Let $\Omega$ be an element of
$\mathfrak{gl}(2,\mathbb{C})\otimes \Omega^1_{U_t} (D+D_{\mathrm{App}}) \otimes 
\mathcal{O}_{\widehat{\mathcal{M}}_{\boldsymbol{t}_{\text{ra}}}\times T_{\boldsymbol{\theta}}} $,
which has an expansion at $t$ as follows:
\begin{equation*}
\Omega = \frac{\Omega_0dx_t}{x_t^{n_t}} + \frac{\Omega_1dx_t}{x_t^{n_t-1}}+
\cdots + \frac{\Omega_{n_t-1}dx_t}{x_t} + \cdots ,
\end{equation*}
where $\Omega_k \in  
\mathfrak{gl}(2,\mathbb{C})\otimes\mathcal{O}_{ 
\widehat{\mathcal{M}}_{\boldsymbol{t}_{\text{ra}}}\times T_{\boldsymbol{\theta}}}$.
Let $g$ be an element of 
$\mathcal{E}nd(\mathcal{O}^{\oplus 2}_{U_t\times 
(\widehat{\mathcal{M}}_{\boldsymbol{t}_{\text{ra}}}\times T_{\boldsymbol{\theta}})}  )  $
such that $g$ has an expansion at $t$ as follows:
\begin{equation}\label{2020.2.16.21.26}
g = g_0 + g_1 x_t + \cdots+ g_{n_t} x_t^{n_t}+ \cdots ,
\end{equation}
where $g_k \in  \mathcal{E}nd(\mathcal{O}^{\oplus 2}_{ 
\widehat{\mathcal{M}}_{\boldsymbol{t}_{\text{ra}}}\times T_{\boldsymbol{\theta}}}  ) $ .

\begin{Lem}\label{2020.1.20.21.24}
\textit{
Let $\hat{\delta}_1$ and $\hat{\delta}_2$ be vector fields on an open subset of 
$\widehat{\mathcal{M}}_{\boldsymbol{t}_{\text{ra}}}\times T_{\boldsymbol{\theta}} $
such that $\Omega$ and $g$ are defined on this open subset.
Let $\psi$  be a formal solution of $d +\Omega=0$ at $t$.
We assume that 
\begin{itemize}
\item $\Omega_k\in  \mathfrak{gl}(2,\mathbb{C})\otimes 
 \pi_{\boldsymbol{t}_{\mathrm{ra}}, \boldsymbol{\theta}_0}^{-1} 
(\mathcal{O}_{(T_{\boldsymbol{t}})_{\boldsymbol{t}_{\text{ra}}} \times T_{\boldsymbol{\theta}}} )$
for $k=0,1,\ldots , n_t -1$,
\item 
$ g_k \in \mathcal{E}nd( 
\pi_{\boldsymbol{t}_{\mathrm{ra}}, \boldsymbol{\theta}_0}^{-1} 
(\mathcal{O}_{(T_{\boldsymbol{t}})_{\boldsymbol{t}_{\text{ra}}} 
\times T_{\boldsymbol{\theta}}} )^{\oplus 2})$
for $k=0,1,\ldots , n_t -1$,
\item we can define the inverse matrix $ g_0^{-1} \in \mathcal{E}nd( 
\pi_{\boldsymbol{t}_{\mathrm{ra}}, \boldsymbol{\theta}_0}^{-1} 
(\mathcal{O}_{(T_{\boldsymbol{t}})_{\boldsymbol{t}_{\text{ra}}} 
\times T_{\boldsymbol{\theta}}} )^{\oplus 2})$ of $g_0$,
\item the $(1,1)$-entry and the $(2,2)$-entry of $g_0^{-1}g_{n_t}$ vanish,
\item $g_0^{-1} \Omega_0 g_0$ is a diagonal matrix, and 
\item $\hat{\delta}_1 (\psi) \psi^{-1}$ and $\hat{\delta}_2 (\psi) \psi^{-1}$ are formally meromorphic at $t$.
\end{itemize}
If we set $\Omega'= g^{-1} dg +g^{-1} \Omega g$ and $\psi'=g^{-1} \psi$, then
the difference 
\begin{equation*}
\mathrm{res}_{x=t}\mathrm{Tr} \left( \hat{\delta} (\Omega')
\wedge \hat{\delta}(\psi') (\psi')^{-1}   \right)
-\mathrm{res}_{x=t}\mathrm{Tr} \left( \hat{\delta} (\Omega)
\wedge \hat{\delta}(\psi) (\psi)^{-1}  
 \right)
\end{equation*}
is a section of $\pi_{\boldsymbol{t}_{\mathrm{ra}}, \boldsymbol{\theta}_0}^{-1}
(\Omega^2_{(T_{\boldsymbol{t}})_{\boldsymbol{t}_{\text{ra}}} \times T_{\boldsymbol{\theta}}} )$.
}
\end{Lem}

\begin{proof}
Put $\hat{\delta}(\Omega) \wedge \hat{\delta}(\psi) \psi^{-1}:= \hat{\delta}_1 (\Omega)
\hat{\delta}_2(\psi) (\psi)^{-1}  
-\hat{\delta}_1(\psi) (\psi)^{-1} \hat{\delta}_2 (\Omega)$.
By Lemma \ref{2020.12.24.17.31}, we have the following equality
\begin{equation}\label{2021.1.1.17.38}
\begin{aligned}
\mathrm{Tr} &\left( \hat{\delta}(\Omega') \wedge \hat{\delta}(\psi') (\psi')^{-1}\right)- 
\mathrm{Tr} \left( \hat{\delta}(\Omega) \wedge \hat{\delta}(\psi) \psi^{-1} \right)\\
=&-\mathrm{Tr} \left( \hat{\delta}_1(\Omega') \tilde{u}^{(2)} -
\tilde{u}^{(1)}\hat{\delta}_2(\Omega')  \right)
-\mathrm{Tr} \left( \hat{\delta}_1(\Omega)   u^{(2)} -
u^{(1)}\hat{\delta}_2(\Omega)   \right)  \\
&+ \mathrm{Tr} \left(  d (\psi^{-1} u^{(1)} \hat{\delta}_2(\psi) 
-\psi^{-1} u^{(2)} \hat{\delta}_1(\psi) ) \right),
\end{aligned}
\end{equation}
where $u^{(i)}:= \hat{\delta}_i(g)g^{-1}$ 
and $\tilde{u}^{(i)}:= g^{-1}\hat{\delta}_i(g)$ for $i \in \{ 1,2\}$.
We set $\hat{\delta}_1= \delta_1 \in 
\Theta_{(\widehat{\mathcal{M}}_{\boldsymbol{t}_{\text{ra}}}\times T_{\boldsymbol{\theta}}) 
/((T_{\boldsymbol{t}})_{\boldsymbol{t}_{\text{ra}}} \times T_{\boldsymbol{\theta}})}$.
We will show that the residue of \eqref{2021.1.1.17.38} at $t$ vanishes.
We consider the residues of
 the first term and the second term of the right hand side of \eqref{2021.1.1.17.38}.
 We calculate $\hat{\delta}_2 (\Omega')$ as follows:
 \begin{equation*}
 \begin{aligned}
\hat{\delta}_2 (\Omega') &= \hat{\delta}_2 \left((g^{-1}_0\Omega_0g_0)(x-t)^{-n_t} \right)dx + \cdots  \\
&= g^{-1}_0\Omega_0g_0 \cdot \hat{\delta}_2 \left((x-t)^{-n_t} \right)dx
+\hat{\delta}_2 \left(g^{-1}_0\Omega_0g_0\right)(x-t)^{-n_t} dx  
+ \cdots  \\
&= n_{t_i}\hat{\delta}_2 (t) \cdot \ g^{-1}_0\Omega_0g_0 \cdot (x-t)^{-n_t-1} dx
+\hat{\delta}_2 \left(g^{-1}_0\Omega_0g_0\right)(x-t)^{-n_t} dx  
+ \cdots.
\end{aligned}
\end{equation*}
 So we have that
the variation
 $\hat{\delta}_2(\Omega')$
has a pole of order $n_t+1$ at $t$
and
the leading coefficient of $\hat{\delta}_2(\Omega')$ is
the diagonal matrix $n_t \hat{\delta}_2 (t) \cdot g_0^{-1} \Omega_0 g_0 $.
Since $\delta_1$ is an element of
$\Theta_{(\widehat{\mathcal{M}}_{\boldsymbol{t}_{\text{ra}}}\times T_{\boldsymbol{\theta}}) 
/((T_{\boldsymbol{t}})_{\boldsymbol{t}_{\text{ra}}} \times T_{\boldsymbol{\theta}})}$,
we may check that
 $\delta_1(\Omega') \tilde{u}^{(2)}$ 
is holomorphic at $t$.
We calculate the residue of the first term as follows:
\begin{equation}\label{2021.1.1.17.51}
\begin{aligned}
-\mathrm{res}_{x=t}\mathrm{Tr} \left( \delta_1(\Omega') \tilde{u}^{(2)} -
\tilde{u}^{(1)}\hat{\delta}_2(\Omega')  \right)
&=\mathrm{res}_{x=t}\mathrm{Tr} \left( 
\tilde{u}^{(1)}\hat{\delta}_2(\Omega')  \right)\\
&=n_t \hat{\delta}_2(t) \cdot \mathrm{Tr} \left( 
g_0^{-1}\delta_1(g_{n_t}) g_0^{-1} \Omega_0 g_0 \right) \\
&=n_t \hat{\delta}_2(t) \cdot \mathrm{Tr} \left( 
\delta_1(g_0^{-1}g_{n_t}) g_0^{-1} \Omega_0 g_0 \right).
\end{aligned}
\end{equation}
Since the diagonal entries of $g_0^{-1}g_{n_t}$ vanish 
and $g_0^{-1} \Omega_0 g_0$ is a diagonal matrix,
the residue
\eqref{2021.1.1.17.51} is zero.
Next we calculate the residue of the second term as follows:
\begin{equation*}
\begin{aligned}
-\mathrm{res}_{x=t}\mathrm{Tr} \left( \delta_1(\Omega)   u^{(2)} -
u^{(1)}\hat{\delta}_2(\Omega)   \right) 
&= \mathrm{res}_{x=t} \mathrm{Tr} \left( 
\delta_1(g) g^{-1} \hat{\delta}_2(\Omega)   \right)  \\
&= \mathrm{res}_{x=t} \mathrm{Tr} \left( 
g^{-1}\delta_1(g) \left(
 \hat{\delta}_2(g^{-1} \Omega g) 
-\hat{\delta}_2(g^{-1}) \Omega g
-g^{-1} \Omega  \hat{\delta}_2(g) \right)
 \right)\\
&=n_t \hat{\delta}_2(t) \cdot \mathrm{Tr} \left( 
g_0^{-1}\delta_1(g_{n_t}) g_0^{-1} \Omega_0 g_0 \right)=0.
\end{aligned}
\end{equation*}
Here remark that $\delta_1(g) \in O(x_t^{n_t})$.
Finally, the residues of the third term of \eqref{2021.1.1.17.38}
at $t$ is zero,
since $\mathrm{Tr}(\psi^{-1} u^{(1)} \hat{\delta}_2(\psi) 
-\psi^{-1} u^{(2)} \hat{\delta}_1(\psi))$ is formally meromorphic at $t$.
Then we have that the residue of \eqref{2020.1.8.12.55} at $t$ vanishes.
We obtain the assertion of this lemma.
\end{proof}

\begin{Thm}\label{2020.1.21.13.48}
\textit{
Set $P(x;\boldsymbol{t}):=\prod_{i=1}^{\nu} (x-t_i)^{n_i}$ and  
$D_i(x;\boldsymbol{t}, \boldsymbol{\theta} ):=D_i(x)$ for $i \in I$.
We put 
\begin{equation}\label{2020.1.21.13.49}
\begin{aligned}
\hat{\omega}':=&\ 
 \sum_{j=1}^{n-3} d \left(\frac{p_j}{P(q_j;\boldsymbol{t})} 
 - \sum_{i=1}^{\nu} \frac{D_i(q_j;\boldsymbol{t}, \boldsymbol{\theta})}{(q_j -t_i)^{n_i}} 
 - D_{\infty}(q_j;\boldsymbol{t}, \boldsymbol{\theta})\right)
  \wedge dq_j\\
&\ + \sum_{i\in I_{\mathrm{un}}} \sum_{l=0}^{n_{i}-2} \left( d  H_{\theta_{l,t_i}^{+}} \wedge d \theta_{l,t_i}^{+}
+ d H_{\theta_{l,t_i}^{-}} \wedge d \theta_{l,t_i}^{-} \right)
+\sum_{i\in \{3,4,\ldots,\nu \}} d H_{t_i} \wedge dt_i.
\end{aligned}
\end{equation}
Then the difference $\hat{\omega}- \hat{\omega}'$ 
is a section of $\pi_{\boldsymbol{t}_{\mathrm{ra}}, \boldsymbol{\theta}_0}^* 
(\Omega^2_{(T_{\boldsymbol{t}})_{\boldsymbol{t}_{\text{ra}}} \times T_{\boldsymbol{\theta}}} )$.
}
\end{Thm}

\begin{proof}
Recall that $\hat{\omega}$ is
\begin{equation*}
\begin{aligned}
\frac{1}{2} \sum_{i\in I} 
\mathrm{res}_{x=t_i} 
\mathrm{Tr} \left( \hat{\delta} (
\widehat\Omega_{(\boldsymbol{t}_{\mathrm{ra}},\boldsymbol{\theta}_0)}^{(n-2)} ) 
\wedge \hat{\delta}(\psi_i)\psi_i^{-1} )   \right)
+\frac{1}{2} \sum_{j=1}^{n-3} \mathrm{res}_{x=q_j} 
\mathrm{Tr} \left( \hat{\delta} 
(\widehat\Omega_{(\boldsymbol{t}_{\mathrm{ra}},\boldsymbol{\theta}_0)}^{(n-2)} )
\wedge \hat{\delta}(\psi_{q_j})\psi_{q_j}^{-1} 
  \right).
\end{aligned}
\end{equation*}
The plan of the proof is as follows.
First we will consider the first term of this formula.
We calculate 
the residue at $t_i$ for some (local) gauge transformation of 
$d+\widehat\Omega_{(\boldsymbol{t}_{\mathrm{ra}},\boldsymbol{\theta}_0)}^{(n-2)}$.
We need to consider the difference between the residue after taking the gauge transformation
and the residue before taking the gauge transformation.
Here the residue before taking the gauge transformation is just
$\mathrm{res}_{x=t_i} \mathrm{Tr} \left( \hat{\delta} (
\widehat\Omega_{(\boldsymbol{t}_{\mathrm{ra}},\boldsymbol{\theta}_0)}^{(n-2)} ) 
\wedge \hat{\delta}(\psi_i)\psi_i^{-1} )   \right)$.
To consider this difference, we will use Lemma \ref{2020.1.20.21.24}.
By calculation of the residue after taking the gauge transformation,
we may derive the second part of \eqref{2020.1.21.13.49}:
$$
\sum_{i\in I_{\mathrm{un}}} \sum_{l=0}^{n_{i}-2} \left( d  H_{\theta_{l,t_i}^{+}} \wedge d \theta_{l,t_i}^{+}
+ d H_{\theta_{l,t_i}^{-}} \wedge d \theta_{l,t_i}^{-} \right)
+\sum_{i\in \{3,4,\ldots,\nu \}} d H_{t_i} \wedge dt_i.
$$
Second we will calculate the second term
(the residue of $\mathrm{Tr} ( \hat{\delta} 
(\widehat\Omega_{(\boldsymbol{t}_{\mathrm{ra}},\boldsymbol{\theta}_0)}^{(n-2)} )
\wedge \hat{\delta}(\psi_{q_j})\psi_{q_j}^{-1} )$ at $x=q_j$)
by the same argument as in the proof of Theorem \ref{2020.1.21.16.17}.
Then we may derive the first part of \eqref{2020.1.21.13.49}:
$$
 \sum_{j=1}^{n-3} d \left(\frac{p_j}{P(q_j;\boldsymbol{t})} 
 - \sum_{i=1}^{\nu} \frac{D_i(q_j;\boldsymbol{t}, \boldsymbol{\theta})}{(q_j -t_i)^{n_i}} 
 - D_{\infty}(q_j;\boldsymbol{t}, \boldsymbol{\theta})\right)
  \wedge dq_j.
$$

First, we consider the residue of 
$\mathrm{Tr} \left( 
\hat\delta(\widehat\Omega_{(\boldsymbol{t}_{\mathrm{ra}},\boldsymbol{\theta}_0)}^{(n-2)})
\wedge \hat{\delta}(\psi_i)\psi_i^{-1}   \right)$
at $t_i$.
Now
we take diagonalizations of 
$d+\widehat\Omega_{(\boldsymbol{t}_{\mathrm{ra}},\boldsymbol{\theta}_0)}^{(n-2)}$ 
until some degree term at each point $t_i$.
For $i \in I$, we put
\begin{equation*}
\Xi_{i}^{\le 2n_i-1} (x_{t_i}):= 
\begin{pmatrix}
1 & 0\\
0 & 1
\end{pmatrix}+
\sum_{s=1}^{2n_i-1} 
\begin{pmatrix}
(\xi^{(i)}_s)_{11} & (\xi^{(i)}_s)_{12} \\
(\xi^{(i)}_s)_{21} & (\xi^{(i)}_s)_{22}
\end{pmatrix}
x_{t_i}^s .
\end{equation*}
Here the coefficient matrices of $\Xi_{i}^{\le 2n_i-1}$ appear in  
Lemma \ref{2019.12.30.22.09}
as the coefficient matrices of $\Xi_{i}$
for the connection 
$d+\widehat\Omega_{(\boldsymbol{t}_{\mathrm{ra}},\boldsymbol{\theta}_0)}^{(n-2)}$.
We put
\begin{equation*}
\begin{aligned}
&\tilde{\Omega}_{i} :=  (\Phi_i\Xi_{i}^{\le 2n_i-1})^{-1} d(\Phi_i\Xi_{i}^{\le 2n_i-1})
+(\Phi_i\Xi_{i}^{\le 2n_i-1})^{-1} 
\widehat\Omega_{(\boldsymbol{t}_{\mathrm{ra}},\boldsymbol{\theta}_0)}^{(n-2)}
 (\Phi_i\Xi_{i}^{\le 2n_i-1}) \quad \text{and} \\
&\tilde{\psi}_{i} := (\Phi_i\Xi_{i}^{\le 2n_i-1})^{-1} 
\psi_{i},
\end{aligned}
\end{equation*}
where $\psi_{i}$ is the formal solution as in Lemma \ref{2019.12.30.22.09}.
We may describe $\tilde\Omega_i$ as follows:
$$
\begin{aligned}
\tilde{\Omega}_{i}&=
\begin{pmatrix}
\theta_{0,t_i}^{+}  & 0 \\
0 & \theta_{0,t_i}^{-}
\end{pmatrix}
\frac{dx}{(x-t_i)^{n_i}}
+ \cdots+
\begin{pmatrix}
\theta_{n_i-1,t_i}^{+}  & 0 \\
0 & \theta_{n_i-1,t_i}^{-}
\end{pmatrix}
 \frac{dx}{x-t_i} \\
&\qquad+\begin{pmatrix}
\theta_{n_i,t_i}^{+}  & 0 \\
0 & \theta_{n_i,t_i}^{-}
\end{pmatrix}dx
+\cdots +
\begin{pmatrix}
\theta_{2n_i-1,t_i}^{+}  & 0 \\
0 & \theta_{2n_i-1,t_i}^{-}
\end{pmatrix} (x-t_i)^{n_i-1} dx +O(x-t_i)^{n_i}.
\end{aligned}
$$
The residue part $\theta_{n_i-1,t_i}^{\pm}$ of $\tilde{\Omega}_{i}$
is constant on $\widehat{\mathcal{M}}_{\boldsymbol{t}_{\text{ra}}}\times T_{\boldsymbol{\theta}}$.
So we have $\hat{\delta}(\theta^{\pm}_{n_i-1,t_i})=0$ for any 
$\hat{\delta} \in \Theta_{\widehat{\mathcal{M}}_{\boldsymbol{t}_{\text{ra}}}\times T_{\boldsymbol{\theta}}}$.
We may check that
the variation $\hat\delta_1(\tilde{\Omega}_{i})$ is equal to
$$
\begin{aligned}
&\begin{pmatrix}
\frac{\hat\delta_1(\theta_{0,t_i}^{+})}{(x-t_i)^{n_i}}  & 0 \\
0 & \frac{\hat\delta_1(\theta_{0,t_i}^{-})}{(x-t_i)^{n_i}}
\end{pmatrix}
dx
+ \cdots+
\begin{pmatrix}
\frac{\hat\delta_1(\theta_{n_i-2,t_i}^{+} )}{(x-t_i)^2} & 0 \\
0 &\frac{ \hat\delta_1(\theta_{n_i-2,t_i}^{-})}{(x-t_i)^2}
\end{pmatrix}dx \\
&\qquad+\begin{pmatrix}
\hat\delta_1(\theta_{n_i,t_i}^{+} ) & 0 \\
0 & \hat\delta_1(\theta_{n_i,t_i}^{-})
\end{pmatrix}dx
+\cdots +
\begin{pmatrix}
\hat\delta_1(\theta_{2n_i-1,t_i}^{+} ) & 0 \\
0 &\hat\delta_1( \theta_{2n_i-1,t_i}^{-})
\end{pmatrix} (x-t_i)^{n_i-1} dx \\
&\quad+n_i\hat\delta_1(t_i)\begin{pmatrix}
\frac{\theta_{0,t_i}^{+} }{(x-t_i)^{n_i+1}} & 0 \\
0 & \frac{\theta_{0,t_i}^{-}}{(x-t_i)^{n_i+1}}
\end{pmatrix}
 dx
+ \cdots+\hat\delta_1(t_i)
\begin{pmatrix}
\frac{\theta_{n_i-1,t_i}^{+}}{(x-t_i)^2}  & 0 \\
0 &\frac{ \theta_{n_i-1,t_i}^{-}}{(x-t_i)^2}
\end{pmatrix}
dx \\
&\qquad
-\hat\delta_1(t_i)\begin{pmatrix}
\theta_{n_i+1,t_i}^{+}  & 0 \\
0 & \theta_{n_i+1,t_i}^{-}
\end{pmatrix}  dx
-\cdots - (n_i-1)\hat\delta_1(t_i)
\begin{pmatrix}
\theta_{2n_i-1,t_i}^{+}  & 0 \\
0 & \theta_{2n_i-1,t_i}^{-}
\end{pmatrix}   (x-t_i)^{n_i-2} dx\\
& \qquad  \quad
-\hat{\delta}_1(t_i)
\begin{pmatrix}
* & * \\
* & * 
\end{pmatrix}
 \cdot (x-t_i)^{n_i-1}dx+O(x-t_i)^{n_i}.
\end{aligned}
$$
We define $\hat{\lambda}_{i,\pm}^{\le 2n_i-1}(x_{t_i})$ as
$$
\begin{aligned}
\hat{\lambda}_{i,\pm}^{\le 2n_i-1}(x_{t_i}) &=
\frac{\theta_{0,t_i}^{\pm}}{-n_i+1}  (x-t_i)^{-n_i+1}
+ \cdots+\frac{\theta_{n_i-2,t_i}^{\pm}}{-1}  (x-t_i)^{-1} + \theta^{\pm}_{n_i-1,t_i} \log (x-t_i)\\
&\qquad+\theta_{n_i,t_i}^{\pm}  (x-t_i)
+\cdots +\frac{\theta_{2n_i-1,t_i}^{\pm}}{n_i}  (x-t_i)^{n_i} .
\end{aligned}
$$
On the other hand, the variation $\hat{\delta}_{2}(\tilde{\psi}_i)\tilde{\psi}_i^{-1}$ is equal to
$$
\begin{aligned}
& \delta((\Phi_i\Xi_{i}^{\le 2n_i-1})^{-1} \Phi_{i} \Xi_{i} )((\Phi_i\Xi_{i}^{\le 2n_i-1})^{-1} \Phi_{i} \Xi_{i})^{-1}\\
 &\qquad + 
((\Phi_i\Xi_{i}^{\le 2n_i-1})^{-1} \Phi_{i} \Xi_{i} ) 
\begin{pmatrix}
-\hat\delta_2 (\hat{\lambda}_i^+(x_{t_i})) & 0 \\
0 &-\hat\delta_2 (\hat{\lambda}_i^-(x_{t_i}))
\end{pmatrix}
 ((\Phi_i\Xi_{i}^{\le 2n_i-1})^{-1} \Phi_{i} \Xi_{i})^{-1}\\
 &=
\begin{pmatrix}
-\hat{\delta}_2 ( \hat{\lambda}^{\le 2n_i-1}_{i,+}(x_{t_i}))& 0 \\
0 &-\hat{\delta}_2 (\hat{\lambda}^{\le 2n_i-1}_{i,-}(x_{t_i}))
\end{pmatrix}  + \hat{\delta}_2(t_i)
\begin{pmatrix}
* & * \\
* & * 
\end{pmatrix}
 x_{t_i}^{n_i} + O(x_{t_i}^{n_{i}+1}).
\end{aligned}
$$
Since $\hat{\delta}_2(\hat{\lambda}_{i,\pm}^{\le 2n_i-1}(x_{t_i}))$ is equal to
$$
\begin{aligned}
&\frac{\hat{\delta}_2(\theta_{0,t_i}^{\pm})}{-n_i+1}  (x-t_i)^{-n_i+1}
+ \cdots+\frac{\hat{\delta}_2(\theta_{n_i-2,t_i}^{\pm})}{-1}  (x-t_i)^{-1}  \\
&\qquad+\hat{\delta}_2(\theta_{n_i,t_i}^{\pm} ) (x-t_i)
+\cdots +\frac{\hat{\delta}_2(\theta_{2n_i-1,t_i}^{\pm})}{n_i}  (x-t_i)^{n_i} \\
&\quad +\theta_{0,t_i}^{\pm} \cdot(-\hat{\delta}_2(t_i)) (x-t_i)^{-n_i}
+ \cdots+\theta_{n_i-2,t_i}^{\pm}\cdot(-\hat{\delta}_2(t_i))  (x-t_i)^{-2} 
+ \theta^{\pm}_{n_i-1,t_i}\cdot(-\hat{\delta}_2(t_i)) (x-t_i)^{-1}\\
&\qquad+\theta_{n_i,t_i}^{\pm} \cdot(-\hat{\delta}_2(t_i)) 
+\cdots +\theta_{2n_i-1,t_i}^{\pm}\cdot (-\hat{\delta}_2(t_i))  (x-t_i)^{n_i-1},
\end{aligned}
$$
we may check that the residue 
$\mathrm{Tr} ( \hat{\delta}_1 (\tilde{\Omega}_i) 
 \hat{\delta}_2(\tilde{\psi}_i)\tilde{\psi}_i^{-1} )$ at $t_i$
coincides with
\begin{equation*}
\begin{aligned}
&\sum_{l\in \{0,1,\ldots,2n_i-2 \} \setminus \{n_i-1 \}}\left( \hat{\delta}_1(\theta_{l,t_i}^{+})\cdot
\frac{\hat{\delta}_2(\theta_{2n_i-l-2,t_i}^{+})}{n_i-l-1}\right)
+\sum_{l\in \{0,1,\ldots,2n_i-2 \} \setminus \{n_i-1 \}}\left( \hat{\delta}_1(\theta_{l,t_i}^{-})\cdot
\frac{\hat{\delta}_2(\theta_{2n_i-l-2,t_i}^{-})}{n_i-l-1}\right)\\
&+\sum_{l\in \{0,1,\ldots,2n_i-1 \} \setminus \{n_i \}}\left(
(n_i-l)\theta_{l,t_i}^{+} \hat{\delta}_1(t_i)\cdot
\frac{\hat{\delta}_2(\theta_{2n_i-l-1,t_i}^{+})}{n_i-l}
-\hat{\delta}_1(\theta_{l,t_i}^{+} )\cdot
\theta_{2n_i-l-1,t_i}^{+}\hat{\delta}_2(t_i)  \right)\\
&+\sum_{l\in \{0,1,\ldots,2n_i-1 \} \setminus \{n_i \}}\left(
(n_i-l)\theta_{l,t_i}^{-} \hat{\delta}_1(t_i)\cdot
\frac{\hat{\delta}_2(\theta_{2n_i-l-1,t_i}^{-})}{n_i-l}
-\hat{\delta}_1(\theta_{l,t_i}^{-} )\cdot
\theta_{2n_i-l-1,t_i}^{-}\hat{\delta}_2(t_i)  \right)\\
&+R_i \cdot \hat{\delta}_1(t_i)
\cdot \hat{\delta}_2(t_i),
\end{aligned}
\end{equation*}
where $R_i$ is a rational functions 
on $\widehat{\mathcal{M}}_{\boldsymbol{t}_{\text{ra}}}\times T_{\boldsymbol{\theta}}$.
We may check that 
$$
\begin{aligned}
&\sum_{l\in \{0,1,\ldots,2n_i-2 \} \setminus \{n_i-1 \}}\left( \hat{\delta}_1(\theta_{l,t_i}^{\pm})\cdot
\frac{\hat{\delta}_2(\theta_{2n_i-l-2,t_i}^{\pm})}{n_i-l-1}\right) \\
&=\sum_{l=0}^{n_i-2 }\hat{\delta}_1(\theta_{l,t_i}^{\pm})\cdot
\frac{\hat{\delta}_2(\theta_{2n_i-l-2,t_i}^{\pm})}{n_i-l-1}
-\sum_{l=0}^{n_i-2 } 
\frac{\hat{\delta}_1(\theta_{2n_i-l-2,t_i}^{\pm})}{n_i-l-1}
\cdot \hat{\delta}_2(\theta_{l,t_i}^{\pm})
\end{aligned}
$$
and
$$
\begin{aligned}
&\sum_{l\in \{0,1,\ldots,2n_i-1 \} \setminus \{n_i \}}\left(
(n_i-l)\theta_{l,t_i}^{\pm} \hat{\delta}_1(t_i)\cdot
\frac{\hat{\delta}_2(\theta_{2n_i-l-1,t_i}^{\pm})}{n_i-l}
-\hat{\delta}_1(\theta_{l,t_i}^{\pm} )\cdot
\theta_{2n_i-l-1,t_i}^{\pm}\hat{\delta}_2(t_i)  \right) \\
&=
\hat{\delta}_1(t_i) \cdot \hat{\delta}_2 \left( 
\sum_{l=0}^{n_i-1}    \theta_{l,t_i}^{\pm} \cdot
\theta_{2n_i-l-1,t_i}^{\pm} \right)
-\hat{\delta}_1 \left( 
\sum_{l=0}^{n_i-1}    \theta_{l,t_i}^{\pm} \cdot
\theta_{2n_i-l-1,t_i}^{\pm} \right) \cdot 
\hat{\delta}_2(t_i),
\end{aligned}
$$
since $\hat{\delta}(\theta_{n_i-1,t_i}^{\pm})=0$.
Remark that $\hat{\delta}_1(t_i) =\hat{\delta}_2(t_i)=0$ for $i=0,1,\infty$.
Then we have 
\begin{equation}\label{2020.1.20.21.33}
\begin{aligned}
&\frac{1}{2}\sum_{i\in I} \mathrm{res}_{x=t_i} 
\mathrm{Tr} \left( \hat{\delta} (\tilde{\Omega}_i) 
\wedge \hat{\delta}(\tilde{\psi}_i)\tilde{\psi}_i^{-1} 
   \right)\\
&=\left(
\sum_{i\in I_{\text{un}}} \sum_{l=0}^{n_{i}-2} \left( d  H_{\theta_{l,t_i}^{+}} \wedge d \theta_{l,t_i}^{+}
+ d H_{\theta_{l,t_i}^{-}} \wedge d \theta_{l,t_i}^{-} \right)
+\sum_{i\in \{3,4,\ldots,\nu \}} d H_{t_i} \wedge dt_i \right) ( \hat{\delta}_1, \hat{\delta}_2).
\end{aligned}
\end{equation}
This is just the second part of \eqref{2020.1.21.13.49}.
We may take $\Phi_i$ and $\Xi_{i}^{\le 2n_i-1} (x_{t_i})$ 
so that $\Phi_i \Xi_{i}^{\le 2n_i-1} (x_{t_i})$ satisfies the assumption of Lemma \ref{2020.1.20.21.24}
(see Remark \ref{2021.4.30.17.15}).
Since $\widehat\Omega_{(\boldsymbol{t}_{\mathrm{ra}},\boldsymbol{\theta}_0)}^{(n-2)}$ 
also satisfies the assumption of Lemma \ref{2020.1.20.21.24},
the difference the residue (\ref{2020.1.20.21.33}) and the residue of 
$\frac{1}{2} \cdot  \mathrm{Tr}
 \left( \hat{\delta}
 (\widehat\Omega_{(\boldsymbol{t}_{\mathrm{ra}},\boldsymbol{\theta}_0)}^{(n-2)} ) 
\wedge \hat{\delta}(\psi_i)\psi_i^{-1}  \right)$
at $t_i$ is zero if $\delta_1 \in 
\Theta_{(\widehat{\mathcal{M}}_{\boldsymbol{t}_{\text{ra}}}\times T_{\boldsymbol{\theta}}) 
/((T_{\boldsymbol{t}})_{\boldsymbol{t}_{\text{ra}}} \times T_{\boldsymbol{\theta}})}$.

Second we calculate 
the residue of 
$\mathrm{Tr} \left( \hat{\delta}
 (\widehat\Omega_{(\boldsymbol{t}_{\mathrm{ra}},\boldsymbol{\theta}_0)}^{(n-2)} ) 
\wedge \hat{\delta}(\psi_{q_j})\psi_{q_j}^{-1}  \right)$ 
at $x=q_j$.
First, $\hat{\delta} 
(\widehat\Omega_{(\boldsymbol{t}_{\mathrm{ra}},\boldsymbol{\theta}_0)}^{(n-2)} )$ 
is described at $x=q_j$ 
as follows:
\begin{equation*}
\hat{\delta} (\widehat\Omega_{(\boldsymbol{t}_{\mathrm{ra}},\boldsymbol{\theta}_0)}^{(n-2)} )=
\begin{pmatrix}
0& \hat{\delta}(1/P(x;\boldsymbol{t})) \\
\hat{\delta}(c_0) & \hat{\delta}(d_0)
\end{pmatrix},
\end{equation*}
where
\begin{equation*}
\begin{aligned}
\hat{\delta}(1/P(x;\boldsymbol{t}))&= \sum_{i=1}^{\nu}
\frac{\partial}{\partial t_i} \left( \frac{1}{P(q_j;\boldsymbol{t})} \right) \hat{\delta}(t_i)
+ O(x-q_j),\\
\hat{\delta}(c_0)&=
\frac{p_j \hat{\delta} (q_j)}{(x-q_j)^2}
+\frac{\hat{\delta}(p_j)}{(x-q_j)}
 + O(x-q_j)^0, \text{ and} \\
\hat{\delta}(d_0)&=
-\frac{\hat{\delta}(q_j)}{(x-q_j)^2}
-\sum_{k\neq j}\frac{\hat{\delta}(q_k)}{(q_j-q_k)^2}
+\sum_{i=3}^{\nu}\left(\frac{\partial d_0}{\partial t_i}(q_j) \hat{\delta} (t_i)\right)
+\sum_{i \in I_{\mathrm{un}}} \sum_{l=0}^{n_i-2} \left(
 \frac{\partial d_0}{\partial \theta_{l,t_i}^{\pm}}(q_j) \hat{\delta} (\theta_{l,t_i}^{\pm})
\right)
+ O(x-q_j).
\end{aligned}
\end{equation*}
Second, we consider $\hat{\delta}(\psi_{q_j})\psi_{q_j}^{-1}$. 
By Lemma \ref{2019.12.30.21.53}, we have
\begin{equation*}
\hat{\delta}(\psi_{q_j})\psi_{q_j}^{-1} 
= \hat{\delta}(\Phi_{q_j} \Xi_{q_j}(x) )(\Phi_{q_j} \Xi_{q_j}(x))^{-1}+ 
(\Phi_{q_j} \Xi_{q_j}(x) ) 
\begin{pmatrix}
0 & 0 \\
0 &\frac{-\hat{\delta}(q_j)}{x-q_j}
\end{pmatrix}
 (\Phi_{q_j} \Xi_{q_j}(x))^{-1}.
\end{equation*}
By using Lemma \ref{2019.12.30.21.53},
we have the following equality:
$$
\begin{aligned}
\Phi_{q_j}\Xi_{q_j}(x)
&= \begin{pmatrix}
1 & 0 \\
p_j & 1 
\end{pmatrix}
+
\begin{pmatrix}
-\frac{p_j}{P(q_j;\boldsymbol{t})} & -\frac{1}{2P(q_j;\boldsymbol{t})} \\
-\frac{p_j^2}{P(q_j;\boldsymbol{t})} & \frac{p_j}{2P(q_j;\boldsymbol{t})} - \sum_{i=1}^{\nu} \frac{D_i(q_j)}{(q_j-t_i)^{n_i}}
+ \sum_{k\neq j} \frac{1}{q_j - q_k}- D_{\infty} (q_j)
\end{pmatrix}(x-q_j) \\
&\qquad + O(x-q_j)^2.
\end{aligned}
$$
By this description of $\Phi_{q_j}\Xi_{q_j}(x)$, we may check that 
the constant term of the expansion of $\delta(\Phi_{q_j} \Xi_{q_j}(x))$ at $q_j$ 
has following description:
$$
\begin{pmatrix}
0 & 0 \\
\delta(p_j) & 0
\end{pmatrix} - \delta(q_j)
\begin{pmatrix}
-\frac{p_j}{P(q_j;\boldsymbol{t})} & -\frac{1}{2P(q_j;\boldsymbol{t})} \\
-\frac{p_j^2}{P(q_j;\boldsymbol{t})} &*
\end{pmatrix}. 
$$
Since $   \sum_{i=1}^{\nu} \frac{D_i(x)}{(x -t_i)^{n_i}} 
 + D_{\infty}(x)=d_0-\sum_{j=1}^{n-3} \frac{-1}{x-q_j}$,
 the coefficient of
the $(x-q_j)$-term of the expansion of $\delta(\Phi_{q_j} \Xi_{q_j}(x))$ 
has following description:
$$
\begin{aligned}
&\begin{pmatrix}
* &0 \\
* & \frac{\delta(p_j)}{2P(q_j;\boldsymbol{t})} + \sum_{k \neq j} \frac{\delta(q_k)}{(q_j - q_k)^2}
\end{pmatrix} - \delta(q_j)
\begin{pmatrix}
* & * \\
* &*
\end{pmatrix}  \\
&\quad +
\begin{pmatrix}
* & -\frac{1}{2} \sum_{i=1}^{\nu}
\frac{\partial}{\partial t_i} \left( \frac{1}{P(q_j;\boldsymbol{t})} \right) \hat{\delta}(t_i) \\
* &\frac{p_j}{2}\sum_{i=3}^{\nu}
\frac{\partial}{\partial t_i} \left( \frac{1}{P(q_j;\boldsymbol{t})} \right) \hat{\delta}(t_i)
 -\sum_{i=3}^{\nu}\left(\frac{\partial d_0}{\partial t_i}(q_j) \hat{\delta} (t_i)\right)
-\sum_{i \in I_{\mathrm{un}}} \sum_{l=0}^{n_i-2} \left(
 \frac{\partial d_0}{\partial \theta_{l,t_i}^{\pm}}(q_j)
  \hat{\delta} (\theta_{l,t_i}^{\pm}) \right)
\end{pmatrix}.
\end{aligned}
$$
Here we put the entries having $\delta(q_j)$
together in the second matrices.
Moreover,
we may check that 
the constant term of the expansion of $(\Phi_{q_j} \Xi_{q_j}(x))^{-1}$ at $q_j$ 
is 
$\begin{pmatrix}
1 & 0 \\
-p_j & 1
\end{pmatrix}$
and the coefficient of
the $(x-q_j)$-term of the expansion of $(\Phi_{q_j} \Xi_{q_j}(x))^{-1}$ 
has following description:
$$
-\begin{pmatrix}
-\frac{p_j}{2P(q_j;\boldsymbol{t})} & -\frac{1}{2P(q_j;\boldsymbol{t})} \\
* & \frac{p_j}{P(q_j;\boldsymbol{t})} - \sum_{i=1}^{\nu} \frac{D_i(q_j)}{(q_j-t_i)^{n_i}}
+ \sum_{k\neq j} \frac{1}{q_j - q_k}- D_{\infty} (q_j)
\end{pmatrix} .
$$
By the calculation of $\delta(\Phi_{q_j} \Xi_{q_j}(x))$ and
$(\Phi_{q_j} \Xi_{q_j}(x))^{-1}$,
we may show that 
$\delta(\Phi_{q_j} \Xi_{q_j}(x)) (\Phi_{q_j} \Xi_{q_j}(x))^{-1}$ is
\begin{equation*}
\begin{pmatrix}
* & \frac{\hat{\delta}(q_j)}{2P(q_j;\boldsymbol{t})} \\
*& *
\end{pmatrix}
+ 
\begin{pmatrix}
* &  x_{12} \\
* &  \frac{\hat{\delta}(p_j)}{P(q_j;\boldsymbol{t})} 
+ \sum_k \frac{ \hat{\delta}(q_k)}{(q_j - q_k)^2}
+x_{22}
\end{pmatrix}(x-q_j)+ O(x-q_j)^2.
\end{equation*}
Here we put
\begin{equation*}
\begin{aligned}
&x_{12}:= f_{12} \hat{\delta}(q_j)-
\frac{1}{2}\sum_{i=3}^{\nu}
\frac{\partial}{\partial t_i} \left( \frac{1}{P(q_j;\boldsymbol{t})} \right)\hat{\delta}(t_i) \text{ and}\\
&x_{22}:=f_{22} \hat{\delta}(q_j)
-\sum_{i=3}^{\nu}\left(\frac{\partial d_0}{\partial t_i}(q_j) \hat{\delta} (t_i) \right)
-\sum_{i \in I_{\mathrm{un}}} \sum_{l=0}^{n_i-2} \left( 
\frac{\partial d_0}{\partial \theta_{l,t_i}^{\pm}}(q_j) \hat{\delta} (\theta_{l,t_i}^{\pm})
\right)
+\frac{p_j}{2}\sum_{i=3}^{\nu}
\frac{\partial}{\partial t_i} \left( \frac{1}{P(q_j;\boldsymbol{t})} \right) \hat{\delta}(t_i),
\end{aligned}
\end{equation*}
where $f_{12}$ and $f_{22}$ are rational functions 
on $\widehat{\mathcal{M}}_{\boldsymbol{t}_{\text{ra}}}\times T_{\boldsymbol{\theta}}$.
Moreover we may show that 
\begin{equation*}
\begin{aligned}
&(\Phi_{q_j} \Xi_{q_j}(x) ) 
\begin{pmatrix}
0 & 0 \\
0 &\frac{-\hat{\delta}(q_j)}{x-q_j}
\end{pmatrix}
 (\Phi_{q_j} \Xi_{q_j}(x))^{-1}\\
 &=
\frac{\begin{pmatrix}
0 & 0 \\
p_j\hat{\delta}(q_j)& -\hat{\delta}(q_j)
\end{pmatrix}}{x-q_j}
+ 
\begin{pmatrix}
* &  \frac{ \hat{\delta}(q_j)}{2P(q_j;\boldsymbol{t})} \\
* &  g^{(0)}_{22} \hat{\delta}(q_j)
\end{pmatrix}+ 
\begin{pmatrix}
* &  g^{(1)}_{12} \hat{\delta}(q_j) \\
* &  g^{(1)}_{22} \hat{\delta}(q_j)
\end{pmatrix}(x-q_j)
+O(x-q_j)^2,
\end{aligned}
\end{equation*}
where $g^{(0)}_{22}$, $g^{(1)}_{12}$ and $g^{(1)}_{22}$ are rational functions 
on $\widehat{\mathcal{M}}_{\boldsymbol{t}_{\text{ra}}}\times T_{\boldsymbol{\theta}}$.
Finally, we have
\begin{equation*}
\begin{aligned}
\frac{1}{2} \cdot \mathrm{res}_{x=q_j} 
\mathrm{Tr} \left( \hat{\delta}
 (\widehat\Omega_{(\boldsymbol{t}_{\mathrm{ra}},\boldsymbol{\theta}_0)}^{(n-2)} ) 
\wedge\hat{\delta}(\psi_{q_j})\psi_{q_j}^{-1} 
\right) 
&= 
\frac{\hat{\delta}_1(p_j)\hat{\delta}_2(q_j)}{P(q_j;\boldsymbol{t})}
-\frac{\hat{\delta}_2(p_j)\hat{\delta}_1(q_j)}{P(q_j;\boldsymbol{t})}
+\sum_{k\neq j}\frac{\hat\delta_1(q_k)\hat\delta_2(q_j)}{(q_j-q_k)^2}
-\sum_{k\neq j}\frac{\hat\delta_1(q_j)\hat\delta_2(q_k)}{(q_j-q_k)^2} \\
&\quad
+p_j\sum_{i=3}^{\nu}\frac{\partial}{\partial t_i} \left( \frac{1}{P(q_j;\boldsymbol{t})} \right)
\left(\hat{\delta}_{1}(t_i)  \hat{\delta}_2(q_j)
-\hat{\delta}_{2}(t_i)  \hat{\delta}_1(q_j) \right)\\
&\qquad
-  \sum_{i=3}^{\nu}\frac{\partial d_0}{\partial t_i}(q_j)
\left(  \hat{\delta}_1 (t_i) \hat{\delta}_2 (q_j) -\hat{\delta}_2 (t_i) \hat{\delta}_1 (q_j) \right)\\
&\qquad \quad
-  \sum_{i \in I_{\mathrm{un}}} \sum_{l=0}^{n_i-2} 
 \frac{\partial d_0}{\partial \theta_{l,t_i}^{\pm}}(q_j)
\left(  \hat{\delta}_1 (\theta_{l,t_i}^{\pm}) \hat{\delta}_2 (q_j) 
-\hat{\delta}_2 (\theta_{l,t_i}^{\pm}) \hat{\delta}_1 (q_j) \right).
\end{aligned}
\end{equation*}
Moreover, this is equal to
$$
\begin{aligned}
& \left( \frac{ \hat{\delta}_1(p_j)}{P(q_j;\boldsymbol{t}) }
+ p_j \sum_{i=3}^{\nu}\frac{\partial}{\partial t_i} \left( \frac{1}{P(q_j;\boldsymbol{t})} \right) \hat\delta_1(t_i)
-\sum_{i=3}^{\nu}\frac{\partial d_0}{\partial t_i}(q_j) \hat\delta_1(t_i)
- \sum_{i \in I_{\mathrm{un}}} \sum_{l=0}^{n_i-2} 
 \frac{\partial d_0}{\partial \theta_{l,t_i}^{\pm}}(q_j) \hat\delta_1(\theta_{l,t_i}^{\pm})
 \right)\hat{\delta}_2(q_j)\\
& - \left( \frac{ \hat{\delta}_2(p_j)}{P(q_j;\boldsymbol{t}) }
+ p_j \sum_{i=3}^{\nu}\frac{\partial}{\partial t_i} \left( \frac{1}{P(q_j;\boldsymbol{t})} \right) \hat\delta_2(t_i)
-\sum_{i=3}^{\nu}\frac{\partial d_0}{\partial t_i}(q_j) \hat\delta_2(t_i)
- \sum_{i \in I_{\mathrm{un}}} \sum_{l=0}^{n_i-2} 
 \frac{\partial d_0}{\partial \theta_{l,t_i}^{\pm}}(q_j) \hat\delta_2(\theta_{l,t_i}^{\pm})
 \right)\hat{\delta}_1(q_j)\\
 &\quad-\sum_{k\neq j}\frac{\hat\delta_1(q_k)\hat\delta_2(q_j)
 -\hat\delta_1(q_j)\hat\delta_2(q_k)}{(q_j-q_k)^2}.
\end{aligned}
$$
In the first term and the second term of this formula,
the exterior derivative of $\frac{p_j}{P(q_j;\boldsymbol{t})} 
  - \sum_{i=1}^{\nu} \frac{D_i(q_j;\boldsymbol{t}, \boldsymbol{\theta})}{(q_j -t_i)^{n_i}} 
  - D_{\infty}(q_j;\boldsymbol{t}, \boldsymbol{\theta})$ on the extended moduli space 
  $\widehat{\mathcal{M}}_{\boldsymbol{t}_{\text{ra}}}\times T_{\boldsymbol{\theta}} $
  appears.
  Here remark that $ 
   \sum_{i=1}^{\nu} \frac{D_i(x;\boldsymbol{t}, \boldsymbol{\theta})}{(x -t_i)^{n_i}} 
 + D_{\infty}(x;\boldsymbol{t}, \boldsymbol{\theta})=
d_0-\sum_{j=1}^{n-3} \frac{-1}{x-q_j} $
and the coefficients of $D_i(x;\boldsymbol{t}, \boldsymbol{\theta})$ and 
$D_{\infty}(x;\boldsymbol{t}, \boldsymbol{\theta})$ are 
independent of the parameters $\{(q_j,p_j) \}_{j=1,2,\ldots,n-3}$.
Since $\sum_{j=1}^{n-3}\sum_{k\neq j}\frac{\hat\delta_1(q_k)\hat\delta_2(q_j)
 -\hat\delta_1(q_j)\hat\delta_2(q_k)}{(q_j-q_k)^2}=0$,
we have
\begin{equation*}
\begin{aligned}
&\frac{1}{2}\sum_{j=1}^{n-3}
 \mathrm{res}_{x=q_j} 
\mathrm{Tr} \left( \hat{\delta}
 (\widehat\Omega_{(\boldsymbol{t}_{\mathrm{ra}},\boldsymbol{\theta}_0)}^{(n-2)}) 
\wedge \hat{\delta}(\psi_{q_j})\psi_{q_j}^{-1} 
   \right)\\
&=\left(
 \sum_{j=1}^{n-3} d \left(\frac{p_j}{P(q_j;\boldsymbol{t})} 
  - \sum_{i=1}^{\nu} \frac{D_i(q_j;\boldsymbol{t}, \boldsymbol{\theta})}{(q_j -t_i)^{n_i}} 
  - D_{\infty}(q_j;\boldsymbol{t}, \boldsymbol{\theta})\right) \wedge dq_j \right) (\hat{\delta}_1,\hat{\delta}_2). 
\end{aligned}
\end{equation*}
We obtain the assertion of this theorem.
\end{proof}

\begin{Cor}\label{2020.1.21.15.11}
\textit{
Set $\eta_j:= \frac{p_j}{P(q_j;\boldsymbol{t})}
 - \sum_{i=1}^{\nu} \frac{D_i(q_j;\boldsymbol{t}, \boldsymbol{\theta})}{(q_j -t_i)^{n_i}} 
 - D_{\infty}(q_j;\boldsymbol{t}, \boldsymbol{\theta})$.
The vector fields $\delta_{\theta_{l,t_i}^{\pm}}^{\mathrm{IMD}}$ 
$(i\in I_{\mathrm{un}}$ and $l=0,1,\ldots,n_i-2)$
and
$\delta_{t_i}^{\mathrm{IMD}}$ $(i=3,4,\ldots,\nu)$
have the following
hamiltonian description}: \textit{
\begin{equation}\label{2020.1.21.14.00}
\begin{aligned}
\delta_{\theta_{l,t_i}^{\pm}}^{\mathrm{IMD}}&= 
\frac{\partial}{\partial \theta_{l,t_i}^{\pm}}
-\sum_{j=1}^{n-3} \left( 
 \frac{\partial H_{\theta_{l,t_i}^{\pm}}}{\partial \eta_j} 
\frac{\partial}{\partial q_j} 
- \frac{\partial H_{\theta_{l,t_i}^{\pm}}}{\partial q_j} 
\frac{\partial}{\partial \eta_j}
\right)
\text{ and} \\
\delta_{t_i}^{\mathrm{IMD}}&= 
\frac{\partial}{\partial t_i}
-\sum_{j=1}^{n-3} \left( 
 \frac{\partial  H_{t_i}}{\partial \eta_j} 
\frac{\partial}{\partial q_j} 
- \frac{\partial H_{t_i}}{\partial q_j} 
\frac{\partial}{\partial \eta_j}
\right),
\end{aligned}
\end{equation}
respectively.
}
\end{Cor}

\begin{proof}
We can put
\begin{equation*}
\begin{aligned}
\delta_{\theta_{l,t_i}^{\pm}}^{\mathrm{IMD}}= 
\frac{\partial}{\partial \theta_{l,t_i}^{\pm}}
+\sum_{j=1}^{n-3} \left( 
X^j_{\theta_{l,t_i}^{\pm}}
\frac{\partial}{\partial q_j} 
+Y^j_{\theta_{l,t_i}^{\pm}}
\frac{\partial}{\partial \eta_j}
\right)
\text{ and } 
\delta_{t_i}^{\mathrm{IMD}}= 
\frac{\partial}{\partial t_i}
+\sum_{j=1}^{n-3} \left( 
X_{t_i}^j
\frac{\partial}{\partial q_j} 
+Y_{t_i}^j
\frac{\partial}{\partial \eta_j}
\right).
\end{aligned}
\end{equation*}
By Theorem \ref{2020.1.21.13.48},
the terms of $dq_j, d\eta_j$ ($j=1,2,\ldots , n-3$) of the $1$-forms
$\hat{\omega}(\delta_{\theta_{l,t_i}^{\pm}}^{\mathrm{IMD}}, *)$
and
$\hat{\omega}(\delta_{t_i}^{\mathrm{IMD}}, *)$
are
\begin{equation*}
\begin{aligned}
 & \sum_{j=1}^{n-3} 
\left( -X^j_{\theta_{l,t_i}^{\pm}} d\eta_j 
+Y^j_{\theta_{l,t_i}^{\pm}} dq_j  \right)  
-\sum_{j=1}^{n-3} 
\left( \frac{\partial H_{\theta_{l,t_i}^{\pm}}}{\partial \eta_j} d\eta_j
+  \frac{\partial H_{\theta_{l,t_i}^{\pm}}}{\partial q_j} dq_j \right) 
 \text{ and}\\ 
&\sum_{j=1}^{n-3} 
\left( -X^j_{t_i} d\eta_j 
+Y^j_{t_i} dq_j  \right)  
-\sum_{j=1}^{n-3} 
\left( \frac{\partial H_{t_i}}{\partial \eta_j} d\eta_j
+  \frac{\partial H_{t_i}}{\partial q_j} dq_j \right),
\end{aligned}
\end{equation*}
respectively.
By Theorem \ref{2020.1.21.13.51}, we have the following equalities
\begin{equation*}
\begin{aligned}
X^j_{\theta_{l,t_i}^{\pm}}&= -\frac{\partial H_{\theta_{l,t_i}^{\pm}}}{\partial \eta_j},  &
Y^j_{\theta_{l,t_i}^{\pm}} &= \frac{\partial H_{\theta_{l,t_i}^{\pm}}}{\partial q_j},  \\
X^j_{t_i} &=-\frac{\partial H_{t_i}}{\partial \eta_j}, \text{ and}  &
Y^j_{\theta_{l,t_i}^{\pm}}&= \frac{\partial H_{t_i}}{\partial q_j}.
\end{aligned}
\end{equation*}
Then we have the hamiltonian description (\ref{2020.1.21.14.00}). 
\end{proof}

\begin{Rem}
In \cite{Kimura}, 
the Hamiltonian systems of the two-dimensional (degenerated) Garnier systems
have been described
by using the coordinates $(q_j,\tilde{\eta}_j)_{0\le j \le 2}$, 
where $\tilde{\eta}_j:= -\frac{p_j}{P(q_j;\boldsymbol{t})}$.
In these cases, the $2$-form
$d \left(  - \sum_{i=1}^{\nu} \frac{D_i(q_j;\boldsymbol{t}, \boldsymbol{\theta})}{(q_j -t_i)^{n_i}} 
- D_{\infty}(q_j;\boldsymbol{t}, \boldsymbol{\theta})\right)
  \wedge dq_j$,
which comes from the residue at an apparent singularity,
is canceled by some terms of the $2$-form 
$\sum_{i\in I } \sum_{l=0}^{n_{i}-2} \left( d  H_{\theta_{l,t_i}^{+}} \wedge d \theta_{l,t_i}^{+}
+ d H_{\theta_{l,t_i}^{-}} \wedge d \theta_{l,t_i}^{-} \right)$,
which comes from the residues at unramified irregular singular points.
\end{Rem}

\section{Ramified irregular singularities}

In this section,
we assume that $ I_{\text{ra}}\neq \emptyset$.
For $i\in I_{\text{ra}}$,
the leading coefficient $\Omega_{t_i}(0)$ is a non trivial Jordan block.
In Section \ref{2020.2.8.10.32}, we define a $2$-form on the fiber
$\mathcal{M}_{\boldsymbol{t}_{0},\boldsymbol{t}_{\text{ra}}}$
by Krichever's formula \cite[Section 5]{Krichever}.
Remark that $\mathcal{M}_{\boldsymbol{t}_{0},\boldsymbol{t}_{\text{ra}}}$
is isomorphic to the moduli space
$\mathfrak{Conn}_{(\boldsymbol{t}_0,\boldsymbol{\theta},\boldsymbol{\theta}_0)}$.
We show that this $2$-form coincides with the symplectic form \eqref{2020.2.7.11.24}.
In Section \ref{2020.2.7.11.34},
we will construct horizontal lifts of $\tilde{\nabla}_{\mathrm{DL, ext}}^{(1)}$.
Let
$\partial / \partial \theta_{l',t_i}$ ($i \in I_{\mathrm{ra}}$, $l'=0,1,\ldots,2n_i-3$)
be the vector fields on $(T_{\boldsymbol{t}})_{\boldsymbol{t}_{\text{ra}}} \times T_{\boldsymbol{\theta}}$.
By the construction of the horizontal lifts,
we have the vector field
$\delta^{\mathrm{IMD}}_{\theta_{l',t_i}}$
on
$\widehat{\mathcal{M}}_{\boldsymbol{t}_{\text{ra}}}\times T_{\boldsymbol{\theta}}$
determined by the integrable deformations
with respect to $\partial / \partial \theta_{l',t_i}$.
Remark that 
$\widehat{\mathcal{M}}_{\boldsymbol{t}_{\text{ra}}}\times T_{\boldsymbol{\theta}}$ 
is isomorphic to the extended moduli space
$\widehat{\mathfrak{Conn}}_{(\boldsymbol{t}_{\mathrm{ra}},\boldsymbol{\theta}_0)}$.
In Section \ref{2020.2.8.10.33},
we define a $2$-form on 
$\widehat{\mathcal{M}}_{\boldsymbol{t}_{\text{ra}}}\times T_{\boldsymbol{\theta}}$
by Krichever's formula.
We show that this $2$-form is the isomonodromy $2$-form.
In Section \ref{2020.2.8.10.35},
we calculate this $2$-form on 
$\widehat{\mathcal{M}}_{\boldsymbol{t}_{\text{ra}}}\times T_{\boldsymbol{\theta}}$
by using Diarra--Loray's global normal form.
Then we obtain an explicit formula of this $2$-form.

For $i \in I \setminus I_{\text{ra}}$,
we fix a compatible framing $\Phi_i$ and
take $\Xi_{i}(x_{t_i})$ as in Lemma \ref{2019.12.30.22.09}.
For each $i \in I_{\text{ra}}$,
we consider the leading coefficient of 
$\Omega^{(n-2)}
_{(\boldsymbol{t}_{0},\boldsymbol{\theta},\boldsymbol{\theta}_0)}$
at $t_i$:
$$
\Omega^{(n-2)}
_{(\boldsymbol{t}_{0},\boldsymbol{\theta},\boldsymbol{\theta}_0)}
= \begin{pmatrix}
0 & \prod_{j \neq i}(t_i - t_j)^{-n_j}  \\
-\frac{\theta_{0,t_i}^2}{4} \cdot \prod_{j \neq i}(t_i - t_j)^{n_j} & \theta_{0,t_i}
\end{pmatrix}\frac{ dx_{t_i}}{x_{t_i}^{n_i}}
 + [\text{ higher order terms }].
$$
This leading coefficient at $t_i$ is independent of $\{ (q_j,p_j) \}_{j=1,2,\ldots,n-3}$.
We fix $\Phi_i \in \mathrm{GL}(2,\mathbb{C})$ so that 
\begin{equation}\label{2021.4.21.10.45}
\Phi_i^{-1} 
\Omega^{(n-2)}
_{(\boldsymbol{t}_{0},\boldsymbol{\theta},\boldsymbol{\theta}_0)}
 \Phi_i
= \begin{pmatrix}
\frac{\theta_{0,t_i}}{2}& \frac{\theta_{1,t_i}}{2} \\
0 & \frac{\theta_{0,t_i}}{2}
\end{pmatrix}\frac{dx_{t_i}}{x_{t_i}^{n_i}}
 + [\text{ higher order terms }].
\end{equation}
We call the matrix $\Phi_i$ a {\it compatible framing at }$t_i$.
If we have another $\Phi_i'$ 
such that the leading coefficient matrix of
 $(\Phi_i')^{-1} \Omega^{(n-2)}
_{(\boldsymbol{t}_{0},\boldsymbol{\theta},\boldsymbol{\theta}_0)} \Phi'_i$
is an upper triangular matrix as in \eqref{2021.4.21.10.45},
then there exists an upper triangular matrix
$$
C_{t_i} = 
\begin{pmatrix}
c_{t_i,11} & c_{t_i,12}\\
0 & c_{t_i,11}
\end{pmatrix}
$$ 
such that
$\Phi_i' =\Phi_i C_{t_i}$.
We define $\zeta_i$ as $x_{t_i}=\zeta_i^2$.
Let $M_{\zeta_i}$ be the matrix \eqref{2020.1.10.13.58}.
For the compatible framing $\Phi_i$,
there exist unique
\begin{itemize}
\item formal transformation 
\begin{equation}\label{2021.4.23.18.23}
\Xi_i(x_{t_i}) := 
\begin{pmatrix}
1 & 0\\
0 & 1
\end{pmatrix}+
\sum_{s=1}^{\infty} 
\begin{pmatrix}
(\xi^{(i)}_s)_{11} & (\xi^{(i)}_s)_{12} \\
(\xi^{(i)}_s)_{21} & (\xi^{(i)}_s)_{22}
\end{pmatrix}
x_{t_i}^s , \text{ and}
\end{equation}

\item $\theta_{l',t_i} \in 
\Gamma(\mathcal{M}_{\boldsymbol{t}_{0},\boldsymbol{t}_{\text{ra}}} ,
\mathcal{O}_{\mathcal{M}_{\boldsymbol{t}_{0},\boldsymbol{t}_{\text{ra}}}})$ 
($l'\geq 2n_i -2$ and $i \in I_{\text{ra}}$)
\end{itemize}
such that
\begin{itemize}
\item[(1)] we have the following equality  
\begin{equation}\label{2021.4.23.12.09}
(\Phi_i\Xi_i (x_{t_i}))^{-1} 
d(\Phi_i \Xi_i (x_{t_i})) +(\Phi_i\Xi_i (x_{t_i}))^{-1} 
\Omega^{(n-2)}
_{(\boldsymbol{t}_{0},\boldsymbol{\theta},\boldsymbol{\theta}_0)}
(\Phi_i\Xi_i (x_{t_i})) =
\begin{pmatrix}
\alpha_i & \beta_i \\
x_{t_i}\beta_i & \alpha_i -\frac{dx_{t_i}}{2x_{t_i}}
\end{pmatrix},
\end{equation}
where we set
\begin{equation*}
\begin{cases}
\alpha_i:= \frac{\theta_{0,t_i}}{2}\frac{dx_{t_i}}{x_{t_i}^{n_i}}+ \cdots 
+\frac{\theta_{2l,t_i}}{2} \frac{dx_{t_i}}{x_{t_i}^{n_i-l}}+\cdots 
+\frac{\theta_{2n_i-2,t_i}}{2} \frac{dx_{t_i}}{x_{t_i}} + \cdots \\
\beta_i:= \frac{\theta_{1,t_i}}{2}\frac{dx_{t_i}}{x_{t_i}^{n_i}}+ \cdots 
+ \frac{\theta_{2l+1,t_i}}{2}\frac{dx_{t_i}}{x_{t_i}^{n_i-l}}+\cdots 
+\frac{\theta_{2n_i-3,t_i}}{2}\frac{dx_{t_i}}{x^2_{t_i}}+\cdots,
\end{cases}
\end{equation*}

\item[(2)] there exists a formal power series $\xi(\zeta_i) \in 
\Gamma(\mathcal{M}_{\boldsymbol{t}_{0},\boldsymbol{t}_{\text{ra}}} ,
\mathcal{O}_{\mathcal{M}_{\boldsymbol{t}_{0},\boldsymbol{t}_{\text{ra}}}})[[ \zeta_i]]$ such that 
\begin{equation}\label{2022.4.9.15.44}
M_{\zeta_i}^{-1} \Xi_i (\zeta_{i}^2) M_{\zeta_i} = 
\begin{pmatrix}
1 & \zeta_i \cdot \xi(\zeta_i)  \\
-\zeta_i \cdot \xi(-\zeta_i) & 1
\end{pmatrix}.
\end{equation}
\end{itemize}
Indeed, the $\zeta_i^{-2n_i+1}$-term
and the $\zeta_i^{-2n_i+2}$-term of the expansion of 
\begin{equation}\label{2022.4.10.10.26}
(\Phi_iM_{\zeta_i})^{-1} 
d(\Phi_i M_{\zeta_i}) +(\Phi_iM_{\zeta_i})^{-1} 
\Omega^{(n-2)}
_{(\boldsymbol{t}_{0},\boldsymbol{\theta},\boldsymbol{\theta}_0)}(\zeta_i^2)
(\Phi_iM_{\zeta_i})
\end{equation}
at $\zeta_i=0$ are diagonal.
The eigenvalues of the $\zeta_i^{-2n_i+2}$-term are distinct.
After the $\zeta_i^{-2n_i+2}$-term, 
we can diagonalize \eqref{2022.4.10.10.26}
by a matrix as in Lemma \ref{2019.12.30.22.09}.
By the argument as in the proof of \cite[Proposition 10]{DF},
we may check that, in this situation, 
this matrix has a form as in the right hand side of \eqref{2022.4.9.15.44}.
Moreover by the gauge transformation of this diagonal matrix by $M_{\zeta_i}^{-1}$,
we have the right hand side of \eqref{2021.4.23.12.09}.
We may check that 
$$
 M_{\zeta_i} 
 \begin{pmatrix}
1 & \zeta_i \cdot \xi(\zeta_i)  \\
-\zeta_i \cdot \xi(-\zeta_i) & 1
\end{pmatrix}
M_{\zeta_i}^{-1}
$$
is invariant under replacing $\zeta_i$ with $-\zeta_i$.
So we have $\Xi_i(x_{t_i})$.
By Lemma \ref{2019.12.30.22.09}, 
such $\Xi_i(x_{t_i})$ is unique.

Let $U_{t_i}$ be an affine open subset on $\mathbb{P}^1$ for $i \in I_{\text{ra}}$
so that $x_{t_i}$ is a coordinate on $U_{t_i}$.
Let $U_{\zeta_i}$ be the inverse image of 
$U_{t_i}$ under the map 
$\mathrm{Spec}\, \mathbb{C}[\zeta_i] \rightarrow
\mathrm{Spec}\, \mathbb{C}[x_{t_i}]$ by $x_{t_i}=\zeta_i^2$.
Let
$$
f_{\zeta_i} \colon U_{\zeta_i} 
\times\mathcal{M}_{\boldsymbol{t}_{0},\boldsymbol{t}_{\text{ra}}}
\longrightarrow U_{t_i}
\times \mathcal{M}_{\boldsymbol{t}_{0},\boldsymbol{t}_{\text{ra}}}
$$
be the map induced by 
$\mathrm{Spec}\, \mathbb{C}[\zeta_i] \rightarrow
\mathrm{Spec}\, \mathbb{C}[x_{t_i}]$.
We consider the pull-back
$$
 (f_{\zeta_i}^* E_{n-2}
 |_{U_{t_i}\times \mathcal{M}_{\boldsymbol{t}_{0},\boldsymbol{t}_{\text{ra}}}}, 
 f^*_{\zeta_i} \tilde\nabla^{(n-1)}_{\text{DL}} 
  |_{U_{t_i}\times \mathcal{M}_{\boldsymbol{t}_{0},\boldsymbol{t}_{\text{ra}}}}
  ).
$$
Let 
$\Omega^{(n-2)}
_{(\boldsymbol{t}_{0},\boldsymbol{\theta},\boldsymbol{\theta}_0)}(\zeta^2_i)$
be the pull-back of the connection matrix $\Omega^{(n-2)}
_{(\boldsymbol{t}_{0},\boldsymbol{\theta},\boldsymbol{\theta}_0)}|_{U_{t_i}
\times \mathcal{M}_{\boldsymbol{t}_{0},\boldsymbol{t}_{\text{ra}}}}$
under the map $ f^*_{\zeta_i}$.
We take a formal fundamental matrix solution of 
$d+\Omega^{(n-2)}
_{(\boldsymbol{t}_{0},\boldsymbol{\theta},\boldsymbol{\theta}_0)}(\zeta^2_i)
=0$ as follows.
We have the following diagonalization:
\begin{equation}\label{2021.5.2.18.41}
\begin{aligned}
& M_{\zeta_i}^{-1} d  M_{\zeta_i} +
 M_{\zeta_i}^{-1}
\begin{pmatrix}
\alpha_i & \beta_i \\
x_{t_i}\beta_i & \alpha_i -\frac{dx_{t_i}}{2x_{t_i}}
\end{pmatrix} M_{\zeta_i} \\
&=\sum_{l=0,1,\ldots }
\begin{pmatrix}
\frac{\theta_{2l,t_i}d\zeta_i}{\zeta_i^{2(n_i-l) -1}} & 0 \\
0 & \frac{\theta_{2l,t_i}d\zeta_i}{\zeta_i^{2(n_i-l) -1}} 
\end{pmatrix}+
\sum_{l=0,1,\ldots}
\begin{pmatrix}
\frac{\theta_{2l+1,t_i}d\zeta_i}{\zeta_i^{2(n_i-l) -2}} & 0 \\
0 &- \frac{\theta_{2l+1,t_i}d\zeta_i}{\zeta_i^{2(n_i-l) -2}}  
\end{pmatrix}.
\end{aligned}
\end{equation}
We set 
$\hat{\lambda}_{i,\pm}(\zeta_i):=
\sum_{l'=0}^{\infty}(\pm 1)^{l'} \theta_{l',t_i} \int \zeta_i^{-2n_i+l'+1} d\zeta_i$,
\begin{equation}\label{2020.1.10.18.22}
\begin{aligned}
\Lambda_{i} (\zeta_{i}) :=& \begin{pmatrix}
 \hat{\lambda}_{i,+}(\zeta_i) & 0 \\
0 &\hat{\lambda}_{i,-}(\zeta_i)
\end{pmatrix},\text{ and } \\
\psi_{\zeta_i}:=&\  \Phi_i \Xi_i(\zeta_{i}^2) M_{\zeta_i} \mathrm{exp}\left(-
\Lambda_{i} (\zeta_{i})
 \right).
\end{aligned}
\end{equation}
Then $\psi_{\zeta_i}$ is a formal matrix solution of 
$d+ 
\Omega^{(n-2)}
_{(\boldsymbol{t}_{0},\boldsymbol{\theta},\boldsymbol{\theta}_0)}(\zeta^2_i) =0$,
that is, $(d+ 
\Omega^{(n-2)}
_{(\boldsymbol{t}_{0},\boldsymbol{\theta},\boldsymbol{\theta}_0)}(\zeta^2_i) 
)\psi_{\zeta_i}=0$.

For $i \in I_{\text{ra}}$,
we take a tuple $(\Phi_i , \Xi_i(x_{t_i}))$
of a compatible framing and a formal transformation as above.
We may give another formal fundamental matrix solution $\psi'_{\zeta_i}$ as follows.
If we set
$$ 
\tilde{C}_{t_i} (x_{t_i}) := 
\begin{pmatrix}
c_{t_i,\text{odd}} (x_{t_i}) & c_{t_i,\text{even}} (x_{t_i})\\
x_{t_i} c_{t_i,\text{even}} (x_{t_i})  & c_{t_i,\text{odd}} (x_{t_i})
\end{pmatrix}
=
\begin{pmatrix}
c_{t_i,1} & c_{t_i,2}\\
 0 & c_{t_i,1}
\end{pmatrix}
+\begin{pmatrix}
c_{t_i,3} & c_{t_i,4}\\
c_{t_i,2} & c_{t_i,3}
\end{pmatrix} x_{t_i} + \cdots  ,
$$ 
then 
$$
\begin{aligned}
&(\Phi_i \Xi_i(\zeta_{i}^2) \tilde{C}_{t_i} (x_{t_i}) M_{\zeta_i})^{-1}
d\Big( 
(\Phi_i \Xi_i(\zeta_{i}^2) \tilde{C}_{t_i} (x_{t_i}) M_{\zeta_i})\Big)\\
&\qquad +(\Phi_i \Xi_i(\zeta_{i}^2) \tilde{C}_{t_i} (x_{t_i}) M_{\zeta_i})^{-1}
\Omega^{(n-2)}
_{(\boldsymbol{t}_{0},\boldsymbol{\theta},\boldsymbol{\theta}_0)}
(\Phi_i \Xi_i(\zeta_{i}^2) \tilde{C}_{t_i} (x_{t_i}) M_{\zeta_i})
\end{aligned}
$$
is also diagonal, since $M_{\zeta_i}^{-1}\tilde{C}_{t_i} (\zeta_i^2) M_{\zeta_i}$ is diagonal.
Let $\tilde{c}_{t_i,11}(\zeta_i)$
and $\tilde{c}_{t_i,22}(\zeta_i)$
be the formal power series such that 
$$
M_{\zeta_i}^{-1}\tilde{C}_{t_i} (\zeta_i^2) M_{\zeta_i}  =  
\begin{pmatrix} 
\tilde{c}_{t_i,11}(\zeta_i)  & 0 \\
0 & \tilde{c}_{t_i,22} (\zeta_i)
\end{pmatrix}.
$$
We define $\Lambda'_{i'} (\zeta_{i'})$ by 
$$
\Lambda'_{i'} (\zeta_{i'})
= \Lambda_{i'} (\zeta_{i'}) + 
\begin{pmatrix}
 \int \tilde{c}_{t_i,11}(\zeta_i)^{-1} d(\tilde{c}_{t_i,11}(\zeta_i))  & 0 \\
0 &  \int \tilde{c}_{t_i,22}(\zeta_i)^{-1} d(\tilde{c}_{t_i,22}(\zeta_i))
\end{pmatrix}.
$$
If we set 
\begin{equation}\label{2021.4.26.21.53}
\psi'_{\zeta_i}:=  \Phi_i \Xi_i(\zeta_{i}^2) \tilde{C}_{t_i} (\zeta_i^2)
 M_{\zeta_i} \mathrm{exp}\left(-
\Lambda'_{i} (\zeta_{i})
 \right),
\end{equation}
we have another formal solution 
$(d+ 
\Omega^{(n-2)}
_{(\boldsymbol{t}_{0},\boldsymbol{\theta},\boldsymbol{\theta}_0)}(\zeta^2_i) 
)\psi'_{\zeta_i}=0$.
There exists a diagonal matrix $\tilde{C}_{t_i}$ such that
$\psi'_{\zeta_i} = \psi_{\zeta_i}\tilde{C}_{t_i}$ 
and $\tilde{C}_{t_i}$ is independent of $\zeta_i$.

\subsection{Symplectic structure}\label{2020.2.8.10.32}

\begin{Def}\label{2022.4.9.10.41}
Let $\delta_1$ and $\delta_2$ be vector fields on 
$\mathcal{M}_{\boldsymbol{t}_{0},\boldsymbol{t}_{\mathrm{ra}}}\subset
\mathrm{Sym}^{(n-3)}(\mathrm{Tot}(\Omega^1_{\mathbb{P}^1}(D)))$,
which is isomorphic to the moduli space
$\mathfrak{Conn}_{(\boldsymbol{t}_0,\boldsymbol{\theta},\boldsymbol{\theta}_0)}$.
We fix a (formal) fundamental matrix solution $\psi_i$ of 
$(d+\Omega^{(n-2)}
_{(\boldsymbol{t}_{0},\boldsymbol{\theta},\boldsymbol{\theta}_0)})\psi_i=0$ at $x=t_i$
for $i\in I \setminus I_{\mathrm{ra}}$
as in Lemma \ref{2019.12.30.22.09}.
We fix a fundamental matrix solution $\psi_{q_j}$ of 
$(d+\Omega^{(n-2)}
_{(\boldsymbol{t}_{0},\boldsymbol{\theta},\boldsymbol{\theta}_0)})\psi_{q_j}=0$ at $x=q_j$
as in Lemma \ref{2019.12.30.21.53}.
Moreover we take $\psi_{\zeta_i}$ defined in (\ref{2020.1.10.18.22}).
We define a $2$-form $\omega$ on
$\mathcal{M}_{\boldsymbol{t}_{0},\boldsymbol{t}_{\mathrm{ra}}}$ as
\begin{equation}\label{2021.4.26.21.57}
\begin{aligned}
\omega (\delta_1, \delta_2) :=\ &  
 \frac{1}{2} \sum_{i \in I\setminus I_{\mathrm{ra}}}   \mathrm{res}_{x=t_i} 
\mathrm{Tr} \left( 
\delta (\Omega^{(n-2)}
_{(\boldsymbol{t}_{0},\boldsymbol{\theta},\boldsymbol{\theta}_0)}) 
\wedge\delta(\psi_i)\psi_i^{-1}    \right)
 +\frac{1}{4} \sum_{i \in I_{\mathrm{ra}}}   \mathrm{res}_{\zeta_i=0} 
\mathrm{Tr} \left( \delta 
(\Omega^{(n-2)}
_{(\boldsymbol{t}_{0},\boldsymbol{\theta},\boldsymbol{\theta}_0)}(\zeta_i^2) )
\wedge\delta(\psi_{\zeta_i})\psi_{\zeta_i}^{-1}   \right)\\
&\qquad +\frac{1}{2} \sum_{j=1}^{n-3} \mathrm{res}_{x=q_j} 
\mathrm{Tr} \left( \delta (
\Omega^{(n-2)}
_{(\boldsymbol{t}_{0},\boldsymbol{\theta},\boldsymbol{\theta}_0)}) 
\wedge \delta(\psi_{q_j})\psi_{q_j}^{-1} 
  \right).
\end{aligned}
\end{equation}
\end{Def}

As in Section \ref{2020.2.8.10.24},
we may check that
the residue of 
$\delta 
(\Omega^{(n-2)}
_{(\boldsymbol{t}_{0},\boldsymbol{\theta},\boldsymbol{\theta}_0)}(\zeta_i^2) )
\wedge\delta(\psi_{\zeta_i})\psi_{\zeta_i}^{-1} $ at $\zeta_i=0$
is well-defined.
Moreover, we may also check that
 the right hand side of \eqref{2021.4.26.21.57} is independent of the choice of 
$\psi_{q_j}$ and a formal solution $\psi_{i}$ ($i \in I\setminus I_{\text{ra}}$).
Let $\psi'_{\zeta_i}$ be another fundamental matrix solution \eqref{2021.4.26.21.53}.
There exists a diagonal matrix $\tilde{C}_{t_i}$ such that
$\psi'_{\zeta_i} = \psi_{\zeta_i}\tilde{C}_{t_i}$ 
and $\tilde{C}_{t_i}$ is independent of $\zeta_i$.
By the same argument to check the independency of the choice of $\psi_{i}$,
we may check that the residue of $\mathrm{Tr} \left( \delta 
(\Omega^{(n-2)}
_{(\boldsymbol{t}_{0},\boldsymbol{\theta},\boldsymbol{\theta}_0)}(\zeta_i^2) )
\wedge\delta(\psi_{\zeta_i})\psi_{\zeta_i}^{-1}   \right)$ at $\zeta_i=0$
is equal to 
$\mathrm{Tr} \left( \delta 
(\Omega^{(n-2)}
_{(\boldsymbol{t}_{0},\boldsymbol{\theta},\boldsymbol{\theta}_0)}(\zeta_i^2) )
\wedge\delta(\psi'_{\zeta_i})(\psi'_{\zeta_i})^{-1}   \right)$ at $\zeta_i=0$.

\begin{Lem}\label{2021.4.26.22.22}
{\it 
For any vector field $\delta$ on 
$\mathcal{M}_{\boldsymbol{t}_{0},\boldsymbol{t}_{\mathrm{ra}}}$,
the formal series $\delta(\psi_{\zeta_i})\psi_{\zeta_i}^{-1} $
descends under the ramification $x_{t_i}= \zeta_i^2$.} 
\end{Lem}

\begin{proof}
We consider $\psi_{\zeta_i}M_{\zeta_i}^{-1}$. 
We decompose $\Lambda_{i} (\zeta_{i})$
into the odd degree part, the even degree part, and the logarithmic term:
$$
\Lambda_{i} (\zeta_{i}) = 
\begin{pmatrix}
 \hat{\lambda}_{\text{odd}}(\zeta_i) & 0 \\
0 &-\hat{\lambda}_{\text{odd}}(\zeta_i)
\end{pmatrix}
+
\begin{pmatrix}
 \hat{\lambda}_{\text{even}}(\zeta_i) & 0 \\
0 &\hat{\lambda}_{\text{even}}(\zeta_i)
\end{pmatrix}
+
\begin{pmatrix}
\theta_{2n_i-2,t_i} \log(\zeta_i) & 0 \\
0 &\theta_{2n_i-2,t_i} \log(\zeta_i)
\end{pmatrix}.
$$
Since 
$$ 
\begin{aligned}
&M_{\zeta_i}  \Lambda_{i} (\zeta_{i}) M_{\zeta_i}^{-1}\\
&= 
\begin{pmatrix}
0 &  \hat{\lambda}_{\text{odd}}(\zeta_i)/\zeta_i \\
\zeta_i \hat{\lambda}_{\text{odd}}(\zeta_i)& 0
\end{pmatrix}
+\begin{pmatrix}
 \hat{\lambda}_{\text{even}}(\zeta_i) & 0 \\
0 &\hat{\lambda}_{\text{even}}(\zeta_i)
\end{pmatrix}
+\frac{1}{2}\begin{pmatrix}
\theta_{2n_i-2,t_i} \log(\zeta_i^2) & 0 \\
0 &\theta_{2n_i-2,t_i} \log(\zeta_i^2)
\end{pmatrix},
\end{aligned}
$$
we may check that
$M_{\zeta_i}  \Lambda_{i} (\zeta_{i}) M_{\zeta_i}^{-1}$ 
descends under the ramification $x_{t_i}= \zeta_i^2$.
Since
$$
\psi_{\zeta_i} M_{\zeta_i}^{-1} =  \Phi_i \Xi_i(\zeta_{i}^2) M_{\zeta_i} \mathrm{exp}\left(-
\Lambda_{i'} (\zeta_{i'}) \right)M_{\zeta_i}^{-1}
=\Phi_i \Xi_i(\zeta_{i}^2)  \mathrm{exp}\left(-M_{\zeta_i}
\Lambda_{i'} (\zeta_{i'}) M_{\zeta_i}^{-1} \right),
$$
we have that $\psi_{\zeta_i}M_{\zeta_i}^{-1}$ 
descends under the ramification $x_{t_i}= \zeta_i^2$.
Since $M_{\zeta_i}$ is independent of the parameters of 
$\mathcal{M}_{\boldsymbol{t}_{0},\boldsymbol{t}_{\mathrm{ra}}}$,
we have $\delta(M_{\zeta_i})=0$.
Then we have 
$\delta (\psi_{\zeta_i}M_{\zeta_i}^{-1}) (\psi_{\zeta_i}M_{\zeta_i}^{-1})^{-1}
=\delta(\psi_{\zeta_i})\psi_{\zeta_i}^{-1}$.
Finally we have that $\delta(\psi_{\zeta_i})\psi_{\zeta_i}^{-1}$
descends under the ramification $x_{t_i}= \zeta_i^2$.
\end{proof}

\begin{Rem}
In the proof of this lemma, we check that 
 $\psi_{\zeta_i}M_{\zeta_i}^{-1}$
descends under the ramification $x_{t_i}= \zeta_i^2$.
If we set $\psi_{i} := \psi_{\zeta_i}M_{\zeta_i}^{-1}$ for $i \in I_{\mathrm{ra}}$,
we have that 
$$
\begin{aligned}
&\frac{1}{4}  \cdot \mathrm{res}_{\zeta_i=0} 
\mathrm{Tr} \left( \delta 
(\Omega^{(n-2)}
_{(\boldsymbol{t}_{0},\boldsymbol{\theta},\boldsymbol{\theta}_0)}(\zeta_i^2) )
\wedge\delta(\psi_{\zeta_i})\psi_{\zeta_i}^{-1}   \right)
=\frac{1}{2} \cdot  \mathrm{res}_{x=t_i} 
\mathrm{Tr} \left( \delta 
(\Omega^{(n-2)}
_{(\boldsymbol{t}_{0},\boldsymbol{\theta},\boldsymbol{\theta}_0)} )
\wedge\delta(\psi_{i})\psi_{i}^{-1}   \right).
\end{aligned}
$$
So we may define the $2$-form $\omega$ as 
\begin{equation*}
\begin{aligned}
\omega (\delta_1, \delta_2) :=\ &  
 \frac{1}{2} \sum_{i \in I}   \mathrm{res}_{x=t_i} 
\mathrm{Tr} \left( 
\delta (\Omega^{(n-2)}
_{(\boldsymbol{t}_{0},\boldsymbol{\theta},\boldsymbol{\theta}_0)}) 
\wedge\delta(\psi_i)\psi_i^{-1}    \right)
+\frac{1}{2} \sum_{j=1}^{n-3} \mathrm{res}_{x=q_j} 
\mathrm{Tr} \left( \delta (
\Omega^{(n-2)}
_{(\boldsymbol{t}_{0},\boldsymbol{\theta},\boldsymbol{\theta}_0)}) 
\wedge \delta(\psi_{q_j})\psi_{q_j}^{-1} 
  \right).
\end{aligned}
\end{equation*}
In this definition of $\omega$,
the variable $\zeta_i$ disappears.
\end{Rem}

\begin{Thm}
\textit{
Let $\omega$ be the $2$-form
on $\mathcal{M}_{\boldsymbol{t}_{0},\boldsymbol{t}_{\text{ra}}}$ 
defined by \eqref{2021.4.26.21.57} in Definition \ref{2022.4.9.10.41}.
The $2$-form $\omega$ coincides with
\begin{equation*}
\sum_{j=1}^{n-3} d \left(\frac{p_j}{P(q_j)} \right) \wedge dq_j.
\end{equation*}
}
\end{Thm}

\begin{proof}
Let $\delta$ be a vector field on 
$\mathcal{M}_{\boldsymbol{t}_{0},\boldsymbol{t}_{\text{ra}}}$.
Since $t_i$ ($i \in I$), 
$\theta_{l,t_i}^{\pm}$ ($i \in I_{\mathrm{un}}, 0\leq l \leq n_i -1$), 
and $\theta_{l',t_i}$ ($i \in I_{\mathrm{ra}}, 0\leq l' \leq 2n_i -2$)
are constants on $\mathcal{M}_{\boldsymbol{t}_{0},\boldsymbol{t}_{\text{ra}}}$,
we have $\delta(t_i)=0$, $\delta(\theta_{l,t_i}^{\pm})=0$, 
and $\delta(\theta_{l',t_i})=0$.
By the equalities \eqref{2021.4.30.13.59} and \eqref{2021.4.30.13.59.2},
we have $\delta(C_i)=\delta(D_i)=0$ 
for $i=1,2,\ldots,\nu, \infty$.
Here $C_i$ and $D_i$ are the polynomials in \eqref{2020.12.26.21.03}.
We compute the residue of the trace of
$ \delta 
(\Omega^{(n-2)}
_{(\boldsymbol{t}_{0},\boldsymbol{\theta},\boldsymbol{\theta}_0)}(\zeta_i^2)) 
\wedge \delta(\psi_{\zeta_i})\psi_{\zeta_i}^{-1} $ 
at $\zeta_i=0$. 
First, we consider the expansion of $\delta
(\Omega^{(n-2)}
_{(\boldsymbol{t}_{0},\boldsymbol{\theta},\boldsymbol{\theta}_0)}(\zeta_i^2)) $ 
at $\zeta_i=0$.
Since $\delta(C_i)=\delta(D_i)=0$ 
for $i=1,2,\ldots,\nu, \infty$, we have
$\delta(c_2)= O(x_{t_i}^0)$ and 
$\delta(d_2)= O(x_{t_i}^0)$.
Second, we consider $\delta(\psi_{\zeta_i})\psi_{\zeta_i}^{-1} $. 
By the definition (\ref{2020.1.10.18.22}), we have that
$\delta(\psi_{\zeta_i})\psi_{\zeta_i}^{-1} $ coincides with
\begin{equation*}
\begin{aligned}
&\Phi_{i} \Xi_{i}(\zeta_i^2) \delta(M_{\zeta_i})M_{\zeta_i}^{-1} (\Phi_{i} \Xi_{i}(\zeta_i^2))^{-1}+
\delta(\Phi_{i} \Xi_{i}(\zeta_i^2) )(\Phi_{i} \Xi_{i}(\zeta_i^2))^{-1} \\
&+  (\Phi_{i} \Xi_{i}(\zeta_i^2) M_{\zeta_i}) 
\begin{pmatrix}
-\delta (\hat{\lambda}_{i,+}(\zeta_i)) & 0 \\
0 &-\delta (\hat{\lambda}_{i,-}(\zeta_i))
\end{pmatrix}
 (\Phi_{i} \Xi_{i}(\zeta_i^2)M_{\zeta_i})^{-1}.
\end{aligned}
\end{equation*}
Since $\delta(M_{\zeta_i})=0$ and $\delta (\hat{\lambda}_{i,\pm}(\zeta_i))=O(\zeta_i)$,
we have 
\begin{equation*}
\mathrm{res}_{\zeta_i=0} 
\mathrm{Tr} \left( \delta 
(\Omega^{(n-2)}
_{(\boldsymbol{t}_{0},\boldsymbol{\theta},\boldsymbol{\theta}_0)}
(\zeta^2_i))
\wedge \delta(\psi_{\zeta_i})\psi_{\zeta_i}^{-1}   \right)=0.
\end{equation*}
 The remained residues are calculated as in the proof of Theorem \ref{2020.1.21.16.17}.
Then we obtain
\begin{equation*}
\omega (\delta_1, \delta_2) = 
\sum_{j=1}^{n-3}
\left( \frac{\delta_1(p_j)\delta_2(q_j)}{P(q_j)}
-\frac{\delta_2(p_j)\delta_1(q_j)}{P(q_j)} \right),
\end{equation*}
which means that 
$\omega$ 
coincides with $\sum_{j=1}^{n-3} d\left( \frac{p_j}{P(q_j)} \right) \wedge dq_j$.
\end{proof}

\subsection{Integrable deformations 
associated to $\theta_{l',t_i}$ for $i\in I_{\mathrm{ra}}$}\label{2020.2.7.11.34}

First we fix $i \in I_{\text{ra}}$ and $l' \in \{ 0,1,2,\ldots, 2n_i -3\}$.
Let 
$$
\tilde{\nabla}_{\mathrm{DL, ext}}^{(1)}=
\begin{cases}
d+ \widehat\Omega_{(\boldsymbol{t}_{\mathrm{ra}},\boldsymbol{\theta}_0)}^{(1)} 
& \text{ on $U_0\times(\widehat{\mathcal{M}}_{\boldsymbol{t}_{\text{ra}}} \times
 T_{\boldsymbol{\theta}})$} \\
d + G_{1}^{-1} dG_{1}+G_{1}^{-1} 
\widehat\Omega_{(\boldsymbol{t}_{\mathrm{ra}},\boldsymbol{\theta}_0)}^{(1)} \, G_{1}
 & \text{ on $U_\infty \times (\widehat{\mathcal{M}}_{\boldsymbol{t}_{\text{ra}}} 
 \times T_{\boldsymbol{\theta}})$}
\end{cases}
$$
be the family \eqref{2021.4.14.10.36}.
Let $\theta_{l',t_i}$ be the natural coordinate of 
$(T_{\boldsymbol{t}})_{\boldsymbol{t}_{\text{ra}}} \times T_{\boldsymbol{\theta}}$
and $\partial/\partial \theta_{l',t_i}$ be the vector field on
$(T_{\boldsymbol{t}})_{\boldsymbol{t}_{\text{ra}}} \times T_{\boldsymbol{\theta}}$
associated to $\theta_{l',t_i}$.
We will construct a horizontal lift of $\tilde{\nabla}_{\mathrm{DL, ext}}^{(1)}$
with respect to $\partial/\partial \theta_{l',t_i}$.

By using the explicit form of $\tilde{\nabla}_{\mathrm{DL, ext}}^{(n-2)}$,
we take a family of compatible framings of 
$\tilde{\nabla}_{\mathrm{DL, ext}}^{(n-2)}$
at $\tilde{t}_{i'}$ for $i' \in I_{\text{ra}}$.
We denote by $\Phi_{i'}$ 
this family of compatible framings at $\tilde{t}_{i'}$ for $i' \in I_{\text{ra}}$.
Let $\Xi_{i'}(x_{t_{i'}})$ be the formal transformation \eqref{2021.4.23.18.23} for $i' \in I_{\text{ra}}$.
If $i' \in I \setminus I_{\text{ra}}$,  
let $\Phi_{i'}$ and $\Xi_{i'}(x_{t_{i'}})$ be a compatible framing at $\tilde{t}_{i'}$
and the (formal) transformation with respect to $\Phi_{i'}$
appeared in Lemma \ref{2019.12.30.22.09}, respectively.
Let $\tilde{G}$ be the matrix defined in \eqref{2021.4.2.20.06}.
For each $i' \in I$,
we denote 
formal expansion of $\tilde{G}^{-1}\Phi_{i'}\Xi_{i'}(x_{t_{i'}})$ at $x_{t_{i'}}=0$ by\
\begin{equation*}
\tilde{G}^{-1}\Phi_{i'}\Xi_{i'}(x_{t_{i'}})
= P_{i',0} + P_{i',1} x_{t_{i'}} + P_{i',2} x^2_{t_{i'}} +\cdots.
\end{equation*}
Set 
\begin{equation*}
\begin{aligned}
P_{i'} :=&\  P_{i',0} + P_{i',1} x_{t_{i'}} + P_{i',2} x^2_{t_{i'}} +\cdots +
 P_{i',2n_{i'}-1} x^{2n_{i'}-1}_{t_{i'}}
\ \  \text{(for $i'\in I$)}, \ \   \text{and} \\ 
P_{\nu+1} := &\ \mathrm{id}.
\end{aligned}
\end{equation*} 
We take an affine open covering $\{ \hat{U}_{i'} \}_{i' \in I\cup \{\nu+1\}}$
of $\mathbb{P}^1\times 
(\widehat{\mathcal{M}}_{\boldsymbol{t}_{\text{ra}}} \times T_{\boldsymbol{\theta}})$
as in Section \ref{2020.1.24.22.25}.
By using the matrices $P_{i'}$,
we define a new trivialization $\hat\varphi_{i'}$ of $\hat{E}_1$ on 
 $\hat{U}_{i'}$ for each $i' \in I\cup \{\nu+1\}$
as in Section \ref{2020.1.24.22.25}.
Let $\hat{\Omega}_{i'}$ be
the connection matrix 
of $ \tilde{\nabla}_{\mathrm{DL, ext}}^{(1)}$
under the new trivialization $\hat\varphi_{i'}$:
\begin{equation*}
\begin{aligned}
\hat{\Omega}_{i'}&=
P_{i'}^{-1} dP_{i'}+P_{i'}^{-1} 
\widehat\Omega_{(\boldsymbol{t}_{\mathrm{ra}},\boldsymbol{\theta}_0)}^{(1)} |_{\hat{U}_{i'}}
P_{i'} \quad \text{ for $i' \in (I\cup \{\nu+1\}) \setminus \{ \infty\} $, and} \\
\hat{\Omega}_{\infty}&=
(G_1P_{\infty})^{-1} d(G_1 P_{\infty})+(G_1P_{\infty})^{-1} 
\widehat\Omega_{(\boldsymbol{t}_{\mathrm{ra}},\boldsymbol{\theta}_0)}^{(1)} |_{\hat{U}_{\infty}}
(G_1P_{\infty} ).
\end{aligned}
\end{equation*}
Remark that
$M_{\zeta_{i'}}^{-1} d M_{\zeta_{i'}} + M_{\zeta_{i'}}^{-1}\hat{\Omega}_{i'}M_{\zeta_{i'}}$ is diagonal 
until the $\zeta_{i'}^{2n_{i'}-3}$-terms for each $i' \in I_{\text{ra}}$.

For the fixed indices $i \in I_{\mathrm{ra}}$ and $l' \in \{ 0,1,\ldots,2n_i -3\}$, 
we define a matrix $B_{\theta_{l',t_i}}$
by
\begin{equation*}
B_{\theta_{l',t_i}}=
\begin{cases} 
\frac{-1}{2(n_i-l-1)}\cdot \frac{\delta(\theta_{2l,t_i})}{\zeta^{2(n_i-1-l)}_i}
\begin{pmatrix}
1& 0 \\
0 &  1
\end{pmatrix} & \text{ if $l'=2l$} \\
\frac{-1}{2(n_i-l)-3}\cdot \frac{\delta(\theta_{2l+1,t_i})}{\zeta^{2(n_i-1 -l) -1}_i}
\begin{pmatrix}
1& 0 \\
0 & - 1
\end{pmatrix} & \text{ if $l'=2l+1$}.
\end{cases}
\end{equation*}
We may check that
\begin{equation*}
M_{\zeta_i}B_{\theta_{l',t_i}}M_{\zeta_i}^{-1}=
\begin{cases} 
\frac{-1}{2(n_i-l-1)}\cdot \delta(\theta_{2l,t_i})
\begin{pmatrix}
x_{t_i}^{-n_i+l+1}& 0 \\
0 &  x_{t_i}^{-n_i+l+1}
\end{pmatrix} & \text{ if $l'=2l$} \\
\frac{-1}{2(n_i-l)-3}\cdot  \delta(\theta_{2l+1,t_i}) 
\begin{pmatrix}
0& x_{t_i}^{-n_i +l+1  } \\
x_{t_i}^{-n_i +l+2 } & 0
\end{pmatrix} & \text{ if $l'=2l+1$}.
\end{cases}
\end{equation*}
In particular, $M_{\zeta_i}B_{\theta_{l',t_i}}M_{\zeta_i}^{-1}$
descends under the ramification $x_{t_i}=\zeta_i^2$.
We define
a vector bundle $(\widehat{E}_1)_{\theta_{l',t_i}}^{\epsilon}$ on 
$\mathbb{P}^1\times(\widehat{\mathcal{M}}_{\boldsymbol{t}_{\text{ra}}} \times
 T_{\boldsymbol{\theta} } ) \times \mathrm{Spec}\, \mathbb{C}[\epsilon]$
by the same argument as in the construction of $(\widehat{E}_1)_{\theta_{l,t_i}^{\pm}}^{\epsilon}$.
That is, we replace $B_{\theta_{l,t_i}^{\pm}}(x_{t_i})$ in \eqref{2021.4.16.15.17}
with $M_{\zeta_i}B_{\theta_{l',t_i}}M_{\zeta_i}^{-1}$.
We define a morphism
$$
\nabla^{\epsilon}_{\partial/\partial \theta_{l',t_i}}  \colon 
(\widehat{E}_1)_{\theta_{l',t_i}}^{\epsilon}
\longrightarrow
(\widehat{E}_1)_{\theta_{l',t_i}}^{\epsilon} \otimes
\tilde\Omega^1_{\partial/\partial t_i}
$$
by the same argument as in the construction of $\nabla^{\epsilon}_{\partial/\partial \theta_{l,t_i}^{\pm}} $
in Section \ref{2020.1.24.22.25}.
That is, we replace $B_{\theta_{l,t_i}^{\pm}}(x_{t_i})$ in \eqref{2021.4.16.15.34}
with $M_{\zeta_i}B_{\theta_{l',t_i}}M_{\zeta_i}^{-1}$.
The $\epsilon$-term of $\nabla^{\epsilon}_{\partial/\partial \theta_{l',t_i}}|_{\hat{U}_{i}^{\epsilon}}$ 
for fixed $i \in I_{\text{ra}}$
is
\begin{equation*}
\begin{aligned}
&d (M_{\zeta_i}B_{\theta_{l',t_i}} M_{\zeta_i}^{-1} )
+ [ \hat{\Omega}_{i}
, M_{\zeta_i}B_{\theta_{l',t_i}} M_{\zeta_i}^{-1} ] \\
&=
M_{\zeta_i}\left( \frac{\partial}{\partial \zeta_i} (B_{\theta_{l',t_i}} )d\zeta_i
+ [M_{\zeta_i}^{-1} d M_{\zeta_i}+ M_{\zeta_i}^{-1}\hat{\Omega}_{i}M_{\zeta_i}
, B_{\theta_{l',t_i}}] \right)M_{\zeta_i}^{-1} \\
&=
\begin{cases}
 \frac{\delta(\theta_{2l,t_i})}{2}
\begin{pmatrix}
x_{t_i}^{-n_i+l}& 0 \\
0 &  x_{t_i}^{-n_i+l}
\end{pmatrix} dx_{t_i} 
+M_{\zeta_i}[M_{\zeta_i}^{-1} d M_{\zeta_i}+ M_{\zeta_i}^{-1}\hat{\Omega}_{i}M_{\zeta_i}
, B_{\theta_{l',t_i}}] M_{\zeta_i}^{-1}
& \text{ if $l'=2l$} \\
\frac{\delta(\theta_{2l+1,t_i}) }{2}
\begin{pmatrix}
0& x_{t_i}^{-n_i +l  } \\
x_{t_i}^{-n_i +l+1 } & 0
\end{pmatrix}dx_{t_i}
+M_{\zeta_i}[M_{\zeta_i}^{-1} d M_{\zeta_i}+ M_{\zeta_i}^{-1}\hat{\Omega}_{i}M_{\zeta_i}
, B_{\theta_{l',t_i}}] M_{\zeta_i}^{-1} & \text{ if $l'=2l+1$}.
\end{cases}
\end{aligned}
\end{equation*}
Since $M_{\zeta_i}^{-1} d M_{\zeta_i}+ M_{\zeta_i}^{-1}\hat{\Omega}_{i}M_{\zeta_i}$
 and $B_{\theta_{l',t_i}}$ are diagonal 
until the $\zeta_i^{2n_i-3}$-terms, 
the negative parts of the relative connection
$\overline{\nabla}^{\epsilon}_{\partial/\partial \theta_{l',t_i}}$
 along the divisor $[x_{t_i}=0]$ is
\begin{equation*}
\begin{pmatrix}
\alpha_i & \beta_i \\
x_{t_i}\beta_i & \alpha_i -\frac{dx_{t_i}}{2x_{t_i}}
\end{pmatrix}.
\end{equation*}
Here if $l'=2l$, the entry $\alpha_i$ is the following
\begin{equation*}
\frac{\theta_{0,t_i}}{2}\frac{dx_{t_i}}{x_{t_i}^{n_i}}+ \cdots 
+\frac{\theta_{2l,t_i} + \epsilon \delta (\theta_{2l,t_i})}{2} \frac{dx_{t_i}}{x_{t_i}^{n_i-l}}+\cdots 
+\frac{\theta_{2n_i-2,t_i}}{2} \frac{dx_{t_i}}{x_{t_i}} 
\end{equation*}
and the coefficients of $\beta_i$ are independent of $\epsilon$ until the $x_{t_i}^{-2}$-term.
If $l'= 2l+1$, the expansion of $\beta_i$ at $t_i$ until the $x_{t_i}^{-2}$-term is the following
\begin{equation*}
 \frac{\theta_{1,t_i}}{2}\frac{dx_{t_i}}{x_{t_i}^{n_i}}+ \cdots 
+ \frac{\theta_{2l+1,t_i} + \epsilon \delta(\theta_{2l+1,t_i})}{2}\frac{dx_{t_i}}{x_{t_i}^{n_i-l}}+\cdots 
+\frac{\theta_{2n_i-3,t_i}}{2}\frac{dx_{t_i}}{x^2_{t_i}} 
\end{equation*}
and the coefficients of $\alpha_i$ are independent of $\epsilon$ until the $x_{t_i}^{-1}$-term.

As in Section \ref{2020.1.24.22.25},
$\widehat{E}_1^{\epsilon}  \cong (\widehat{E}_1)_{\theta_{l',t_i}}^{\epsilon}$.
If we consider the pull-back of $\nabla^{\epsilon}_{\partial/\partial \theta_{l',t_i}}$
under this isomorphism, 
then we have a horizontal lift of $\tilde{\nabla}_{\mathrm{DL, ext}}^{(1)}$
with respect to $\partial/\partial \theta_{l',t_i}$.
If we take a relativization of this horizontal lift,
we have a family of connections parametrized by
$( \widehat{\mathcal{M}}_{\boldsymbol{t}_{\text{ra}}}\times T_{\boldsymbol{\theta}})
\times \mathrm{Spec}\, \mathbb{C}[\epsilon]$.
This family gives a map from the base space
$( \widehat{\mathcal{M}}_{\boldsymbol{t}_{\text{ra}}}\times T_{\boldsymbol{\theta}})
\times \mathrm{Spec}\, \mathbb{C}[\epsilon]$
to the moduli space $\widehat{\mathfrak{Conn}}_{(\boldsymbol{t}_{\mathrm{ra}},\boldsymbol{\theta}_0)}$.
By taking composition with $\widehat{\mathrm{App}}$ defined in \eqref{2020.1.6.13.00},
we have a map 
\begin{equation}\label{2021.4.24.10.40}
( \widehat{\mathcal{M}}_{\boldsymbol{t}_{\text{ra}}}\times T_{\boldsymbol{\theta}})
\times \mathrm{Spec}\, \mathbb{C}[\epsilon]
\longrightarrow
\widehat{\mathcal{M}}_{\boldsymbol{t}_{\text{ra}}} \times T_{\boldsymbol{\theta}}.
\end{equation}

\begin{Def}
Then 
we may define the vector field 
on $\widehat{\mathcal{M}}_{\boldsymbol{t}_{\text{ra}}}\times T_{\boldsymbol{\theta}}$
associated 
to the map \eqref{2021.4.24.10.40}.
We denote by $\delta_{\theta_{l',t_i}}^{\mathrm{IMD}}$ this vector field 
on $\widehat{\mathcal{M}}_{\boldsymbol{t}_{\text{ra}}}\times T_{\boldsymbol{\theta}}$.
\end{Def}

Let $f_{\theta_{l',t_i}}^{\mathrm{IMD}} \colon 
( \widehat{\mathcal{M}}_{\boldsymbol{t}_{\text{ra}}}\times T_{\boldsymbol{\theta}})
\times \mathrm{Spec}\, \mathbb{C}[\epsilon]
\longrightarrow
\widehat{\mathcal{M}}_{\boldsymbol{t}_{\text{ra}}} \times T_{\boldsymbol{\theta}}$
be the map induced by the vector field
$\delta_{\theta_{l',t_i}}^{\mathrm{IMD}}$.
We have $\widehat{E}^{\epsilon}_1=
(\mathrm{id} \times f_{\theta_{l',t_i}}^{\mathrm{IMD}})^* \widehat{E}_1$.
We denote by
\begin{equation}\label{2021.4.16.12.11.ra}
\begin{cases}
d+ \widehat\Omega_{(\boldsymbol{t}_{\mathrm{ra}},\boldsymbol{\theta}_0)}^{(1)} 
+\epsilon \delta_{\theta_{l',t_i}}^{\mathrm{IMD}}
(\widehat\Omega_{(\boldsymbol{t}_{\mathrm{ra}},\boldsymbol{\theta}_0)}^{(1)} ) \\
\qquad \qquad 
\text{ on $U_0 \times(\widehat{\mathcal{M}}_{\boldsymbol{t}_{\text{ra}}}\times T_{\boldsymbol{\theta}} )
\times \mathrm{Spec}\, \mathbb{C} [\epsilon]$ } \\
d + G_{1}^{-1} dG_{1}+G_{1}^{-1} 
\widehat\Omega_{(\boldsymbol{t}_{\mathrm{ra}},\boldsymbol{\theta}_0)}^{(1)} 
 \, G_{1}
+\epsilon G_{1}^{-1} 
\delta_{\theta_{l',t_i}}^{\mathrm{IMD}}
(\widehat\Omega_{(\boldsymbol{t}_{\mathrm{ra}},\boldsymbol{\theta}_0)}^{(1)} )
G_{1}  \\
\qquad\qquad 
\text{ on $U_\infty \times(\widehat{\mathcal{M}}_{\boldsymbol{t}_{\text{ra}}}\times T_{\boldsymbol{\theta}} )
\times \mathrm{Spec}\, \mathbb{C} [\epsilon]$ }
\end{cases}
\end{equation}
the pull-back of $\tilde{\nabla}_{\mathrm{DL, ext}}^{(1)}$ 
under the map $\mathrm{id} \times f_{\theta_{l',t_i}}^{\mathrm{IMD}}$.
As in Section \ref{2020.1.24.22.25},
we have a lift of 
$(\mathrm{id} \times f_{\theta_{l',t_i}}^{\mathrm{IMD}})^*\tilde{\nabla}_{\mathrm{DL, ext}}^{(1)}$:
\begin{equation}\label{2021.4.16.12.12.ra}
\begin{cases}
\hat{d}+ \widehat\Omega_{(\boldsymbol{t}_{\mathrm{ra}},\boldsymbol{\theta}_0)}^{(1)} 
+\epsilon \delta_{\theta_{l',t_i}}^{\mathrm{IMD}}
(\widehat\Omega_{(\boldsymbol{t}_{\mathrm{ra}},\boldsymbol{\theta}_0)}^{(1)} )
+ \Upsilon_{\theta_{l',t_i}}^{\mathrm{IMD}}
 d\epsilon\\
\qquad \qquad 
\text{ on $U_0 \times(\widehat{\mathcal{M}}_{\boldsymbol{t}_{\text{ra}}}\times T_{\boldsymbol{\theta}} )
\times \mathrm{Spec}\, \mathbb{C} [\epsilon]$ } \\
\hat{d} + G_{1}^{-1} dG_{1}+G_{1}^{-1} 
\widehat\Omega_{(\boldsymbol{t}_{\mathrm{ra}},\boldsymbol{\theta}_0)}^{(1)} 
 \, G_{1}
+\epsilon G_{1}^{-1} 
\delta_{\theta_{l',t_i}}^{\mathrm{IMD}}
(\widehat\Omega_{(\boldsymbol{t}_{\mathrm{ra}},\boldsymbol{\theta}_0)}^{(1)} )
G_{1} 
+ G_1^{-1} \Upsilon_{\theta_{l',t_i}}^{\mathrm{IMD}}
 G_1 d\epsilon \\
\qquad\qquad 
\text{ on $U_\infty \times(\widehat{\mathcal{M}}_{\boldsymbol{t}_{\text{ra}}}\times T_{\boldsymbol{\theta}} )
\times \mathrm{Spec}\, \mathbb{C} [\epsilon]$ },
\end{cases}
\end{equation}
which is 
a morphism $\widehat{E}^{\epsilon}_1 \rightarrow 
\widehat{E}^{\epsilon}_1 \otimes
\tilde\Omega^1_{\partial/\partial t_i}$
with the Leibniz rule
and the following equality
\begin{equation}\label{2021.4.27.14.24}
 \delta_{\theta_{l',t_i}}^{\mathrm{IMD}}
(\widehat\Omega_{(\boldsymbol{t}_{\mathrm{ra}},\boldsymbol{\theta}_0)}^{(1)} )
 = 
d  \Upsilon_{\theta_{l',t_i}}^{\mathrm{IMD}} 
+[\widehat\Omega_{(\boldsymbol{t}_{\mathrm{ra}},\boldsymbol{\theta}_0)}^{(1)}  ,
\Upsilon_{\theta_{l',t_i}}^{\mathrm{IMD}}],
\end{equation}
which means the integrable condition.

\begin{Rem}
In this paper, we consider only rank 2 connections on $\mathbb{P}^1$.
Moreover we impose some Zariski open conditions, for example 
the underlying vector bundles are isomorphic to 
$\mathcal{O}_{\mathbb{P}^1} \oplus \mathcal{O}_{\mathbb{P}^1}(1)$.
That is, we consider only a Zariski open subset of the moduli space
of connections constructed in \cite{I1}.
Horizontal lifts for more general situations are constructed by Inaba
in \cite[Section 9]{I2}.
\end{Rem}

\subsection{Isomonodromy $2$-form}\label{2020.2.8.10.33}

\begin{Def}
Let $\hat{\delta}_1$ and $\hat{\delta}_2$ be vector fields on 
$\widehat{\mathcal{M}}_{\boldsymbol{t}_{\text{ra}}}\times T_{\boldsymbol{\theta}}$,
which is isomorphic to the extended moduli space
$\widehat{\mathfrak{Conn}}_{(\boldsymbol{t}_{\mathrm{ra}},\boldsymbol{\theta}_0)}$.
For each $i \in I_{\text{ra}} $,
we fix a formal fundamental matrix solution $\psi_{\zeta_i}$ of 
$(d+\widehat\Omega_{(\boldsymbol{t}_{\mathrm{ra}},\boldsymbol{\theta}_0)}^{(n-2)} (\zeta^2_i))
\psi_{\zeta_i}=0$ at $x=t_i$ defined in \eqref{2020.1.10.18.22}.
For each $i \in I_{\text{reg}} \cup I_{\text{un}}$,
we fix a formal fundamental matrix solution $\psi_i$ of 
$(d+\widehat\Omega_{(\boldsymbol{t}_{\mathrm{ra}},\boldsymbol{\theta}_0)}^{(n-2)})
\psi_i=0$ at $x=t_i$
as in Lemma \ref{2019.12.30.22.09}.
We take a fundamental matrix solution $\psi_{q_j}$ of 
$(d+\widehat\Omega_{(\boldsymbol{t}_{\mathrm{ra}},\boldsymbol{\theta}_0)}^{(n-2)})
\psi_{q_j}=0$ at $x=q_j$ as in Lemma \ref{2019.12.30.21.53}.
We define a $2$-form $\hat{\omega}$ 
on $\widehat{\mathcal{M}}_{\boldsymbol{t}_{\text{ra}}}\times T_{\boldsymbol{\theta}}$ as
\begin{equation*}
\begin{aligned}
\hat{\omega} (\hat{\delta}_1,\hat{\delta}_2) :=\ &  \frac{1}{2} \sum_{i\in I \setminus I_{\mathrm{ra}}} 
\mathrm{res}_{x=t_i} 
\mathrm{Tr} \left( \hat{\delta}
(\widehat\Omega_{(\boldsymbol{t}_{\mathrm{ra}},\boldsymbol{\theta}_0)}^{(n-2)} ) 
\wedge \hat{\delta}(\psi_i)\psi_i^{-1}  \right) 
+\frac{1}{4} \sum_{i\in I_{\mathrm{ra}}} 
\mathrm{res}_{\zeta_i=0} 
\mathrm{Tr} \left( \hat{\delta} 
(\widehat\Omega_{(\boldsymbol{t}_{\mathrm{ra}},\boldsymbol{\theta}_0)}^{(n-2)}(\zeta^2_i) ) 
\wedge \hat{\delta}(\psi_{\zeta_i})\psi_{\zeta_i}^{-1} \right) \\
&+\frac{1}{2} \sum_{j=1}^{n-3} \mathrm{res}_{x=q_j} 
\mathrm{Tr} \left( \hat{\delta}
 (\widehat\Omega_{(\boldsymbol{t}_{\mathrm{ra}},\boldsymbol{\theta}_0)}^{(n-2)} )
\wedge \hat{\delta}(\psi_{q_j})\psi_{q_j}^{-1}  \right).
\end{aligned}
\end{equation*}
\end{Def}

By the same argument as in Section \ref{2020.2.8.10.32},
We may check that the residue of 
$\hat{\delta} (\widehat\Omega_{(\boldsymbol{t}_{\mathrm{ra}},\boldsymbol{\theta}_0)}^{(n-2)}) 
 \wedge \hat{\delta}(\psi_i)\psi_i^{-1}$ at $\tilde{t}_i$ (for $i \in I \setminus I_{\text{ra}}$)
 and 
$ \hat{\delta} 
(\widehat\Omega_{(\boldsymbol{t}_{\mathrm{ra}},\boldsymbol{\theta}_0)}^{(n-2)}(\zeta^2_i) ) 
\wedge \hat{\delta}(\psi_{\zeta_i})\psi_{\zeta_i}^{-1}$
at $\zeta_i=0$ (for $i \in  I_{\text{ra}}$)
are well-defined.
Moreeover,
the right hand side of \eqref{2020.1.7.17.07} is independent of the choice of 
$\psi_{q_j}$ and $\psi_{i}$ ($i \in I \setminus I_{\text{ra}}$) and $\psi_{\zeta_i}$
($i \in  I_{\text{ra}}$).
By the same argument as in the proof of
Lemma \ref{2021.4.26.22.22},
we may check that 
the formal series $\hat{\delta}(\psi_{\zeta_i})\psi_{\zeta_i}^{-1} $
descends under the ramification $x_{t_i}= \zeta_i^2$
for any vector field $\hat{\delta}$ on 
$\widehat{\mathcal{M}}_{\boldsymbol{t}_{\text{ra}}}\times T_{\boldsymbol{\theta}}$.

\begin{Prop}
{\it 
Let
$\tilde{G}$ and
$\tilde{G}_{\infty}$
be the matrices defined in
\eqref{2021.4.2.20.06}.
If $\infty$ is an element of $I_{\mathrm{reg}} \cup I_{\mathrm{un}}$,
we set $\tilde\psi_i := \tilde{G}^{-1} \psi_i$ for any 
$i \in (I_{\mathrm{reg}} \cup I_{\mathrm{un}})  \setminus  \{\infty \}$,
 $\tilde\psi_{\zeta_i} := \tilde{G}^{-1} \psi_{\zeta_i}$ for any 
$i \in I_{\mathrm{ra}} $,
and
$\tilde\psi_{\infty} := \tilde{G}^{-1}_{\infty} \psi_\infty$.
If $\infty$ is an element of $I_{\mathrm{ra}} $,
we set $\tilde\psi_i := \tilde{G}^{-1} \psi_i$ for any 
$i \in I_{\mathrm{reg}} \cup I_{\mathrm{un}}$,
 $\tilde\psi_{\zeta_i} := \tilde{G}^{-1} \psi_{\zeta_i}$ for any 
$i \in I_{\mathrm{ra}} \setminus \{ \infty\} $,
and
$\tilde\psi_{\zeta_{\infty}} := \tilde{G}^{-1}_{\infty} \psi_{\zeta_{\infty}}$.
We have the following equality: 
\begin{equation}\label{eq:2020.11.19.16.08.ra}
\begin{aligned}
\hat\omega (\hat\delta_1, \hat\delta_2) &=
 \frac{1}{2} \sum_{i \in I \setminus I_{\mathrm{ra}}}   \mathrm{res}_{x=t_i} 
\mathrm{Tr} \left( \hat\delta
(\widehat\Omega_{(\boldsymbol{t}_{\mathrm{ra}},\boldsymbol{\theta}_0)}^{(1)}) 
\wedge\hat\delta(\tilde\psi_i)(\tilde\psi_i)^{-1}  \right) \\
&\qquad +\frac{1}{4} \sum_{i\in I_{\mathrm{ra}}} 
\mathrm{res}_{\zeta_i=0} 
\mathrm{Tr} \left( \hat{\delta} 
(\widehat\Omega_{(\boldsymbol{t}_{\mathrm{ra}},\boldsymbol{\theta}_0)}^{(1)}(\zeta^2_i) ) 
\wedge \hat{\delta}(\tilde\psi_{\zeta_i})\tilde\psi_{\zeta_i}^{-1} \right).
\end{aligned}
\end{equation}}
\end{Prop}

\begin{proof}
Since $\hat{\delta}(\psi_{\zeta_i})\psi_{\zeta_i}^{-1} $
descends under the ramification $x_{t_i}= \zeta_i^2$,
we have 
$$
\frac{1}{4} \cdot 
\mathrm{res}_{\zeta_i=0} 
\mathrm{Tr} \left( \hat{\delta} 
(\widehat\Omega_{(\boldsymbol{t}_{\mathrm{ra}},\boldsymbol{\theta}_0)}^{(1)}(\zeta^2_i) ) 
\wedge \hat{\delta}(\tilde\psi_{\zeta_i})\tilde\psi_{\zeta_i}^{-1} \right)
=\frac{1}{2}  \cdot 
\mathrm{res}_{x=t_i} 
\mathrm{Tr} \left( \hat{\delta} 
(\widehat\Omega_{(\boldsymbol{t}_{\mathrm{ra}},\boldsymbol{\theta}_0)}^{(1)} ) 
\wedge \hat{\delta}(\tilde\psi_{\zeta_i})\tilde\psi_{\zeta_i}^{-1} \right).
$$
By this equality and the same argument as in Proposition \ref{2021.4.26.22.19},
we may check the equality \eqref{eq:2020.11.19.16.08.ra}.
\end{proof}

\begin{Thm}\label{2020.1.31.9.16}
\textit{
For the vector field $\delta_{\theta_{l',t_i}}^{\mathrm{IMD}}$, we have
$\hat{\omega}(\delta_{\theta_{l',t_i}}^{\mathrm{IMD}}, \hat{\delta}) =0$ for any
vector field $\hat{\delta}
 \in \Theta_{\widehat{\mathcal{M}}_{\boldsymbol{t}_{\text{ra}}}\times T_{\boldsymbol{\theta}}}$.
}
\end{Thm}

\begin{proof} 
By the equality \eqref{eq:2020.11.19.16.08.ra}, we have 
\begin{equation*}
\begin{aligned}
\hat\omega (\delta_{\theta_{l',t_i}}^{\mathrm{IMD}}, \hat\delta) &=
 \frac{1}{2} \sum_{i' \in I \setminus I_{\mathrm{ra}} }   \mathrm{res}_{x=t_{i'}} 
\mathrm{Tr} \left( \delta_{\theta_{l',t_i}}^{\mathrm{IMD}}
(\widehat\Omega_{(\boldsymbol{t}_{\mathrm{ra}},\boldsymbol{\theta}_0)}^{(1)}) 
\hat\delta(\tilde\psi_{i'})\tilde\psi_{i'}^{-1} - 
\delta_{\theta_{l',t_i}}^{\mathrm{IMD}}(\tilde\psi_{i'}) \tilde\psi_{i'}^{-1}
\hat\delta (\widehat\Omega_{(\boldsymbol{t}_{\mathrm{ra}},\boldsymbol{\theta}_0)}^{(1)}) 
  \right)  \\
  &\quad +
   \frac{1}{4} \sum_{i' \in I_{\text{ra}}}   \mathrm{res}_{\zeta_{i'}=0} 
\mathrm{Tr} \left( \delta_{\theta_{l',t_i}}^{\mathrm{IMD}}
(\widehat\Omega_{(\boldsymbol{t}_{\mathrm{ra}},\boldsymbol{\theta}_0)}^{(1)} (\zeta^2_{i'})) 
\hat\delta(\tilde\psi_{\zeta_{i'}})\tilde\psi_{\zeta_{i'}}^{-1} - 
\delta_{\theta_{l',t_i}}^{\mathrm{IMD}}(\tilde\psi_{\zeta_{i'}}) \tilde\psi_{\zeta_{i'}}^{-1}
\hat\delta (\widehat\Omega_{(\boldsymbol{t}_{\mathrm{ra}},\boldsymbol{\theta}_0)}^{(1)} (\zeta^2_{i'})) 
  \right).
\end{aligned}
\end{equation*}
We take $i' \in I \setminus I_{\text{ra}}$.
 We may show that 
$$
 \mathrm{res}_{x=t_{i'}} 
\mathrm{Tr} \left( \delta_{\theta_{l',t_i}}^{\mathrm{IMD}}
(\widehat\Omega_{(\boldsymbol{t}_{\mathrm{ra}},\boldsymbol{\theta}_0)}^{(1)}) 
\hat\delta(\tilde\psi_{i'})(\tilde\psi_{i'})^{-1} - 
\delta_{\theta_{l',t_i}}^{\mathrm{IMD}}(\tilde\psi_{i'}) (\tilde\psi_{i'})^{-1}
\hat\delta (\widehat\Omega_{(\boldsymbol{t}_{\mathrm{ra}},\boldsymbol{\theta}_0)}^{(1)}) 
  \right)=0
$$
by the same argument as in the proof of Theorem \ref{2020.1.21.13.51}.
So we have 
\begin{equation*}
\hat\omega (\delta_{\theta_{l',t_i}}^{\mathrm{IMD}}, \hat\delta) =
  \frac{1}{4} \sum_{i' \in I_{\text{ra}}}   \mathrm{res}_{\zeta_{i'}=0} 
\mathrm{Tr} \left( \delta_{\theta_{l',t_i}}^{\mathrm{IMD}}
(\widehat\Omega_{(\boldsymbol{t}_{\mathrm{ra}},\boldsymbol{\theta}_0)}^{(1)} (\zeta^2_{i'})) 
\hat\delta(\tilde\psi_{\zeta_{i'}})\tilde\psi_{\zeta_{i'}}^{-1} - 
\delta_{\theta_{l',t_i}}^{\mathrm{IMD}}(\tilde\psi_{\zeta_{i'}}) \tilde\psi_{\zeta_{i'}}^{-1}
\hat\delta (\widehat\Omega_{(\boldsymbol{t}_{\mathrm{ra}},\boldsymbol{\theta}_0)}^{(1)} (\zeta^2_{i'})) 
  \right).
\end{equation*}

We take $i' \in I_{\text{ra}}$.
Let $U_{t_{i'}}$ be an affine open subset on $\mathbb{P}^1$
so that $x_{t_{i'}}$ is a coordinate on $U_{t_{i'}}$.
We define $\zeta_{i'}$ as $x_{t_{i'}}=\zeta_{i'}^2$.
Let $U_{\zeta_{i'}}$ be the inverse image of 
$U_{t_{i'}}$ under the map 
$\mathrm{Spec}\, \mathbb{C}[\zeta_{i'}] \rightarrow
\mathrm{Spec}\, \mathbb{C}[x_{t_{i'}}]$ by $x_{t_{i'}}=\zeta_{i'}^2$.
We take an analytic open subset $V$ of 
$\widehat{\mathcal{M}}_{\boldsymbol{t}_{\text{ra}}}\times T_{\boldsymbol{\theta}}$.
We take an inverse image of $V$ under the composition of
 $$
\hat{f}_{\zeta_{i'}} \colon
 U_{\zeta_{i'}}  \times
( \widehat{\mathcal{M}}_{\boldsymbol{t}_{\text{ra}}}\times T_{\boldsymbol{\theta}})
\longrightarrow
U_{t_{i'}} \times
( \widehat{\mathcal{M}}_{\boldsymbol{t}_{\text{ra}}}\times T_{\boldsymbol{\theta}})
 $$
 and 
 the projection 
 $p_{U_{t_{i'}}} \colon U_{t_{i'}} \times
( \widehat{\mathcal{M}}_{\boldsymbol{t}_{\text{ra}}}\times T_{\boldsymbol{\theta}})
\rightarrow
\widehat{\mathcal{M}}_{\boldsymbol{t}_{\text{ra}}}\times T_{\boldsymbol{\theta}}$.
Here $\hat{f}_{\zeta_{i'}}$ is defined by $x_{t_i}=\zeta_i^2$.
Let $\widehat{\Delta}_{i'}^{\text{an}}$ ($i' \in I$) be an analytic open
subset of the inverse image $( p_{U_{t_{i'}}} \circ \hat{f}_{\zeta_{i'}}  )^{-1} (V)$
such that $[\zeta_{i'}=0] \cap 
( p_{U_{t_{i'}}} \circ \hat{f}_{\zeta_{i'}}  )^{-1} (V) \subset \widehat{\Delta}_{i'}^{\text{an}}$ 
and the fibers of 
$(p_{U_{t_{i'}}} \circ \hat{f}_{\zeta_{i'}} )|_{\widehat{\Delta}_{i'}^{\text{an}}}
\colon \widehat{\Delta}_{i'}^{\text{an}} \rightarrow V$ 
for each point of $V$
are unit disks such that $\zeta_{i'}$ gives a coordinate 
of the unit disks.
We denote by
$\widehat\Omega_{(\boldsymbol{t}_{\mathrm{ra}},\boldsymbol{\theta}_0)}^{(1)} (\zeta^2_{i'})$
the pull-back of the connection matrix 
$\widehat\Omega_{(\boldsymbol{t}_{\mathrm{ra}},\boldsymbol{\theta}_0)}^{(1)}
|_{U_{t_{i'}} \times
( \widehat{\mathcal{M}}_{\boldsymbol{t}_{\text{ra}}}\times T_{\boldsymbol{\theta}})}$
under the map $\hat{f}_{\zeta_{i'}}$.
We define a matrix $S(\zeta_{i'})$ on 
$ U_{\zeta_{i'}}  \times
( \widehat{\mathcal{M}}_{\boldsymbol{t}_{\text{ra}}}\times T_{\boldsymbol{\theta}})$ as follows:
$$
S(\zeta_{i'}) :=
\tilde{G}^{-1}\Phi_{i'}
\begin{pmatrix}
1 & 0 \\
0 & \zeta_{i'}
\end{pmatrix}.
$$
Remark that $\tilde{G}^{-1}\Phi_{i'}$ is a compatible framing of 
$\widehat\Omega_{(\boldsymbol{t}_{\mathrm{ra}},\boldsymbol{\theta}_0)}^{(1)}
|_{U_{t_{i'}} \times
( \widehat{\mathcal{M}}_{\boldsymbol{t}_{\text{ra}}}\times T_{\boldsymbol{\theta}})}$
at $t_{i'}$.
So we can define the (local) elementary transformation of 
$\widehat\Omega_{(\boldsymbol{t}_{\mathrm{ra}},\boldsymbol{\theta}_0)}^{(1)} (\zeta^2_{i'})$
by $S(\zeta_{i'})$. We denote by 
${}'\widehat\Omega_{(\boldsymbol{t}_{\mathrm{ra}},\boldsymbol{\theta}_0)}^{(1)} (\zeta^2_{i'})$
the elementary transformation.
That is,
$$
{}'\widehat\Omega_{(\boldsymbol{t}_{\mathrm{ra}},\boldsymbol{\theta}_0)}^{(1)} (\zeta^2_{i'}) 
= S(\zeta_{i'})^{-1} d S(\zeta_{i'})
+ S(\zeta_{i'})^{-1}
\widehat\Omega_{(\boldsymbol{t}_{\mathrm{ra}},\boldsymbol{\theta}_0)}^{(1)} (\zeta^2_{i'})
S(\zeta_{i'}).
$$
Let
${}'\widehat\Omega_{(\boldsymbol{t}_{\mathrm{ra}},\boldsymbol{\theta}_0)}^{(1)}
(\zeta_{i'}^2,t)$
be a connection matrix
such that
${}'\widehat\Omega_{(\boldsymbol{t}_{\mathrm{ra}},\boldsymbol{\theta}_0)}^{(1)}
(\zeta_{i'}^2,0)=
{}'\widehat\Omega_{(\boldsymbol{t}_{\mathrm{ra}},\boldsymbol{\theta}_0)}^{(1)}(\zeta^2_{i'})
|_{\widehat{\Delta}_{i'}^{\text{an}} }$
and
$\frac{\partial}{\partial t}
{}'\widehat\Omega_{(\boldsymbol{t}_{\mathrm{ra}},\boldsymbol{\theta}_0)}^{(1)}
(\zeta_{i'}^2,t)|_{t=0}
= \delta_{\theta_{l',t_i}}^{\mathrm{IMD}}
({}'\widehat\Omega_{(\boldsymbol{t}_{\mathrm{ra}},\boldsymbol{\theta}_0)}^{(1)} (\zeta^2_{i'}))
|_{\widehat{\Delta}_{i'}^{\text{an}}}$.
Let $\widehat{\Sigma} \subset \widehat{\Delta}^{\text{an}}_{{i'}}$ 
be a family of sufficiently small sectors in 
$\widehat{\Delta}^{\text{an}}_{{i'}}
\rightarrow 
\widehat{\mathcal{M}}_{\boldsymbol{t}_{\text{ra}}}\times T_{\boldsymbol{\theta}}$. 
By \cite[Theorem 12.1]{Wasow},
we may take a fundamental matrix
solution ${}'\Psi_{\widehat{\Sigma}}(\zeta_{i'},t)$ 
on $ \widehat{\Sigma}  \times U_t^{\text{an}}$
of the differential equation
$$
d({}'\Psi_{\widehat\Sigma}(\zeta_{i'},t))
+{}'\widehat\Omega_{(\boldsymbol{t}_{\mathrm{ra}},\boldsymbol{\theta}_0)}^{(1)}
(\zeta^2_{i'},t)
({}'\Psi_{\widehat\Sigma}(\zeta_{i'},t)) =0
$$
with uniform asymptotic relation 
$$
{}'\Psi_{\widehat\Sigma}(\zeta_{i'},t) \exp\left( \Lambda^-_{i'}(\zeta_{i'} ,t )\right) \sim
\widehat{P}_{i'}(\zeta_{i'},t)
\quad (\zeta_{i'} \rightarrow 0, \zeta_{i'} \in \widehat\Sigma).
$$
Here 
we set
$$
\begin{aligned}
&\Lambda^-_{i'}(\zeta_{i'}) :=\sum_{l'=0}^{2n_{i'}-2}
\begin{pmatrix}
\theta_{l',t_{i'}}  \int \zeta_{i'}^{-2n_{i'}+l'+1} d\zeta_{i'}  & 0 \\
0 &  (-1)^{l'}\theta_{l',t_{i'}}  \int \zeta_{i'}^{-2n_{i'}+l'+1} d\zeta_{i'}
\end{pmatrix} \\
&\widehat{P}_{i'}(\zeta_{i'}) :=
\begin{pmatrix}
1 & 0 \\
0 & \frac{1}{\zeta_{i'}}
\end{pmatrix}
\Xi_{i'} (\zeta^2_{i'}) M_{\zeta_{i'}}
\end{aligned}
$$
and we take
$$
\begin{aligned}
&\widehat{P}_{i'}(\zeta_{i'},t) =\widehat{P}_{i',0}(t) + \widehat{P}_{i',1}(t) \zeta_{i'}+\cdots
\quad \text{and} \quad \\
&\Lambda^-_{i'}(\zeta_{i'} ,t )= \sum_{l'=0}^{2n_{i'}-2}
\begin{pmatrix}
\theta_{l',t_{i'}} (t) \int \zeta_{i'}^{-2n_{i'}+l'+1} d\zeta_{i'}  & 0 \\
0 &  (-1)^{l'}\theta_{l',t_{i'}}(t) \int \zeta_{i'}^{-2n_{i'}+l'+1} d\zeta_{i'}
\end{pmatrix}
\end{aligned}
$$
so that the expansions of $\Lambda^-_{i'}(\zeta_{i'} ,\epsilon )$
and $\widehat{P}_{i'}(\zeta_{i'},\epsilon)$ 
with respect to $\epsilon$ are the following:
$$
\begin{aligned}
\Lambda^-_{i'}(\zeta_{i'} ,\epsilon ) &= \Lambda^-_{i'}(\zeta_{i'}) 
+\epsilon \cdot \delta_{\theta_{l',t_i}}^{\mathrm{IMD}} (\Lambda^-_{i'}(\zeta_{i'})) \\
\widehat{P}_{i'}(\zeta_{i'},\epsilon)
&=\widehat{P}_{i'}(\zeta_{i'})
+\epsilon \cdot \delta_{\theta_{l',t_i}}^{\mathrm{IMD}}
(\widehat{P}_{i'}(\zeta_{i'})).
\end{aligned}
$$
Remark that $\widehat{P}_{i'}(\zeta_{i'})^{-1}$ has no pole at $\zeta_{i'}=0$ and
$S(\zeta_{i'}) \widehat{P}_{i'}(\zeta_{i'}) \exp(-\Lambda^-_{i'}(\zeta_{i'})) = \tilde\psi_{\zeta_{i'}}$.
We set 
$$
{}'\tilde\psi_{\zeta_{i'}} := \widehat{P}_{i'}(\zeta_{i'}) \exp(-\Lambda^-_{i'}(\zeta_{i'})) .
$$
Let $\Upsilon_{\theta_{l',t_{i}}}^{\mathrm{IMD}} (\zeta^2_{i'})$ be the 
pull-back of $\Upsilon_{\theta_{l',t_{i}}}^{\mathrm{IMD}}|_{U_{t_{i'}} \times
( \widehat{\mathcal{M}}_{\boldsymbol{t}_{\text{ra}}}\times T_{\boldsymbol{\theta}})}$
under the map $\hat{f}_{\zeta_{i'}}$.
We set 
$$
{}'\Upsilon_{\theta_{l',t_{i}}}^{\mathrm{IMD}} (\zeta^2_{i'})
:= 
S(\zeta_{t_{i'}})^{-1}\Upsilon_{\theta_{l',t_{i}}}^{\mathrm{IMD}} (\zeta^2_{i'})S(\zeta_{t_{i'}})
- S(\zeta_{i'})^{-1}\delta_{\theta_{l',t_{i}}}^{\mathrm{IMD}}(S(\zeta_{i'})) .
$$
By the same argument as the verification of \eqref{2021.4.18.17.41_1},
we may check the following asymptotic relation:
\begin{equation}\label{2021.11.11.14.25}
 \delta_{\theta_{l',t_{i}}}^{\mathrm{IMD}}({}'\tilde\psi_{\zeta_{i'}}) ({}'\tilde\psi_{\zeta_{i'}})^{-1}
\sim 
{}'\Upsilon_{\theta_{l',t_{i}}}^{\mathrm{IMD}} (\zeta^2_{i'})
-\widehat{P}_{i'}(\zeta_{t_{i'}}) \tilde{C}^{\text{diag}}_{{i'}}\widehat{P}_{i'}(\zeta_{t_{i'}})^{-1}.
\end{equation}
Here $\tilde{C}^{\text{diag}}_{t_{i'}}$ is a diagonal matrix such that 
$\tilde{C}^{\text{diag}}_{t_{i'}}$ is independent of $\zeta_{i'}$.
By the asymptotic relation \eqref{2021.11.11.14.25} and 
the definition of ${}'\Upsilon_{\theta_{l',t_{i}}}^{\mathrm{IMD}} (\zeta^2_{i'})$,
we have 
\begin{equation*}
 \delta_{\theta_{l',t_{i}}}^{\mathrm{IMD}}(\tilde\psi_{\zeta_{i'}}) (\tilde\psi_{\zeta_{i'}})^{-1}
\sim \Upsilon_{\theta_{l',t_{i}}}^{\mathrm{IMD}} (\zeta^2_{i'})
-S(\zeta_{t_{i'}})\widehat{P}_{i'}(\zeta_{t_{i'}}) 
\tilde{C}^{\text{diag}}_{{i'}}
\widehat{P}_{i'}(\zeta_{t_{i'}})^{-1}
S(\zeta_{t_{i'}})^{-1}.
\end{equation*}
By this asymptotic relation and
the same calculation as in the proof of Theorem \ref{2020.1.21.13.51},
we may check the following equalities:
\begin{equation}\label{2021.4.27.14.30}
\begin{aligned}
& \mathrm{res}_{\zeta_{i'}=0} 
\mathrm{Tr} \left( \delta_{\theta_{l',t_i}}^{\mathrm{IMD}}
(\widehat\Omega_{(\boldsymbol{t}_{\mathrm{ra}},\boldsymbol{\theta}_0)}^{(1)}(\zeta^2_{i'})) 
\hat\delta(\tilde\psi_{\zeta_{i'}})(\tilde\psi_{\zeta_{i'}})^{-1} - 
\delta_{\theta_{l',t_i}}^{\mathrm{IMD}}(\tilde\psi_{\zeta_{i'}}) (\tilde\psi_{\zeta_{i'}})^{-1}
\hat\delta (\widehat\Omega_{(\boldsymbol{t}_{\mathrm{ra}},\boldsymbol{\theta}_0)}^{(1)}(\zeta^2_{i'})) 
  \right) \\
&= \mathrm{res}_{\zeta_{i'}=0} 
\mathrm{Tr} \left( \delta_{\theta_{l',t_i}}^{\mathrm{IMD}}
(\widehat\Omega_{(\boldsymbol{t}_{\mathrm{ra}},\boldsymbol{\theta}_0)}^{(1)}(\zeta^2_{i'})) 
\hat\delta(\tilde\psi_{\zeta_{i'}})(\tilde\psi_{\zeta_{i'}})^{-1} - 
\Upsilon_{\theta_{l',t_{i}}}^{\mathrm{IMD}} (\zeta^2_{i'})
\hat\delta (\widehat\Omega_{(\boldsymbol{t}_{\mathrm{ra}},\boldsymbol{\theta}_0)}^{(1)}(\zeta^2_{i'})) 
  \right) \\
  &\qquad +
\mathrm{res}_{\zeta_{i'}=0}\mathrm{Tr} \Big(
\tilde{C}^{\text{diag}}_{t_{i'}} d\Big(
\widehat{P}_{i'}(\zeta_{i'},0)^{-1} S(\zeta_{t_{i'}})^{-1}
\hat\delta(
S(\zeta_{t_{i'}})\widehat{P}_{i'}(\zeta_{i'},0))  \Big)
\Big)  \\
&\qquad \qquad +\mathrm{res}_{\zeta_{i'}=0}\mathrm{Tr} \Big(\tilde{C}^{\text{diag}}_{t_{i'}}
d\Big( \hat\delta(  -\Lambda^-_{i'}(\zeta_{i'}))\Big)
\Big)\\
&=\mathrm{res}_{\zeta_{i'}=0} 
\mathrm{Tr} \left( \delta_{\theta_{l',t_i}}^{\mathrm{IMD}}
(\widehat\Omega_{(\boldsymbol{t}_{\mathrm{ra}},\boldsymbol{\theta}_0)}^{(1)}(\zeta^2_{i'})) 
\hat\delta(\tilde\psi_{\zeta_{i'}})(\tilde\psi_{\zeta_{i'}})^{-1} - 
\Upsilon_{\theta_{l',t_{i}}}^{\mathrm{IMD}} (\zeta^2_{i'})
\hat\delta (\widehat\Omega_{(\boldsymbol{t}_{\mathrm{ra}},\boldsymbol{\theta}_0)}^{(1)}(\zeta^2_{i'})) 
  \right).
  \end{aligned}
\end{equation}
By the integrable condition \eqref{2021.4.27.14.24},
we have the equality 
$$
 \delta_{\theta_{l',t_i}}^{\mathrm{IMD}} 
(\widehat\Omega_{(\boldsymbol{t}_{\mathrm{ra}},\boldsymbol{\theta}_0)}^{(1)}  (\zeta_{i'}^2))
 = 
d ( \Upsilon_{\theta_{l',t_i}}^{\mathrm{IMD}} (\zeta_{i'}^2))
+[\widehat\Omega_{(\boldsymbol{t}_{\mathrm{ra}},\boldsymbol{\theta}_0)}^{(1)}(\zeta_{i'}^2)  ,
\Upsilon_{\theta_{l',t_i}}^{\mathrm{IMD}}(\zeta_{i'}^2)].
$$
We may check that 
\begin{equation}\label{2021.4.27.14.31}
\begin{aligned}
&\mathrm{res}_{\zeta_{i'}=0} 
\mathrm{Tr} \left( \delta_{\theta_{l',t_i}}^{\mathrm{IMD}}
(\widehat\Omega_{(\boldsymbol{t}_{\mathrm{ra}},\boldsymbol{\theta}_0)}^{(1)}(\zeta^2_{i'})) 
\hat\delta(\tilde\psi_{\zeta_{i'}})(\tilde\psi_{\zeta_{i'}})^{-1} - 
\Upsilon_{\theta_{l',t_{i}}}^{\mathrm{IMD}} (\zeta^2_{i'})
\hat\delta (\widehat\Omega_{(\boldsymbol{t}_{\mathrm{ra}},\boldsymbol{\theta}_0)}^{(1)}(\zeta^2_{i'})) 
  \right)\\
  &= \mathrm{res}_{\zeta_{i'}=0} 
\mathrm{Tr} \left(  
d \left( \Upsilon_{\theta_{l',t_i}}^{\mathrm{IMD}} (\zeta_{i'}^2)
\hat\delta(\tilde\psi_{\zeta_{i'}})(\tilde\psi_{\zeta_{i'}})^{-1}
\right)
\right)=0
  \end{aligned}
\end{equation}
by the same calculation as in the proof of Theorem \ref{2020.1.21.13.51}.
By combining \eqref{2021.4.27.14.30} and \eqref{2021.4.27.14.31},
we obtain that 
\begin{equation*}
\hat\omega (\delta_{\theta_{l,t_i}^{\pm}}^{\mathrm{IMD}}, \hat\delta) =
  \frac{1}{4} \sum_{i' \in I_{\text{ra}}}   \mathrm{res}_{\zeta_{i'}=0} 
\mathrm{Tr} \left( \delta_{\theta_{l',t_i}}^{\mathrm{IMD}}
(\widehat\Omega_{(\boldsymbol{t}_{\mathrm{ra}},\boldsymbol{\theta}_0)}^{(1)} (\zeta^2_{i'})) 
\hat\delta(\tilde\psi_{\zeta_{i'}})\tilde\psi_{\zeta_{i'}}^{-1} - 
\delta_{\theta_{l',t_i}}^{\mathrm{IMD}}(\tilde\psi_{\zeta_{i'}}) \tilde\psi_{\zeta_{i'}}^{-1}
\hat\delta (\widehat\Omega_{(\boldsymbol{t}_{\mathrm{ra}},\boldsymbol{\theta}_0)}^{(1)} (\zeta^2_{i'})) 
  \right)=0.
\end{equation*}
That is, we have the assertion of this theorem.
\end{proof}

As in Section \ref{2020.1.24.22.25} and Section \ref{2020.1.25.8.11},
we may define vector fields 
$\delta_{\theta_{l,t_i}^{\pm}}^{\mathrm{IMD}}$ ($i \in I_{\mathrm{un}}$, $0\le l\le n_i-2$)
and $\delta_{t_i}^{\mathrm{IMD}}$ 
($i \in  \{ 3,4,\ldots,\nu\} \cap (I_{\mathrm{reg}} \cup I_{\mathrm{un}})$), respectively.
By the same argument as in the proof of 
Theorem \ref{2020.1.21.13.51} and Theorem \ref{2020.1.31.9.16},
we obtain the following theorem.

\begin{Thm}\label{2020.2.2.19.05}
\textit{
Let $i \in I_{\mathrm{un}}$.
For the vector field $\delta_{\theta_{l,t_i}^{\pm}}^{\mathrm{IMD}}$, we have
$\hat{\omega}(\delta_{\theta_{l,t_i}^{\pm}}^{\mathrm{IMD}},\hat{\delta}) =0$ for any
vector field $\hat{\delta} 
\in \Theta_{\widehat{\mathcal{M}}_{\boldsymbol{t}_{\text{ra}}}\times T_{\boldsymbol{\theta}}}$.
Moreover, 
for the vector field $\delta_{t_i}^{\mathrm{IMD}}$, we have
$\hat{\omega}(\delta_{t_i}^{\mathrm{IMD}}, \hat{\delta}) =0$ for any
vector field $\hat{\delta} \in 
\Theta_{\widehat{\mathcal{M}}_{\boldsymbol{t}_{\text{ra}}}\times T_{\boldsymbol{\theta}}}$.
}
\end{Thm}

Then we have that the $2$-form $\hat{\omega}$ is the isomonodromy $2$-form.

\subsection{Hamiltonian systems}\label{2020.2.8.10.35}

We have the following diagonalization:
$$
\begin{aligned}
\Omega_{t_i}^{\text{diag}} (\zeta_i^2):=&\ 
(\Phi_i \Xi_iM_{\zeta_i})^{-1} d(\Phi_i \Xi_iM_{\zeta_i}) +(\Phi_i \Xi_iM_{\zeta_i})^{-1} 
\widehat\Omega_{(\boldsymbol{t}_{\mathrm{ra}},\boldsymbol{\theta}_0)}^{(n-2)} (\zeta_i^2)
(\Phi_i \Xi_iM_{\zeta_i})\\
\Omega_{t_i}^{\text{diag}} (\zeta_i^2)=&\ 
\begin{pmatrix}
\theta_{0,t_i}  & 0 \\
0 &  \theta_{0,t_i}
\end{pmatrix}
\frac{d\zeta_i}{\zeta_i^{2n_i-1}}
+ \cdots+
\begin{pmatrix}
\theta_{2n_i-2,t_i}  & 0 \\
0 & \theta_{2n_i-2,t_i}
\end{pmatrix}
 \frac{d\zeta_i}{\zeta_i} \\
&\qquad+\begin{pmatrix}
\theta_{2n_i-1,t_i}  & 0 \\
0 & -\theta_{2n_i-1,t_i}
\end{pmatrix}d\zeta_i
+\cdots +
\begin{pmatrix}
\theta_{4n_i-4,t_i}  & 0 \\
0 & \theta_{4n_i-4,t_i}
\end{pmatrix} \zeta_i^{2n_i-3} dx
+\cdots.
\end{aligned}
$$
 Remark that
 we have an equation $(d +\Omega_{t_i}^{\text{diag}} (\zeta_i^2)) \exp(-\Lambda_i( \zeta_i))=0$.
We set 
\begin{equation*}
\Lambda^+_i( \zeta_i)= 
\begin{pmatrix}
\sum_{l'=2n_i-1}^{\infty} \theta_{l',t_i} \int \zeta_i^{-2n_i+l'+1} d\zeta_i & 0 \\
0 &\sum_{l'=2n_i-1}^{\infty}(- 1)^{l'} \theta_{l',t_i} \int \zeta_i^{-2n_i+l'+1} d\zeta_i
\end{pmatrix} .
\end{equation*}

\begin{Def}\label{2020.2.7.12.01}
For each $t_i $ ($i \in I_{\text{ra}}$) and 
each $l'$ ($0\le l' \le 2n_i-3$), we define rational function 
$H_{\theta_{l' , t_i}}$ on 
$\widehat{\mathcal{M}}_{\boldsymbol{t}_{\text{ra}}}\times T_{\boldsymbol{\theta}}$ as
\begin{equation*}
\begin{aligned}
H_{\theta_{l' , t_i}}&= 
-[\text{ the coefficient of the $\zeta_i^{2(n_i-1) - l'}$-term of
the $(1,1)$-entry of $\Lambda^+_i( \zeta_i)$ }] \\
&=-\frac{\theta_{4(n_i-1)- l',t_i}}{2(n_i-1)-l'}.
\end{aligned}
\end{equation*}
We call $H_{\theta_{l' , t_i}}$ the {\it Hamiltonian associated to $\theta_{l',t_i}$}.
\end{Def}

We define $H_{\theta_{l,t_i}^{\pm}}$ 
($i\in I_{\mathrm{un}}$ and $l=0,1,\ldots, n_i-2$)
and $H_{t_i}$ 
($i\in \{ 3,4,\ldots,\nu\} \cap (I_{\mathrm{reg}} \cup I_{\mathrm{un}})$)
as in Definition \ref{2020.2.7.11.59} and Definition \ref{2020.2.7.12.00}.

\begin{Thm}\label{2020.1.29.11.46}
\textit{
Set $P(x;\boldsymbol{t}):=\prod_{i=1}^{\nu} (x-t_i)^{n_i}$ and  
$D_i(x;\boldsymbol{t}, \boldsymbol{\theta} ):=D_i(x)$ for $i \in I$.
We put
\begin{equation*}
\begin{aligned}
\hat{\omega}' :=&\ \sum_{j=1}^{n-3} d \left(\frac{p_j}{P(q_j;\boldsymbol{t})} 
 - \sum_{i=1}^{\nu} \frac{D_i(q_j;\boldsymbol{t}, \boldsymbol{\theta})}{(q_j -t_i)^{n_i}} 
 - D_{\infty}(q_j;\boldsymbol{t}, \boldsymbol{\theta})\right) \wedge dq_j \\
&+ \sum_{i\in I_{\mathrm{un}}} 
\sum_{l=0}^{n_{i}-2} \left( d  H_{\theta_{l,t_i}^{+}} \wedge d \theta_{l,t_i}^{+}
+ d H_{\theta_{l,t_i}^{-}} \wedge d \theta_{l,t_i}^{-} \right) \\
&+ \sum_{i \in I_{\mathrm{ra}}}\sum_{l'=0}^{2n_i-3} d H_{\theta_{l' , t_i}} \wedge d\theta_{l' , t_i}
+\sum_{i\in \{3,4,\ldots,\nu \}\cap (I_{\mathrm{reg}} \cup I_{\mathrm{un}})} d H_{t_i} \wedge dt_i.
\end{aligned}
\end{equation*}
Then the difference $\hat{\omega}- \hat{\omega}'$ 
is a section of $\pi_{\boldsymbol{t}_{\mathrm{ra}}, \boldsymbol{\theta}_0}^* 
(\Omega^2_{(T_{\boldsymbol{t}})_{\boldsymbol{t}_{\text{ra}}} \times T_{\boldsymbol{\theta}}} )$.
}
\end{Thm}

\begin{proof}
Recall that $\hat{\omega} (\hat{\delta}_1,\hat{\delta}_2)$ is
\begin{equation*}
\begin{aligned}
 &  \frac{1}{2} \sum_{i\in I \setminus I_{\mathrm{ra}}} 
\mathrm{res}_{x=t_i} 
\mathrm{Tr} \left( \hat{\delta}
(\widehat\Omega_{(\boldsymbol{t}_{\mathrm{ra}},\boldsymbol{\theta}_0)}^{(n-2)} ) 
\wedge \hat{\delta}(\psi_i)\psi_i^{-1}  \right) 
+\frac{1}{4} \sum_{i\in I_{\mathrm{ra}}} 
\mathrm{res}_{\zeta_i=0} 
\mathrm{Tr} \left( \hat{\delta} 
(\widehat\Omega_{(\boldsymbol{t}_{\mathrm{ra}},\boldsymbol{\theta}_0)}^{(n-2)}(\zeta^2_i) ) 
\wedge \hat{\delta}(\psi_{\zeta_i})\psi_{\zeta_i}^{-1} \right) \\
&+\frac{1}{2} \sum_{j=1}^{n-3} \mathrm{res}_{x=q_j} 
\mathrm{Tr} \left( \hat{\delta}
 (\widehat\Omega_{(\boldsymbol{t}_{\mathrm{ra}},\boldsymbol{\theta}_0)}^{(n-2)} )
\wedge \hat{\delta}(\psi_{q_j})\psi_{q_j}^{-1}  \right).
\end{aligned}
\end{equation*}
The plan of this proof is as follows.
Since the calculation of the first term and the third term in this formula 
is the same as in the proof of Theorem \ref{2020.1.21.13.48}.
So we omit the calculation of these terms.
Now we consider only the second term
$$
\frac{1}{4} \sum_{i\in I_{\mathrm{ra}}} 
\mathrm{res}_{\zeta_i=0} 
\mathrm{Tr} \left( \hat{\delta} 
(\widehat\Omega_{(\boldsymbol{t}_{\mathrm{ra}},\boldsymbol{\theta}_0)}^{(n-2)}(\zeta^2_i) ) 
\wedge \hat{\delta}(\psi_{\zeta_i})\psi_{\zeta_i}^{-1} \right).
$$
We calculate 
the residue at $\zeta_i=0$ for some (local) gauge transformation of 
$d+\widehat\Omega_{(\boldsymbol{t}_{\mathrm{ra}},\boldsymbol{\theta}_0)}^{(n-2)} (\zeta_i^2)$.
We need to consider the difference between the residue after taking the gauge transformation
and the residue before taking the gauge transformation.
Here the residue before taking the gauge transformation is just
the residue of $\mathrm{Tr} ( \hat{\delta} 
(\widehat\Omega_{(\boldsymbol{t}_{\mathrm{ra}},\boldsymbol{\theta}_0)}^{(n-2)}(\zeta^2_i) ) 
\wedge \hat{\delta}(\psi_{\zeta_i})\psi_{\zeta_i}^{-1} )$ at $\zeta_i=0$.
This difference is more complicated than the unramified irregular singular cases,
since the $(1,2)$-entry of the residue part of \eqref{2021.4.23.12.09}
is not a piece of the local formal data.
We consider this difference by using Lemma \ref{2020.12.24.17.31}.

In fact, for each $i \in I_{\mathrm{ra}}$,
we consider the residue 
\begin{equation}\label{2020.1.24.11.09}
\frac{1}{4}\cdot \mathrm{res}_{\zeta_i=0}
\mathrm{Tr} \left( \hat{\delta}
(\widehat\Omega_{(\boldsymbol{t}_{\mathrm{ra}},\boldsymbol{\theta}_0)}^{(n-2)} (\zeta^2_{i'}) ) 
\wedge \hat{\delta}(\psi_{\zeta_i})\psi_{\zeta_i}^{-1} 
\right).
\end{equation}
Now
we take a diagonalization of 
$d+\widehat\Omega_{(\boldsymbol{t}_{\mathrm{ra}},\boldsymbol{\theta}_0)}^{(n-2)} (\zeta^2_{i'})$ 
until some degree term at $\zeta_i=0$.
We put
\begin{equation*}
\Xi_{i}^{\le 2n_i-1}(x_{t_i}) := 
\begin{pmatrix}
1 & 0\\
0 & 1
\end{pmatrix}+
\sum_{s=1}^{2n_i-1} 
\begin{pmatrix}
(\xi^{(i)}_s)_{11} & (\xi^{(i)}_s)_{12} \\
(\xi^{(i)}_s)_{21} & (\xi^{(i)}_s)_{22}
\end{pmatrix}
x_{t_i}^s 
\end{equation*}
for $i\in I$.
Here the coefficient matrices of $\Xi_{i}^{\le 2n_i-1}(x_{t_i})$ appear in  
as the coefficient matrices of $\Xi_{i}(x_{t_i})$ defined by \eqref{2021.4.23.18.23}.
Moreover, we put
\begin{equation*}
\begin{aligned}
&\tilde\Omega_{i} :=  (\Phi_i\Xi_{i}^{\le 2n_i-1})^{-1} d(\Phi_i\Xi_{i}^{\le 2n_i-1})
+(\Phi_i\Xi_{i}^{\le 2n_i-1})^{-1} \widehat\Omega_{(\boldsymbol{t}_{\mathrm{ra}},\boldsymbol{\theta}_0)}^{(n-2)} 
 (\Phi_i\Xi_{i}^{\le 2n_i-1})  , \\
&\tilde\Omega'_{\zeta_i}:= M_{\zeta_i}^{-1} dM_{\zeta_i} 
+M_{\zeta_i}^{-1} \tilde{\Omega}_i  M_{\zeta_i},\\
&\tilde{\psi}_{\zeta_i} := (\Phi_i\Xi_{i}^{\le 2n_i-1})^{-1} 
\psi_{\zeta_i},
\quad \text{and} \quad 
\tilde{\psi}_{\zeta_i}':=M_{\zeta_i}^{-1} \tilde{\psi}_{\zeta_i},
\end{aligned}
\end{equation*}
where $\psi_{\zeta_i}$ is the formal solution (\ref{2020.1.10.18.22}).
We may describe $\tilde\Omega_i$ as follows:
$$
\begin{aligned}
\tilde{\Omega}_{\zeta_i}'&=
\begin{pmatrix}
\theta_{0,t_i}  & 0 \\
0 &  \theta_{0,t_i}
\end{pmatrix}
\frac{d\zeta_i}{\zeta_i^{2n_i-1}}
+ \cdots+
\begin{pmatrix}
\theta_{2n_i-2,t_i}  & 0 \\
0 & \theta_{2n_i-2,t_i}
\end{pmatrix}
 \frac{d\zeta_i}{\zeta_i} \\
&\qquad+\begin{pmatrix}
\theta_{2n_i-1,t_i}  & 0 \\
0 & -\theta_{2n_i-1,t_i}
\end{pmatrix}d\zeta_i
+\cdots +
\begin{pmatrix}
\theta_{4n_i-4,t_i}  & 0 \\
0 & \theta_{4n_i-4,t_i}
\end{pmatrix} \zeta_i^{2n_i-3} dx +O(\zeta_i^{2n_i-2}).
\end{aligned}
$$
The residue part $\theta_{2n_i-2,t_i}$ of $\tilde{\Omega}_{\zeta_i}'$
is constant on $\widehat{\mathcal{M}}_{\boldsymbol{t}_{\text{ra}}}\times T_{\boldsymbol{\theta}}$.
So we have $\hat{\delta}(\theta_{2n_i-2,t_i})=0$ for any 
$\hat{\delta} \in \Theta_{\widehat{\mathcal{M}}_{\boldsymbol{t}_{\text{ra}}}\times T_{\boldsymbol{\theta}}}$.
Then the variation $\hat{\delta}_1(\tilde{\Omega}_{\zeta_i}')$ is equal to
$$
\begin{aligned}
&\begin{pmatrix}
\hat{\delta}_1(\theta_{0,t_i})  & 0 \\
0 & \hat{\delta}_1( \theta_{0,t_i})
\end{pmatrix}
\frac{d\zeta_i}{\zeta_i^{2n_i-1}}
+ \cdots+
\begin{pmatrix}
\hat{\delta}_1 (\theta_{2n_i-3,t_i})  & 0 \\
0 & -\hat{\delta}_1( \theta_{2n_i-3,t_i})
\end{pmatrix}
 \frac{d\zeta_i}{\zeta_i^2} \\
&\qquad+\begin{pmatrix}
\hat{\delta}_1(\theta_{2n_i-1,t_i} ) & 0 \\
0 & - \hat{\delta}_1(\theta_{2n_i-1,t_i})
\end{pmatrix}d\zeta_i
+\cdots +
\begin{pmatrix}
\hat{\delta}_1(\theta_{4n_i-4,t_i} ) & 0 \\
0 &  \hat{\delta}_1(\theta_{4n_i-4,t_i})
\end{pmatrix} \zeta_i^{2n_i-3} dx +O(\zeta_i^{2n_i-2}).
\end{aligned}
$$
We define $\hat{\lambda}^{\le 4n_i-4}_{i,\pm}(\zeta_i)$ as
$$
\begin{aligned}
\hat{\lambda}^{\le 4n_i-4}_{i,\pm}(\zeta_i)&=
(\pm 1)^{0}  \frac{\theta_{0,t_i}}{-2n_i+2}  \zeta_i^{-2n_i+2} 
+\cdots +
(\pm 1)^{2n_i-3}  \frac{\theta_{2n_i-3,t_i}}{-1}  \zeta_i^{-1} + \theta_{2n_i-2,t_i} \log \zeta_i \\
&\quad +(\pm 1)^{2n_i-1} \theta_{2n_i-1,t_i}  \zeta_i
+\cdots
+(\pm 1)^{4n_i-4}  \frac{\theta_{4n_i-4,t_i}}{2n_i-2}  \zeta_i^{2n_i-2}.
\end{aligned}
$$
On the other hand, the variation $\hat{\delta}_{2}(\tilde{\psi}'_{\zeta_i})(\tilde{\psi}'_{\zeta_i})^{-1}$ is equal to
$$
\begin{aligned}
& \delta((\Phi_i\Xi_{i}^{\le 2n_i-1}M_{\zeta_i})^{-1} \Phi_{i} \Xi_{i} M_{\zeta_i})
 ((\Phi_i\Xi_{i}^{\le 2n_i-1}M_{\zeta_i})^{-1} \Phi_{i} \Xi_{i}M_{\zeta_i})^{-1}\\
 &\qquad + 
((\Phi_i\Xi_{i}^{\le 2n_i-1}M_{\zeta_i})^{-1} \Phi_{i} \Xi_{i} M_{\zeta_i}) 
\begin{pmatrix}
-\hat\delta_2 (\hat{\lambda}_{i,+}(\zeta_i)) & 0 \\
0 &-\hat\delta_2 (\hat{\lambda}_{i,-}(\zeta_i))
\end{pmatrix}
 ((\Phi_i\Xi_{i}^{\le 2n_i-1}M_{\zeta_i})^{-1} \Phi_{i} \Xi_{i}M_{\zeta_i})^{-1}\\
 &=
\begin{pmatrix}
-\hat{\delta}_2 ( \hat{\lambda}^{\le 4n_i-4}_{i,+}(\zeta_i))& 0 \\
0 &-\hat{\delta}_2 (\hat{\lambda}^{\le 4n_i-4}_{i,-}(\zeta_i))
\end{pmatrix}  + O(\zeta_i^{2n_{i}-1}).
\end{aligned}
$$
Since
$$
\begin{aligned}
\hat{\delta}_2(\hat{\lambda}^{\le 4n_i-4}_{i,\pm}(\zeta_i))&=
(\pm 1)^{0}  \frac{\hat{\delta}_2(\theta_{0,t_i})}{-2n_i+2}  \zeta_i^{-2n_i+2} 
+\cdots +
(\pm 1)^{2n_i-3}  \frac{\hat{\delta}_2(\theta_{2n_i-3,t_i})}{-1}  \zeta_i^{-1}  \\
&\quad +(\pm 1)^{2n_i-1}\hat{\delta}_2( \theta_{2n_i-1,t_i}  )\zeta_i
+\cdots
+(\pm 1)^{4n_i-4}  \frac{\hat{\delta}_2(\theta_{4n_i-4,t_i})}{2n_i-2}  \zeta_i^{2n_i-2},
\end{aligned}
$$
we may check that the residue of 
$\mathrm{Tr} ( \hat{\delta}_1 (\tilde{\Omega}'_{\zeta_i}) 
 \hat{\delta}_2(\tilde{\psi}'_{\zeta_i})(\tilde{\psi}_{\zeta_i}')^{-1}   ) $
at $\zeta_i=0$ is equal to 
\begin{equation*}
\begin{aligned}
& \sum_{l'\in \{0,1,\ldots,4n_i-4 \} \setminus \{2n_i-2 \}} 2 \cdot \left( \hat{\delta}_1(\theta_{l',t_i})\cdot
\frac{\hat{\delta}_2(\theta_{4(n_i-1)-l',t_i})}{2(n_i-1)-l'}\right)\\
&=2 \cdot \left(\sum_{l'=0}^{2n_i-3}   \hat{\delta}_1(\theta_{l',t_i})\cdot
\frac{\hat{\delta}_2(\theta_{4(n_i-1)-l',t_i})}{2(n_i-1)-l'}
-\sum_{l'=0}^{2n_i-3} 
\frac{\hat{\delta}_1(\theta_{4(n_i-1)-l',t_i})}{2(n_i-1)-l'}
\cdot \hat{\delta}_2(\theta_{l',t_i}) \right).
\end{aligned}
\end{equation*}
Then we have 
\begin{equation}\label{2020.1.24.11.10}
\frac{1}{4}\sum_{i \in I_{\mathrm{ra}}}\mathrm{res}_{\zeta_i=0} 
\mathrm{Tr} \left( \hat{\delta} (\tilde{\Omega}'_{\zeta_i}) 
\wedge \hat{\delta}(\tilde{\psi}'_{\zeta_i})(\tilde{\psi}_{\zeta_i}')^{-1}   \right) 
=\left(
\sum_{i \in I_{\mathrm{ra}}} 
\sum_{l'=0}^{2n_{i}-3}  d  H_{\theta_{l',t_i}} \wedge d \theta_{l',t_i}
  \right) (\hat{\delta}_1 , \hat{\delta}_2).
\end{equation}

Now we consider the difference between the residue after taking the gauge transformation
and the residue before taking the gauge transformation.
 Put $g_{\zeta_i}:=\Phi_i\Xi_{i}^{\le 2n_i-1} M_{\zeta_i}$ 
 and $g:= \Phi_i\Xi_{i}^{\le 2n_i-1}$.
We consider
the difference between 
the residue (\ref{2020.1.24.11.09}) and the residue (\ref{2020.1.24.11.10})
when $\hat{\delta}_1 = \delta_1 \in 
\Theta_{(\widehat{\mathcal{M}}_{\boldsymbol{t}_{\text{ra}}}\times T_{\boldsymbol{\theta}}) 
/((T_{\boldsymbol{t}})_{\boldsymbol{t}_{\text{ra}}} \times T_{\boldsymbol{\theta}})}$:
\begin{equation}\label{2021.4.30.18.44}
\begin{aligned}
&\mathrm{Tr} \left( \hat{\delta} (\tilde{\Omega}'_{\zeta_i}) 
\wedge \hat{\delta}(\tilde{\psi}'_{\zeta_i})(\tilde{\psi}_{\zeta_i}')^{-1}   \right) 
-\mathrm{Tr} \left( \hat{\delta} 
(\widehat\Omega_{(\boldsymbol{t}_{\mathrm{ra}},\boldsymbol{\theta}_0)}^{(n-2)}(\zeta^2_i) ) 
\wedge \hat{\delta}(\psi_{\zeta_i})\psi_{\zeta_i}^{-1} \right) \\
&=
-\mathrm{Tr} ( \delta_1(\tilde{\Omega}_{\zeta_i}') \tilde{u}^{(2)} -
\tilde{u}^{(1)}\hat{\delta}_2(\tilde{\Omega}_{\zeta_i}')  )
-\mathrm{Tr} \left( 
\delta_1(\widehat\Omega_{(\boldsymbol{t}_{\mathrm{ra}},\boldsymbol{\theta}_0)}^{(n-2)} (\zeta^2_{i'}))
u^{(2)} -u^{(1)}
\hat{\delta}_2(\widehat\Omega_{(\boldsymbol{t}_{\mathrm{ra}},\boldsymbol{\theta}_0)}^{(n-2)} (\zeta^2_{i'}))
  \right)\\
  &\qquad + \mathrm{Tr}\left( d( \psi_{\zeta_i}^{-1} u^{(1)} \hat\delta_2(\psi_{\zeta_i})
  -\psi_{\zeta_i}^{-1} u^{(2)} \delta_1(\psi_{\zeta_i}))
  \right).
\end{aligned}
\end{equation}
Here we set $u^{(k)}:= \hat{\delta}_{k} ( g_{\zeta_i}) g_{\zeta_i}^{-1}$ 
and $\tilde{u}^{(k)}:= g_{\zeta_i}^{-1} \hat{\delta}_{k} ( g_{\zeta_i}) $ for $k=1,2$.
We calculate the residue of the second term of the right hand side of \eqref{2021.4.30.18.44}
at $\zeta_i=0$.
Since $t_i$ ($i\in I_{\mathrm{ra}}$) is not a deformation parameter,
$\hat{\delta}_{k} (t_i)=\hat{\delta}_{k} (\zeta_i)=0$ for $k=1,2$.
Then $u^{(k)}$ coincides with
$\hat{\delta}_{k} (g) g^{-1}$.
We expand $g$ as \eqref{2020.2.16.21.26}.
Since $g_0,\ldots, g_{n_i-2}$ are parametrized by only 
$(T_{\boldsymbol{t}})_{\boldsymbol{t}_{\text{ra}}} \times T_{\boldsymbol{\theta}}$
and $\delta_1 \in 
\Theta_{(\widehat{\mathcal{M}}_{\boldsymbol{t}_{\text{ra}}}\times T_{\boldsymbol{\theta}}) 
/((T_{\boldsymbol{t}})_{\boldsymbol{t}_{\text{ra}}} \times T_{\boldsymbol{\theta}})}$,
the variations $\delta_1(g_0), \ldots, \delta_1(g_{n_i-2})$ vanish.

We will calculate the variation $\delta_1(g_{n_i-1})$.
We consider the following gauge transformation 
\begin{equation}\label{2022.4.9.15.55}
(g_0 M_{\zeta_i})^{-1}d(g_0 M_{\zeta_i})
+(g_0 M_{\zeta_i})^{-1}
\widehat\Omega_{(\boldsymbol{t}_{\mathrm{ra}},\boldsymbol{\theta}_0)}^{(n-2)}(\zeta_i^2)
(g_0 M_{\zeta_i}).
\end{equation}
The $\zeta_i^{-2n_i+1}$-term
and the $\zeta_i^{-2n_i+2}$-term 
of the expansion of \eqref{2022.4.9.15.55}
at $\zeta_i=0$
are diagonal.
The eigenvalues of the $\zeta_i^{-2n_i+2}$-term are distinct.
The terms of this expansion after the $\zeta_i^{-2n_i+2}$-term are 
diagonalized by
the gauge transformation
by the right hand side of \eqref{2022.4.9.15.44}. 
Since this negative part 
of \eqref{2022.4.9.15.55}
is independent of 
$(q_j,p_j)_{1\le j\le n-3}$,
we have that
the coefficients of
 $\zeta_i \cdot \xi(\zeta_i)$ in the right hand side of \eqref{2022.4.9.15.44}
are independent of 
$(q_j,p_j)_{1\le j\le n-3}$
until the $\zeta_i^{2n_i-3}$-term.
So the $(2,1)$-entry of $g_0^{-1} g_{n_i-1}$ is independent of 
$(q_j,p_j)_{1\le j\le n-3}$.
Since $\delta_1 \in 
\Theta_{(\widehat{\mathcal{M}}_{\boldsymbol{t}_{\text{ra}}}\times T_{\boldsymbol{\theta}}) 
/((T_{\boldsymbol{t}})_{\boldsymbol{t}_{\text{ra}}} \times T_{\boldsymbol{\theta}})}$
and the $(2,1)$-entry of $g_0^{-1} g_{n_i-1}$ is independent of 
$(q_j,p_j)_{1\le j\le n-3}$, 
we have that
the $(2,1)$-entry of $g_0^{-1}\delta_1(g_{n_i-1})$ vanishes. 
Moreover we have that $\mathrm{Tr}(g_0^{-1} g_{n_i-1})$ is constant.
By comparing the $x_{t_i}^{-1}$-terms of the expansions of the both sides of  
$$
g_0^{-1}g\,  \tilde{\Omega}_i 
= g_0^{-1}dg+ 
g_0^{-1} \widehat\Omega_{(\boldsymbol{t}_{\mathrm{ra}},\boldsymbol{\theta}_0)}^{(n-2)}
g_0(g_0^{-1}  g),$$ 
we have the following equality:
\begin{equation}\label{2021.1.2.11.06}
\begin{aligned}
&g_0^{-1}g_{n_i-1}
\begin{pmatrix}
\theta_{0,t_i}/2 & \theta_{1,t_i}/2 \\
0 & \theta_{0,t_i}/2
\end{pmatrix}+ \cdots + g_0^{-1} g_0 \begin{pmatrix}
\theta_{2n_i-2,t_i}/2 & \theta_{2n_i-1,t_i}/2 \\
\theta_{2n_i-3,t_i}/2 & \theta_{2n_i-2,t_i}/2
\end{pmatrix}\\
&= \begin{pmatrix}
\theta_{0,t_i}/2 & \theta_{1,t_i}/2 \\
0 & \theta_{0,t_i}/2
\end{pmatrix} g_0^{-1}g_{n_i-1} + \cdots 
+ g_0^{-1}
(\widehat\Omega_{(\boldsymbol{t}_{\mathrm{ra}},\boldsymbol{\theta}_0)}^{(n-2)})_{n_i-1}
g_0g_0^{-1} g_0
\end{aligned}.
\end{equation}
We consider the variations of the $x_{t_i}^{-1}$-terms of the both sides
of \eqref{2021.1.2.11.06}.
In particular, we focus on the $(1,2)$-entries of the both sides.
Since $\delta_1 \in 
\Theta_{(\widehat{\mathcal{M}}_{\boldsymbol{t}_{\text{ra}}}\times T_{\boldsymbol{\theta}}) 
/((T_{\boldsymbol{t}})_{\boldsymbol{t}_{\text{ra}}} \times T_{\boldsymbol{\theta}})}$
and $\mathrm{Tr}(g_0^{-1} g_{n_i-1})$ is constant,
we have explicit descriptions of
the $(1,1)$-entry and the $(2,2)$-entry of $g_0^{-1}\delta_1(g_{n_i-1})$.
So we have 
\begin{equation}\label{2020.3.5.21.08}
g_0^{-1}\delta_1(g_{n_i-1})=
\begin{pmatrix}
-\frac{\delta_1(\theta_{2n_i-1,t_i})}{2\theta_{1,t_i}} & *   \\
0 &   \frac{\delta_1(\theta_{2n_i-1,t_i})}{2\theta_{1,t_i}}
\end{pmatrix}.
\end{equation}
Since $\delta_1(g_0), \ldots, \delta_1(g_{n_i-2})$ vanish, 
we have
\begin{equation}\label{2021.11.10.17.18}
 g_0^{-1}\delta_1(g)=
\begin{pmatrix}
-\frac{\delta_1(\theta_{2n_i-1,t_i})}{2\theta_{1,t_i}} & *   \\
0 &   \frac{\delta_1(\theta_{2n_i-1,t_i})}{2\theta_{1,t_i}}
\end{pmatrix} (x-t_i)^{n_i-1} + O((x-t_i)^{n_i}).
\end{equation}

On the other hand, 
$\hat{\delta}_2(\widehat\Omega_{(\boldsymbol{t}_{\mathrm{ra}},\boldsymbol{\theta}_0)}^{(n-2)})$
has the following expansion:
$$
\hat{\delta}_2(\widehat\Omega_{(\boldsymbol{t}_{\mathrm{ra}},\boldsymbol{\theta}_0)}^{(n-2)})
= \begin{pmatrix}
0 & - \frac{\hat{\delta}_2(\prod_{j \neq i}(t_i - t_j)^{n_j})}{\prod_{j \neq i}(t_i - t_j)^{2n_j}} \\
-\frac{1}{4} \cdot \hat{\delta}_2(\theta^2_{0,t_i} \prod_{j \neq i}(t_i - t_j)^{n_j}) & 
\hat{\delta}_2(\theta_{0,t_i})
\end{pmatrix} \frac{1}{(x-t_i)^{n_i}}+
[\text{ higher order terms }].
$$
Now we take a compatible framing $g_0$ as 
\begin{equation}\label{2021.1.2.18.21}
g_0 = 
\begin{pmatrix}
1 &\frac{\theta_{0,t_i}}{2} \prod_{j \neq i}(t_i - t_j)^{n_j} \\
\frac{\theta_{0,t_i}}{2} \prod_{j \neq i}(t_i - t_j)^{n_j} 
&(\frac{\theta_{0,t_i}^2}{4} \prod_{j \neq i}(t_i - t_j)^{n_j}+
\frac{\theta_{1,t_i}}{2} ) \prod_{j \neq i}(t_i - t_j)^{n_j}
\end{pmatrix}.
\end{equation}
Then the leading term of
$g_0^{-1}
\widehat\Omega_{(\boldsymbol{t}_{\mathrm{ra}},\boldsymbol{\theta}_0)}^{(n-2)} g_0$
coincides with the leading term of the right hand side of \eqref{2021.4.21.10.45}.
We have the following expansion of 
$g_0^{-1}\hat{\delta}_2(
\widehat\Omega_{(\boldsymbol{t}_{\mathrm{ra}},\boldsymbol{\theta}_0)}^{(n-2)}) g_0$ 
at $t_i$:
\begin{equation}\label{2021.11.10.17.18.2}
\begin{pmatrix}
-\frac{\theta_{0,t_i}}{2} \cdot  
\frac{\hat{\delta}_2(\prod_{j \neq i}(t_i - t_j)^{n_j})}{\prod_{j \neq i}(t_i - t_j)^{n_j}} & *   \\
0 & 
\frac{\theta_{0,t_i}}{2} \cdot  
\frac{\hat{\delta}_2(\prod_{j \neq i}(t_i - t_j)^{n_j})}{\prod_{j \neq i}(t_i - t_j)^{n_j}}
+\hat{\delta}_2(\theta_{0,t_i})
\end{pmatrix}\frac{1}{(x-t_i)^{n_i}}+
[\text{ higher order terms }].
\end{equation}
Since $\delta_1(\widehat\Omega_{(\boldsymbol{t}_{\mathrm{ra}},\boldsymbol{\theta}_0)}^{(n-2)} ) 
\hat{\delta}_2(g)g^{-1}$
is holomorphic at $t_i$, we have the following equalities:
\begin{equation}\label{2021.1.2.19.08}
\begin{aligned}
&-\frac{1}{4} \cdot
\mathrm{res}_{\zeta_i=0}
\mathrm{Tr} \left( 
\delta_1(\widehat\Omega_{(\boldsymbol{t}_{\mathrm{ra}},\boldsymbol{\theta}_0)}^{(n-2)}(\zeta_i^2) )
u^{(2)} -u^{(1)}
\hat{\delta}_2(
\widehat\Omega_{(\boldsymbol{t}_{\mathrm{ra}},\boldsymbol{\theta}_0)}^{(n-2)}(\zeta_i^2)) \right) \\
&=-\frac{1}{2} \cdot
\mathrm{res}_{x=t_i} \mathrm{Tr}
 \left( \delta_1(\widehat\Omega_{(\boldsymbol{t}_{\mathrm{ra}},\boldsymbol{\theta}_0)}^{(n-2)}) 
\hat{\delta}_2(g)g^{-1} - \delta_1(g)g^{-1}
\hat{\delta}_2(\widehat\Omega_{(\boldsymbol{t}_{\mathrm{ra}},\boldsymbol{\theta}_0)}^{(n-2)} )  \right) \\
&=\frac{1}{2} \cdot
\mathrm{res}_{x=t_i} \mathrm{Tr} \left( \delta_1(g) g^{-1}
\hat{\delta}_2(\widehat\Omega_{(\boldsymbol{t}_{\mathrm{ra}},\boldsymbol{\theta}_0)}^{(n-2)} ) 
 \right)\\
 &=\frac{1}{2} \cdot
\mathrm{res}_{x=t_i} \mathrm{Tr} \left( \left( g_0^{-1} \delta_1(g) \right) g^{-1} g_0
\left(g_0^{-1}
\hat{\delta}_2(\widehat\Omega_{(\boldsymbol{t}_{\mathrm{ra}},\boldsymbol{\theta}_0)}^{(n-2)} )
g_0 \right) 
 \right)\\
&= \frac{\delta_1\left( \theta_{2n_i-1,t_i}\right)}{4\theta_{1,t_i}}  \hat{\delta}_2(\theta_{0,t_i})
+2\cdot \frac{\delta_1(\theta_{2n_i-1,t_i})}{8} \cdot 
\frac{\theta_{0,t_i}}{\theta_{1,t_i}} \cdot 
 \frac{\hat{\delta}_2(\prod_{j \neq i}(t_i - t_j)^{n_j})}{\prod_{j \neq i}(t_i - t_j)^{n_j}}\\
 &= \frac{\delta_1(\theta_{2n_i-1,t_i})}{4} \cdot 
\frac{\hat{\delta}_{2}(\theta_{0,t_i}  \prod_{j \neq i}(t_i - t_j)^{n_j})}
{\theta_{1,t_i} \prod_{j \neq i}(t_i - t_j)^{n_j}}.
\end{aligned}
\end{equation}
Remark that 
$\delta_1 \in 
\Theta_{(\widehat{\mathcal{M}}_{\boldsymbol{t}_{\text{ra}}}\times T_{\boldsymbol{\theta}}) 
/((T_{\boldsymbol{t}})_{\boldsymbol{t}_{\text{ra}}} \times T_{\boldsymbol{\theta}})}$.
The fourth equality follows from the equalities \eqref{2021.11.10.17.18} and 
\eqref{2021.11.10.17.18.2}.

Next we calculate the residue of the
first term of the right hand side of \eqref{2021.4.30.18.44} at $\zeta_i=0$.
We calculate this residue as follows:
$$
\begin{aligned}
&-\mathrm{res}_{\zeta_i=0} \mathrm{Tr} \left( \delta_1(\tilde{\Omega}_{\zeta_i}') \tilde{u}^{(2)} -
\tilde{u}^{(1)}\hat{\delta}_2(\tilde{\Omega}_{\zeta_i}')  \right)\\
&= -\mathrm{res}_{\zeta_i=0}
\mathrm{Tr} \left( \delta_1(\tilde{\Omega}_{\zeta_i}') M_{\zeta_i}^{-1}g^{-1} \hat{\delta}_{2} ( gM_{\zeta_i} ) -
M_{\zeta_i}^{-1}g^{-1}\hat{\delta}_{2} ( gM_{\zeta_i}) \hat{\delta}_2(\tilde{\Omega}_{\zeta_i}')  \right)\\
&= -\mathrm{res}_{\zeta_i=0}
\mathrm{Tr} \left( \delta_1(M_{\zeta_i}\tilde{\Omega}_{\zeta_i}' M_{\zeta_i}^{-1}) g^{-1} \hat{\delta}_{2} ( g ) -
g^{-1}\hat{\delta}_{2} ( g) \hat{\delta}_2(M_{\zeta_i}\tilde{\Omega}_{\zeta_i}'M_{\zeta_i}^{-1})  \right)\\
&=-2\cdot \mathrm{res}_{x=t_i}
\mathrm{Tr} \left( \delta_1(\tilde{\Omega}_{i}) g^{-1}\hat{\delta}_{2} ( g) -
g^{-1}\delta_{1} ( g)\hat{\delta}_2(\tilde{\Omega}_i)  \right).
\end{aligned}
$$
Here the last equality follows from $\delta_1(M^{-1}_{\zeta_i}dM_{\zeta_i})=
\hat{\delta}_2(M^{-1}_{\zeta_i}dM_{\zeta_i})=0$.
The coefficients of the expansion of $\delta_1(\tilde{\Omega}_{i})$ at $x=t_i$
vanish until the $x_{t_i}^{-2}$-term.
The $(1,2)$-entry of the $x_{t_i}^{-1}$-term of $\tilde{\Omega}_i$ 
depends on $(q_j,p_j)_{1\le j\le n-3}$ 
and the other entries of the $x_{t_i}^{-1}$-term of $\tilde{\Omega}_i$ 
are independent of $(q_j,p_j)_{1\le j\le n-3}$.
The $(1,2)$-entry of the $x_{t_i}^{-1}$-term of $\delta_1(\tilde{\Omega}_{i})$ 
is $\delta_1(\theta_{2n_i-1,t_i})/2$ and the other entries are zero.
On the other hand, the $(2,1)$-entry of the 
leading coefficient of $g^{-1}\hat{\delta}_{2} ( g)$ is 
$\hat{\delta}_{2}(\theta_{0,t_i}  \prod_{j \neq i}(t_i - t_j)^{n_j})/
(\theta_{1,t_i} \prod_{j \neq i}(t_i - t_j)^{n_j})$, 
since we set $g_0$ as \eqref{2021.1.2.18.21}.
Then the residue of $\mathrm{Tr}( \delta_1(\tilde{\Omega}_{i}) g^{-1}\hat{\delta}_{2} ( g))dx$ at $t_i$
is 
\begin{equation*}
\begin{aligned}
&\frac{\delta_1(\theta_{2n_i-1,t_i})}{2} \cdot 
\frac{\hat{\delta}_{2}(\theta_{0,t_i}  \prod_{j \neq i}(t_i - t_j)^{n_j})}
{\theta_{1,t_i} \prod_{j \neq i}(t_i - t_j)^{n_j}}.
\end{aligned}
\end{equation*}
We consider the residue of
$\mathrm{Tr}(g^{-1}\delta_{1} ( g)\hat{\delta}_2(\tilde{\Omega}_i)  )$ at $t_i$.
By \eqref{2020.3.5.21.08}, the residue of
$\mathrm{Tr}(g^{-1}\delta_{1} ( g)\hat{\delta}_2(\tilde{\Omega}_i)  )$ at $t_i$
vanishes.
Then we have
\begin{equation}\label{2021.4.30.20.14}
\begin{aligned}
&-\frac{1}{4} \cdot \mathrm{res}_{\zeta_i=0} 
\mathrm{Tr} ( \delta_1(\tilde{\Omega}_{\zeta_i}') \tilde{u}^{(2)} -
\tilde{u}^{(1)}\hat{\delta}_2(\tilde{\Omega}_{\zeta_i}')  ) \\
&=-\frac{1}{2} \cdot \mathrm{res}_{x=t_i} \mathrm{Tr} 
( \delta_1(\tilde{\Omega}_{i}) g^{-1}\hat{\delta}_{2} ( g) -
g^{-1}\delta_{1} ( g)\hat{\delta}_2(\tilde{\Omega}_i)  ) \\
&=
-\frac{\delta_1(\theta_{2n_i-1,t_i})}{4} \cdot 
\frac{\hat{\delta}_{2}(\theta_{0,t_i}  \prod_{j \neq i}(t_i - t_j)^{n_j})}
{\theta_{1,t_i} \prod_{j \neq i}(t_i - t_j)^{n_j}}.
\end{aligned}
\end{equation}

By combining \eqref{2021.4.30.18.44}, \eqref{2021.1.2.19.08} and \eqref{2021.4.30.20.14}, we have 
\begin{equation*}
\begin{aligned}
&\frac{1}{4} \cdot \mathrm{res}_{\zeta_i=0} \mathrm{Tr} \left( \hat{\delta} 
(\widehat\Omega_{(\boldsymbol{t}_{\mathrm{ra}},\boldsymbol{\theta}_0)}^{(n-2)}(\zeta^2_i) ) 
\wedge \hat{\delta}(\psi_{\zeta_i})\psi_{\zeta_i}^{-1} \right) \\
&=\frac{1}{4} \cdot \mathrm{res}_{\zeta_i=0}\mathrm{Tr} \left( \hat{\delta} (\tilde{\Omega}'_{\zeta_i}) 
\wedge \hat{\delta}(\tilde{\psi}'_{\zeta_i})(\tilde{\psi}_{\zeta_i}')^{-1}   \right)   
 -\frac{1}{4} \cdot \mathrm{res}_{\zeta_i=0} \mathrm{Tr}\left( d( \psi_{\zeta_i}^{-1} u^{(1)} \hat\delta_2(\psi_{\zeta_i})
  -\psi_{\zeta_i}^{-1} u^{(2)} \delta_1(\psi_{\zeta_i}))
  \right)\\
  &= \frac{1}{4} \cdot \mathrm{res}_{\zeta_i=0}\mathrm{Tr} \left( \hat{\delta} (\tilde{\Omega}'_{\zeta_i}) 
\wedge \hat{\delta}(\tilde{\psi}'_{\zeta_i})(\tilde{\psi}_{\zeta_i}')^{-1}   \right) .
\end{aligned}
\end{equation*}
By the equality \eqref{2020.1.24.11.10}, we have
\begin{equation*}
\frac{1}{4} \cdot \mathrm{res}_{\zeta_i=0} \mathrm{Tr} \left( \hat{\delta} 
(\widehat\Omega_{(\boldsymbol{t}_{\mathrm{ra}},\boldsymbol{\theta}_0)}^{(n-2)}(\zeta^2_i) ) 
\wedge \hat{\delta}(\psi_{\zeta_i})\psi_{\zeta_i}^{-1} \right) 
=\left(
\sum_{i \in I_{\mathrm{ra}}} 
\sum_{l'=0}^{2n_{i}-3}  d  H_{\theta_{l',t_i}} \wedge d \theta_{l',t_i}
\right) (\hat{\delta}_1 , \hat{\delta}_2)
\end{equation*}
when $\delta_1 \in 
\Theta_{(\widehat{\mathcal{M}}_{\boldsymbol{t}_{\text{ra}}}\times T_{\boldsymbol{\theta}}) 
/((T_{\boldsymbol{t}})_{\boldsymbol{t}_{\text{ra}}} \times T_{\boldsymbol{\theta}})}$.
So we have $\hat{\omega} (\delta_1,\hat{\delta}_2)- \hat{\omega}'(\delta_1,\hat{\delta}_2)=0$
when $\delta_1 \in 
\Theta_{(\widehat{\mathcal{M}}_{\boldsymbol{t}_{\text{ra}}}\times T_{\boldsymbol{\theta}}) 
/((T_{\boldsymbol{t}})_{\boldsymbol{t}_{\text{ra}}} \times T_{\boldsymbol{\theta}})}$.
Then we obtain the assertion of this theorem.
\end{proof}

By Theorem \ref{2020.1.31.9.16}, Theorem \ref{2020.2.2.19.05} and Theorem \ref{2020.1.29.11.46},
we obtain the following corollary:
\begin{Cor}\label{2020.2.2.19.01}
\textit{
Set $\eta_j:= \frac{p_j}{P(q_j;\boldsymbol{t})}
 - \sum_{i=1}^{\nu} \frac{D_i(q_j;\boldsymbol{t}, \boldsymbol{\theta})}{(q_j -t_i)^{n_i}} 
 - D_{\infty}(q_j;\boldsymbol{t}, \boldsymbol{\theta})$.
The vector fields $\delta_{\theta_{l,t_i}^{\pm}}^{\mathrm{IMD}}$ 
$(i\in I_{\mathrm{un}}$ and $l=0,1,\ldots,n_i-2)$,
$\delta_{t_i}^{\mathrm{IMD}}$ $(i\in \{ 3,4,\ldots,\nu\} \cap (I_{\mathrm{reg}} \cup I_{\mathrm{un}}))$,
and $\delta_{\theta_{l',t_i}}$ $(i\in I_{\mathrm{ra}}$ and $l'=0,1,\ldots,2n_i-3)$
have the following
hamiltonian description}: \textit{
\begin{equation*}
\begin{aligned}
\delta_{\theta_{l,t_i}^{\pm}}^{\mathrm{IMD}}&= 
\frac{\partial}{\partial \theta_{l,t_i}^{\pm}}
-\sum_{j=1}^{n-3} \left( 
 \frac{\partial H_{\theta_{l,t_i}^{\pm}}}{\partial \eta_j} 
\frac{\partial}{\partial q_j} 
- \frac{\partial H_{\theta_{l,t_i}^{\pm}}}{\partial q_j} 
\frac{\partial}{\partial \eta_j}
\right), \\
\delta_{t_i}^{\mathrm{IMD}}&= 
\frac{\partial}{\partial t_i}
-\sum_{j=1}^{n-3} \left( 
 \frac{\partial  H_{t_i}}{\partial \eta_j} 
\frac{\partial}{\partial q_j} 
- \frac{\partial H_{t_i}}{\partial q_j} 
\frac{\partial}{\partial \eta_j}
\right),\text{ and} \\
\delta_{\theta_{l',t_i}}^{\mathrm{IMD}}&= 
\frac{\partial}{\partial \theta_{l',t_i}}
-\sum_{j=1}^{n-3} \left( 
 \frac{\partial H_{\theta_{l',t_i}}}{\partial \eta_j} 
\frac{\partial}{\partial q_j} 
- \frac{\partial H_{\theta_{l',t_i}}}{\partial q_j} 
\frac{\partial}{\partial \eta_j}
\right),
\end{aligned}
\end{equation*}
respectively.
}
\end{Cor}

\section{Examples}

\subsection{Example ($\nu=2, n_1=n_2=n_{\infty}=2$)}

We consider the connection $d + \Omega$ 
on $\mathcal{O}\oplus \mathcal{O}(4)$ with the following connection matrix:
\begin{equation*}
\Omega=
\begin{pmatrix}
0& \frac{1}{P(x)} \\
c_0(x) & d_0(x)
\end{pmatrix}dx.
\end{equation*}
Here we put
$P(x):=x^2(x-1)^{2}$,
\begin{equation*}
\begin{aligned}
c_0(x)&:=\frac{C_0^{(0)} + C_0^{(1)}x}{x^2}+ \frac{C_1^{(0)} + C_1^{(1)}(x-1)}{(x-1)^{2}}
+\sum_{j=1}^{3}\frac{p_j}{x-q_j} \\
&\qquad+\tilde{C}^{(0)}+\tilde{C}^{(1)}x+\tilde{C}^{(2)}x^2+ C^{(0)}_{\infty}x^3
+C^{(1)}_{\infty}x^4   
, \text{ and} \\
d_0(z)&:=\frac{D_0^{(0)} + D_0^{(1)}x}{x^2}
+ \frac{D_1^{(0)} + D_1^{(1)}(x-1)}{(x-1)^{2}}
+\sum_{j=1}^{n-3}\frac{-1}{x-q_j} +D_{\infty}^{(0)}.
\end{aligned}
\end{equation*}
We set $t_1:=0$, $t_2:=1$, and $t_{\infty}:=\infty$.
The polar divisor of the connection $d+ \Omega$ is $2\cdot t_1+2\cdot t_2 + 2 \cdot t_{\infty}
+q_1+q_2+q_3$.
We assume that
the leading coefficients $\Omega_{t_i}(0)$ are semi-simple for $i=1,2,\infty$.
We put $x_{t_i}:=x-t_i$ for $i=1,2$ and $x_{t_{\infty}}=w$.
We fix the formal type of the negative part of $d+\Omega$ for each $t_i$.
That is, we fix $\theta_{l,t_i}^{\pm}$ for $l=0,1$ and $i=1,2,\infty$,
and the negative part of $d+\Omega$ for each $t_i$ 
is diagonalizable as
\begin{equation*}
\frac{\begin{pmatrix}
\theta_{0,t_i}^{+} & 0\\
0 & \theta_{0,t_i}^{-} 
\end{pmatrix}}{x_{t_i}^2}
+\frac{\begin{pmatrix}
\theta_{1,t_i}^{+} & 0\\
0 & \theta_{1,t_i}^{-} 
\end{pmatrix}}{x_{t_i}}
\end{equation*}
by a formal transformation (see Section \ref{2020.1.21.14.24}).
Then the coefficients of $c_2$ and $d_2$ are determined as follows. 
\begin{equation}\label{2020.1.16.14.21}
\begin{cases}
C_0^{(0)}+C_0^{(1)}x =  - \theta_{0,0}^+\theta_{0,0}^- 
+(2\theta_{0,0}^+\theta_{0,0}^-  - \theta_{0,0}^+\theta_{1,0}^-  - \theta_{0,0}^-\theta_{1,0}^+ )x \\
D_0^{(0)}+D_0^{(1)}x =\theta_{0,0}^+ + \theta_{0,0}^- 
+(\theta_{1,0}^+ + \theta_{1,0}^- )x,
\end{cases}
\end{equation}
\begin{equation}\label{2020.1.16.14.22}
\begin{cases}
C_1^{(0)}+C_1^{(1)}(x-1) =  - \theta_{0,1}^+\theta_{0,1}^- 
-(2\theta_{0,1}^+\theta_{0,1}^-  + \theta_{0,1}^+\theta_{1,1}^-  + \theta_{0,1}^-\theta_{1,1}^+ )(x - 1)\\
D_1^{(0)}+D_1^{(1)}(x-1) =\theta_{0,1}^+ + \theta_{0,1}^- 
+(\theta_{1,1}^+ + \theta_{1,1}^- )(x - 1),
\end{cases}
\end{equation}
and
\begin{equation}\label{2020.1.16.14.23}
\begin{cases}
C_{\infty}^{(0)}+C_{\infty}^{(1)}x =  2\theta_{0,\infty}^+ \theta_{0,\infty}^- 
- \theta_{0,\infty}^- \theta_{1,\infty}^+ 
- \theta_{0,\infty}^+ \theta_{1,\infty}^-
- (\theta_{0,\infty}^+\theta_{0,\infty}^-)x  \\
D_{\infty}^{(0)} = - \theta_{0,\infty}^+-\theta_{0,\infty}^-.
\end{cases}
\end{equation}
Moreover we assume that $q_1$, $q_2$, and $q_3$ are apparent singularities.
We define $\tilde{C}_{q_j}$ for $j=1,2,3$ so that 
$\tilde{C}^{(0)}+\tilde{C}^{(1)}x+\tilde{C}^{(2)}x^2$
is equal to
\begin{equation}\label{2020.1.16.14.24}
\tilde{C}_{q_1}(x-q_2)(x-q_3)+ \tilde{C}_{q_2}(x-q_1)(x-q_3)+\tilde{C}_{q_3}(x-q_1)(x-q_2).
\end{equation}
Since $q_1$, $q_2$, and $q_3$ are apparent singularities, we have
\begin{equation}\label{2020.1.16.14.25}
\tilde{C}_{q_j}=\frac{1}{Q'(q_j)}\left(
\frac{p_j^2}{q_j^2(q_j-1)^2} - \sum_{i=1,2}
\frac{ D_i(q_j)p_j + C_i(q_j)}{(q_j-t_i)^2} +\sum_{k\in \{1,2,3\} \setminus\{j\}} \frac{p_j-p_k}{q_j-q_k} 
- D_{\infty}^{(0)} p_j- C_{\infty}^{(0)}q_j^3 - C_{\infty}^{(1)} q_j^4 \right)
\end{equation}
for $j=1,2,3$, 
where we put $Q(x):=(x-q_1)(x-q_2)(x-q_3)$.
We determine the matrices $\Phi_i$ and $\Xi_i$ 
as in Lemma \ref{2019.12.30.22.09} as follows.
\begin{equation*}
\begin{aligned}
\Phi_0= 
\begin{pmatrix}
1 & \frac{1}{\theta_{0,0}^-} \\
\theta_{0,0}^+ & 1
\end{pmatrix}, 
\Xi_1^{(0)} &=
\begin{pmatrix}
0 & -\frac{2\theta_{0,0}^--\theta_{1,0}^-}{(\theta_{0,0}^+-\theta_{0,0}^-)\theta_{0,0}^-} \\
\frac{(2\theta_{0,0}^+-\theta_{1,0}^+)\theta_{0,0}^-}{\theta_{0,0}^+-\theta_{0,0}^-} & 0
\end{pmatrix},
\Xi_2^{(0)} =
\begin{pmatrix}
0 & (\xi_0^{(2)})_{12} \\
(\xi_0^{(2)})_{21} & 0
\end{pmatrix}\\
\Phi_1= 
\begin{pmatrix}
1 & \frac{1}{\theta_{0,1}^-} \\
\theta_{0,1}^+ & 1
\end{pmatrix}, 
\Xi_1^{(1)} &=
\begin{pmatrix}
0 & \frac{2\theta_{0,1}^-+\theta_{1,1}^-}{(\theta_{0,1}^+-\theta_{0,1}^-)\theta_{0,1}^-} \\
-\frac{(2\theta_{0,1}^++\theta_{1,1}^+)\theta_{0,1}^-}{\theta_{0,1}^+-\theta_{0,1}^-} & 0
\end{pmatrix},
\Xi_2^{(1)} =
\begin{pmatrix}
0 & (\xi_1^{(2)})_{12} \\
(\xi^{(2)}_1)_{21} & 0
\end{pmatrix},\text{ and}\\
\Phi_{\infty}= 
\begin{pmatrix}
1 &- \frac{1}{\theta_{0,\infty}^-} \\
-\theta_{0,\infty}^+ & 1
\end{pmatrix}, 
\Xi_1^{(\infty)} &=
\begin{pmatrix}
0 & 
\frac{2\theta_{0,\infty}^--\theta_{1,\infty}^-}{(\theta_{0,\infty}^+-\theta_{0,\infty}^-)\theta_{0,\infty}^-} \\
-\frac{(2\theta_{0,\infty}^+-\theta_{1,\infty}^+)\theta_{0,\infty}^-}{\theta_{0,\infty}^+-\theta_{0,\infty}^-} & 0
\end{pmatrix},
\Xi_2^{(\infty)} =
\begin{pmatrix}
0 & (\xi_1^{(2)})_{12} \\
(\xi^{(2)}_1)_{21} & 0
\end{pmatrix}.
\end{aligned}
\end{equation*}
Here the descriptions of $ (\xi_i^{(2)})_{12} $ and $(\xi_i^{(2)})_{21}$ are omitted.
Set $\Xi^{\le 2}_i:= \mathrm{id}+ \Xi^{(i)}_1 x_{t_i} + \Xi^{(i)}_2 x^2_{t_i}$.
Let $\theta_{2,t_i}^{\pm}$ ($i=1,2,\infty$) be the coefficient as in Lemma \ref{2019.12.30.22.09}.
That is,
\begin{equation*}
\begin{aligned}
&(\Phi_i \Xi^{\le 2}_i )^{-1}d(\Phi_i \Xi^{\le 2}_i )+
(\Phi_i \Xi^{\le 2}_i )^{-1} \Omega (\Phi_i \Xi^{\le 2}_i ) \\
&=\frac{\begin{pmatrix}
\theta_{0,t_i}^{+} & 0\\
0 & \theta_{0,t_i}^{-} 
\end{pmatrix}}{x_{t_i}^2}
+\frac{\begin{pmatrix}
\theta_{1,t_i}^{+} & 0\\
0 & \theta_{1,t_i}^{-} 
\end{pmatrix}}{x_{t_i}}
+\begin{pmatrix}
\theta_{2,t_i}^{+} & 0\\
0 & \theta_{2,t_i}^{-} 
\end{pmatrix}+ O(x_{t_i})
\end{aligned}
\end{equation*}
for $i=0,1,\infty$.
Remark that $\Xi^{\le 2}_i$ is degree $2$ in $x_{t_i}$.
This degree is sufficient to define Hamiltonians
since there is no parameter corresponding to the positions of irregular singularities.
By the equations
(\ref{2020.1.16.14.21}),
(\ref{2020.1.16.14.22}),
(\ref{2020.1.16.14.23}),
(\ref{2020.1.16.14.24}),
and (\ref{2020.1.16.14.25}), 
we can determine the Hamiltonians $H_{\theta_{0,0}^{\pm}}$, 
$H_{\theta_{0,1}^{\pm}}$, and $H_{\theta_{0,\infty}^{\pm}}$ as follows.
\begin{equation*}
\begin{aligned}
H_{\theta_{0,0}^{\pm}}:=-\theta_{2,0}^{\pm}=\frac{-1}{\theta_{0,0}^{\pm}-\theta_{0,0}^{\mp}}
&\Big( \theta_{0,0}^{\pm}\theta_{0,0}^{\mp}+ \theta_{1,0}^{\pm}\theta_{1,0}^{\mp}
-2(\theta_{0,0}^{\pm}\theta_{1,0}^{\mp}+\theta_{1,0}^{\pm}\theta_{0,0}^{\mp})\\
& + \tilde{C}^{(0)} +(C_1^{(0)} - C_1^{(1)})
 + (D_1^{(0)} - D_1^{(1)} + D_{\infty}^{(0)} )\theta_{0,0}^{\pm}\\
&-\frac{ p_1 - \theta_{0,0}^{\pm}}{q_1}
-\frac{  p_2 - \theta_{0,0}^{\pm} }{q_2}
-\frac{p_3 - \theta_{0,0}^{\pm} }{q_3} \Big),\\
H_{\theta_{0,1}^{\pm}}:=-\theta_{2,1}^{\pm}=\frac{-1}{\theta_{0,1}^{\pm}-\theta_{0,1}^{\mp}}
&\Big( \theta_{0,1}^{\pm}\theta_{0,1}^{\mp}+ \theta_{1,1}^{\pm}\theta_{1,1}^{\mp}
+2(\theta_{0,1}^{\pm}\theta_{1,1}^{\mp}+\theta_{1,1}^{\pm}\theta_{0,1}^{\mp})  \\
&+
( \tilde{C}^{(0)} +\tilde{C}^{(1)} +\tilde{C}^{(2)} )
+( C_0^{(0)} + C_0^{(1)}) + (C_{\infty}^{(0)}  + C_{\infty}^{(1)} )
+ (D_0^{(0)} + D_0^{(1)} + D_{\infty}^{(0)})\theta_{0,1}^{\pm}\\
&-\frac{p_1 - \theta_{0,1}^{\pm} }{q_1-1}
-\frac{ p_2 - \theta_{0,1}^{\pm} }{q_2-1}
-\frac{ p_3 - \theta_{0,1}^{\pm}}{q_3-1} \Big), \text{ and}\\
H_{\theta_{0,\infty}^{\pm}}:=-\theta_{2,\infty}^{\pm}=\frac{-1}{\theta_{0,\infty}^{\pm}-\theta_{0,\infty}^{\mp}} 
&\Big(
\theta_{0,\infty}^{\pm} \theta_{0,\infty}^{\mp} + \theta_{1,\infty}^{\pm} \theta_{1,\infty}^{\mp}
 - 2(\theta_{0,\infty}^{\pm}\theta_{1,\infty}^{\mp}  +\theta_{1,\infty}^{\pm}\theta_{0,\infty}^{\mp})\\
&+\tilde{C}^{(2)}-(D_0^{(0)}  + D_1^{(0)} +D_1^{(1)}) \theta_{0,\infty}^{\pm}
+ (q_1 + q_2 + q_3) \theta_{0,\infty}^{\pm} \Big).
\end{aligned}
\end{equation*}
Set $\eta_j:=\frac{p_j}{q_j^2(q_j-1)^2}-
\frac{D_0^{(0)} + D_0^{(1)}q_j}{q_j^2}
- \frac{D_1^{(0)} + D_1^{(1)}(q_j-1)}{(q_j-1)^{2}}
-D_{\infty}^{(0)}$ for $j=1,2,3$.
By Corollary \ref{2020.1.21.15.11}, the vector field determined by the 
generalized isomonodromic deformations is described as
\begin{equation*}
\frac{\partial}{\partial \theta_{0,t_i}^{\pm}}
-\sum_{j=1}^{3} \left( 
 \frac{\partial H_{\theta_{0,t_i}^{\pm}}}{\partial \eta_j} 
\frac{\partial}{\partial q_j} 
- \frac{\partial H_{\theta_{0,t_i}^{\pm}}}{\partial q_j} 
\frac{\partial}{\partial \eta_j}
\right).
\end{equation*}

\subsection{Example corresponding to Kimura's $L(9/2;2)$}

In this section, we consider Kimura's family $L(9/2;2)$ of rank $2$ linear differential equations 
in \cite[p.37]{Kimura}. 
We describe the corresponding global normal form (see \cite[Section 6]{DF})
and consider the integrable deformations of the family given by 
the global normal form.
Then we can reproduce Kimura's Hamiltonian $H(9/2)$
in \cite[p.40]{Kimura}.

Let $D$ be the effective divisor defined as 
$D=5 \cdot \infty$.
We consider the connection $d+\Omega^{(\infty)}$ 
on $\mathcal{O}_{\mathbb{P}^1} \oplus \mathcal{O}_{\mathbb{P}^1}(3)$
with
\begin{equation*}
\Omega^{(\infty)}=
\begin{pmatrix}
0& -\frac{1}{w^5} \\
c^{(\infty)}_0(w) & d^{(\infty)}_0(w)
\end{pmatrix}dw
\end{equation*}
where 
\begin{equation}\label{2020.1.29.11.49}
\begin{aligned}
c^{(\infty)}_0(w)&:=
 -\frac{9}{w^4} 
 - \frac{9 t_1}{w^2}
 -\frac{3 t_2}{w} 
 -3 K_2 
 -3 K_1 w
-\sum_{i=1}^2\frac{p_iw^2}{1-q_iw}, \text{ and} \\
d^{(\infty)}_0(w)&:=\sum_{i=1}^2\frac{1}{w(1-q_iw)} - \frac{3}{w}
\end{aligned}
\end{equation}
(see \cite[Section 6]{DF}).
The polar divisor is $D + q_1+q_2$.
Assume that $w=1/q_1$ and $w=1/q_2$ are apparent singularities.
Then we can determine $K_1$ and $K_2$ as
rational functions whose variables are $t_1,t_2,q_1,q_2,p_1,p_2$.

If we set $\Phi_{\infty}:= 
\begin{pmatrix}
1 & 0 \\
0 & -3
\end{pmatrix}$ and
\begin{equation}\label{2020.1.29.11.50}
\begin{aligned}
\Xi_{\infty}^{\le 6}:=&\ \mathrm{id} 
+\begin{pmatrix}
-\frac{t_1}{4} & 0 \\
0& \frac{t_1}{4}
\end{pmatrix}w^2
+\begin{pmatrix}
-\frac{t_2}{12} & \frac{1}{24} \\
0& \frac{t_2}{12}
\end{pmatrix}w^3
+\begin{pmatrix}
\frac{t_1^2}{8} - \frac{K_2}{12} & - \frac{q_1+q_2}{12} \\
-\frac{1}{24}& -\frac{t_1^2}{8} + \frac{K_2}{12}
\end{pmatrix}w^4 \\
&\qquad +\begin{pmatrix}
 \frac{t_1t_2}{12} -  \frac{K_1}{12} & 0 \\
\frac{q_1+q_2}{12}
&  -\frac{t_1t_2}{12} +  \frac{K_1}{12}
\end{pmatrix}w^5
+\begin{pmatrix}
* & * \\
* & *
\end{pmatrix}w^6,
\end{aligned}
\end{equation}
then we have  
\begin{equation}\label{2020.1.29.11.51}
\begin{aligned}
\tilde{\Omega}_\infty=&
\frac{\begin{pmatrix}
0& 3 \\
0 & 0
\end{pmatrix}}{w^5}
+
\frac{\begin{pmatrix}
0&  0 \\
3 & 0
\end{pmatrix}}{w^4}
+
\frac{\begin{pmatrix}
0 & \frac{3t_1}{2} \\
0 &   0
\end{pmatrix}}{w^3}
+
\frac{\begin{pmatrix}
0 & \frac{t_2}{2} \\
\frac{3t_1}{2} & 0
\end{pmatrix}}{w^2} 
+
\frac{\begin{pmatrix}
-\frac{1}{4} & b_3 \\
\frac{t_2}{2} & -\frac{1}{4} -\frac{1}{2}
\end{pmatrix}}{w} +
\begin{pmatrix}
a_1 & b_4 \\
b_3 & a_1
\end{pmatrix}
+ O(w)^2,
\end{aligned}
\end{equation}
where $\tilde{\Omega}_\infty dw:=
 (\Phi_{\infty} \Xi^{\le 6}_{\infty})^{-1} \Omega^{(\infty)} (\Phi_{\infty} \Xi^{\le 6}_{\infty})
+ (\Phi_{\infty} \Xi^{\le 6}_{\infty})^{-1}d(\Phi_{\infty} \Xi^{\le 6}_{\infty}) $.
We have
\begin{equation*}
a_1=  \frac{q_1}{2}+\frac{q_2}{2},\ 
b_3= -\frac{3t_1^2}{8} + \frac{K_2}{2}, \text{ and } 
b_4= - \frac{t_1 t_2}{4}  + \frac{K_1}{2}.
\end{equation*}
After ramification
$w=\zeta^2$
and the following transformation of $\tilde{\Omega}_\infty$
\begin{equation*}
\tilde{\Omega}'_{\zeta} d \zeta :=M_{\zeta}^{-1} (\tilde{\Omega}_\infty dw) M_{\zeta} 
+M^{-1}_{\zeta} d M_{\zeta},
\text{ where }
M_{\zeta} = 
\begin{pmatrix}
1 & 1 \\
\zeta & -\zeta
\end{pmatrix},
\end{equation*}
we have an unramified irregular singular point with matrix connection
\begin{equation*}
\begin{aligned}
\tilde{\Omega}'_{\zeta}=&
\frac{\begin{pmatrix}
6& 0 \\
0 & -6
\end{pmatrix}}{\zeta^8}
+
\frac{\begin{pmatrix}
3t_1 & 0 \\
0 &   -3t_1
\end{pmatrix}}{\zeta^4}
+
\frac{\begin{pmatrix}
t_2 & 0 \\
0 & -t_2
\end{pmatrix}}{\zeta^2} 
+
\frac{\begin{pmatrix}
-\frac{1}{2} & 0 \\
0 & -\frac{1}{2} 
\end{pmatrix}}{\zeta} \\
&+
\begin{pmatrix}
2b_3 & 0 \\
0 & -2b_3
\end{pmatrix}
+
\begin{pmatrix}
2a_1 & 0 \\
0 & 2a_1
\end{pmatrix}\zeta
+
\begin{pmatrix}
2b_4 & 0 \\
0 & -2b_4
\end{pmatrix}\zeta^2
+O(\zeta)^3.
\end{aligned}
\end{equation*}
We define Hamiltonians 
\begin{equation*}
\begin{aligned}
H_1:=&\ -[\text{ the coefficient of the $\zeta^{3}$-term of 
$\sum_{k=9}^{\infty} \theta_{k,\infty} \int \zeta^{k-9} d\zeta$  }]\\
=&\ -\frac{2 b_4}{3} =
-\frac{K_1}{3}+\frac{t_1 t_2}{6} \text{ and } \\
H_2:=&\ - [\text{ the coefficient of the $\zeta$-term of 
$\sum_{k=9}^{\infty} \theta_{k,\infty} \int \zeta^{k-9} d\zeta$ }] \\
=&\  - 2 b_3 =  - K_2 + \frac{3t_1^2}{4}.
\end{aligned}
\end{equation*}
Then the $2$-form $\hat{\omega}'$ defined in Theorem \ref{2020.1.29.11.46} is described as 
\begin{equation}\label{2020.2.2.10.22}
\begin{aligned}
\hat{\omega}'&= \sum_{i=1,2} dp_i \wedge dq_i
+ dH_1 \wedge d(3t_1) + dH_2 \wedge d(t_2) \\
&=  - \left(\sum_{i=1,2} d\eta_i \wedge dq_i
- d (3 H_1) \wedge dt_1 - dH_2 \wedge dt_2 \right)\\
&=  - \left(\sum_{i=1,2} d\eta_i \wedge dq_i
+ dK_1 \wedge dt_1 + dK_2 \wedge dt_2 -  t_1 dt_1 \wedge dt_2 \right),
\end{aligned}
\end{equation}
where $\eta_i := -p_i$ for $i=1,2$.
By Theorem \ref{2020.1.31.9.16} and Theorem \ref{2020.1.29.11.46},
the vector field determined by the 
integrable deformations is described as 
\begin{equation*}
\frac{\partial}{\partial t_i}
-\sum_{j=1}^{2} \left( 
 \frac{\partial K_i}{\partial \eta_j} 
\frac{\partial}{\partial q_j} 
- \frac{\partial K_i}{\partial q_j} 
\frac{\partial}{\partial \eta_j}
\right)
\end{equation*}
for $i=1,2$. 
This description is given in \cite{Kimura}.

\begin{Rem}
We may check that $\hat{\omega}=\hat{\omega}'$
by the calculation of the right hand side of (\ref{2020.1.8.12.55}) for
(\ref{2020.1.29.11.49}), (\ref{2020.1.29.11.50}), and (\ref{2020.1.29.11.51}).
Then the $2$-form \eqref{2020.2.2.10.22} is the isomonodromy $2$-form.
In fact, we may check the equality
\begin{equation*}
\frac{\partial (3 H_1)}{\partial t_2} - \frac{\partial H_2}{\partial t_1} 
- \sum_{i=1,2}\left(  
\frac{\partial (3 H_1)}{\partial p_i} \frac{\partial  H_2}{\partial q_i}
-\frac{\partial (3 H_1)}{\partial q_i} \frac{\partial H_2}{\partial p_i}
\right)=0
\end{equation*}
directly.
\end{Rem}


\noindent
{\bf Acknowledgments.}
The author would like to thank Professor Frank Loray for leading him 
to the subject treated in this paper and also for valuable discussions.
The author would like to thank Professor Masa-Hiko Saito
and Professor Michi-aki Inaba
for valuable comments and for warm encouragement.
In particular, Professor Inaba's comments have been very helpful
for the improvement of this paper.
He is supported by Japan Society for the Promotion of Science
KAKENHI Grant Numbers
17H06127,
18J00245, 
19K14506, and
22H00094.
He is very grateful to the anonymous referee's
insightful suggestions which helped to significantly improve the paper.


\end{document}